\documentclass[12pt,a4paper,twoside]{article}
\usepackage{pdfsync}
\usepackage[english]{babel}
\usepackage[utf8]{inputenc}
\usepackage{graphicx}
\graphicspath{{./figuras/}{./figuras/figuras/}{./}}
\usepackage{epsfig}
\usepackage[dvipsnames]{xcolor}
\usepackage[hypertex]{hyperref}
\usepackage{amsthm, amsmath,amssymb,color,srcltx}   
\input{xy}
\xyoption{all}
\usepackage{fancyhdr}

\def\cocoa{{\hbox{\rm C\kern-.13em o\kern-.07em C\kern-.13em o\kern-.15em A}}}
\def\deg{{\rm deg}}

\def\qed{\hfill  \framebox(5,5){}}
\def\det{{\rm det}}

\def\dim{{\rm dim}}

\def\Res{{\rm Res}}

\def\Sing{{\rm Sing}}
\def\Iso{{\rm Iso}}

\newcommand{\mult}{\mathop{\mathrm{mult}}\nolimits}

\newcommand{\Con}{\mathop{\mathrm{Con}}\nolimits}
\def\normal{\mathrm{nor}}

\newcommand{\nin}{\noindent}

\newtheorem{Theorem}{{\bf Theorem}}[section]
\newtheorem{Remark}[Theorem]{{\bf Remark}}
\newtheorem{Definition}[Theorem]{{\bf Definition}}
\newtheorem{Corollary}[Theorem]{{\bf Corollary}}
\newtheorem{Proposition}[Theorem]{{\bf Proposition}}

\newtheorem{Lemma}[Theorem]{{\bf Lemma}}
\newtheorem{Assumptions}[Theorem]{{\bf Assumptions}}

\newtheorem{Example}[Theorem]{{\bf Example}}
\newtheorem*{ExampleIntro}{{\bf Example}}

\def \C{\mathbb{C}}

\def \N{\mathbb{N}}
\def \P{\mathbb{P}}

\def \R{\mathbb{R}}

\def \Z{\mathbb{Z}}
\def \K{\mathbb{K}}

\def\bV {{\mathbf V}}

\def\cA{{\mathcal A}}
\def\cB{{\mathcal B}}
\def\cC{{\mathcal C}}

\def\cF {{\mathcal F}}
\def\cG {{\mathcal G}}

\def\cI {{\mathcal I}}
\def\cJ {{\mathcal J}}

\def\cL {{\mathcal L}}
\def\cM {{\mathcal M}}

\def\cO {{\mathcal O}}
\def\cP {{\mathcal P}}
\def\cQ {{\mathcal Q}}

\def\cT {{\mathcal T}}

\def\cV {{\mathcal V}}
\def\cW {{\mathcal W}}

\def\cZ {{\mathcal Z}}

\def\mfF {{\mathfrak  F}}
\def\mfG {{\mathfrak  G}}

\def\mfQ {{\mathfrak  Q}}
\def\mfR {{\mathfrak  R}}
\def\mfS {{\mathfrak  S}}

\usepackage{makeidx}
\makeindex

\title{Total Degree Formula for the Generic Offset to a Parametric Surface\footnote{\small Preprint of an article to be published at the International Journal of Algebra and Computation, World Scientific Publishing, DOI:10.1142/S0218196711006807}}
\author{San Segundo, F. \and Sendra, J.R.}

\begin{document}
\maketitle

\tableofcontents

\section*{Introduction}\label{sec:01-Introduction}

\nin This paper focuses on the study of the total degree w.r.t. the spatial variables of the multivariate polynomial defining the
generic offset to a rational surface in three-dimensional space.
So, before continuing with this introduction, let us describe, at least informally, the offsetting construction; for a detailed explanation see
Section \ref{sec:GenericOffsetDegreeProblem} and, more specifically, for the concept of generic offset, see Definition \ref{def:ch1:GenericOffset} (page \pageref{def:ch1:GenericOffset}).

\begin{figure}
\begin{center}
\includegraphics[width=10cm]{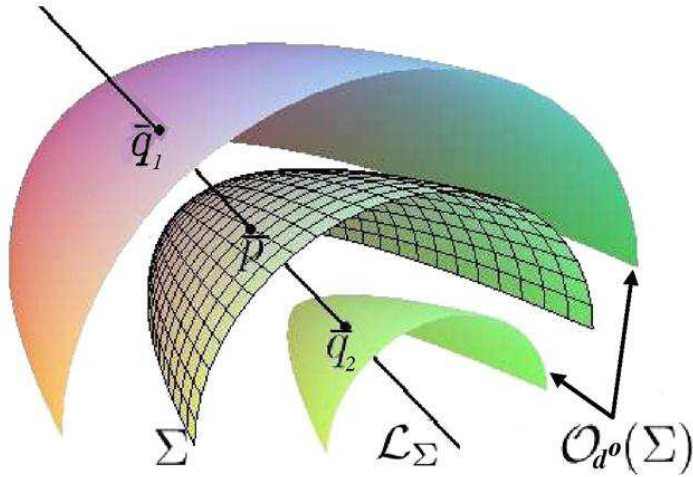}
\end{center}
\caption{Informal Definition of Offset to a Generating Surface \label{Figure:InformalDefinitionOffsetSurface}}
\end{figure}

\noindent Let $\Sigma$ be a surface in three-dimensional space. At each point $\bar p$ of $\Sigma$, consider the normal line
$\cL_{\Sigma}$ to the surface (assume, for this informal introduction, that the normal line to $\Sigma$ at $\bar p$
is well defined). Let $\bar q$ be a point on that line, at a distance $d^o$ of $\bar p$ (there are two such
points $\bar q$); equivalently, we consider the intersection points of $\cL_{\Sigma}$ with a sphere centered at
$\bar p$ and with radius $d^o$. The offset surface to $\Sigma$ at distance value $d^o$, is the
set $\cO_{d^o}(\Sigma)$ of all the points $\bar q$ obtained by this geometric construction,
illustrated in Figure \ref{Figure:InformalDefinitionOffsetSurface}.
$\Sigma$ is said to be the {\sf generating surface} of $\cO_{d^o}(\Sigma)$.

\nin The classical offset construction for algebraic hypersurfaces has been, and still is, an active research subject of scientific interest. Even though the historical origins of the study of offset curves can be traced back to the work of classical geometers (\cite{Leibniz1692},\cite{loria1902},\cite{salmon-treatise}), often under the denomination of {\em parallel curves}, the subject received increased attention when the technological advance in the fields of Computer Assisted Design and Computer Assisted Manufacturing (CAD/CAM) resulted in a strong demand of effective algorithms for the manipulation of curves and surfaces. We quote the following from one of the two seminal papers (\cite{Farouki1990},\cite{Farouki1990a}) by Farouki and Neff: {\em ``Apart from numerical-control machining, offset curves arise in a variety of practical applications such as tolerance analysis, geometric optics, robot path-planning, and in the formulation of simple geometric procedures (growing/shrinking, blending, filleting, etc.)''}. To the applications listed by these authors we should add here some recent ones, e.g. the connection with the medial axis transform; for these, and related applications see Chapter 11 in \cite{patrikalakis2002sic}, and the references contained therein.

\nin As a result of this interest coming from the applications, many new methods and algorithms have been developed by engineers and
mathematicians, and many geometric and algebraic properties of the offset construction have been studied in recent years; see, e.g. the
references \cite{alcazar2008},\cite{alcazar2008localsurfaces},\cite{alcazarLocalCurves2009},\cite{alcazar2007local},\cite{anton2005offset},\cite{Arrondo1997},\cite{Arrondo1999},\cite{Farouki1990},\cite{Farouki1990a},\cite{Hoffmann1989},\cite{hoschek1993fundamentals},\cite{Lu1995},
\cite{Lu95TR},\cite{Pottmann1995},\cite{pottmann1996rational},\cite{pottmann1998laguerre},\cite{SS05},\cite{SSS09},\cite{Sendra2000},\cite{Sendra2000a}.
In addition to these references, we also refer to the theses \cite{Alcazar2007}, \cite{PhdSanSegundo2010} and \cite{Sendra1999}, developed  within the
research group of Prof.\,J.R. Sendra. In \cite{Sendra1999} the fundamental algebraic properties of offsets to hypersurfaces are deduced, the
unirationality of the offset components are characterized, and  the genus problem (for the curve case) is studied. In \cite{Alcazar2007}
the topological behavior of the offset curve is analyzed. The present paper contains the results about the total degree of the generic offset to a rational surface in Chapter 4 of the Ph.D. Thesis \cite{PhdSanSegundo2010}
.

\nin With the exception of certain degenerated situations, that are indeed well known, the offset to an algebraic hypersurface is again a
hypersurface (see \cite{Sendra2000}). Thus, one might answer all the problems mentioned above (parametrization expressions, genus computation,
topologic types determination, degree analysis, etc), by applying the available algorithms to the resulting (offset) hypersurface. However,
in most cases, this strategy results unfeasible. The reason is that the offsetting process generates a huge size increment of the data
defining the offset in comparison to the data of the original variety.  The challenge, therefore, is to derive information (say algebraic or
geometric properties) of the offset hypersurface from the information that could be easily derived from the original (in general much simpler)
hypersurface.

\nin Framed in the above philosophy, the goal of this paper is to provide an efficient formula for the total degree of the generic offset to a rational
surface. More precisely: let $f(y_1,y_2,y_3)$ be the defining polynomial of  $\Sigma$, and let us treat $d$ as variable. Then, we introduce a
new polynomial $g(d,x_1,x_2,x_n)$ such that for almost all non-zero values $d^o$ of $d$ the specialization $g(d^o,x_1,x_2,x_3)$ defines the
offset to $\Sigma$ at distance $d^0$. Such a polynomial is called the generic offset polynomial
(see Definition \ref{def:ch1:GenericOffsetEquation}, page \pageref{def:ch1:GenericOffsetEquation}), and the hypersurface that it defines
(in four-dimensional space) is called the generic offset of $\Sigma$ (see Def \ref{def:ch1:GenericOffset}, page \pageref{def:ch1:GenericOffset}).
In this situation, the goal of this paper is to describe an effective solution for the problem of computing the total degree
in $\{x_1,x_2,x_3\}$ of $g$.

\nin There exist some (few) contributions in the literature concerning the degree problem for offset curves and surfaces. To our knowledge,
the first attempt to provide a degree formula for offset curves was given by Salmon, in \cite{salmon-treatise}. This formula was proved
wrong in the already mentioned paper by Farouki and Neff \cite{Farouki1990}. In this paper, the authors provide a degree formula for rational
curves given parametrically. They also deal separately with the case of polynomial parametrizations. Our papers \cite{SS05} and \cite{SSS09}
provide a complete and efficient solution for the problem of computing the degree structure (that is total and partial degrees, and degree w.r.t. the distance) for the case of plane algebraic curves, both in the implicit and the parametric case. Finally, let us mention that Anton
et al. provide in \cite{anton2005offset} an alternative formula (to those presented in \cite{SS05}) for computing the total degree of
the offsets to an algebraic curve.

\nin In contrast with the case of curves, even in the case of a generating parametric surface, there are, up to our knowledge, no available
results for the offset degree problem in the scientific literature. In this paper, concretely in Theorem \ref{thm:ch4:DegreeFormula} (page
\pageref{thm:ch4:DegreeFormula}) we provide a formula for the total degree of the generic offset to a rational surface, given
in parametric form, provided that the Assumption \ref{rem:ch4:NotInfinitelyManyOffsetsThroughOrigin}  holds (see page \pageref{rem:ch4:NotInfinitelyManyOffsetsThroughOrigin}). {\sf The parametrization of the surface is not assumed to be proper}, and the formula in fact provides the product
of the total offset degree times the tracing index of the parametrization. However, since there are available efficient algorithms for
computing the tracing index of a surface parametrization (see \cite{SendraPerez2004JPAA-DegreeSurfaceParam}) this assumption does not limit
the applicability of the formula.

\nin The strategy for this offset degree problem is, as in our previous papers, based in the analysis of the intersection between the generic offset and a pencil of lines through the origin. The restriction to the rational case, combined with this strategy, results in a reduction in the dimension of the space needed to study of the intersection problem. Thus, we  are led to consider again an intersection problem of plane curves. The auxiliary curves involved in this case are obtained by eliminating the variables corresponding to a point in the generating surface from the offset-line intersection system. The main technical differences between this paper and our previous ones (see \cite{SS05} and \cite{SSS09}) are that:
\begin{itemize}
 \item Here we need to consider more than two intersection curves. Thus, the total degree formula is expressed as a generalized resultant of the equations of these auxiliary curves.
 \item Furthermore, all the curves involved in the intersection problem depend on parameters. Thus, the notion of fake point and their characterization is technically more demanding.
\end{itemize}
Generally speaking, the dimensional advantage gained by working with a parametric representation is partially compensated by the fact that we
are not dealing directly with the points of the surface but with their parametric representation, and thus we are losing some geometric
intuition. In the general situation of an implicitly given generating surface, if one were to apply a similar strategy to the offset
degree problem, we believe that one is bound to consider a surface intersection problem, instead of the simpler curve intersection problem
used here. However, in this paper we do not address the offset degree problem in that general situation.

\nin As those skilled in the art know, going from the curve to the surface case usually implies a huge step in the difficulty of the proofs.
For us, this has indeed been the case. As a result, some of the proofs in this paper are rather technical. And in one particular case,
we have not been able to extend to the surface case the proof of a result that we obtained for plane curves. Specifically, in
Lemma 4  of our paper \cite{SS05} we proved that there are only finitely many distance values $d^o$ for which the origin belongs to
$\cO_{d^o}(\cC)$. Our conjecture is that a similar property holds for all algebraic surfaces. However, as we said, we have not been able
to provide a proof. The following example shows that, even for curves, this property does not hold if we consider the analytic case

\vspace{3mm}
\begin{ExampleIntro}
We want to emphasize that this proposition does not hold in a non-algebraic context. For example, the ``offset'' to the analytic curve with
implicit equation $y^3-\sin(x)=0$, passes through the origin for infinitely many values of $d$. In fact, for this curve, all the offsets
with values of $d$ equal to $k\pi$, for $k\in\Z$, pass through the origin. This is illustrated in Figure
\ref{Figure:SmoothCurveInfiniteOffsetssPassingThroughOrigin}, where the curve (in red) and some offsets passing through the origin (for $k=1,2,3$) are shown.
\begin{figure}[h]
\begin{center}
\includegraphics[width=10cm]{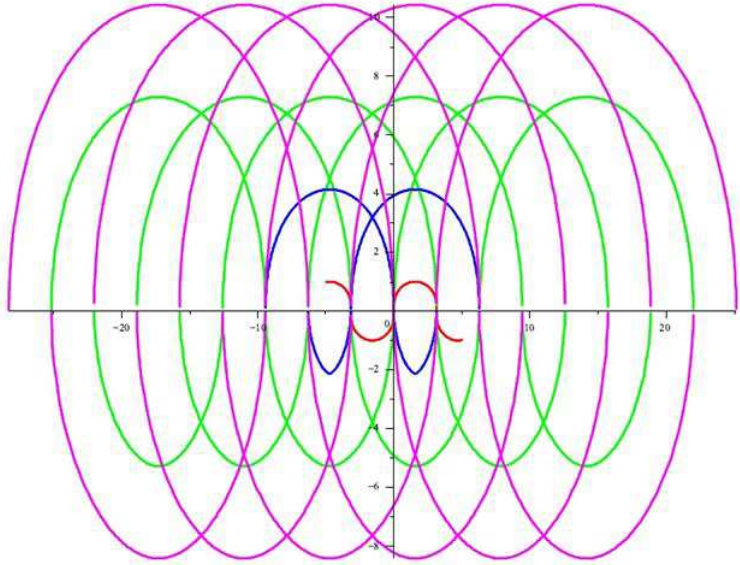}
\end{center}
\caption{A smooth curve with infinitely many offsets through the origin}\label{Figure:SmoothCurveInfiniteOffsetssPassingThroughOrigin}
\end{figure}
\end{ExampleIntro}
\begin{center}
\rule{2cm}{0.5pt}
\quad\\
\end{center}

\nin Besides, because of its own nature, our strategy fails in the case of some simple surfaces.
 We have met similar situations in our previous papers (\cite{SS05} and \cite{SSS09}), where we needed to exclude circles centered at the origin and lines through the origin from our considerations. Correspondingly,  in this paper we need to exclude the case in which the generating surface is a sphere centered at the origin. In this case, however, the generic offset degree (in fact the generic offset equation) is known beforehand. Therefore, excluding it does not really affect the generality of the degree formula that we present here. The above observations are the reason for the {\sf following assumptions:}

\vspace{3mm}
\noindent{\bf Assumptions \ref{rem:ch4:NotInfinitelyManyOffsetsThroughOrigin}} (page \pageref{rem:ch4:NotInfinitelyManyOffsetsThroughOrigin})
Let $\Sigma$ denote the generating surface. In this paper, we assume that:
\begin{enumerate}
  \item[(1)] {\sf There exists a finite subset $\Delta^1$ of $\C$ such that, for $d^o\not\in\Delta^1$ the origin does not belong to $\cO_{d^o}({\Sigma})$.}
 \item[(2)] {\sf $\Sigma$ is not a sphere centered at the origin.}
\end{enumerate}

\begin{center}
\subsubsection*{Structure of the Paper}\label{subsubsec:ch1:StructurePaper}
\end{center}

\subsubsection*{}
\vspace{-18mm} In Subsection \ref{subsec:ch1:FormalDefinitionPropertiesGenericOffset} of {\sf Section \ref{sec:GenericOffsetDegreeProblem}} (page \pageref{sec:GenericOffsetDegreeProblem}) we begin by reviewing the fundamental concepts related to the classical offset construction, and its basic properties that will be used in the sequel. Then  we  introduce the notion of generic offset, which can be considered as a hypersurface that collects as level curves all the classical offsets to a given hypersurface. This notion is the natural generalization of the classical concept, by considering the distance as a new variable. The defining polynomial of this object is the {generic offset polynomial}; see Definition \ref{def:ch1:GenericOffsetEquation}  (page \pageref{def:ch1:GenericOffsetEquation}).  We establish the fundamental specialization property of the generic offset polynomial in Theorem \ref{thm:Ch1:GenericOffsetSpecialization} (page \pageref{thm:Ch1:GenericOffsetSpecialization}). In Subsection \ref{subsec:ch4:SurfaceParametrizationsAndtheirAssociatedNormalVector}, we recall some basic notions on parametric algebraic surfaces, and some technical lemmas about them. We introduce the notion of associated normal vector, and we also construct a parametric analogous of the Generic Offset System. In Subsection \ref{sec:ch1:IntersectionCurvesResultants} (page \pageref{sec:ch1:IntersectionCurvesResultants}) we present two technical lemmas about the use of univariate resultants to study the problem of the intersection of plane algebraic curves.

\noindent {\sf Section \ref{sec:ch4:OffsetLineIntersectionforRationalSurfaces}} (page \pageref{sec:ch4:OffsetLineIntersectionforRationalSurfaces}) describes the theoretical foundation of the strategy. Subsection \ref{subsec:ch4:IntersectionWithLines} contains the analysis of the intersection between the generic offset and a pencil of lines through the origin. In Subsection \ref{sec:ch4:AuxiliaryCurvesForRationalSurfaces} we will see that, when elimination techniques are brought into our strategy, the dimension of the space in which we count the points in ${\mathcal O}_d(\Sigma)\cap{\mathcal L}_{\bar k}$ is reduced, and we arrive at an intersection problem between projective plane curves. Then we begin the analysis of that problem. Specifically, in Subsection \ref{subsec:ch4:EliminationAndAuxiliaryPolynomials} we describe the auxiliary polynomials obtained by using elimination techniques in the Parametric Offset-Line System, and we introduce the Auxiliary System \ref{sys:ch4:AuxiliaryCurvesSystem} (page \pageref{sys:ch4:AuxiliaryCurvesSystem}), denoted by ${{\mathfrak S}^P_3}(d,\bar k)$. Some geometric properties of the solutions of ${{\mathfrak S}^P_3}(d,\bar k)$ (see  Proposition \ref{prop:ch4:ExtendableSolutions}, page \pageref{prop:ch4:ExtendableSolutions} and Lemma \ref{lem:ch4:SignOfLambdaAndOffsetting}, page \pageref{lem:ch4:SignOfLambdaAndOffsetting}) will be used in the sequel to study the relation between the solution sets of Systems $\mathfrak S^P_2(d,\bar k)$ and ${{\mathfrak S}^P_3}(d,\bar k)$.  In Subsection \ref{subsec:ch4:FakePoints} (page \pageref{subsec:ch4:FakePoints}) we define the corresponding notion of fake points and invariant points for the Affine Auxiliary System ${{\mathfrak S}^P_3}(d,\bar k)$. The relation between these two notions is then shown in Proposition \ref{prop:ch4:FakePointsAndInvariantSolutionsCoincide} (page \pageref{prop:ch4:FakePointsAndInvariantSolutionsCoincide}).

The statement and proof of the degree formula appear in {\sf Section \ref{sec:ch4:TotalDegreeFormulaForRationalSurfaces}} (page \pageref{sec:ch4:TotalDegreeFormulaForRationalSurfaces}). This section is structured into four subsections as follows. In Subsection \ref{subsec:ProjectivizationParametrizationSurface} we study the projective version of the auxiliary curves introduced in the preceding section, and we introduce the Projective Auxiliary System \ref{sys:ch4:AuxiliaryCurvesSystem-ProjectiveAndPrimitive} (page \pageref{sys:ch4:AuxiliaryCurvesSystem-ProjectiveAndPrimitive}). The polynomials that define this system are the basic ingredients of the degree formula. Subsection \ref{subsec:ch4:InvariantSolutionsOfS5}. (page \pageref{subsec:ch4:InvariantSolutionsOfS5}) deals with the invariant solutions of the Projective Auxiliary System. In Subsection \ref{subsec:ch4:MultiplicityOfIntersectionAtNon-FakePoints} (page \pageref{subsec:ch4:MultiplicityOfIntersectionAtNon-FakePoints}) we will prove that the value of the multiplicity of intersection of the auxiliary curves at their non-invariant points of intersection equals one (in Proposition \ref{prop:ch4:MultiplicityAtNonFakePoints}, page \pageref{prop:ch4:MultiplicityAtNonFakePoints}). Subsection \ref{subsec:ch4:DegreeFormula} (page \pageref{subsec:ch4:DegreeFormula}) contains the statement and proof of the degree formula, in Theorem \ref{thm:ch4:DegreeFormula} (page \pageref{thm:ch4:DegreeFormula}).

\begin{center}
\subsubsection*{Notation and Terminology}\label{subsubsec:ch1:NotationTerminology}
\end{center}
\nin For ease of reference, we introduce and collect here the main {\sf notation and terminology} that will be used throughout this paper.
\begin{itemize}
 \item As usual, $\mathbb C$ and $\mathbb R$ correspond to the fields of complex and real numbers, respectively. The $n$-dimensional affine space
 is the set ${\C}^3$, and the associated projective space will be denoted by $\mathbb P^3$.
 \item We will use  $(y_1,y_2,y_3)$ for the affine coordinates in ${\C}^3$, and $(y_0:y_1:y_2:y_3)$ for the projective coordinates in $\mathbb P^3$,
  as well as  the abbreviations:
    \[\bar y=(y_1,y_2,y_3),\quad \bar y_h=(y_0:y_1:y_2:y_3).\]
 \item  In order to distinguish offset surfaces from  their generating surfaces, we will also use $\bar x=(x_1,x_2,x_3),\quad
 \bar x_h=(x_0:x_1:x_2:x_3)$ to refer to the affine and projective coordinates of a point in the offset,
  and $\bar y, \bar y_h$ as above for the original surface.
 \item A point in ${\C}^3$ will be denoted by
 \[\bar y^o=(y_1^o,y_2^o,y_3^o)\]
 and, correspondingly, a point in $\mathbb P^3$ will be denoted by
 \[\bar y^o_h=(y_0^o:y_2^o:y_3^o)\]
 Throughout this work, we will frequently use this $^o$ superscript to indicate a particular value of a variable.
 \item The Zariski closure of a  set $A\subset {\C}^n$ will be denoted by $A^*$. The projective closure of an algebraic set $A$
 will be denoted by  $\overline A$.
 \item Let $A$ be an algebraic set. We denote by $\Sing_a(A)$ the affine singular locus of $A$, and by $\Sing(A)$ the projective singular locus of $A$; i.e. the singular locus of $\overline{A}$.
 \item If $I$ is a polynomial ideal, $\mathbf V(I)$ denotes the affine algebraic set defined by $I$; that is,
 \[\mathbf V(I)=\{\bar x^o\in{\C}^n/\forall f\in I, f(\bar x^o)=0\}\]
 \item When we homogenize a polynomial $g\in{\C}[y_1,y_2,y_3]$,  we will use capital letters, as in $G(y_0,y_1,y_2,y_3)$, to denote the
  homogenization of $g$ w.r.t. $y_0$. Also, by abuse of notation, we will write $g(\bar y)$, $G(\bar y_h), {\C}[\bar y], {\C}[\bar y_h]$.
 \item The partial derivatives w.r.t. $y_i$ of $g(y_1,y_2,y_3)$ and of its homogenization $G(y_0,y_1,y_2,y_3)$ will be denoted
 $g_i$ and $G_i$ respectively, for $i=0,\ldots,3$. The symbol $\nabla g$ (resp. $\nabla G$) denotes the gradient vector of partial derivatives,
 i.e.:
 \[\nabla g(\bar y)=\left(g_1,g_2,g_3\right)(\bar y),\quad(\mbox{ resp. }\nabla G(\bar y_h)=\left(G_1,\ldots,G_3\right)(\bar y_h) )\]
 \item The symbol $\Sigma$ denotes a {\sf rational algebraic surface} defined over $\C$ by the irreducible polynomial $f(\bar y)\in\C[\bar y]$.
 \item We assume that we are given a {\sf non-necessarily proper rational parametrization} of $\Sigma$:
 \[
 {P}(\bar t)=\left(
 \dfrac{P_{1}(\bar t)}{P_{0}(\bar t)},
 \dfrac{P_{2}(\bar t)}{P_{0}(\bar t)},
 \dfrac{P_{3}(\bar t)}{P_{0}(\bar t)}
 \right).
 \]
 Here $\bar t=(t_1,t_2)$, and $P_0,\ldots,P_3\in\C[\bar t]$ with $\gcd(P_0,\ldots,P_3)=1$.

\item The {\sf projectivization} $P_h$ of ${P}$ is obtained by homogenizing the components of $P$ w.r.t. a new variable $t_0$, multiplying both the numerators and denominators if necessary by a suitable power of $t_0$.  It will be denoted by
\[
P_h(\bar t_h)=
\left(
\dfrac{X(\bar t_h)}{W(\bar t_h)},
\dfrac{Y(\bar t_h)}{W(\bar t_h)},
\dfrac{Z(\bar t_h)}{W(\bar t_h)}
\right)
\]
where $\bar t_h=(t_0:t_1:t_2)$, and $X, Y, Z, W\in\C[\bar t_h]$ are homogeneous polynomials of the same degree $d_{P}$, for which $\gcd(X,Y,Z,W)=1$ holds.
 \item The (classical) offset at distance $d^o\in {\C}^\times$ to ${\Sigma}$ is denoted by ${\cal O}_{d^o}({\Sigma})$ and the generic (classical) offset to ${\Sigma}$ by  ${\cal O}_{d}({\Sigma})$ (see Definition \ref{def:ch1:GenericOffset} in page \pageref{def:ch1:GenericOffset}). In this work, the
     variable $d$ always represents the distance values.
 \item We denote by $g\in{\C}[d, \bar x]$  the generic offset equation for ${\Sigma}$ (see Definition \ref{def:ch1:GenericOffsetEquation} in
 page \pageref{def:ch1:GenericOffsetEquation}).
\item $\delta$ is the total degree of $g$ w.r.t. $\bar x$; i.e. $\delta=\deg_{\bar x}(g)$.
 \item Given  $\phi(\bar y), \psi(\bar y)\in{\C}[\bar y]$ we denote by $\operatorname{Res}_{y_i}(\phi,\psi)$ the
     univariate resultant of $\phi$ and $\psi$ w.r.t. $y_i$, for $i=0,\ldots,n$. And if $A$ is a subset of the set of variables $\{y_0,\ldots,y_n\}$, we denote by $\operatorname{PP}_{A}(\phi)$ (resp. $\Con_{A}(\phi)$) the primitive part (resp. the content) of the polynomial $\phi$ w.r.t. $A$.
\item\label{def:ch0:NotationForGenericLine} $\cL_{\bar k}$ denotes a generic  line through the origin, whose direction is determined by the values of a variable $\bar k=(k_1,k_2,k_3)$.
 More precisely, for a particular value of $\bar k$, denoted by $\bar k^o$, the parametric equations of $\cL_{\bar k}$ are
\[
\ell_i(\bar k,l,\bar x):\,{ x_1}-{ k_1}\,l=0, \mbox{ for }i=1,2,3,\\
\]
where $l$ is the parameter.
 \item We will keep the convention of always using the letter $\Delta$ to indicate a finite subset of values of the variable $d$.
  Accordingly, the letter $\Theta$ denotes a Zariski closed set of values of $\bar k$.
 A Zariski open subset of ${\C}\times{\C}^{3}$, formed by pairs of values of $(d,\bar k)$, will be denoted by $\Omega$. In some proofs,
 an open set $\Omega$ will be constructed in several steps. In these cases we will use a superscript to indicate the step in the construction.
 Thus,  $\Omega_1^0, \Omega_1^1, \Omega_1^2$, etc. are open sets, defined in sucessive steps in the construction of $\Omega_1$.

\item A similar convention will be used for systems of equations and their solutions. A system of equations will be denoted by $\mfS$, with
 sub and superscripts to distinguish between systems, and the set of solutions of the system will be denoted by $\Psi$, with the same choice
of sub and superscripts.
 \item We will also need to consider local parametrizations of algebraic varieties. To distinguish local from (global) rational parametrizations, we will use calligraphic typeface for local  parametrizations. Thus, a local parametrization will be denoted by
\[\cP(\bar t)=\left(\cP_1(\bar t),\ldots,\cP_n(\bar t)\right)\]
\item Let $V$ be an $n$-dimensional vector space over ${\C}$. Two vectors $\bar v=(v_1,\ldots,v_n)$ and $\bar w=(w_1,\ldots,w_n)$ are said to be {\sf parallel} if and only if
\[v_iw_j-v_jw_i=0\mbox{, for }i,j=1,\ldots,n\]
In this case we write $\bar v\parallel\bar w$.
\item Given two vectors $(a_1,b_1,c_1), (a_2,b_2,c_2)\in{\C}^3$, their {\sf cross product} is defined as
$$(a_1,b_1,c_1)\wedge(a_2,b_2,c_2)=(b_1c_2-b_2c_1,a_2c_1-a_1c_2,b_1c_2-b_2c_1)$$
\item We consider in $V$  the symmetric bilinear form defined by:
\[\Xi(\bar v,\bar w)=\sum_{i=1}^nv_iw_i\]
which induces in $V$ a metric vector space  (see \cite{Roman2005}, \cite{Snapper1971}) with light cone of isotropy given by:
\[
L_{\Xi}=\{(x_1,\ldots,x_n)\in V/x_1^2+\cdots+x_n^2=0\}
\]
Note that, when ${\C}=\C$, this
is not the usual unitary space $\mathbb C^n$. On the other hand, when we consider the field $\R$, it is the usual Euclidean metric space, thus it preserves the usefulness of our results for applications.
In this work, the {\sf norm} $\|\bar v\|$ of a vector $\bar v\in V$ denotes a square root of $\bar v\cdot\bar v$, that is
\[\|\bar v\|^2=\Xi(\bar v,\bar v)=\sum_{i=1}^nv_i^2\]
Moreover, a vector $\bar v\in V$ is {\sf isotropic} if $\bar v\in L_{\Xi}$ (equivalently if $\|\bar v\|=0$). Note that for a non-isotropic vectors there are precisely two choices of norm, which differ only by multiplication by $-1$.
\item We denote
\[\normal_{(i,j)}(\bar x,\bar y)=f_i(\bar y)(x_j-y_j)-f_j(\bar y)(x_i-y_i).\]
Recall that $f(\bar y)$ is the defining polynomial of the irreducible hypersurface ${\Sigma}$, and that $f_i(\bar y)$ denotes its partial derivative w.r.t. $y_i$.
 \item \label{notation:ch1:Hodograph} The {\sf affine normal-hodograph} of $f$ is the polynomial:
 \begin{equation*}\label{def:ch1:ImplicitHodographNormalIsotropic}
  h(\bar y)=f_1^2(\bar y)+f_2^2(\bar y)+f_3^2(\bar y)
 \end{equation*}
 Following our convention of notation, $H$ is the homogenization of $h$ w.r.t. $y_0$; that is:
 \[H(\bar y_h)=F_1^2(\bar y_h)+F_2^2(\bar y_h)+F_3^2(\bar y_h)\]
 $H$ is called the {\sf projective normal-hodograph} of ${\Sigma}$. Moreover, a point $\bar y^o_h\in\overline{{\Sigma}}$ (resp. $\bar y^o\in{{\Sigma}}$) is called {\sf normal-isotropic} if $H(\bar y^o_h)=0$ (resp. $h(\bar y^o)=0$).
 \item We denote by ${\Sigma}_o$ the set of non normal-isotropic affine points of ${\Sigma}$; that is:
 \[{\Sigma}_o=\{\bar y^o\in{\Sigma}/h(\bar y^o)\neq 0\}\]
 In the rest of this work {\sf we will  assume} that the Zariski-open subset ${\Sigma}_o$ is non-empty. In \cite{Sendra1999}, Proposition 2, it is proved that this is equivalent to $H$ not being a multiple of $F$. We denote by $\Iso({\Sigma})$ the closed set of affine normal-isotropic points of ${\Sigma}$. Note that $\Sing_a({\Sigma})\subset\Iso({\Sigma})$.
\item If $K$ is an irreducible component of an algebraic set $A$, and $K\subset\Iso(A)$ we will say that $K$ is normal-isotropic.
\end{itemize}

\section{The Generic Offset}\label{sec:GenericOffsetDegreeProblem}

\begin{itemize}
\item In {\sf Subsection \ref{subsec:ch1:FormalDefinitionPropertiesGenericOffset}},
we first review  the concept of classical offset and the related properties that will be used in the sequel.
In order to do this, we follow the incidence diagram formalism in \cite{Arrondo1997} and \cite{Sendra2000} (see Diagram \ref{sys:ch1:IncidenceDiagramOffsetFixedDistance}  in page \pageref{sys:ch1:IncidenceDiagramOffsetFixedDistance};
the reader may also see our previous papers \cite{SS06} and \cite{SSS09})\footnote{The refered works deal
with general hypersurfaces. We consider here the special case of three dimensional surfaces.}.
Besides, the relation of this notion with Elimination Theory is established.
We also collect several fundamental properties of the classical offset construction that we
will refer to in the sequel. Next, the notions of generic offset and generic offset polynomial are introduced; see Definition
\ref{def:ch1:GenericOffsetEquation}  (page \pageref{def:ch1:GenericOffsetEquation}). The fundamental specialization property of the generic offset polynomial is established in Theorem \ref{thm:Ch1:GenericOffsetSpecialization} (page \pageref{thm:Ch1:GenericOffsetSpecialization}).

\item In {\sf Subsection \ref{subsec:ch4:SurfaceParametrizationsAndtheirAssociatedNormalVector}}, we recall some basic notions on parametric algebraic surfaces, and some technical lemmas about them. We also introduce the notion of associated normal vector, and we review some of its properties. Besides, we construct a parametric analogous of the Generic Offset System (see System \ref{sys:ch4:ParametricOffsetSystem}).

\item {\sf Subsection \ref{sec:ch1:IntersectionCurvesResultants}} (page \pageref{sec:ch1:IntersectionCurvesResultants}) contains some technical results about the use of univariate resultants to study the problem of the intersection of plane algebraic curves. The classical setting for the computation of the intersection points of two plane curves by means of resultants is well known (see for instance \cite{Brieskorn1986},  \cite{Walker1950} and \cite{Sendra2007}). This requires in general a linear change  of coordinates. However, in this work,  we need to analyze the behavior of the resultant when some of the standard requirements are not satisfied. This is the content of Lemma \ref{lem:ch1:MultiplicityUsingResultants} (page \pageref{lem:ch1:MultiplicityUsingResultants}). Similarly, we also need to analyze the case when more than two curves are involved, by using generalized resultants. This is done in Lemma \ref{lem:ch1:GeneralizedResultants} (page \pageref{lem:ch1:GeneralizedResultants}).

\end{itemize}

\subsection{Formal definition and basic properties of the generic offset}\label{subsec:ch1:FormalDefinitionPropertiesGenericOffset}

We begin by recalling the formal definition of classical offset, that can be found  in \cite{Arrondo1997} and \cite{Sendra1999}. With the notation introduced above, let $d^o\in{\C}^\times$ be a fixed value, and let
$\Psi_{d^o}(\Sigma)\subset{\C}^{7}$ be the set of solutions (in the variables $(\bar x,\bar y,u)$) of the following polynomial system:
\begin{equation}\label{sys:Ch1:OffsetSystemFixed_d}
\left.\begin{array}{lr}
&f(\bar y)=0\\
\normal_{(i,j)}(\bar x,\bar y): &f_i(\bar y)(x_j-y_j)-f_j(\bar y)(x_i-y_i)=0\\
(\mbox{for }i,j=1,\ldots,3; i<j)&\\
b_{d^o}(\bar x, \bar y ):& (x_1-y_1)^2+(x_2-y_2)^2+(x_3-y_3)^2-(d^o)^2=0\\
w(\bar y,u):& u\cdot(\|\nabla f(\bar y)\|^2)-1=0
\end{array}\right\}
\equiv\mathfrak{S}_1(d^o)
\end{equation}
\noindent Let us consider the following:
\begin{center}
\vspace{-5mm}
\begin{equation}\label{sys:ch1:IncidenceDiagramOffsetFixedDistance}
\begin{array}{c}
\mbox{\sf Offset Incidence Diagram}\\[1mm]
\xymatrix{
&\Psi_{d^o}({\Sigma})\subset {\C}^3\times{\C}^3\times{\C}
\ar[ld]^{\pi_1}\ar[rd]_{\pi_2}&\\
{\mathcal A}_{d^o}({\Sigma})\subset {\C}^3&&{\Sigma}\subset {\C}^3}
\end{array}
\end{equation}
\end{center}
\noindent where
\[
\begin{array}{ccc}
\begin{cases}
\pi_1:{\C}^{7}\mapsto{\C}^{3}\\
\pi_1(\bar x,\bar y,u)=\bar x
\end{cases}
&
\mbox{and}
&
\begin{cases}
\pi_2:{\C}^{7}\mapsto{\C}^{3}\\
\pi_2(\bar x,\bar y,u)=\bar y
\end{cases}
\end{array}
\]
and $\cA_{d^o}({\Sigma})=\pi_1(\Psi_{d^o}({{\Sigma}}))$.

\medskip
\vspace{3mm}\begin{Definition}\label{Def:ClassicalOffset}
The {\sf (classical) offset} to $\Sigma$ at distance $d^o$ is the algebraic set $\cA_{d^o}(\Sigma)^*$ (recall that the asterisk indicates
Zariski closure). It will be denoted by ${\cO}_{d^o}(\Sigma)$.
\end{Definition}
\vspace{3mm}
\begin{Remark}\label{rem:ch1:CommentariesToFormalDefinitionClassicalOffset}
\begin{enumerate}
 \item[]
 \item If there is a solution of the system \ref{sys:Ch1:OffsetSystemFixed_d} of the form $\bar p^o=(\bar x^o,\bar y^o,u^o)$, then we say that
the point $\bar y^o\in{\Sigma}$ and the point $\bar x^o\in\cO_{d^o}(\Sigma)$ are {\sf associated points}.
 \item Let $I(d^o)\subset{\C}[\bar x, \bar y, u]$ be the ideal generated by the polynomials in ${\mathfrak S}_1(d^o)$; that is:
 \[I(d^o)=<f(\bar y),b_{d^o}(\bar x, \bar y ),\normal_{(1,2)}(\bar x, \bar y ),\normal_{(1,3)}(\bar x, \bar y ),\normal_{(2,3)}(\bar x, \bar y ),w(\bar y,u)>.\]
 This means that
 \[\Psi_{d^o}({\Sigma})=\mathbf V(I(d^o))\]
 is the affine algebraic set defined by $I(d^o)$, and
 \[{\cO}_{d^o}({\Sigma})=\mathbf V(\tilde I(d^o))\]
 where $\tilde I(d^o)=I(d^o)\cap{{\C}}[\bar x]$ is the $(\bar y,u)$-elimination ideal
 of $I(d^o)$ (see \cite{Cox1997}, Closure Theorem, p. 122).
 In particular, this means that the offset can be computed by elimination techniques, such as Gr\"{o}bner bases, resultants, characteristic sets, etc.
\end{enumerate}
\end{Remark}
\nin Next, we will refer to some properties of the classical offset construction, that we collect here for the reader's convenience. We start with a very important geometric property regarding the normal vector of the classical offset construction.
\vspace{3mm}
\begin{Proposition}[\sf Fundamental Property of the Classical Offsets]\label{prop:ch1:FundamentalPropertyOffsets}
Let $\bar y^o\in\cV_o$, and let $\bar x^o\in\cO_{d^o}(\cV)$ be a point associated to $\bar y^o$. Then the normal line to $\cV$ at $\bar y^o$ is also normal to $\cO_{d^o}(\cV)$ at $\bar x^o$.
\end{Proposition}
\begin{proof}
See \cite{Sendra1999}, Theorem 16.
\end{proof}

\nin In order to prove that we can avoid degenerated situations, we will sometimes need information about the dimension of certain sets of points. The basic tools for doing this will be the incidence diagrams, analogous to \ref{sys:ch1:IncidenceDiagramOffsetFixedDistance} (page \pageref{sys:ch1:IncidenceDiagramOffsetFixedDistance}) and well known results about the dimension of the fiber of a regular map. For ease of reference, we include here a statement of one such result, in a form that meets our needs. The proof can be found in \cite{Harris1992}
\vspace{3mm}
\begin{Lemma}\label{lem:ch1:FiberDimension}
Let $A$ be an affine algebraic set, and let $f:A\mapsto\C^p$ be a regular map. Let us denote $B=f(A)^*$. For $\bar a^o\in A$, let $\mu(\bar a^o)=\dim(f^{-1}(f(a^o)))$
Then, if $A_o\subset A$ is any irreducible component, $B_o=f(A_o)$ its image, and $\mu^o$ is the minimum value of $\mu(\bar a^o)$ for $\bar a^o\in A_o$, we have
\[\dim(A_o)=\dim(B_o)+\mu^o\]
In particular, if there exists $a^o\in A_o$ for which $\dim(f^{-1}(f(a^o)))=0$, then $\dim(A_o)=\dim(B_o)$.
\end{Lemma}

\nin We next analyze the number and dimension of the irreducible components of the offset.
\vspace{3mm}
\begin{Proposition}\label{prop:ch1:OffsetHasAtMostTwoComponents}
$\cO_{d^o}(\cV)$ has at most two irreducible components.
\end{Proposition}
\begin{proof}
See \cite{Sendra1999}, Theorem 1.
\end{proof}
\vspace{3mm}
\begin{Proposition}\label{prop:ch1:ComponentsPsiHaveSameDimensionAsV}
The irreducible components of $\Psi_{d^o}(\cV)$ have the same dimension as $\cV$.
\end{Proposition}
\begin{proof}
In \cite{Sendra1999}, Lemma 1, this is proved using local parametrizations. However, since, as we have seen above, $\pi_2$ is a $2:1$ map, this can also be considered a straightforward application of the preceding Lemma \ref{lem:ch1:FiberDimension}.
\end{proof}
\vspace{3mm}
\begin{Remark}
This implies immediately that $\cO_{d^o}(\cV)$ has at most two irreducible components, whose dimension is less or equal than $\dim(\cV)$.
\end{Remark}

\noindent To present the next two results, we recall some of the terminology introduced in \cite{Sendra1999}:
\vspace{3mm}\begin{Definition}\label{def:ch1:DegenerationAndSimpleAndSpecialComponents}
\begin{enumerate}
 \item[]
 \item The offset $\cO_{d^o}(\cV)$ is called {\sf degenerated} if at least one of its components is not a hypersurface.
 \item A component $\cM\subset\cO_{d^o}(\cV)$ is said to be a {\sf simple component} if there exists a non-empty Zariski dense subset $\cM_1\subset\cM$ such that every point of $\cM_1$ is associated to exactly one point of $\cV$. Otherwise, $\cM$ is called a {\sf special component} of the offset. Furthermore, we say that $\cO_{d^o}(\cV)$ is {\sf simple} if all its components are simple, and {\sf special} if it has at least a special component (in this case, it has precisely one special component, see Proposition \ref{prop:ch1:OffsetPropertiesRegardingSimpleAndSpecialComponents} below).
\end{enumerate}
\end{Definition}
\nin Note that if the offset is degenerated, and taking into account Lemma \ref{lem:ch1:FiberDimension}, the map $\pi_1$ must have a non-zero dimensional fiber for some point in $\Psi_{d^o}(\cV)$.

\noindent The following two results tell us that degeneration and special components are very infrequent phenomena.
\vspace{3mm}
\begin{Proposition}\label{prop:ch1:OffsetDoesNotDegenerateForAlmostAllDistances}
There is a finite set $\Delta_0\subset\K$ such that if $d^o\not\in\Delta_0$, then $\cO_{d^o}(\cV)$ is not degenerated.
\end{Proposition}
\begin{proof}
See \cite{Sendra1999}, Theorem 2.
\end{proof}
\quad\\
\vspace{3mm}
\begin{Proposition}\label{prop:ch1:OffsetPropertiesRegardingSimpleAndSpecialComponents}
\begin{enumerate}
 \item[]
 \item Let $\cM$ be an irreducible and non-degenerated component of $\cO_{d^o}(\cV)$. Then $\cM$ is special if and only if $\cO_{d^o}(\cM)=\cV$.
 \item  $\cO_{d^o}(\cV)$ has at least a simple component.
 \item If $\cO_{d^o}(\cV)$ is irreducible, then it is simple.
 \item There is a finite set $\Delta_1\subset\C$ such that, if $d^o\not\in\Delta_1$ then $\cO_{d^o}(\cV)$ is simple, and  the irreducible components of $\cO_{d^o}(\cV)$ are not contained in $\Iso(\cV)$.
\end{enumerate}
\end{Proposition}
\begin{proof}
See \cite{Sendra1999}, Theorems 7, 8 and Corollary 6.
\end{proof}
\nin The next result shows that --as expected, being a metric construction-- the offset construction is invariant under rigid motions of the affine space.
\vspace{3mm}
\begin{Proposition}\label{prop:ch1:OffsetInvarianceRigidMovements}
Let $\cT$ be a rigid motion of the affine space $\C^3$. Then
\[\cT(\cO_{d^o}(\cV))=\cO_{d^o}(\cT(\cV))\]
\end{Proposition}
\begin{proof}
See \cite{Sendra1999}, Lemma 2.5 in Chapter 2.
\end{proof}

\begin{center}
\subsubsection*{Formal definition of the generic offset}
\end{center}

The concept of generic offset to an algebraic surface was formally introduced in our previous paper \cite{SSS09}.
The motivation for this concept is the following: As the distance value $d^o$ varies, different
offset varieties are obtained. The idea is to have a global expression of the offset for all (or almost all)
distance values. This motivates the concept of {\sf generic polynomial of the offset to ${\Sigma}$}. This is a
polynomial, depending on the distance variable $d$, such that for every (or almost every, see the examples below)
non-zero value $d^o$, the polynomial specializes to the defining polynomial of the offset at that particular distance.
Let us see a couple of examples that give some insight into the situation.
\vspace{3mm}
\begin{Example}
\begin{itemize}
\item[]
\item[(a)] Using this informal definition of generic offset polynomial, and using Gr\"{o}bner basis techniques, one can see that if $\cC$ is the parabola of equation $y_2-y_1^2=0$, the generic polynomial of its offset is:\\
 \vspace{1mm}
 \noindent
 $g(d,x_1,x_2)=-48\,{d}^{2}{x_1}^{4}-32\,{d}^{2}{x_1}^{2}{x_2}^{2}+48
 \,{d}^{4}{x_1}^{2}+16\,{x_1}^{6}+16\,{x_2}^{2}{x_1 }^{4}+16\,{d}^{4}{x_2}^{2}
 -16\,{d}^{6}-40\,x_2\,{x_1}^{ 4}-32\,{x_1}^{2}{x_2}^{3}+8\,{d}^{2}x_2\,{x_1}^{2}
 -32\,{d}^{2}{x_2}^{3}+32\,{d}^{4}x_2+{x_1}^{4}+32\,{{ \it
 x1}}^{2}{x_2}^{2}+16\,{x_2}^{4}-20\,{d}^{2}{x_1}^{2
 }-8\,{d}^{2}{x_2}^{2}-8\,{d}^{4}-2\,x_2\,{x_1}^{2}-8\,{
 x_2}^{3}+8\,x_2\,{d}^{2}+{x_2}^{2}-{d}^{2}.$\\
 \noindent In addition, and using again Gr\"{o}bner basis techniques, one may check that for every distance the generic offset polynomial specializes properly (see Example \ref{exm:ch1:StandardParabolaGenericOffset} in page \pageref{exm:ch1:StandardParabolaGenericOffset} below, for a detailed description of this example and the preceding claims).
 \item[(b)] On the other hand, the generic offset polynomial of  the circle of equation $y_1^2+y_2^2-1=0$ factors as the product of two circles of radius $1+d$ and $1-d$; that is:
 \[g(d,x_1,x_2)= \left( x_1^2+x_2^2-(1+d)^2\right) \left(  x_1^2+x_2^2-(1-d)^2\right).
 \]
 Now, observe that for $d^o=1$, this generic polynomial gives
 \[g(1,x_1,x_2)= \left( x_1^2+x_2^2-2^2\right) \left( x_1^2+x_2^2\right)=
 \left( x_1^2+x_2^2-2^2\right)\left( x_1+i x_2\right)\left( x_1-i x_2\right)
 \]
 which describes the union of a circle of radius $2$, and two complex lines. This is not a correct representation of the offset at distance $1$ to $\cC$, which consists of the union of the circle of radius $2$ and a point (the origin). In fact, using Gr\"obner basis techniques, one has that the
elimination ideal $\tilde I(1)$ (see Remark \ref{rem:ch1:CommentariesToFormalDefinitionClassicalOffset}(2), page \pageref{rem:ch1:CommentariesToFormalDefinitionClassicalOffset}) is:
\[\tilde I(1)=<x_2(x_1^2+x_2^2-4), x_1(x_1^2+x_2^2-4)>.\]
\nin Thus, in this example we see that the generic offset polynomial does not specialize properly for $d^o=1$. Nevertheless, for every other value of $d^o$ the specialization is correct.
\end{itemize}
\end{Example}
\begin{center}
\rule{2cm}{0.5pt}
\quad\\
\end{center}
\nin After these examples, we describe the formal definition of generic offset and generic offset polynomial, adapted to the case of surfaces in three
dimensional space. We start by considering the following system of equations:
\begin{equation}\label{sys:Ch1:GenericOffsetSystem}
\left.\begin{array}{lr}
&f(\bar y)=0\\
\normal_{(i,j)}(\bar x,\bar y): &f_i(\bar y)(x_j-y_j)-f_j(\bar y)(x_i-y_i)=0\\
(\mbox{for }i,j=1,\ldots,3; i<j)&\\
b(d,\bar x, \bar y ):& (x_1-y_1)^2+(x_2-y_2)^2+(x_3-y_3)^2-d^2=0\\
w(\bar y,u):& u\cdot(\|\nabla f(\bar y)\|^2)-1=0
\end{array}\right\}
\equiv\mathfrak{S}_1(d)
\end{equation}
The above system will be called the {\sf Generic Offset System}. In this system, we consider $d$ as a variable,
so that $b\in{\C}[d, \bar x, \bar y]$. A solution of this system is thus a point of the form
$(d^o, \bar x^o, \bar y^o, u^o)\in{\C}^8$.

\noindent Let $\Psi({\Sigma})\subset{\C}\times{\C}^3\times{\C}^3\times{\C}$ be the set of solutions of $\mathfrak{S}_1(d)$, and
consider the following:\\
\begin{center}
\vspace{-10mm}
\begin{equation}\label{sys:ch1:IncidenceDiagramGenericOffset}
\begin{array}{c}
\mbox{\sf Generic Offset Incidence Diagram}\\[3mm]
\xymatrix{
&\Psi({\Sigma})\subset {\C}\times{\C}^3\times{\C}^3\times{\C}
\ar[ld]^{\pi_1}\ar[rd]_{\pi_2}&\\
{\mathcal A}({\Sigma})\subset {\C}^{4}&&{\C}\times{\Sigma}\subset {\C}^{4}
}
\end{array}
\end{equation}
\end{center}
\noindent where
\[
\begin{array}{ccc}
\begin{cases}
\pi_1:{\C}^{8}\mapsto{\C}^{4}\\
\pi_1(d, \bar x, \bar y, u)=(d, \bar x)
\end{cases}
&
\mbox{and}
&
\begin{cases}
\pi_2:{\C}^{8}\mapsto{\C}^{4}\\
\pi_2(d, \bar x, \bar y, u)=(d, \bar y)
\end{cases}
\end{array}
\]
and $\cA({\Sigma})=\pi_1(\Psi({{\Sigma}}))$.

\nin Then one has the following definition (recall that  the asterisk denotes the Zariski closure of a set):
\vspace{3mm}\begin{Definition}\label{def:ch1:GenericOffset}
The {\sf generic offset} to ${\Sigma}$ is
\[{\cO}_{d}({\Sigma})=\cA({\Sigma})^*=\pi_1(\Psi({{\Sigma}}))^*\subset{\C}^{4}\]
\end{Definition}

\vspace{3mm}
\begin{Remark}\label{rem:ch1:EliminationIdealGenericOffset}
\begin{enumerate}
 \item[]
 \item Let
 \[I(d)=<f(\bar y),b(d, \bar x, \bar y),\normal_{(1,2)}(\bar x, \bar y ),\normal_{(1,3)}(\bar x, \bar y ),\normal_{(2,3)}(\bar x, \bar y ),w(\bar y,u)>\]
 be the ideal in ${\C}[d, \bar x, \bar y, u]$ generated by the polynomials in System \ref{sys:Ch1:GenericOffsetSystem}. Note that the above definition implies that
 \[{\cO}_{d}({\Sigma})=\mathbf V(\tilde I(d))\]
 where $\tilde I(d)=I(d)\cap{\C}[d, \bar x]$ is the $(\bar y,u)$-elimination ideal of $I(d)$.
\item The Closure Theorem from Elimination Theory (see e.g. Theorem 3 in page 122 of \cite{Cox1997}) implies that the dimension of the set
\[\cO_d({\Sigma})\setminus\pi_1(\Psi_1({\Sigma}))\]
is smaller than the dimension of $\cO_d({\Sigma})$. This is the set of points of the generic offset associated with singular or normal-isotropic points of ${\Sigma}$.
\end{enumerate}

\end{Remark}

\begin{center}
\subsubsection*{Basic properties of the generic offset}\label{subsec:ch1:BasicPropertiesGenericOffset}
\end{center}

\nin In the following Proposition we will see that the properties of the offset at a fixed distance,
regarding its dimension and number of components (see Lemma 1, Theorem 1 and Theorem 2 in \cite{Sendra1999}),
are reflected in the generic offset. In particular, this Proposition shows that the generic offset is a
hypersurface, and thus guarantees the existence of the generic polynomial
(see below, Definition \ref{def:ch1:GenericOffsetEquation}).
\vspace{3mm}
\begin{Proposition}\label{prop:ch1:GenericOffsetHypersurface}
\begin{enumerate}
 \item[]
 \item $\cO_d({\Sigma})$ has at most two components.
 \item Each component of $\cO_d({\Sigma})$ is a hypersurface in ${\C}^{4}$.
\end{enumerate}
\end{Proposition}
\begin{proof} (Adapted from Lemma 1, Theorem 1 and Theorem 2 in \cite{Sendra1999}).
We begin by showing that if $K$ is a component of $\Psi_1({\Sigma})$, then $\dim(K)=3$. Thus
\begin{equation}\label{eq:ch1:DimensionPsiVEquals_n}
\dim(\Psi_1({\Sigma}))=3
\end{equation}
Let $\psi^o=(d^o, \bar x^o, \bar y^o, u^o)\in K$. Then, $\bar y^o\in{\Sigma}$ is a regular point of ${\Sigma}$. Let $\cP(\bar t)$,
with $\bar t=(t_1,\ldots,t_{n-1})$, be a local parametrization of ${\Sigma}$ at $\bar y^o$, with $\cP(\bar t^o)=\bar y^0$. Then, it holds
that one of the local parametrizations defined by:
\[
\cP^\pm(d, \bar t)=\left(d,\, \cP(\bar t)\pm d\dfrac{\nabla f(\cP(\bar t))}{\|\nabla f(\cP(\bar t))\|},\, \cP(\bar t),\,\dfrac{1}{\|\nabla f(\cP(\bar t))\|^2}\right)
\]
parametrizes $\Psi_1({\Sigma})$ locally at $\psi^o$ (we choose sign so that $\cP^\pm(\bar t^o)=\psi^o$). Since $(d,\cP(\bar t))$ parametrizes
${\C}\times{\Sigma}$, we get that $(d, \bar t)$ are algebraically independent, and so $\dim(K)=3$.

\nin Now we can prove the first statement of the proposition. Since the number of components of
$\cO_{d}({\Sigma})$ is at most the number of components of $\Psi_1({\Sigma})$,
one just only has to prove that  $\Psi_1({\Sigma})$ has at most two components. Let us suppose that
$\Gamma_{1}, \Gamma_{2}$ y $\Gamma_{3}$  are three different components of $\Psi_1({\Sigma})$ and let
$Z=\pi_{2}(\Gamma_{1})\cap \pi_{2}(\Gamma_{2})\cap\pi_{2}(\Gamma_{3})$, where  $\pi_{2}$ is the  projection
of the incidence diagram \ref{sys:ch1:IncidenceDiagramGenericOffset}. Then, it holds that  $\dim(Z)=3$. Observe
that if $\dim(Z) <3$ then $\dim({\Sigma}\setminus Z)=\dim(\bigcup^{3}_{i=1}({\Sigma}\setminus \pi_{2}(\Gamma_{i})))=3$. Hence,
at least one of the sets ${\Sigma}\setminus \pi_{2}(\Gamma_{i})$ is of dimension $3$, which is impossible since  $\pi_{2}(\Gamma_{i})$
are constructible sets of dimension $3$. On the other hand, it holds that $\dim(Z\cap \pi_{2}(\Gamma_{i}\cap \Gamma_{j}))< 3$ for $i<j$. Then
\[Z\setminus \bigcup_{i \neq j}(Z\cap \pi_{2}(\Gamma_{i}\cap \Gamma_{j})) \neq \emptyset. \]
Now, take $\bar p=(d^o, \bar y^o)\in Z\setminus\bigcup_{i \neq j}(Z \cap \pi_{2}(\Gamma_{i}\cap
\Gamma_{j}))$, then  $\pi_{2}^{-1}(\bar p)=\{\bar q_{1}, \bar q_{2}, \bar q_{3} \}$
where $\bar q_{i} \neq \bar q_{j}$ for $i<j$, which is impossible since the mapping
$\pi_{2}$ is $(2:1)$ on  $\pi_2(\Psi_1({\Sigma}))$.

\nin Finally we can prove statement 2 in the proposition. We analyze the dimension of the tangent space to a component of the generic offset.
Let  $(d^o, \bar y^o)\in\pi_2(\Psi_1({\Sigma}))$, such that the two points $(d^o,\bar x^o_1),(d^o,\bar x^o_2)\in\cO_d({\Sigma})$ generated
by $(d^o, \bar y^o)$ satisfy that the dimension of their tangent spaces is the dimension of the  corresponding component of
$\cO_d({\Sigma})$. Let $u^o=\dfrac{1}{\|\nabla f(\bar y^o)\|^2}$,  and let $\cP(\bar t)$ be a local parametrization of  ${\Sigma}$ at
$\bar x^o$. Then, it holds that :
\[
{\tilde\cP}^+(d, \bar t)=
\left(d,\, \cP(\bar t)+d\dfrac{\nabla f(\cP(\bar t))}{\|\nabla f(\cP(\bar t))\|}\right)
 \quad\mbox{ and }\quad
{\tilde\cP}^-(d, \bar t)=
\left(d,\, \cP(\bar t)-d\dfrac{\nabla f(\cP(\bar t))}{\|\nabla f(\cP(\bar t))\|}\right)
\]
parametrize locally $\cO_d({\Sigma})$ at $(d^o,\bar x^o_1)$, and $(d^o,\bar x^o_2)$. In this situation, let $Q^\pm$ be as above, and
consider the following map:
\[ \begin{array}{cccccccc}
\psi^{+}: & {{\C}}^{3} & \longrightarrow &\Psi_1({\Sigma}) & \stackrel{\varphi^{+}}{\longrightarrow}
 & \cA({\Sigma}) &\stackrel{i}{\hookrightarrow} & {{\C}}^{4} \\
 & (d,\bar{t}) & \longrightarrow & \cP^{+}(d,\bar{t}\,) &
\longrightarrow &\tilde\cP^{+}(d, \bar t)& \longrightarrow & \tilde\cP^{+}(d, \bar t).
\end{array} \]
Similarly, we define $\psi^{-}$ and $\varphi^{-}$. Now consider the following homomorphism, defined by the differential $d\psi^{+}$
(similarly for $d\psi^{-}$), between the tangent space to  $\Psi_1({\Sigma})_{1}$ at $(d^o,\bar x^o_{1},\bar y^o,u^o)$ and  the
tangent space  ${\cT}_{(d^o,\bar x^o_1)}$ to $\cA({\Sigma})_{1}$ at $(d^o,\bar x^o_1)$, where $\Psi_1({\Sigma})_{1}$ and $\cA({\Sigma})_{1}$
denote the  component of $\Psi_1({\Sigma})$ and $\cA({\Sigma})$ containing the points $(d^o,\bar x^o_{1},\bar y^o,u^o)$ and
$(d^o,\bar x^o_{1})$, respectively. Then one has that
\[\dim(\cA({\Sigma})_{1})\geq \dim({\cT}_{(d^o,\bar x^o_1)})\geq \dim(\operatorname{Im}(d\varphi^{+}))
=\operatorname{rank}({\cal J}_{\varphi^{+}}),\]
where ${\cal J}_{\varphi^{+}}$ denotes the jacobian matrix of
$\varphi^{+}$. Furthermore, by Equation \ref{eq:ch1:DimensionPsiVEquals_n} at the beginning of this proof, one has  that
\[ 3=\dim(\Psi_1({\Sigma})_{1})\geq \dim(\cA({\Sigma})_{1})\geq \operatorname{rank}({\cal J}_{\varphi^{+}}).\]
On the other hand, if we take any point of the form $(0,\bar t^o)\in{{\C}}^{n}$, that is, with $d^o=0$, we must
get $\operatorname{rank}({\cal J}_{\varphi^{+}})=3$ at that point; otherwise, one would conclude that the rank of the  jacobian of
${\cal P}(\,\bar{t}\,)$ is smaller than   $2$, which is impossible since  ${\Sigma}$ is a surface.
\end{proof}

\noindent As a first consequence of this Proposition, ${\cO}_{d}({\Sigma})$ is defined by a polynomial $g(d, \bar x)\in{\C}[d, \bar x]$ (see \cite{Shafarevich1994}, p.69, Theorem 3). Thus, we arrive at the following definition:
\vspace{3mm}
\begin{Definition}\label{def:ch1:GenericOffsetEquation}
The {\sf generic offset polynomial} is the defining polynomial of the hypersurface ${\cO}_{d}({\Sigma})$. In the sequel, {\sf we denote} by  $g(d, \bar x)$ the generic
offset equation.
\end{Definition}

\nin The first property of the generic offset polynomial that we study regards its factorization:
\vspace{3mm}
\begin{Lemma}\label{lem:ch1:GenericOffsetIsPrimitiveIn_x}
The generic offset polynomial is primitive w.r.t. $\bar x$
\end{Lemma}
\begin{proof}
Suppose, on the contrary, that $g(d,\bar x)$ has a non-constant factor in ${\C}[d]$. That is
\[g(d,\bar x)=A(d)\tilde g(d,\bar x)\]
Let $d^o\neq 0$ be any root of $A(d)$. Then the hypersurface $\cZ$ in ${\C}\times{\C}^3$ defined by $d=d^o$ is contained in
$\cO_d({\Sigma})$. Taking Remark \ref{rem:ch1:EliminationIdealGenericOffset}(2) (page \pageref{rem:ch1:EliminationIdealGenericOffset}) into
account, one has that there is an open non-empty subset of $\cZ$ contained in $\pi_1(\Psi_1({\Sigma}))$. This in turn implies that
there is an open subset $\tilde\cZ$ of ${\C}^3$ such that if $\bar x^o\in\tilde\cZ$, then $\bar x^o\in\cO_{d^o}(\Sigma)$
(the classical offset at distance $d^o$). This is a contradiction, since we know that $\cO_{d^o}(\Sigma)$ has dimension less or equal
to $2$. Thus, we are left with the case when $A(d)$ is a power of $d$. The argument must be different in this case,
since the classical offset is only defined for $d^o\in{\C}^\times$. However, the reasoning is similar: we conclude that there is an open
non-empty subset $\tilde\cZ_0$ of ${\C}^3$ such that if $\bar x^o\in\tilde\cZ_0$, then the system (``classical offset system for distance $0$'')
\[
\left.\begin{array}{lr}
f(\bar y)=0\\
f_i(\bar y)(x^o_j-y_j)-f_j(\bar y)(x^o_i-y_i)=0\\
(\mbox{for }i,j=1,\ldots,3; i<j)&\\
(x^o_1-y_1)^2+(x^o_2-y_2)^2+(x^o_3-y_3)^2=0\\
u\cdot(\|\nabla f(\bar y)\|^2)-1=0
\end{array}\right\}
\equiv\mathfrak{S}_1(d)
\]
has solutions. Now, if $(\bar y^o,u^o)$ is a solution of this system then
\begin{enumerate}
 \item $\nabla f (\bar y^o)$ is not isotropic,
 \item $\bar y^o-\bar x^o$  is isotropic,
 \item and $\nabla f (\bar y^o)$ is  parallel to $\bar y^o-\bar x^o$,
\end{enumerate}
Thus one has that  $\bar y^o-\bar x^o=0$. But since $\bar x^o$ runs through an open subset of ${\C}^3$,
this contradicts the fact that ${\Sigma}$ is a surface.
\end{proof}
\vspace{3mm}
\begin{Remark}\label{rem:ch1:GenericOffsetEqSqfreeAndHasAtMostTwoFactors}
\begin{enumerate}
 \item[]
 \item Observe that the polynomial $g$ may be reducible (recall the example of the circle) but by construction it is always square-free. Moreover, by Proposition \ref{prop:ch1:GenericOffsetHypersurface} (page \pageref{prop:ch1:GenericOffsetHypersurface}) and Lemma \ref{lem:ch1:GenericOffsetIsPrimitiveIn_x}, $g$ is either irreducible or factors into two irreducible factors not depending only on $d$.
 \item We will also call $g(d,\bar x)=0$ the {\sf generic offset equation} of ${\Sigma}$.
\end{enumerate}
\end{Remark}
\noindent The following theorem gives the fundamental property of the generic offset.
\vspace{3mm}
\begin{Theorem}\label{thm:Ch1:GenericOffsetSpecialization}
For all but finitely many exceptions, the generic offset polynomial specializes properly.
That is, there exists a finite (possibly empty) set $\Delta_2\subset{\C}$ such
that if $d^o\not\in\Delta_2$, then
\[g(d^o, \bar x)=0\]
is the equation of ${\cO}_{d^o}({\Sigma})$.
\end{Theorem}
\begin{proof}
 Let $G(d)$ be a reduced Gr\"{o}bner basis of $I(d)$ w.r.t. an elimination ordering that eliminates $(\bar y,u)$. Then, up to multiplication by a non-zero constant, $G(d)\cap{\C}[d, \bar x]$ is a Gr\"{o}bner basis of $\tilde I(d)$. Proposition \ref{prop:ch1:GenericOffsetHypersurface} above shows that $\tilde G(d)=G(d)\cap{\C}[d, \bar x]=<\nu(d)g(d,\bar x)>$, where $\nu(d)$ is a non-zero polynomial, depending only on $d$ (see the Remark preceding this proof).  But then (see \cite{Cox1997}, exercise 7, page 283) there is a finite (possibly empty) set $\Delta_2^1\subset{\C}$ such that for $d^o\not\in\Delta_2^1$, $G(d)$ specializes well to a Gr\"{o}bner basis of $I(d^o)$ (defined in Remark \ref{rem:ch1:CommentariesToFormalDefinitionClassicalOffset}, page \pageref{rem:ch1:CommentariesToFormalDefinitionClassicalOffset}). It follows that, since $\tilde I(d^o)=I(d^o)\cap{\C}[\bar x]$, then $\tilde G(d^o)=\{\nu(d^o)g(d^o, \bar x)\}$ is a Gr\"{o}bner basis of $\tilde I(d^o)$. In particular, if $\Delta_2^2$ is the finite set of zeros of $\nu(d)$, then for $d^o\not\in\Delta_2=\Delta_2^1\cup\Delta_2^2$, and $d^o\neq 0$, one has that $g(d^o, \bar x)$ is the equation for $\cO_{d^o}(\Sigma)$.
\end{proof}
\nin For future reference, we collect in the following corollary all the information about the --finite-- set of {\em bad} distances that appear in the offsetting construction.
\vspace{3mm}
\begin{Corollary}\label{cor:ch1:BadDistancesFiniteSet}
There is a finite set $\Delta\subset{\C}^{\times}$ such that for $d^o\not\in\Delta$, the following hold:
\begin{enumerate}
 \item[(1)] (non degeneracy): ${\cO}_{d^o}({\Sigma})$ is not degenerated.
 \item[(2)] (simplicity): ${\cO}_{d^o}({\Sigma})$ is simple.
 \item[(3)] (good specialization): if $g(d, \bar x)=0$ is the generic offset polynomial, $g(d^o, \bar x)=0$ is the equation of ${\cO}_{d^o}({\Sigma})$.
 \item[(4)] (degree invariance):
 \[\deg_{\bar x}(\cO_{d}({\Sigma}))=\deg_{\bar x}(\cO_{d^o}(\Sigma)),\quad \deg_{x_i}(\cO_{d}({\Sigma}))=\deg_{x_i}(\cO_{d^o}(\Sigma))\mbox{ for }i=1,\ldots,3.\]
\end{enumerate}
\end{Corollary}
\begin{proof}
Take $\Delta^1=\Delta_0\cup\Delta_1\cup\Delta_2$, with $\Delta_0$ as in Proposition \ref{prop:ch1:OffsetDoesNotDegenerateForAlmostAllDistances},
$\Delta_1$ as in Proposition \ref{prop:ch1:OffsetPropertiesRegardingSimpleAndSpecialComponents}(4) and $\Delta_2$ as in Theorem
\ref{thm:Ch1:GenericOffsetSpecialization} above. Furthermore, let $p(d)\bar x^{\mu}$ be a term of $g(d,\bar x)$ of maximal degree
w.r.t. $\bar x$. That is, $\bar\mu=(\mu_1,\mu_2,\mu_3)\in\N^3$, with $\sum \mu_i=\deg_{\bar x}(g)$, where $p(d)\in\C[d]$ is a non-zero polynomial. Then take:
\[\Delta^{\bar x}=\Delta\cup\{d^o\in\C\,|\,p(d^o)=0\},\]
and similarly, for $i=1,2,3$ construct $\Delta^{\bar x_i}$, by considering a term of $g(d,\bar x)$ of maximal degree w.r.t $x_i$. Finally, taking
\[\Delta=\Delta^1\cup\Delta^{\bar x}\cup\Delta^{\bar x_1}\cup\Delta^{\bar x_2}\cup\Delta^{\bar x_3},\]
our claim holds.
\end{proof}
\noindent Let us see a first example of a generic offset polynomial.
\vspace{3mm}
\begin{Example}\label{exm:ch1:StandardParabolaGenericOffset}
For the parabola $\cC$ with defining polynomial $f(y_1,y_2)=y_2-y_1^2$, the generic offset system turns into:
\[
\left.\begin{array}{lr}
&f(\bar y)=y_2-y_1^2\\
\normal_{(1,2)}(\bar x,\bar y): &-2y_1(x_2-y_2)-(x_1-y_1)=0\\
b(d, \bar x, \bar y ):& (x_1-y_1)^2+(x_2-y_2)^2-d^2=0\\
w(\bar y,u):& u\cdot(\|4y_1^2+1\|^2)-1=0
\end{array}\right\}
\]
Computing a Gr\"{o}bner elimination basis of $I(d)=<f,\normal_{(1,2)},b,w>$, we obtain (with the notation in the proof of Theorem \ref{thm:Ch1:GenericOffsetSpecialization}):
\[G(d)=\{g(d, \bar x),\chi_1(d, \bar x),\ldots,\chi_8(d, \bar x)\}\]
where:\\
$g(d, \bar x)=
16\,{x_{{1}}}^{6}+16\,{x_{{1}}}^{4}{x_{{2}}}^{2}-40\,{x_{{1}}}^{4}x_{{
2}}-32\,{x_{{1}}}^{2}{x_{{2}}}^{3}+ \left( -48\,{d}^{2}+1 \right) {x_{
{1}}}^{4}+ \left( -32\,{d}^{2}+32 \right) {x_{{1}}}^{2}{x_{{2}}}^{2}+
16\,{x_{{2}}}^{4}+ \left( 8\,{d}^{2}-2 \right) {x_{{1}}}^{2}x_{{2}}+
 \left( -32\,{d}^{2}-8 \right) {x_{{2}}}^{3}+ \left( 48\,{d}^{4}-20\,{
d}^{2} \right) {x_{{1}}}^{2}+ \left( 16\,{d}^{4}-8\,{d}^{2}+1 \right)
{x_{{2}}}^{2}+ \left( 32\,{d}^{4}+8\,{d}^{2} \right) x_{{2}}-16\,{d}^{
6}-8\,{d}^{4}-{d}^{2}
$\\
and\\
$\chi_1(d, \bar x)=12\,{d}^{2}u{x_{{1}}}^{2}+16\,{d}^{2}u{x_{{2}}}^{2}-4\,{d}^{2}ux_{{2}}
+ \left( -12\,{d}^{4}+{d}^{2} \right) u-4\,{x_{{1}}}^{4}+8\,{x_{{1}}}^
{2}x_{{2}}+8\,{d}^{2}{x_{{1}}}^{2}-4\,{x_{{2}}}^{2}+4\,{d}^{2}x_{{2}}-
4\,{d}^{4}+3\,{d}^{2}
$\\
$\chi_2(d, \bar x)=64\,{d}^{2}u{x_{{2}}}^{3}-48\,{d}^{2}u{x_{{1}}}^{2}-16\,{d}^{2}u{x_{{2
}}}^{2}+28\,{d}^{2}ux_{{2}}+ \left( -60\,{d}^{4}-3\,{d}^{2} \right) u-
64\,{x_{{1}}}^{4}x_{{2}}+128\,{x_{{1}}}^{2}{x_{{2}}}^{2}+128\,{d}^{2}{
x_{{1}}}^{2}x_{{2}}-64\,{x_{{2}}}^{3}+36\,{d}^{2}{x_{{1}}}^{2}+112\,{d
}^{2}{x_{{2}}}^{2}+ \left( -64\,{d}^{4}+36\,{d}^{2} \right) x_{{2}}-36
\,{d}^{4}+3\,{d}^{2}
$\\
$\chi_3(d, \bar x)=12\,y_{{2}}-16\,u{x_{{1}}}^{2}-16\,u{x_{{2}}}^{2}+8\,ux_{{2}}+ \left(
16\,{d}^{2}-1 \right) u-8\,x_{{2}}+1
$\\
$\chi_4(d, \bar x)=12\,y_{{1}}x_{{2}}+ \left( -12\,{d}^{2}-3 \right) y_{{1}}-8\,y_{{2}}x_
{{1}}x_{{2}}-14\,y_{{2}}x_{{1}}+8\,{x_{{1}}}^{3}+8\,x_{{1}}{x_{{2}}}^{
2}-6\,x_{{1}}x_{{2}}+ \left( -8\,{d}^{2}+3 \right) x_{{1}}
$\\
$\chi_5(d, \bar x)=3\,y_{{1}}x_{{1}}+2\,y_{{2}}x_{{2}}-y_{{2}}-2\,{x_{{1}}}^{2}-2\,{x_{{2
}}}^{2}+2\,{d}^{2}
$\\
$\chi_6(d, \bar x)=12\,{d}^{2}{u}^{2}-4\,u{x_{{1}}}^{2}-16\,u{x_{{2}}}^{2}-4\,ux_{{2}}+
 \left( 4\,{d}^{2}-1 \right) u+4\,x_{{2}}+1
$\\
$\chi_7(d, \bar x)=12\,{d}^{2}y_{{1}}u+8\,y_{{2}}ux_{{1}}x_{{2}}+14\,y_{{2}}ux_{{1}}-3\,y
_{{1}}-8\,u{x_{{1}}}^{3}-8\,ux_{{1}}{x_{{2}}}^{2}+6\,ux_{{1}}x_{{2}}+
 \left( 8\,{d}^{2}+3 \right) ux_{{1}}
$\\
$\chi_8(d, \bar x)={y_{{1}}}^{2}+{y_{{2}}}^{2}-2\,y_{{1}}x_{{1}}-2\,y_{{2}}x_{{2}}+{x_{{1
}}}^{2}+{x_{{2}}}^{2}-{d}^{2}.$\\
In particular,
\[G(d)\cap{\C}[d, \bar x]=<g(d, \bar x)>.\]
And so $g(d, \bar x)$ is the generic offset polynomial for the parabola $\cC$.\\
This Gr\"{o}bner basis has been computed considering the generators of $I(d)$ as polynomials in ${\C}(d)[\bar x,\bar y,u]$. This means that we have relationships of the form:
\[g(d, \bar x)=
a_1(d)f(\bar x)+a_2(d)\normal_{(1,2)}(\bar x,\bar y)+a_3(d)b(d, \bar x, \bar y)+a_4(d)w(\bar y,u)\]
and for $i=1,\ldots,8$:
\[
\chi_i(d, \bar x)=
b_{i1}(d)f(\bar x)+b_{i2}(d)\normal_{(1,2)}(\bar x,\bar y)
+b_{i3}(d)b(d, \bar x, \bar y)+b_{i4}(d)w(\bar y,u)
\]
where $a_1,\ldots,a_4,b_{11},\ldots,b_{84}\in{\C}(d)$. The result in Exercise 7, in page 283 of \cite{Cox1997} indicates that the Gr\"{o}bner basis specializes well for all values $d^o$ such that none of the denominators of the rational functions $a_i$ and $b_{ij}$ vanish at $d^o$. In this particular example, one may compute these rational functions and check that they are all constant. Therefore, specializing $g(d, \bar x)$ provides the offset equation for every non-zero value of $d$.
The computations in this example were obtained with the computer algebra system Singular (see \cite{GPS}). We do not include here the details of the computations, because of obvious space limitations.
\end{Example}
\begin{center}
\rule{2cm}{0.5pt}
\quad\\
\end{center}
\nin The following result, about the dependence on $d$ of the generic offset polynomial, is an easy consequence of Theorem \ref{thm:Ch1:GenericOffsetSpecialization} above. In fact, Theorem \ref{thm:Ch1:GenericOffsetSpecialization} implies that there are infinitely many values $d^o$ such that $g(d^o, \bar x)$ is the polynomial of $\cO_{d^o}(\Sigma)$ and, simultaneously, $g(-d^o,\bar x)$ is the polynomial of $\cO_{-d^o}(\Sigma)$. But, because of the symmetry in the construction, the offsets $\cO_{d^o}(\Sigma)$ and $\cO_{-d^o}({\Sigma})$ are exactly the same algebraic set. Thus, it follows that for infinitely many values of $d^o$ it
holds that up to multiplication by a non-zero constant:
\[g(d^o, \bar x)=g(-d^o,\bar x).\]
Hence, we have proved the following proposition:
\vspace{3mm}
\begin{Proposition}\label{prop:ch1:OffsetEquationIsQuadraticInd}
The generic offset polynomial belongs to ${\C}[\bar x][d^2]$. That is, it only contains even powers of $d$.
\end{Proposition}

\nin Now we can describe precisely the central problem of this work.
\vspace{3mm}\begin{Remark}
The {\sf total degree problem} for $\Sigma$ consists of finding (efficient) formulae to compute the total degree of $g$ in the variables $\bar x$.  We denote
this total degree by $\delta$.
\end{Remark}

\subsection{Surface Parametrizations and Parametric Offset System}
\label{subsec:ch4:SurfaceParametrizationsAndtheirAssociatedNormalVector}

Since $\Sigma$ is an algebraic surface over ${\C}$, all of its irreducible components have dimension $2$ over ${\C}$. Besides, assuming that the
surface $\Sigma$ is {\sf unirational (or parametric)} means that there exists a rational map $P:{\C}^2\mapsto\Sigma$ such that the image of
${P}$ is dense in $\Sigma$ w.r.t. the Zariski topology. The map ${P}$ is called an {\sf (affine) parametrization} of $\Sigma$. If ${P}$
is a birational map, then $\Sigma$  is called a {\sf rational surface},  and ${P}$ is called a {\sf proper parametrization} of $\Sigma$.
{\sf In this paper  we will not assume that ${P}$ is proper} (see Lemma \ref{lem:ch4:PropertiesSurfaceParametrization} in
page \pageref{lem:ch4:PropertiesSurfaceParametrization}, and the observations preceding it). It is well known that a rational surface is
always irreducible.

\nin Thus, a parametrization ${P}$ of $\Sigma$ is given through a non-constant triplet of rational functions in two parameters. We will use $\bar t=(t_1,t_2)$ for the parameters of ${P}$ and, as usual, $\bar t^o=(t_1^o,t_2^o)$ stands for a particular value in ${\C}^2$ of the pair of parameters. By a simple algebraic manipulation, we can assume that the three components of ${P}$ have a common denominator. Thus, we can write:
\begin{equation}\label{eq:ch4:SurfaceAffineParametrization}
{P}(\bar t)=\left(
\dfrac{P_{1}(\bar t)}{P_{0}(\bar t)},
\dfrac{P_{2}(\bar t)}{P_{0}(\bar t)},
\dfrac{P_{3}(\bar t)}{P_{0}(\bar t)}
\right)
\end{equation}
where $P_0,\ldots,P_3\in\C[\bar t]$ and $\gcd(P_0,\ldots,P_3)=1$. The number
\[d_{{P}}=\max_{i=0,\ldots,3}\left(\{
\deg_{\bar t}(P_{i})\}\right)\]
is then called the {\sf degree} of ${P}$.

\nin Over the algebraically closed field ${\C}$, the notions of rational and parametric surface
are equivalent (see the Castelnuovo Theorem \cite{Castelnuovo1939}). Furthermore, there exists an algorithm by
Schicho (see \cite{schicho1998rps}) to obtain a proper parametrization of a rational surface given by its implicit equation.
Thus, in principle, given a non-proper parametrization of a surface, it is possible (though computationally very expensive)
to implicitize, and then apply Schicho's algorithm to obtain a proper parametrization. In addition,
\cite{perezdiaz2006ppr} shows how to properly reparametrize certain special families of rational surfaces. However,
{\sf in this paper  we will not assume that ${P}$ is proper}, and the degree formulas below take this fact into account.

\noindent The parametrization $P$ has two {\sf associated tangent vectors}, denoted by
\begin{equation}\label{def:ch4:SurfaceAffineParametrizationAssociatedTangentVectors}
\dfrac{\partial {P}(\bar t)}{\partial t_1}\mbox{ and }\dfrac{\partial {P}(\bar t)}{\partial t_2}.
\end{equation}
That is:
 \[
 \dfrac{\partial {P}}{\partial t_i}=
\left(\frac{P_{1,i}P_{0}-P_{1}P_{0,i}}{(P_{0})^2},
\frac{P_{2,i}P_{0}-P_{2}P_{0,i}}{(P_{0})^2},
\frac{P_{3,i}P_{0}-P_{3}P_{0,i}}{(P_{0})^2}\right)
 \]
 where $P_{j,i}$ denotes the partial derivative of $P_{j}$ w.r.t. $t_i$, for $j=0,\ldots,2$ and $i=1,2$.

\nin The following Lemma states those properties of the surface parametrization $P$ that we will need in the sequel.
\vspace{3mm}
\begin{Lemma}\label{lem:ch4:PropertiesSurfaceParametrization}
There are non-empty Zariski open subsets $\Upsilon_1\subset\C^2$ and $\Upsilon_2\subset\Sigma$ such that:
\[{P}:\Upsilon_1\mapsto\Upsilon_2\]
is a surjective regular application of degree $m$. In particular, this means that ${P}$ defines a $m:1$ correspondence between $\Upsilon_1$ and $\Upsilon_2$. Thus, given $\bar y^o\in\Upsilon_2$, there are precisely $m$ different values $\bar t_1^o,\ldots,\bar t_m^o$ of the parameter $\bar t$ such that $P(\bar t_i^o)=\bar y^o$ for $i=1,\ldots,m$. Furthermore, if $\bar t^o\in\Upsilon_1$, the rank of the Jacobian matrix $\left(\dfrac{\partial{P}}{\partial\bar t}\right)$ evaluated at $\bar t^o$  is two.
\end{Lemma}
\begin{proof}
See e.g. \cite{Perez-Diaz2002}.
\end{proof}

\begin{Remark}\label{rem:ch4:TracingIndex_m}
\begin{enumerate}
 \item[]
 \item The number $m$ is also called, as in the case of curves, the {\sf tracing index} of ${P}$. See \cite{Sendra2007} for an algorithm to compute $m$. In the sequel, we will {\sf denote} by $m$ the tracing index of $P$.
 \item As a consequence of this lemma, the part of the surface $\Sigma$ not covered by the image of $P$ is a proper closed subset (i.e. a finite collection of curves and points).
\end{enumerate}
\end{Remark}

\nin Starting with the parametrization $P$ of $\Sigma$ as in (\ref{eq:ch4:SurfaceAffineParametrization}) above, we will construct a polynomial normal vector to $\Sigma$, that will be used in the statements of the degree formulas for rational surfaces. This particular choice of normal vector will be called in the sequel the {\sf associated normal vector of $P$}\index{normal vector $\bar n$ to a rational surface, associated to $P$}, and it will be {\sf denoted} by ${\bar n}(\bar t)$.

\noindent To construct ${\bar n}(\bar t)$, we first take the cross product of the associated tangent vectors introduced in \ref{def:ch4:SurfaceAffineParametrizationAssociatedTangentVectors}, page \pageref{def:ch4:SurfaceAffineParametrizationAssociatedTangentVectors}. Let us denote:
\[V(\bar t)=\dfrac{\partial {P}(\bar t)}{\partial t_1}\wedge\dfrac{\partial {P}(\bar t)}{\partial t_2}\]
This vector $V(\bar t)$ has the following form:
\[
V(\bar t)=\left( \dfrac{A_1(\bar t)}{A_0(\bar t)}, \dfrac{A_2(\bar t)}{A_0(\bar t)},
\dfrac{A_3(\bar t)}{A_0(\bar t)}\right)
\]
where $A_i\in{\C}[\bar t]$. Let $G(\bar t)=\gcd(A_1,A_2,A_3)$.

\vspace{3mm}
\begin{Definition}\label{def:ch4:AffineAssociatedNormalVector}
With the above notation, the {\sf associated normal vector} ${\bar n}=(n_1,n_2,n_3)$ to $P$ is the vector whose components are the polynomials:
\[n_i(\bar t)=\dfrac{A_i(\bar t)}{G(\bar t)}\mbox{ for }i=1,2,3.\]
\end{Definition}

\begin{Remark}\label{rem:ch4:RelationBetweenImplicitAndParametricNormalVector}
\begin{enumerate}
\item[]
\item Note that ${\bar n}$ is a normal vector to $\Sigma$ at $P(\bar t)$, vanishing at most at a finite set of points in the $\bar t$ plane. To see this observe that, because of their construction,  $n_1, n_2, n_3$ have no common factors. Besides, at most one of the polynomials $n_i$ is constant (otherwise the surface is a plane). Thus, the non constant components of ${\bar n}$ define a system of at least two plane curves without common components.
\item In particular,  there are some $\mu\in\N$ and $\beta(\bar t)\in\C[\bar t]$, with $\gcd(\beta,P_0)=1$, such that
\begin{equation}\label{eq:ch4:RelationBetweenImplicitAndParametricNormalVectors}
f_i(P(\bar t))=\dfrac{\beta(\bar t)}{P_0(\bar t)^{\mu}}n_i(\bar t)\mbox{ for }i=1,2,3.
\end{equation}
That is:
\[\nabla f(P(\bar t))=\dfrac{\beta(\bar t)}{P_0(\bar t)^{\mu}}\cdot\bar n(\bar  t)\]
\item Note that the polynomial $\beta(\bar t)$ introduced above is not identically zero. Otherwise, one has $f_i(P(\bar t))=0$ for $i=1,2,3$, and this implies that $f(\bar y)$ is a constant polynomial, which is a contradiction.
\end{enumerate}
\end{Remark}

\vspace{3mm}
\begin{Definition}
The polynomial $h\in\C[\bar t]$ defined as
\[h(\bar t)=n_1(\bar t)^2+n_2(\bar t)^2+n_3(\bar t)^2\]
is called the {\sf parametric (affine) normal-hodograph\index{hodograph, affine parametric}\index{normal-hodograph, affine parametric}} of the parametrization $P$.
\end{Definition}
\vspace{3mm}
\begin{Remark}
In the sequel, if we need to refer to the implicit normal-hodograph introduced in page \pageref{notation:ch1:Hodograph}, we will denote it by $H_{\operatorname{imp}}$ in the projective case, resp. $h_{\operatorname{imp}}$ in the affine case.
\end{Remark}
\nin The following lemma will be used below to exclude from our discussion certain pathological cases, associated to some particular parameter values.
\vspace{3mm}
\begin{Lemma}\label{lem:ch4:ExcludeBadParameterValues}
The sets $\Upsilon_1$ and $\Upsilon_2$ in Lemma \ref{lem:ch4:PropertiesSurfaceParametrization} (page \pageref{lem:ch4:PropertiesSurfaceParametrization}) can be chosen so that if $\bar t^o\in\Upsilon_1$, then
\[P_0(\bar t^o)h(\bar t^o)\beta(\bar t^o)\neq 0.\]
In particular, ${\bar n}(\bar t^o)\neq0$.
\end{Lemma}
\begin{proof}
Note that $P_0$, $h$ and $\beta$ are non-zero polynomials. Thus, the equation:
\[P_0(\bar t)h(\bar t)\beta(\bar t)=0\]
defines an algebraic curve. Let us call it $\cC$. Then it suffices to replace $\Upsilon_1$ (resp. $\Upsilon_2$) in Lemma \ref{lem:ch4:PropertiesSurfaceParametrization} with  $\Upsilon_1\setminus\cC$ (resp. $\Upsilon_2\setminus P(\cC)$).
\end{proof}

\begin{center}
\subsubsection*{Parametric system for the generic offset}
\label{subsec:ch4:ParametricSystemforGenericOffset}
\end{center}

\noindent Let $\Sigma$ and ${P}$ be as above. In order to describe $\cO_d(\Sigma)$ from a parametric point of view, we introduce the following system, to be called the {\sf parametric system for the generic offset}:

\vspace{-7mm}
\begin{equation}\label{sys:ch4:ParametricOffsetSystem}
\begin{minipage}{14cm}
\[\hspace{3mm}\mathfrak S_1^{{P}}(d)\equiv\begin{cases}
b^{P}(d,\bar t,\bar x):
\left (P_0{x_1}-P_1\right)^{_2}+\left (P_0{ x_2}-P_2\right )^{_2}+\left ( P_0{ x_3}-P_3\right)^{_2}-{d}^{_2}{P_0}^{_2}=0\\
\normal^{P}_{(1,2)}(\bar t,\bar x):\,n_1\cdot (P_0x_2-P_2)-n_2\cdot (P_0x_1-P_1)=0\\
\normal^{P}_{(1,3)}(\bar t,\bar x):\,n_1\cdot (P_0x_3-P_3)-n_3\cdot (P_0x_1-P_1)=0\\
\normal^{P}_{(2,3)}(\bar t,\bar x):\,n_2\cdot (P_0x_3-P_3)-n_3\cdot (P_0x_2-P_2)=0\\
w^{P}(r,\bar t):\,r\cdot P_0\cdot h\cdot \beta-1=0\\
\end{cases}\]
\end{minipage}
\end{equation}
Our first result will show that this system provides an alternative description for the generic offset. To state this, we will introduce some additional notation. Let
\[\Psi^{P}\subset\C\times\C\times\C^2\times\C^3\]
be the set of solutions, in the variables $(d,r, \bar t,\bar x)$, of the system  $\mathfrak S_1^{{P}}(d)$. We also consider the projection maps
\[
\begin{array}{ccc}
\begin{cases}
\pi_1^{{P}}:\C\times\C\times\C^2\times\C^3\mapsto\C\times\C^{3}\\
\pi_1^{{P}}(d,r, \bar t,\bar x)=(d,\bar x)
\end{cases}
&
\mbox{and}
&
\begin{cases}
\pi_2^{{P}}:\C\times\C\times\C^2\times\C^3\mapsto\C\times\C^{2}\\
\pi_2^{{P}}(d,r, \bar t,\bar x)=(r,\bar t)
\end{cases}
\end{array}
\]
and we define $\cA^{{P}}=\pi_1^P(\Psi^{P})$. Recall that $\left(\cA^{{P}}\right)^*$ denotes the Zariski closure of $\cA^{{P}}$.

\vspace{3mm}
\begin{Proposition}\label{prop:ch4:ParametricDescriptionOffset}
\[
{\cO}_{d}(\Sigma)=\left(\cA^{{P}}\right)^*.
\]
\end{Proposition}
\begin{proof}
With the notation introduced in Definition \ref{def:ch1:GenericOffset}, page \pageref{def:ch1:GenericOffset}, recall that
\[{\cO}_{d}(\Sigma)=\cA(\Sigma)^*=\pi_1(\Psi({\Sigma}))^*.\]
Note that in this proof we use $\pi_1, \pi_2$ as in page \pageref{def:ch1:GenericOffset}, to be distinguished from $\pi^P, \pi_2^P$ introduced above. Let $\Upsilon_1, \Upsilon_2$ be as in Lemma \ref{lem:ch4:ExcludeBadParameterValues}, page \pageref{lem:ch4:ExcludeBadParameterValues}, and let us denote  $$\cB^P_{\Sigma}=\pi_2^{-1}(\C\times\Upsilon_2).$$
$\cB^P_{\Sigma}$ is a non-empty dense subset of $\Psi({\Sigma})$, because $\C\times\Upsilon_2$  is dense in $\C\times\Sigma$. It follows that ${\cO}_{d}(\Sigma)=\pi_1(\cB^P_{\Sigma})^*$. We will show that $\pi_1(\cB^P_{\Sigma})=\cA^P$, thus completing the proof.

\nin If $(d^o, \bar x^o)\in\pi_1(\cB^P_{\Sigma})$, there are $\bar y^o, u^o$ and $\bar t^o\in\Upsilon_1$ such that $(d^o, \bar x^o, \bar y^o, u^o)\in\Psi({\Sigma})$, with $\bar y^o=P(\bar t^o)$. Since $u^o\neq 0$ and also $P_0(\bar t^o)h(\bar t^o)\beta(\bar t^o)\neq 0$, we can define:
\[
r^o=\dfrac{u^o\,\beta(\bar t^o)}{P_0(\bar t^o)^{2\mu+1}}.
\]
where $\mu$ is as in Equation \ref{eq:ch4:RelationBetweenImplicitAndParametricNormalVectors}, page \pageref{eq:ch4:RelationBetweenImplicitAndParametricNormalVectors}. Then, substituting $P(\bar t^o)$  by $\bar y^o$ in System \ref{sys:ch4:ParametricOffsetSystem}, and using also Equation \ref{eq:ch4:RelationBetweenImplicitAndParametricNormalVectors}, one has that:
\begin{equation}\label{sys:EquivalenceBetweenImplicitAndParametricOffsetSystems}
\begin{cases}
b^P(d^o,\bar t^o,\bar x^o)=P_0(\bar t^o)^2 b(d^o,\bar x^o,\bar y^o)=0\\[3mm]
\normal^P_{(i,j)}(\bar t^o,\bar x^o)=\dfrac{P_0(\bar t^o)^{\mu+1}}{\beta(\bar t^o)}\normal_{(i,j)}(\bar x^o,\bar y^o)=0\\[3mm]
w^P(r^o,\bar t^o)=w(u^o,\bar y^o)=0,&\
\end{cases}
\end{equation}
because $(d^o, \bar x^o, \bar y^o, u^o)\in\Psi({\Sigma})$. Therefore, one concludes that $(d^o,r^o,\bar t^o,\bar x^o)\in\Psi^{P}$, and so $(d^o, \bar x^o)\in\cA^P$. This proves that $\pi_1(\cB^P_{\Sigma})\subset\cA^P$.

\nin Conversely, let $(d^o, \bar x^o)\in\cA^P$. Then, there are $\bar t^o, r^o$ such that $(d^o,r^o,\bar t^o,\bar x^o)\in\Psi^{P}$. Since $P_0(\bar t^o)h(\bar t^o)\beta(\bar t^o)\neq 0$,
\[\bar y^o=P(\bar t^o)\quad\mbox{ and }\quad u^o=\dfrac{r^o\,P_0(\bar t^o)^{2\mu+1}}{\beta(\bar t^o)}\]
are well defined. The equations (\ref{sys:EquivalenceBetweenImplicitAndParametricOffsetSystems}) still hold, and in this case, they imply that $(d^o, \bar x^o, \bar y^o, u^o)\in\Psi({\Sigma})$. Besides, $\pi_2(d^o, \bar x^o,\bar y^o, u^o)=(d^o,\bar y^o)\in\C\times\Upsilon_2$, and so $(d^o, \bar x^o)\in\pi_1(\cB^P_{\Sigma})$. This proves that $\cA^P\subset\pi_1(\cB^P_{\Sigma})$, thus finishing the proof.
\end{proof}

\subsection{Intersection of Curves and Resultants}
\label{sec:ch1:IntersectionCurvesResultants}

In the following sections we will show that we can translate the offset degree problem into a suitably constructed planar curves intersection
problem. In this subsection we gather some results about the planar curves intersection problem to be used in the sequel.

\nin It is well known that the intersection points of two plane curves, without common components, as
well as their multiplicity of intersection, can be computed by means of resultants. For this, a
suitable preparatory change of coordinates may be required (see for instance, \cite{Brieskorn1986},  \cite{Walker1950} and,
for a modern treatment of the subject, \cite{Sendra2007}). In this work, for reasons that will turn out to be clear in
subsequent sections, we need to analyze the behavior of the resultant factors, and their correspondence with multiplicities of intersection,
when some of the standard requirements are not satisfied. Similarly, we also need to analyze the case when more than two curves are involved.

More precisely, we will use two technical lemmas. The first one, whose proof can be found in \cite{SS05}.
shows that, under certain conditions, the multiplicity of intersection is reflected in the factors appearing in the resultant, even though
the curves are not properly set. In particular, the requirement that no two intersection points lie on a line through the origin can be relaxed, obtaining in this case the total multiplicity of intersection along that line.

\nin The second lemma, Lemma \ref{lem:ch1:GeneralizedResultants}, is a generalization of Corollary 1 in \cite{Perez-Diaz2002}. It shows that generalized resultants can be used to study the intersection points of a finite family of curves. This lemma will be applied in Section \ref{sec:ch4:TotalDegreeFormulaForRationalSurfaces} to the case of surfaces.

\nin As we have said in the preceding paragraphs, the multiplicity of intersection of two projective plane curves can be read at the resultant of their defining polynomials. In fact, this is often used to define the multiplicity of intersection. More precisely (see \cite{Sendra2007}), let $\cC_1$ and $\cC_2$ be projective plane curves, without common components, such that $(1:0:0)\not\in\cC_1\cup\cC_2$, and $(1:0:0)$ does not belong to any line connecting two points in $\cC_1\cap\cC_2$. Let $F(y_0,y_1,y_2)$, resp. $G(y_0,y_1,y_2)$, be the defining
polynomials of $\cC_1$, resp. $\cC_2$. Let $\bar y^o_h=(y^o_0:y^o_1:y^o_2)\in\cC_1\cap\cC_2$, and let
\[R(y_1,y_2)=\Res_{y_0}(F,G)\]
Then the multiplicity of intersection of $\cC_1$ and $\cC_2$ at $\bar y^o_h$, denoted by $\mult_{\bar y^o_h}(\cC_1,\cC_2)$, equals
the multiplicity of the corresponding factor $(y^o_2y_1-y^o_1y_2)$ in $R(y_1,y_2)$. However, in the following Lemma we see how the
multiplicity of intersection of two curves on a line through the origin can be read in the resultant, under certain circumstances,
even though the curves are not properly set. This lemma can be seen as a generalization of Theorem 5.3, page 111 in \cite{Walker1950}.

\vspace{3mm}
\begin{Lemma}\label{lem:ch1:MultiplicityUsingResultants}
Let ${\cC}_1$ and ${\cC}_2$ be two projective algebraic plane curves without common components, given by the homogeneous polynomials
$F(y_0,y_1,y_2)$ and $G(y_0,y_1,y_2)$, respectively. Let $p_1,\ldots,p_k$ be the intersection points, different from $(1:0:0)$, of
${\cC}_1$ and ${\cC}_2$ lying on the line of equation $\beta y_1-\alpha y_2=0$. Then the factor $(\beta y_1-\alpha y_2)$ appears in
$\Res_{y_0}(F,G)$ with multiplicity equal to
\[\sum_{i=1}^k\mult_{p_i}({\cC}_1,{\cC}_2)\]
\end{Lemma}
\begin{proof}
See Lemma 19 in \cite{SS05}.
\end{proof}
\nin The following Lemma is a generalization of Corollary 1 in \cite{Perez-Diaz2002}. It shows that generalized
resultants can be used to study the intersection points of a finite family of curves.
\vspace{3mm}
\begin{Lemma}\label{lem:ch1:GeneralizedResultants}
Let $\cC_0,\ldots,\cC_m$ be the projective plane curves, defined by the homogeneous polynomials $F_0,\ldots,F_m\in{\C}[\bar t_h]$, respectively. Let us suppose that the following hold:
\begin{enumerate}
 \item[\em (i)] $F_1,\ldots,F_m$ have positive degree in $t_0$.
 \item[\em  (ii)] $\deg_{\bar t_h}(F_1)=\cdots=\deg_{\bar t_h}(F_m)$.
 \item[\em  (iii)] $\gcd(F_1,\ldots,F_m)=1$.
\end{enumerate}
Let us denote:
\[F(\bar c,\bar t_h)=c_1F_1(\bar t^h)+\cdots+c_mF_m(\bar t^h)\]
and let
\[R(\bar c,\bar t)=
 \operatorname{Res}_{t_0}\left(F_0(\bar t^h),F(\bar c,\bar t_h)\right),
\]
(note that by {\em (iii)}, $R(\bar c,\bar t)$ is not identically zero). Finally, let $\operatorname{lc}_{t_0}(F_0)\in\C[\bar t]$ and $\operatorname{lc}_{t_0}(F)\in\C[\bar c, \bar t]$ denote, respectively, the leading coefficients w.r.t. $t_0$ of $F_0$  and $F$.

\nin If $\bar t^o=(t^o_1,t^o_2)\in{\C}^2\setminus\{\bar 0\}$ is such that $\operatorname{Cont}_{\bar c}\left(R\right)(\bar t^o)=0$ and
$$\operatorname{lc}_{t_0}(F_0)(\bar t^o)\cdot\operatorname{lc}_{t_0}(F)(\bar c,\bar t^o)\neq 0,$$
there exists $t^o_0$ such that
$\bar t^o_h=(t^o_0:t^o_1:t^o_2)\in\bigcap_{i=0}^m \cC_i$.
\end{Lemma}
\begin{proof}
First, observe that if $\deg_{t_0}(F_0)=0$, then $\operatorname{lc}_{t_0}(F_0)=F_0$ and $R(\bar c,\bar t)=F_0^{\deg_{\bar t_0}(F_1)}$. Thus, in this case the lemma holds trivially, since there is no $\bar t^o\in{\C}^2\setminus\{\bar 0\}$ satisfying the hypothesis of the lemma. Thus, w.l.o.g., in the rest of the proof, we assume that $\deg_{t_0}(F_0)>0$.

\nin Since $\operatorname{lc}_{t_0}(F)(\bar c,\bar t^o)\neq 0$, there exists an open set
$\Phi\subset{\C}^m$ such that if $\bar c^o=(c^o_1,\ldots,c^o_m)\in\Phi$, the leading
coefficient w.r.t.  $t_0$ of $F(\bar c^o,\bar t^o, t_0)\in\C[t_0]$ is $\operatorname{lc}_{t_0}(F)(\bar c^o,\bar t^o)$, and it is non-zero.  Therefore, by the Extension Theorem, (see \cite{Cox1997}, page 159), there exists
$\zeta(\bar c^o)\in{\C}$ (which, in principle, could depend on $\bar c^o$) such that
\[F_0(\zeta(\bar c^o),t^o_1,t^o_2)=F(\zeta(\bar c^o),t^o_1,t^o_2)=0.\]
We claim that there is $t^o_0\in{\C}$ (not depending on $\bar c^o$), such that
$$F_0(t^o_0,t^o_1,t^o_2)=F_1(t^o_0,t^o_1,t^o_2)=\cdots=F_m(t^o_0,t^o_1,t^o_2)=0.$$
To see this note that, since $\operatorname{lc}_{t_0}(F_0)\neq 0$, there is a non-empty finite
set of solutions of the following equation in $t_0$:
$$F_0(t_0,t^o_1,t^o_2)=0.$$
Let $\zeta_1,\ldots,\zeta_p$ be the solutions. If
$$F_1(\zeta_j,t^o_1,t^o_2)=\cdots=F_m(\zeta_j,t^o_1,t^o_2)=0$$
holds  for some $j=1,\ldots,p$, then it suffices to take
$t^o_0=\zeta_j$. Let us suppose that this is not the case, and we will derive a contradiction. Then
there exists an open set $\Phi_1\subset\Phi$, such that if
$\bar c^o\in\Phi_1$, then
$$F(\zeta_j,t^0_1,t^0_2)= c^o_1F_1(\zeta_j,t^0_1,t^0_2)+\cdots+ c^o_mF_m(\zeta_j,t^0_1,t^0_2)\neq 0$$
for every $j=1,\ldots,p$. This means that, for $\bar c^o\in\Phi_1$,
there is no solution of:
\[
\begin{cases}
F_0(t_0,\bar t^o)=0\\
F(\bar c^o,t_0,\bar t^o)=0
\end{cases}
\]
Since the resultant specializes properly in $\Phi$, this implies that:
$R(\bar c^o,\bar t^o)\neq 0$. But, denoting
\[
M(\bar t)=\operatorname{Cont}_{\bar c}\left(R(\bar c,\bar t)\right),
\quad \mbox{ and }
N(\bar c,\bar t)=\operatorname{PP}_{\bar c}\left(R(\bar c,\bar t)\right),
\]
we have
\[R(\bar c^o,\bar t^o)=M(\bar t^o)N(\bar c^o,\bar t^o)=0
\]
because, by hypothesis $M(\bar t^o)=0$. This contradiction proves the result.
\end{proof}
\color{black}

\section{Offset-Line Intersection for Parametric Surfaces}
\label{sec:ch4:OffsetLineIntersectionforRationalSurfaces}\label{ch04-DegreeFormulaeRationalSurfaces}

As we said in the introduction to this paper, the parametric character of $\Sigma$ results in a reduction of the dimension of the space in which we count the points in ${\mathcal O}_d(\Sigma)\cap{\mathcal L}_{\bar k}$. This is so because, instead of counting directly those points, we count the values of the $\bar t$ parameters that generate them. In this section we will show how, with this approach, we are led to an intersection problem between projective plane curves, and we will analyze that problem.
More precisely:
\begin{itemize}
\item Subsection \ref{subsec:ch4:IntersectionWithLines} (page \pageref{subsec:ch4:IntersectionWithLines}) is devoted to the analysis of the intersection between the generic offset and a pencil of lines through the origin. The results in this subsection (see Theorem \ref{thm:ch4:TheoreticalFoundation}, page \pageref{thm:ch4:TheoreticalFoundation}) constitute the theoretical foundation of the degree formula to be derived in Section \ref{sec:ch4:TotalDegreeFormulaForRationalSurfaces}.
\item In Subsection \ref{subsec:ch4:EliminationAndAuxiliaryPolynomials} we describe the auxiliary polynomials obtained by using elimination techniques in the Parametric Offset-Line System, and we introduce a new auxiliary system,  see System \ref{sys:ch4:AuxiliaryCurvesSystem}. Also, we obtain some geometric properties of the solutions of this new system ${{\mathfrak S}^P_3}(d,\bar k)$ in Proposition \ref{prop:ch4:ExtendableSolutions} (page \pageref{prop:ch4:ExtendableSolutions}) and the subsequent Lemma \ref{lem:ch4:SignOfLambdaAndOffsetting} (page \pageref{lem:ch4:SignOfLambdaAndOffsetting}). These results will be used in the sequel to elucidate the relation between the solution sets of Systems $\mathfrak S^P_2(d,\bar k)$ and ${{\mathfrak S}^P_3}(d,\bar k)$.
\item In Subsection \ref{subsec:ch4:FakePoints} (page \pageref{subsec:ch4:FakePoints}) we define the corresponding notion of fake points and invariant points for the Affine Auxiliary System ${{\mathfrak S}^P_3}(d,\bar k)$. The main result of this subsection is Proposition \ref{prop:ch4:FakePointsAndInvariantSolutionsCoincide} (page \pageref{prop:ch4:FakePointsAndInvariantSolutionsCoincide}), that shows the relation between these two notions.
\end{itemize}

\subsection{Intersection with lines}\label{subsec:ch4:IntersectionWithLines}

As in the case of plane curves (see our paper \cite{SS05}), we will address the degree problem for surfaces
by counting the number of intersection points between ${\mathcal O}_d(\Sigma)$ and a generic line
through the origin. More precisely, let us consider a family of lines through the origin,
{\sf denoted} by $\cL_{\bar k}$, whose direction is determined
by the values of the variable $\bar k=(k_1,k_2,k_3)$. The family $\cL_{\bar k}$ is described by the following set of parametric equations:
\[\cL_{\bar k}\equiv
\begin{cases}
\ell_1(\bar k,l,\bar x):\,{ x_1}-k_1\,l=0\\
\ell_2(\bar k,l,\bar x):\,{ x_2}-k_2\,l=0\\
\ell_3(\bar k,l,\bar x):\,{ x_3}-k_3\,l=0
\end{cases}
\]
A particular line of the family, corresponding to the value $\bar k^o$, will be denoted by $\cL_{\bar k^o}$.
We add the equations $\ell_1, \ell_2, \ell_3$ of $\cL_{\bar k}$ to the equations of the parametric system for the generic offset (System \ref{sys:ch4:ParametricOffsetSystem} in page \pageref{sys:ch4:ParametricOffsetSystem}), and we arrive at the following system:
\begin{equation}\label{sys:ch4:IntersectionOffsetLine}
\hspace{-8mm}\begin{minipage}{14cm}
\[\hspace{3mm}\mathfrak S^P_2(d,\bar k)\equiv
\begin{cases}
b^{P}(d,\bar t,\bar x):
\left (P_0{x_1}-P_1\right)^{_2}+\left (P_0{ x_2}-P_2\right )^{_2}+\left ( P_0{ x_3}-P_3\right)^{_2}-{d}^{_2}{P_0}^{_2}=0\\
\normal^{P}_{(1,2)}(\bar t,\bar x):\,n_1\cdot (P_0x_2-P_2)-n_2\cdot (P_0x_1-P_1)=0\\
\normal^{P}_{(1,3)}(\bar t,\bar x):\,n_1\cdot (P_0x_3-P_3)-n_3\cdot (P_0x_1-P_1)=0\\
\normal^{P}_{(2,3)}(\bar t,\bar x):\,n_2\cdot (P_0x_3-P_3)-n_3\cdot (P_0x_2-P_2)=0\\
w^{P}(r,\bar t):\,r\cdot P_0\cdot \beta\cdot h-1=0\\
\ell_1(\bar k,l,\bar x):\,{ x_1}-k_1\,l=0\\
\ell_2(\bar k,l,\bar x):\,{ x_2}-k_2\,l=0\\
\ell_3(\bar k,l,\bar x):\,{ x_3}-k_3\,l=0
\end{cases}
\]
\end{minipage}
\end{equation}
We will refer to this as the {\sf Parametric Offset-Line System}. The next step is the study of the generic solutions of this system.

\nin We need to exclude certain degenerate situations that arise for a set of values of $(d,\bar k)$. For example, a degenerated situation arises if
the set of points of ${\Sigma}$ where the normal line to ${\Sigma}$ passes through the origin is {\em too big}. The next Lemma says
that this can only happen if ${\Sigma}$ is a sphere centered at the origin.
\vspace{3mm}
\begin{Lemma}\label{lem:Ch1:LineMeetsVarietyNormallyInProperClosedSet}
Let ${\Sigma}_{\bot}\subset{\Sigma}$ denote the set of regular points $\bar y^o\in{\Sigma}$ such that the
normal line to ${\Sigma}$ at $\bar y^o$ is parallel to $\bar y^o$. If ${\Sigma}$ is not a sphere centered at the origin,
then ${\Sigma}_{\bot}^*$ is a proper (possibly empty) closed subset of ${\Sigma}$.
\end{Lemma}
\begin{proof} Let us assume that ${\Sigma}_{\bot}$ is nonempty.
Let, as usual, $f(\bar{y})$ be the irreducible polynomial defining ${\Sigma}$, and let $\tilde{\Sigma}$ be the  algebraic set in ${\C}^3$ defined by:
\[
\begin{cases}
f(\bar{y})=0\\
f_i(\bar y)y_j-f_j(\bar y)y_i=0&(\mbox{for }i,j=1,\ldots,3; i<j).
\end{cases}
\]
Note that this set of equations implies $\bar y^o\parallel\nabla f(\bar y^o)$ for $\bar y^o\in{\Sigma}$. Then
${\Sigma}_{\bot}\subset\tilde{\Sigma}\subset{\Sigma}$. Therefore, it suffices to prove that $\tilde{\Sigma}\neq{\Sigma}$. Let
us suppose that $\tilde{\Sigma}={\Sigma}$. Let
\[K(\bar y)=\bar y\cdot\nabla f(\bar y)=\sum_{j=1}^3y_jf_j(\bar{y}).\]
Then for every $\bar y^o\in{\Sigma}$, using that $f_i(\bar y^o)y_j^o=f_j(\bar y^o)y_i^o$
one has that
\[f_i(\bar y^o)K(\bar y^o)=
\sum_{j=1}^3f_i(\bar y^o)y_j^of_j(\bar y^o)=y_i^o\sum_{j=1}^3f_j(\bar y^o)^2=y_i^oh(\bar y^o),\]
for $i=1,2,3$.
Now let $\bar t=(t_1,t_{2})$ and let $\cQ(\bar t)=(Q_1,Q_2,Q_3)(\bar t)$ be a local parametrization of ${\Sigma}$. Substituting $\cQ$ in the
above relation:
\[ f_i(\cQ(\bar t))K(\cQ(\bar t))=Q_i(\bar t) h(P(\bar t))\]
that is, $K(\cQ(\bar t))\nabla f(\cQ(\bar t))=h(\cQ(\bar t))\cQ(\bar t)$.
Using Prop. 2 in \cite{Sendra1999}, we know that $h(\cQ(\bar t))\neq 0$, and so $K(\cQ(\bar t))\neq 0$. Thus:
\[
\frac{h(\cQ(\bar t))}{K(\cQ(\bar t))}Q_i(\bar t)=f_i(\cQ(\bar t)).
\]
On the other hand, since $f(\cQ(\bar t))=0$, deriving w.r.t. $t_j, (j=1,2)$  one has:
\[
\sum_{i=1}^3f_i(\cQ(\bar t))\frac{\partial Q_i(\bar t)}{\partial t_j}=
\frac{h(\cQ(\bar t))}{K(\cQ(\bar t))}\sum_{i=1}^3Q_i(\bar t)\frac{\partial Q_i(\bar t)}{\partial t_j}=
0\]
From this, one concludes that
\[
\frac{\partial }{\partial t_j}\left(\sum_{i=1}^3Q_i^2(\bar t)\right)=2
\sum_{i=1}^3Q_i(\bar t)\frac{\partial Q_i(\bar t)}{\partial t_j}=0
\]
for $j=1,2$. This means that $\sum_{i=1}^3Q_i^2(\bar t)=c$ for some constant $c\in{\C}$. Since ${\Sigma}$ is assumed not to be normal-isotropic,
one has $c\neq 0$, and since the parametrization converges locally, we conclude that ${\Sigma}$ equals a sphere centered at the origin.
\end{proof}
\vspace{3mm}
\begin{Remark}
Note that, if in Lemma \ref{lem:Ch1:LineMeetsVarietyNormallyInProperClosedSet} we consider those regular points $\bar y^o$ of ${\Sigma}$
such that the normal line to ${\Sigma}$ at $\bar y^o$ is parallel to the vector $\bar y^o-\bar a$ for a fixed $\bar a\in{\C}^3$,
then ${\Sigma}^*_{\bot}$ is a proper (possibly empty) closed subset of ${\Sigma}$, unless ${\Sigma}$ is a sphere centered at $\bar a$.
\end{Remark}

\noindent A closer analysis of the proof of Lemma \ref{lem:Ch1:LineMeetsVarietyNormallyInProperClosedSet} shows that in fact we have also proved
the following:
\vspace{3mm}
\begin{Corollary}\label{cor:ch1:PositiveDimensionComponentsOfBotVareInSpheres}
If $\cW$ is any irreducible component of ${\Sigma}_{\bot}^*$ with $\dim(\cW)>0$ then $\cW$ is contained in a sphere centered
at the origin. That is, there exists $d^o\in{\C}^\times$ such that if $\bar y^o\in\cW$, then
\[(y^0_1)^2+(y^0_2)^2+(y^0_3)^2=(d^o)^2.\]
Since ${\Sigma}_{\bot}^*$ has at most finitely many irreducible components, it follows that there is a finite set of distances $\{d_1^{\bot},\ldots,d_p^{\bot}\}$ such that ${\Sigma}_{\bot}^*$ is contained in the union the spheres centered at the origin and with radius $d_i^{\bot}$ for $i=1,\ldots,p$.
\end{Corollary}
\nin We will use the notation $\Upsilon({\Sigma}_{\bot})=\{d_1^{\bot},\ldots,d_p^{\bot}\}$, and we will say that $\Upsilon({\Sigma}_{\bot})$ is {\sf the set of critical distances of ${\Sigma}$.}

\nin The following lemma is the basic tool to avoid the remaining degenerated situations in the analysis of the offset-line intersection:
for a given proper closed subset $\mfF\subset\Sigma$, it shows
\begin{enumerate}
\item how to avoid the set of values $(d^o,\bar k^o)$ such that  $\cL_{\bar k^o}\setminus\{\bar 0\}$ meets $\Sigma$ in a point $\bar y^o\in\mfF$,
\item and how to avoid the set of values $(d^o,\bar k^o)$ such that $\cL_{\bar k^o}\setminus\{\bar 0\}$ meets $\cO_{d^o}(\Sigma)$ in a point $\bar x^o$ associated to $\bar y^o\in\mfF$.
\end{enumerate}

\nin In the proof of the Lemma we will use the polynomials $f, h, b$ and $\normal_{(i,j)}$ (for $i,j=1,\ldots,3; i<j$), introduced with System $\mfG_1(d)$ in page \pageref{sys:Ch1:GenericOffsetSystem}. For the convenience of the reader we repeat that system here:
\[
\left.\begin{array}{lr}
&f(\bar y)=0\\
\normal_{(i,j)}(\bar x,\bar y): &f_i(\bar y)(x_j-y_j)-f_j(\bar y)(x_i-y_i)=0\\
(\mbox{for }i,j=1,\ldots,3; i<j)&\\
b(d, \bar x, \bar y ):& (x_1-y_1)^2+(x_2-y_2)^2+(x_3-y_3)^2-d^2=0\\
w(\bar y,u):& u\cdot(\|\nabla f(\bar y)\|^2)-1=0
\end{array}\right\}
\equiv\mathfrak{S}_1(d).
\]
and that $h(\bar t)=n_1(\bar t)^2+n_2(\bar t)^2+n_3(\bar t)^2$, while $h_{\operatorname{imp}}(\bar t)=\|\nabla f(\bar y)\|^2$.
\vspace{3mm}
\begin{Lemma}\label{lem:ch4:AvoidClosedSubsetsSigma}
Let $\mfF\subsetneq\Sigma$ be closed. There exists an open $\Omega_{\mfF}\subset\C^4$, such that if $(d^o,\bar k^o)\in \Omega_{\mfF}$, the following hold:
\begin{enumerate}
 \item[(1)] $\cL_{\bar k^o}\cap\left(\mfF\setminus\{\bar 0\}\right)=\emptyset$.
 \item[(2)] If $\bar x^o\in(\cL_{\bar k^o}\cap \cO_{d^o}(\Sigma)){\setminus \{\bar 0\}}$, there is no solution $(d^o,\bar x^o,\bar y^o,u^o)$ of System $\mfG_1(d)$ (System \ref{sys:Ch1:GenericOffsetSystem} in page \pageref{sys:Ch1:GenericOffsetSystem}) with $\bar y^o\in\mfF$.
\end{enumerate}
\end{Lemma}
\nin{\em Proof.}
If $\mfF$ is empty, the result is trivial. Thus, let us assume that $\mfF\neq\emptyset$, and let the defining polynomials of $\mfF$ be $\{\phi_1(\bar y),\ldots,\phi_p(\bar y)\}\subset\C[\bar y]$.  We will show that one may take $\Omega_{\mfF}=\Omega^1_{\mfF}\cap\Omega^2_{\mfF}$, where  $\Omega^1_{\mfF}, \Omega^2_{\mfF}$ are two open sets constructed as follows:
\begin{enumerate}
 \item[(a)]  Let us consider the following ideal in $\C[\bar k,\rho,v,\bar y]$:
\[\cI=<f(\bar y), \phi_1(\bar y),\ldots,\phi_p(\bar y),\bar y-\rho\cdot\bar k,v\cdot\rho-1>,\]
and the projection maps defined in its solution set $\bV(\cI)$ as follows:
\[\pi_{(1,1)}(\bar k,\rho,v,\bar y)=\bar y, \quad \pi_{(1,2)}(\bar k,\rho,v,\bar y)=\bar k\]
We show first that $\pi_{(1,1)}(\bV(\cI))=\mfF$. The inclusion $\pi_{(1,1)}(\bV(\cI))\subset\mfF$ is trivial; and if $\bar y^o\in\mfF$, then since $\cF\subset\Sigma$, $(\bar y^o,1,1,\bar y^o)\in\bV(\cI)$ proves the reversed inclusion. Therefore, since $\cF\subsetneq\Sigma$, $\dim(\pi_{(1,1)}(\bV(\cI)))=\dim(\mathfrak{F})<2$. Besides, for every $\bar y^o\in\pi_{(1,1)}(\bV(\cI))$ one has:
\[\pi_{(1,1)}^{-1}(\bar y^o)=\{(v^o\bar y^o,\frac{1}{v^o},v^o,\bar y^o)\,|\,v^o\in\C^\times \} \]
from where one has that $\dim(\pi_{(1,1)}^{-1}(\bar y^o))=1$. Since the dimension of the fiber does not depend on $\bar y^o$, applying Lemma \ref{lem:ch1:FiberDimension} (page \pageref{lem:ch1:FiberDimension}), we obtain
$\dim(\bV(\cI))<3$. Thus, $\dim(\pi_{(1,2)}(\bV(\cI)))<3$. It follows that $(\pi_{(1,2)}(\bV(\cI)))^*$ is a proper closed subset of $\C^3$. Let $\Theta^1=\C^3\setminus(\pi_{(1,2)}(\bV(\cI)))^*$, and let $\Omega_{\mfF}^1=\C\times\Theta^1$.

\item[(b)] Let us consider the following ideal in $\C[d,\bar k,\rho,v,\bar x,\bar y]$:
\[
\begin{array}{l}
\cJ=<f(\bar y), b(d,\bar x,\bar y),\normal_{(1,2)}(\bar x,\bar y),\normal_{(1,3)}(\bar x,\bar y),\normal_{(2,3)}(\bar x,\bar y),\\
\bar x-\rho\cdot \bar k,v\cdot\rho\cdot d\cdot h_{\operatorname{imp}}(\bar y)-1, \phi_1(\bar y),\ldots,\phi_p(\bar y)>
\end{array}
\]
and the projection maps defined in its solution set $\bV(\cJ)\subset\C^{12}$ as follows:
\[\pi_{(2,1)}(d,\bar k,\rho,v,\bar x,\bar y)=\bar y, \quad \pi_{(2,2)}(d,\bar k,\rho,v,\bar x,\bar y)=(d,\bar k)\]
Then $\pi_{(2,1)}(\bV(\cJ))\subset\mfF$. Therefore $\dim(\pi_{(2,1)}(\bV(\cJ)))\leq 1$. Let $\bar y^o\in\pi_{(2,1)}(\bV(\cJ))$.
Note that then $h_{\operatorname{imp}}(\bar y^o)\neq 0$. We denote $\sigma^o=\sqrt{h_{\operatorname{imp}}(\bar y^o)}$ (a particular choice
of the square root); clearly $\sigma^o\neq 0$. Then, it holds that:
\[ \pi_{(2,1)}^{-1}(\bar y^o)=
\biggl\{
\biggl(d^o,\frac{1}{\rho^o} \left(\bar y^o\pm\frac{d^o}{\sigma^o}\nabla(\bar y^o)\right),\rho^o,\frac{1}{(\sigma^o)^{2}\rho^od^o},
\bar y^o\pm\frac{d^o}{\sigma^o}\nabla(\bar y^o), \bar y^o
\biggr)
\,
\biggl\lvert\, d^o, \rho^o\in {\C}^\times\biggr\}
\]
Therefore $\dim(\pi_{(2,1)}^{-1}(\bar y^o))=2.$ Applying Lemma \ref{lem:ch1:FiberDimension} again, one has
\[\dim(\bV(\cJ))=2+\dim(\pi_{(2,1)}(\bV(\cJ))\leq 3.\]
It follows that $\dim(\pi_{(2,2)}(\bV(\cJ)))\leq 3$. Let us take $\Omega^2_{\mfF}={\C}^{4}\setminus \pi_{(2,2)}({\cal V})^*.$
\end{enumerate}
Let $\Omega_{\mfF}=\Omega^1_{\mfF}\cap\Omega^2_{\mfF}$ and let $(d^o,\bar k^o)\in\Omega_\mfF$.
\begin{enumerate}
\item If $\bar y^o\in\cL_{\bar k^o}\cap\left(\mfF\setminus\{\bar 0\}\right)$, then there is some $\rho^o\in\C^\times$ such that $\bar y^o=\rho^o\bar k^o$. It follows that $(\bar k^o,\rho^o,\dfrac{1}{\rho^o},\bar y^o)\in\bV(\cI)$, and so $\bar k^o\in\pi_{(1,2)}(\bV(\cI))$, contradicting the construction of $\Omega^1_{\mfF}$. This proves statement (1).
\item If $\bar x^o\in(\cL_{\bar k^o}\cap \cO_{d^o}(\Sigma)){\setminus \{\bar 0\}}$, and there is a solution $(d^o,\bar x^o,\bar y^o,u^o)$ of System $\mfG_1(d)$ with $\bar y^o\in\mfF$, then there is some $\rho^o\in\C^\times$ such that $\bar x^o=\rho^o\bar k^o$. It follows that $(d^o,\bar k^o,\rho^o,\dfrac{1}{\rho^o\cdot d^o\cdot h_{\operatorname{imp}}(\bar y^o)},\bar x^o,\bar y^o)\in\bV(\cJ)$. Therefore $(d^o,\bar k^o)\in\pi_{(2,2)}(\bV(\cJ))$, contradicting the construction of $\Omega^2_{\mfF}$. This proves statement (2).\hspace{35mm}$\Box$
\end{enumerate}

\vspace{3mm}
\begin{Remark}
Note that the origin may belong to $\mfF$. In that case, Lemma \ref{lem:ch4:AvoidClosedSubsetsSigma}(1) guarantees that the origin is the only point in $\cL_{\bar k^o}\cap\cF$. Correspondingly, part (2) of the lemma guarantees that the remaining points in $\cL_{\bar k^o}\cap \cO_{d^o}(\Sigma)$ cannot be extended to a solution $(d^o,\bar x^o,\bar y^o,u^o)$ of System $\mfG_1(d)$ with $\bar y^o\in\mfF$.
\end{Remark}

\nin Our next goal is to prove a theorem (Theorem \ref{thm:ch4:TheoreticalFoundation} below), that gives the theoretical foundation for our approach
to the degree problem. Theorem \ref{thm:ch4:TheoreticalFoundation} is the analogous of Theorem 5 in our paper \cite{SS05}. That theorem is
preceded by Lemma 4, that states that for a curve $\cC$, $\bar 0\in\cO_{d^o}(\cC)$ for at most finitely many values $d^o\in\C$.
However, we have not been able to prove a similar result for the case of surfaces: the main difficulty is that a surface can have
infinitely many singular points. Even if we restrict ourselves to the case of parametric surfaces, we still have to take into account the possible
existence of a singular curve contained in $\Sigma$, and not contained in the image of the parametrization.  Besides, in the proof of the
theorem we will use Lemma \ref{lem:Ch1:LineMeetsVarietyNormallyInProperClosedSet} (page \pageref{lem:Ch1:LineMeetsVarietyNormallyInProperClosedSet}),
that does not apply when $\Sigma$ is a sphere centered at the origin. This is the reason for the Assumptions that we announced in the Introduction
of this paper, and that we state formally here. In the sequel, we assume that:

\vspace{3mm}
\begin{Assumptions}\label{rem:ch4:NotInfinitelyManyOffsetsThroughOrigin}
\begin{enumerate}
 \item[]
 \item[(1)] {\sf There exists a finite subset $\Delta^1$ of $\C$ such that, for $d^o\not\in\Delta^1$ the origin does not belong to $\cO_{d^o}({\Sigma})$.}
 \item[(2)] {\sf $\Sigma$ is not a sphere centered at the origin.}
\end{enumerate}
\end{Assumptions}

\nin Before stating the theorem we have to introduce some terminology.
\vspace{3mm}
\begin{Remark}
For $(d^o,\bar k^o)\in\C^4$ we will {\sf denote} by $\Psi_2^P(d^o,\bar k^o)$ the set of solutions of System {$\mathfrak S^P_2(d^o,\bar k^o)$} in the variables $(l,r,\bar t,\bar x)$ (see (\ref{sys:ch4:IntersectionOffsetLine}) in page \pageref{sys:ch4:IntersectionOffsetLine}).
\end{Remark}
\vspace{3mm}
\begin{Theorem}\label{thm:ch4:TheoreticalFoundation}
Let $\Sigma$ satisfy the hypothesis in Remark \ref{rem:ch4:NotInfinitelyManyOffsetsThroughOrigin}. There exists a non-empty Zariski-open subset $\Omega_0\subset\C^4$,
such that if {$(d^o,\bar k^o)\in\Omega_0$}, then
\begin{itemize}
 \item[(a)] if $\bar y^o\in\mathcal L_{\bar k_0}\cap (\Sigma\setminus \{\bar 0\})$, then no normal vector to $\Sigma$ at $\bar y^o$ is parallel to $\bar y^o$.
 \item[(b)] $\Psi_2^P(d^o,\bar k^o)$ has precisely  $m\delta$ elements (recall that $m$ is the tracing index of $P$ and $\delta$ the total degree of the generic offset).
 Besides, the set $\Psi_2^P(d^o,\bar k^o)$ can be partitioned as a disjoint union:
 \[\Psi_2^P(d^o,\bar k^o)=
 \Psi_{2}^1(d^o,\bar k^o)\cup\cdots\cup\Psi_{2}^{\delta}(d^o,\bar k^o),
 \]
 such that:
 \begin{enumerate}
  \item[(b1)] $\#\Psi_{2}^i(d^o,\bar k^o)=m$ for $i=1,\ldots,\delta$.
  \item[(b2)] The $m$ elements of $\#\Psi_{2}^i(d^o,\bar k^o)$ have the same values of the variables $(l,r,\bar x)$, and differ only in the value of $\bar t$. Besides, for $(l^o,r^o,\bar t^o,\bar x^o)\in\Psi_{2}^i(d^o,\bar k^o)$, the point $P(\bar t^o)\in\Sigma$ does not depend on the choice of $\bar t^o$.
  \end{enumerate}
  Let us denote by $(l^o_i,r^o_i,\bar t^o_{h,i},\bar x^o_i)$ an element of $\Psi_{2}^i(d^o,\bar k^o)$. Then
 \begin{enumerate}
   \item[(b3)] The points $\bar x^o_1,\ldots,\bar x^o_{\delta}$ are all different (and different from $\bar 0$), and
  $$\cL_{\bar k^o}\cap\cO_{d^o}(\Sigma)=\{\bar x^o_1,\ldots,\bar x^o_{\delta}\}.
  $$
  Furthermore, $\bar x^o_i$ is non normal-isotropic in $\cO_{d^o}(\Sigma)$, for $i=1,\ldots,\delta$.
  \item[(b4)] The $\delta$ points
  $$\bar y^o_1=P(\bar t^o_{h,1}),\cdots,\bar y^o_{\delta}=P(\bar t^o_{h,_{\delta}})$$
  are affine, distinct and non normal-isotropic points of $\Sigma$.
 \end{enumerate}
\item[(c)] $k_i^o\neq 0$ for $i=1,2,3$.
\end{itemize}
\end{Theorem}

\nin{\em Proof.}\quad
Let $\Delta_0^1=\{d^o\in {\C} \,|\,  g(d^o,\bar 0)\neq 0\}$.
The assumption in Remark \ref{rem:ch4:NotInfinitelyManyOffsetsThroughOrigin} (page \pageref{rem:ch4:NotInfinitelyManyOffsetsThroughOrigin}) implies that $\Delta_0^1$ is an
open non-empty subset of $\C$. Let  $\Delta$ be as in Corollary
\ref{cor:ch1:BadDistancesFiniteSet}, (page \pageref{cor:ch1:BadDistancesFiniteSet}), and let $\Omega^0_0=(\Delta_0^1\cap(\C\setminus\Delta))\times(\C^3\setminus\bigl(\{\bar k^o / k_i^o=0\mbox{ for some }i=1,2,3\}\bigr))$.
Next, let us consider $g(d,\bar x)$ expressed as follows:
\[g(d,\bar x)=\sum_{i=0}^{\delta}g_i(d,\bar x)\]
where $g_i$ is a degree $i$ form in $\bar x$. We consider:
\[\tilde{g}(d,\bar k,\rho)=g(d,\rho\bar k)=\sum_{i=0}^{\delta}g_i(d,\bar k)\rho^i.\]
This polynomial is not identically zero, is primitive w.r.t. $\bar x$ (see Lemma \ref{lem:ch1:GenericOffsetIsPrimitiveIn_x}, page \pageref{lem:ch1:GenericOffsetIsPrimitiveIn_x}), and it is squarefree; note that $g(d,\bar x)$ is square-free by Remark \ref{rem:ch1:GenericOffsetEqSqfreeAndHasAtMostTwoFactors} (page \pageref{rem:ch1:GenericOffsetEqSqfreeAndHasAtMostTwoFactors}), and therefore $\tilde g$ is square-free too. Thus, the discriminant
\[Q(d, \bar k)=\operatorname{Dis}_{\rho}\left(\tilde{g}(d,\bar k,\rho)\right)\]
is not identically zero either.

\nin In this situation, let us take
\[\Omega^0_1=\Omega^0_0\setminus \{(d^o,\bar k^o)\in\C^4 / Q(d^o, \bar k^o)\cdot g_{0}(d^o, \bar k^o)\cdot g_{\delta}(d^o, \bar k^o)=0\}.\]
For $(d^o,\bar k^o)\in\Omega^0_1$, $g(d^o,\rho\bar k_0)$ has $\delta$ different and non-zero roots; say, $\rho_1,\ldots,\rho_{\delta}$.
Therefore, $\cL_{\bar k^o}$ intersects $\cO_{d^o}(\Sigma)$ in $\delta$ different points:
\[\bar x^o_1=\rho_1\bar k^o,\ldots,\bar x^o_{\delta}=\rho_{\delta}\bar k^o,\]
and none of these points is the origin.

\nin We will now construct an open subset $\Omega^0_2\subset\Omega^0_1$ such that for $(d^o,\bar k^o)\in\Omega^0_2$, the points $\bar x^o_1,\ldots,\bar x^o_{\delta}$ are non-normal isotropic points in $\cO_{d^o}(\Sigma)$, and each one of them is associated with a unique non-normal isotropic point of $\Sigma$.  To do this, recall that $\Iso(\Sigma)$ is the closed set of normal-isotropic points of $\Sigma$ (see page \pageref{notation:ch1:Hodograph}), and  let $\Omega_{\Iso(\Sigma)}$ be the set obtained when applying Lemma \ref{lem:ch4:AvoidClosedSubsetsSigma} (page \pageref{lem:ch4:AvoidClosedSubsetsSigma}) to the closed subset $\mathfrak F=\Iso(\Sigma)$. Note that if $(d^o,\bar k^o)\in\Omega^0_1\cap\Omega_{\Iso(\Sigma)}$, then the points $\bar x^o_1,\ldots,\bar x^o_{\delta}$ are not associated with normal-isotropic points of $\Sigma$.  Let us consider the polynomial
\[\Gamma(d,\bar x)=\left(\dfrac{\partial g}{\partial x_1}(d,\bar x)\right)^2+\left(\dfrac{\partial g}{\partial x_2}(d,\bar x)\right)^2+\left(\dfrac{\partial g}{\partial x_3}(d,\bar x)\right)^2.\]
This polynomial is not identically zero, because in that case for every $d^o\not\in\Delta$ all the points in $\cO_{d^o}(\Sigma)$ would be normal-isotropic, contradicting Proposition \ref{prop:ch1:OffsetPropertiesRegardingSimpleAndSpecialComponents}(3) (page \pageref{prop:ch1:OffsetPropertiesRegardingSimpleAndSpecialComponents}). Let then
$\tilde{\Gamma}(d, \bar k, r)=\Gamma(d,r \bar k)$, and consider the resultant:
\[\Phi(d, \bar k)=\Res_r(\tilde{g}(d,\bar k,r),\tilde{\Gamma}(d, \bar k, r))\]
If $\Phi(d, \bar k)\equiv 0$, then  $\tilde{g}(d,\bar k,r)$ y $\tilde{\Gamma}(d,\bar k,r)$ have a common factor of positive degree in $r$. Let us show that this leads to a contradiction. Suppose that
\[
\begin{cases}
\tilde{g}(d,\bar k,r)=M(d, \bar k, r) G(d,\bar k,r), \\
\tilde{\Gamma}(d,\bar k,r)=M(d, \bar k, r) \Gamma^*(d,\bar k,r).
\end{cases}
\]
Then $M$ depends on $\bar k$ (because $\tilde g$ cannot have a non constant factor in $\C[d,r]$). Take  therefore
$r^o\in\C^\times$ such that  $M(d,\frac{\bar k}{r^o}, r^o)$ depends
on $\bar k$. Then:
\[
\begin{cases}
g(d,\bar x)=g(d,r^o\frac{\bar x}{r^o})=
\tilde{g}(d,\frac{\bar x}{r^o},r^o)=M(d,\frac{\bar x}{r^o}, r^o) G(d,\frac{\bar x}{r^o}, r^o) \\[3mm]
\Gamma(d,\bar x)=
\Gamma(d,r^o\frac{\bar x}{r^o} )=\tilde{\Gamma}(d,\frac{\bar x}{r^o},r^o)=
M(d,\frac{\bar x}{r^o}, r^o) \Gamma^*(d,\frac{\bar x}{r^o}, r^o)
\end{cases}
\]
But since $g$ has at most two irreducible components, this would imply that for $d^o\not\in\Delta$, $\cO_{d^o}(\Sigma)$ would have at least a normal-isotropic component, contradicting Proposition \ref{prop:ch1:OffsetPropertiesRegardingSimpleAndSpecialComponents}(3) (page \pageref{prop:ch1:OffsetPropertiesRegardingSimpleAndSpecialComponents}). Therefore, the equation $\Phi(d, \bar k)=0$ defines a proper closed subset of $\C^4$. This shows that we can take:
$$\Omega^0_2=\left(\Omega^0_1\cap\Omega_{\Iso(\Sigma)}\right)\setminus \{(d^o,\bar k^o):\Phi(d^o, \bar k^o)=0 \}$$
Then, for $(d^o,\bar k^o)\in\Omega^0_2$,  each of the points $\bar x^o_i,$ for $i=1,\ldots,\delta$,
is associated with a unique non-normal isotropic point $\bar y^o_i$ of $\Sigma$ (recall that $d^o\in\Delta$, and so the irreducible components of $\cO_{d^o}(\Sigma)$ are simple).

\nin Let $\Omega_{\bot}$ be the open subset of  $\C\times\C^3$ obtained by applying Lemma \ref{lem:ch4:AvoidClosedSubsetsSigma} (page \pageref{lem:ch4:AvoidClosedSubsetsSigma}) to the closed subset $\Sigma_{\bot}$ whose existence is guaranteed by Lemma \ref{lem:Ch1:LineMeetsVarietyNormallyInProperClosedSet} (page \pageref{lem:Ch1:LineMeetsVarietyNormallyInProperClosedSet}). Recall that, by assumption (see Remark \ref{rem:ch4:NotInfinitelyManyOffsetsThroughOrigin}(2), page \pageref{rem:ch4:NotInfinitelyManyOffsetsThroughOrigin}), $\Sigma$ is not a sphere centered at the origin. Besides, let $\Theta=\C^3\setminus {\cal L}_0$, where
\[{\cal L}_0=\begin{cases}
\emptyset\mbox{ if } \bar 0\not\in\Sigma\mbox{ or if }\bar 0\in\Sing(\Sigma)\\
\mbox{ the normal line to }\Sigma\mbox{ at } \bar 0\mbox{ otherwise}.
\end{cases}
\]
and set
\[ \Omega^0_3=\Omega^0_2\cap\Omega_{\bot}\cap(\C\times\Theta). \]
Then for $(d^o,\bar k^o)\in\Omega^o_3$, the points $\bar y^o_i, i=1,\ldots,\delta$, are different. To prove this, note that if  $\bar y_i=\bar y_j$, with $i\neq j$, then $\bar
y^o_i$ generates $\bar x^o_i$ and $\bar x^o_j$. Thus, since $\bar y^o_i,\bar x^o_i,
\bar x^o_j$ are all in the normal line to $\Sigma$ at $\bar y^o_i$ and in $\cL_{\bar k^o}$, it follows that these two lines coincide. This means that $\bar y^o_i\in\Sigma_{\bot}$. Since $(d^o,\bar k^o)\in\Omega_{\bot}$, then (by Lemma \ref{lem:ch4:AvoidClosedSubsetsSigma}) $\bar y^o_i\in {\cL_{\bar k^o}}\cap\Omega_{\bot}$ implies that $\bar y^o_i=\bar 0$, in contradiction with $\bar k^o\in\C\times\Theta$.

\nin We will now show that it is possible to restrict the values of $(d,\bar k)$ so that the points $\bar y^o_i$ belong to the image of the parametrization $P$. Let $\Upsilon_2$ be as in Lemma \ref{lem:ch4:PropertiesSurfaceParametrization} (page \pageref{lem:ch4:PropertiesSurfaceParametrization}), and let $\Omega_{\Upsilon_2}\subset\C\times\C^3$ be the open subset obtained applying Lemma \ref{lem:ch4:AvoidClosedSubsetsSigma} to $\Sigma\setminus\Upsilon_2$. Then take
$\Omega^0_4=\Omega^0_3\cap\Omega_{\Upsilon_2}$.

\nin Let us show that we can take $\Omega_0=\Omega^0_4$. If $(d^o,\bar k^o)\in\Omega^0_4$, then for each of the points $\bar y^o_i$ there are $\mu$ values $\bar t^o_{(i,j)}$ (with $i=1,\ldots,\delta$, $j=1,\ldots,m$) such that $P(\bar t^o_{(i,j)})=\bar y^o_i$. Setting $\Psi_2^i(d^o,\bar k^o)=\{\bar t^o_{(i,j)}\}_{j=1,\ldots,m}$, one has that
\[
\Psi_2^P(d^o,\bar k^o)=\Psi_2^1(d^o,\bar k^o)\cup\cdots\cup\Psi_2^\delta(d^o,\bar k^o)
\]
and so the first part of claim (2) is proved. Furthermore:
\begin{itemize}
 \item claim (a) holds because of the construction of $\Omega^0_3$.
 \item the structure of $\Psi_2(d^o,\bar k^o)$ in claims (b1) and (b2) holds because of the construction of $\Omega^0_4$.
 \item Claims (b3) and (b4) hold because of the construction of $\Omega^0_0, \Omega^0_1$ and $\Omega^0_2$.
 \item Claim (c) follows the construction of $\Omega^0_0$.
\end{itemize}
\qed

\vspace{3mm}
\begin{Remark}\label{rem:ch4:InOmega0GoodSpecialization}
Note that, by the construction of $\Omega_0^0$ in the proof of Theorem \ref{thm:ch4:TheoreticalFoundation} (page \pageref{thm:ch4:TheoreticalFoundation}), if $(d^o,\bar k^o)\in\Omega_0$, then $g(d^o, \bar x)=0$ is the equation of ${\cO}_{d^o}(\Sigma)$.
\end{Remark}

\subsection{Elimination and auxiliary polynomials}\label{subsec:ch4:EliminationAndAuxiliaryPolynomials}\label{sec:ch4:AuxiliaryCurvesForRationalSurfaces}

 To continue with our strategy, we proceed to eliminate the variables $(l,r,\bar x)$ in the Parametric Offset-Line System $\mathfrak S^P_2(d,\bar k)$ (page \pageref{sys:ch4:IntersectionOffsetLine}). This elimination process leads us to consider the following system of equations:
\begin{equation}\label{sys:ch4:AuxiliaryCurvesSystem}
\hspace{-5mm}\begin{minipage}{14cm}
\[\hspace{3mm}
{{\mathfrak S}^P_3}(d,\bar k)\equiv\begin{cases}
s_1(d,\bar k, \bar t):=h(\bar t)(k_2P_3-k_3P_2)^2-d^2P_0(\bar t)^2(k_2n_3-k_3n_2)^2=0,
\\
s_2(d,\bar k, \bar t):=h(\bar t)(k_1P_3-k_3P_1)^2-d^2P_0(\bar t)^2(k_1n_3-k_3n_1)^2=0,
\\
s_3(d,\bar k, \bar t):=h(\bar t)(k_1P_2-k_2P_1)^2-d^2P_0(\bar t)^2(k_1n_2-k_2n_1)^2=0.
\end{cases}
\]
\end{minipage}
\end{equation}
We will refer to this as the {\sf Affine Auxiliary System}.

\nin We recall that $P=\left(\dfrac{P_1}{P_0},\dfrac{P_2}{P_0},\dfrac{P_3}{P_0}\right)$, $\bar k=(k_1,k_2,k_3)$, $\bar n=(n_1,n_2,n_3)$ and $h(\bar t)=n_1(t)^2+n_2(t)^2+n_3(t)^2$.
Along with the polynomials $s_1, s_2, s_3$ introduced in the above system, we will also need to consider the following polynomial:
\[s_0(\bar k, \bar t)=k_1(P_2n_3-P_3n_2)-k_2(P_1n_3-P_3n_1) +k_3(P_1n_2-P_2n_1)\]
The geometrical meaning of $s_0$ is clear when one expresses it as a determinant, as follows:
\begin{equation}\label{eq:ch4:GeometricalInterpretationOfS0}
s_0(\bar k, \bar t)=\det\left(
\begin{array}{ccc}
k_1&k_2&k_3\\
P_1&P_2&P_3\\
n_1&n_2&n_3
\end{array}
\right).
\end{equation}
We will introduce some additional notation to simplify the expression of the polynomials $s_i$ for $i=1,2,3$. More precisely, we {\sf denote:}
\begin{equation}\label{def:ch4:NotationMandGforCurvesSi}
\hspace{-5mm}\begin{minipage}{14cm}
\[\hspace{3mm}
\left\{\begin{array}{lll}
M_1(\bar k,\bar t)=k_2P_3-k_3P_2,&G_1(\bar k,\bar t)=k_2n_3-k_3n_2,\\
M_2(\bar k,\bar t)=k_3P_1-k_1P_3,&G_2(\bar k,\bar t)=k_3n_1-k_1n_3,\\
M_3(\bar k,\bar t)=k_1P_2-k_2P_1,&G_3(\bar k,\bar t)=k_1n_2-k_2n_1.\\
\end{array}\right.
\]
\end{minipage}
\end{equation}
With this notation one has
\[s_i(d,\bar k,\bar t)=h(\bar t)M_i^2(\bar k,\bar t)-d^2P_0(\bar t)^2G_i^2(\bar k,\bar t)\mbox{ for }i=1,2,3.\]
Note also that
\begin{equation}\label{eq:ch4:MandGasVectorProducts}
\begin{cases}
(M_1,M_2,M_3)(\bar k,\bar t)=\bar k\wedge\bigl(P_1(\bar t),P_2(\bar t),P_3(\bar t)\bigr)\\
(G_1,G_2,G_3)(\bar k,\bar t)=\bar k\wedge\bar n(\bar t).
\end{cases}
\end{equation}
\nin Let
 \[
I^P_ 2(d)=<b^{P},\normal^{P}_{(1,2)},\normal^{P}_{(1,3)}\normal^{P}_{(2,3)},w^{P},\ell_1,\ell_2,\ell_3>
\subset \C[d,\bar k,l,r,\bar t,\bar x]
 \]
be the ideal generated by the polynomials that define the Parametric Offset-Line System $\mathfrak S^P_2(d,\bar k)$.  We consider the projection associated with the elimination:
\[\pi_{(2,1)}:\C\times\C^3\times\C\times\C\times\C^2\times\C^3\mapsto\C\times\C^3\times\C^2\]
given by
\[\pi_{(2,1)}(d,\bar k,l,r,\bar t,\bar x)=(d,\bar k,\bar t)\]

\nin The next lemma relates the polynomials $s_0,\ldots,s_3\in\C[d,\bar k,\bar t]$ in System ${\mathfrak S_3}(d,\bar k)$ with the elimination process. We {\sf denote} by $\tilde I^P_ 2(d)$ the elimination ideal $I^P_ 2(d)\cap\C[d,\bar k,\bar t]$. For $(d^o,\bar k^o)\in\C\times\C^3$, the set of solutions of the Parametric Offset-Line system is denoted by $\Psi^P_2(d^o,\bar k^o)$, and the set of solutions of ${\mathfrak S_3^P}(d^o,\bar k^o)$ is denoted by $\Psi^P_3(d^o,\bar k^o)$\label{def:ch4:SolutionSet3}. Note that $\Psi^P_2(d^o,\bar k^o)=\bV(I^P_ 2(d))$.
\vspace{3mm}
\begin{Lemma}\label{lem:ch4:AuxiliaryPolynomialsBelongToEliminationIdeal}
$s_i\in\tilde I^P_ 2(d)$ for $i=0,\ldots,3$. In particular, if $(l^o,r^o,\bar t^o,\bar x^o)\in\Psi^P_2(d^o,\bar k^o)$, then $\bar t^o\in{\Psi_3^P}(d^o,\bar k^o)$.
\end{Lemma}
\begin{proof}
The polynomials $s_i$ can be expressed as follows:
\[
s_i=c^{(i)}_1\,b^{P}+c^{(i)}_2\,\normal^{P}_{(1,2)}+c^{(i)}_3\,\normal^{P}_{(1,3)}+c^{(i)}_4\,\normal^{P}_{(2,3)}
+c^{(i)}_5\,w^{P}+c^{(i)}_6\,\ell_1+c^{(i)}_7\,\ell_2+c^{(i)}_8\,\ell_3
\]
where $c^{(i)}_j\in\C[d,\bar k,l,r,\bar t,\bar x]$ for $i=0,\ldots,3$,  $j=1,\ldots,8$. This polynomials (obtained with the CAS Singular \cite{SingularWeb}) can be found in Appendix \ref{Ap3-ComplementsToSomeProofs} (page \pageref{Ap3-ComplementsToSomeProofs}).
\end{proof}

\nin The next step appears naturally to be the converse analysis: which are the $\bar t^o\in{\mathfrak S_3^P}(d^o,\bar k^o)$ that can be extended to a solution $(l^o,r^o,\bar t^o,\bar x^o)\in\Psi^P_2(d^o,\bar k^o)$?  In order to describe them, we need some notation and a lemma. Let $\cA$ {\sf denote} the set of values $\bar t^o\in\C^2$ such that:
\begin{equation}\label{def:ch4:Set_A_OfExtensibleSolutions}
\left\{\begin{array}{c}
P_0(\bar t^o)h(\bar t^o)\bigl(P_2(\bar t^o)n_3(\bar t^o)-P_3(\bar t^o)n_2(\bar t^o)\bigr)\neq 0\\
\mbox{ or }\\
P_0(\bar t^o)h(\bar t^o)\bigl(P_1(\bar t^o)n_3(\bar t^o)-P_3(\bar t^o)n_1(\bar t^o)\bigr)\neq 0\\
\mbox{ or }\\
P_0(\bar t^o)h(\bar t^o)\bigl(P_1(\bar t^o)n_2(\bar t^o)-P_2(\bar t^o)n_1(\bar t^o)\bigr)\neq 0
\end{array}\right.
\end{equation}
\nin Now we can describe which solutions of $\bar t^o\in{\mathfrak S_3^P}(d^o,\bar k^o)$ can be extended.
\vspace{3mm}
\begin{Proposition}\label{prop:ch4:ExtendableSolutions}
Let $\Omega_0$ be as in Theorem \ref{thm:ch4:TheoreticalFoundation} (page \pageref{thm:ch4:TheoreticalFoundation}), $(d^o,\bar k^o)\in\Omega_0$ and $\bar t^o\in{\Psi_3^P}(d^o,\bar k^o)$.  Then the following holds:
\begin{enumerate}
 \item[(a)] There exists $\lambda^o\in\C^\times$ such that:
 \[\bar k^o\wedge\bigl(P_1(\bar t^o),P_2(\bar t^o),P_3(\bar t^o)\bigr)=\lambda^o\left(\bar k^o\wedge\bar n(\bar t^o)\right).\]
 That is,
 \[M_i(\bar k^o,\bar t^o)=\lambda^o G_i(\bar k^o,\bar t^o)\mbox{ for }i=1,2,3.\]
 \item[(b)] If $\bar t^o\in\cA$, then  $s_0(d^o,k^o,\bar t^o)=0$.
 \item[(c)] $(d^o,\bar k^o,\bar t^o)\in\pi_{(2,1)}({\Psi_2^P}(d^o,\bar k^o))$ if and only if $\bar t^o\in\cA$.
\end{enumerate}
In particular,
 \[m\delta=\#\left(\cA\cap{\Psi_3^P}(d^o,\bar k^o)\right).\]
Recall that $m$ is the tracing index of $P$, and $\delta$ is the total degree w.r.t $\bar x$ of the generic offset equation.
\end{Proposition}
\quad\\
\begin{proof}
\begin{enumerate}
 \item[]
 \item[(a)] To prove the existence of $\lambda^o$, notice that $\bar t^o\in{\Psi_3^P}(d^o,\bar k^o)$ implies:
  \[h(\bar t^o)M_i^2(\bar k^o,\bar t^o)=(d^o)^2P_0(\bar t^o)^2G_i^2(\bar k^o,\bar t^o)\mbox{ for }i=1,2,3.\]
  Since $\bar t^o\in\cA$, $h(\bar t^o)\neq 0$. Therefore one concludes that there exist $\epsilon_i$, with $\epsilon_i^2=1$, such that
  \[M_i(\bar k^o,\bar t^o)=\epsilon_i\dfrac{d^oP_0(\bar t^o)}{\sqrt{h(\bar t^o)}}G_i(\bar k^o,\bar t^o)\mbox{ for }i=1,2,3.\]
   Since there are three of them, two of the $\epsilon_i$ must coincide. We will show that the third one must coincide as well. That is, we will show that either $\epsilon_1=\epsilon_2=\epsilon_3=1$, or $\epsilon_1=\epsilon_2=\epsilon_3=-1$ holds. We will study one particular case, the other possible combinations can be treated similarly. Let us suppose, e.g., that $\epsilon_1=\epsilon_2=1$. Then:
  \[
  \begin{cases}
  k_2^oP_3(\bar t^o)-k_3^oP_2(\bar t^o)=\dfrac{d^oP_0(\bar t^o)}{\sqrt{h(\bar t^o)}}\left(k_2^on_3(\bar t^o)-k_3^on_2(\bar t^o)\right)\\[3mm]
  k_3^oP_1(\bar t^o)-k_1^oP_3(\bar t^o)=\dfrac{d^oP_0(\bar t^o)}{\sqrt{h(\bar t^o)}}\left(k_3^on_1(\bar t^o)-k_1^on_3(\bar t^o)\right)
  \end{cases}\]
  Multiplying the first equation by $k_1^o$ and the second by $k_2^o$, and subtracting one has:
  \[k_3^o\left(k_1^oP_2(\bar t^o)-k_2^oP_1(\bar t^o)\right)=k_3^o\dfrac{d^oP_0(\bar t^o)}{\sqrt{h(\bar t^o)}}\left(k_1^on_2(\bar t^o)-k_2^on_1(\bar t^o)\right)\]
  Since $(d^o,\bar k^o)\in\Omega_0$, we have $k_3^o\neq 0$ (see Theorem \ref{thm:ch4:TheoreticalFoundation}(c), page \pageref{thm:ch4:TheoreticalFoundation}). Thus, we have shown that
\[k_1^oP_2(\bar t^o)-k_2^oP_1(\bar t^o)=\dfrac{d^oP_0(\bar t^o)}{\sqrt{h(\bar t^o)}}\left(k_1^on_2(\bar t^o)-k_2^on_1(\bar t^o)\right)\]
and so $\epsilon_3=1$. Therefore, $\lambda^o=\dfrac{d^oP_0(\bar t^o)}{\sqrt{h(\bar t^o)}}$, and it is non-zero because $\bar t^o\in\cA$.
 \item[(b)] From the identity in (a) it follows immediately that $\bar k^o$, $(P_1(\bar t^o),P_2(\bar t^o),P_3(\bar t^o)\bigr)$ and $\bar n(\bar t^o)$ are coplanar vectors. Thus, $s_0(d^o,\bar k^o,\bar t^o)=0$ (recall the geometric interpretation of $s_0$ in equation \ref{eq:ch4:GeometricalInterpretationOfS0}, page \pageref{eq:ch4:GeometricalInterpretationOfS0}).
 \item[(c)] If $(d^o,\bar k^o,\bar t^o)\in\pi_{(2,1)}({\Psi_2^P}(d^o,\bar k^o))$, then $P_0(\bar t^o)\beta(\bar t^o)h(\bar t^o)\neq 0$ follows from equation $w^{P}$ in the Parametric Offset-Line System \ref{sys:ch4:IntersectionOffsetLine} (page \pageref{sys:ch4:IntersectionOffsetLine}). Besides,
\[
\bigl(P_2n_3-P_3n_2\bigr)(\bar t^o)=\bigl(P_1n_3-P_3n_1\bigr)(\bar t^o)=\bigl(P_1n_2-P_2n_1\bigr)(\bar t^o)=0
\]
is impossible because of Theorem \ref{thm:ch4:TheoreticalFoundation}(a) (page \pageref{thm:ch4:TheoreticalFoundation}).  Thus $\bar t^o\in\cA$.\\
Conversely, let us suppose that $\bar t^o\in\cA$. More precisely, let us suppose w.l.o.g. that
\[P_0(\bar t^o)h(\bar t^o)\bigl(P_2(\bar t^o)n_3(\bar t^o)-P_3(\bar t^o)n_2(\bar t^o)\bigr)\neq 0.\]
The other cases can be proved in a similar way. First we note that
\[G_1(\bar k^o,\bar t^o)=k_3^on_2(\bar t^o)-k_2^on_3(\bar t^o)\neq 0.\]
since, using that $s_1(d^o,\bar k^o,\bar t^o)=0$ and $h(\bar t^o)\neq 0$, one has that
\[k_2^oP_3(\bar t^o)-k_3^oP_2(\bar t^o)=0.\]
Then, from the system:
\[
\begin{cases}
k_2^on_3(\bar t^o)-k_3^on_2(\bar t^o)=0\\
k_2^oP_3(\bar t^o)-k_3^oP_2(\bar t^o)=0
\end{cases}
\]
and the fact that $k_2^ok_3^o\neq 0$ (again, this is Theorem \ref{thm:ch4:TheoreticalFoundation}(c)), one deduces that
\[P_2(\bar t^o)n_3(\bar t^o)-P_3(\bar t^o)n_2(\bar t^o)=0,\]
that is a contradiction. Thus, we can define
\[r^o=\dfrac{1}{P_0(\bar t^o)\beta(\bar t^o)h(\bar t^o)}, \mbox{ and }
l^o=\dfrac{P_3(\bar t^o)n_2(\bar t^o)-P_2(\bar t^o)n_3(\bar t^o)}{-P_0(\bar t^o)G_1(\bar k^o,\bar t^o)}\]
We also define $\bar x^o=l^o\bar k^o$. We claim that $(l^o,r^o,\bar t^o,\bar x^o)\in\Psi^P_2(d^o,\bar k^o)$, and therefore $(d^o,\bar k^o,\bar t^o)\in\pi_{(2,1)}({\Psi_2^P}(d^o,\bar k^o))$. To prove our claim we substitute $(l^o,r^o,\bar t^o,\bar x^o)$ in the equations  of the Parametric Offset-Line System \ref{sys:ch4:IntersectionOffsetLine} (page \pageref{sys:ch4:IntersectionOffsetLine}), and we check that all of them vanish. The vanishing of $w^P(r^o,\bar t^o)$ and $\ell_i(\bar k^o,l^o,\bar x^o)$ for $i=1,2,3$ is a trivial consequence of the definitions. Substitution in $\normal^P_{(2,3)}$ leads to a polynomial whose numerator vanishes immediately. Substituting in $\normal^P_{(1,2)}$ (resp. in $\normal^P_{(1,3)}$) one obtains:
\[
\normal^{P}_{(1,2)}(\bar t^o,\bar x^o)=\dfrac{n_2(\bar t^o)s_0(\bar k^o,\bar t^o)}{n_2\bar t^ok_3^o-n_3(\bar t^o)k_2^o}=0,
\]
(respectively
\[
\normal^{P}_{(1,3)}(\bar t^o,\bar x^o)=\dfrac{n_3(\bar t^o)s_0(\bar k^o,\bar t^o)}{n_2\bar t^ok_3^o-n_3(\bar t^o)k_2^o}=0),
\]
where both equations hold because of part (a). Finally, substituting in $b^P(d^o,\bar t^o,\bar x^o)$ one has:
\begin{equation}\label{eq:ch4:bSubstitutedWithElevationFormulae}
b^P(d^o,\bar t^o,\bar x^o)=\dfrac{s_2(d^o,\bar k^o,\bar t^o)+\phi_1(\bar k^o,\bar t^o)s_0(\bar k^o,\bar t^o)}{(n_2\bar t^ok_3^o-n_3(\bar t^o)k_2^o)^2}=0
\end{equation}
with $\phi_1(\bar k,\bar t)=k_2 n_1P_3 + k_2 n_3 P_1 - k_3 n_1P_2 - k_3 n_2 P_1 - k_1 n_3 P_2 + k_1 n_2 P_3$.
Equation \ref{eq:ch4:bSubstitutedWithElevationFormulae} holds because of part (a) and because $s_2(d^o,\bar k^o,\bar t^o)=0$.
\end{enumerate}
\nin The claim that
 \[m\delta=\#\left(\cA\cap{\Psi_3^P}(d^o,\bar k^o)\right)\]
follows easily from Theorem \ref{thm:ch4:TheoreticalFoundation}(b) (page \pageref{thm:ch4:TheoreticalFoundation}) and the above result (c). This shows that, for $(d^o,\bar k^o)\in\Omega_0$, there is a bijection (under $\pi_{(2,1)}$) between   the points of ${\Psi_2^P}(d^o,\bar k^o)$ and the points in $\cA\cap{\Psi_3^P}(d^o,\bar k^o)$.
This finishes the proof of the proposition.
\end{proof}
\vspace{3mm}
\begin{Remark}\label{rem:ch4:SignOfLambdaAndOffsetting}
In the proof of Proposition \ref{prop:ch4:ExtendableSolutions} (page \ref{prop:ch4:ExtendableSolutions}) we have seen that there is a vector equality:
 \[
\bar M(\bar k^o,\bar t^o)=\epsilon\dfrac{d^oP_0(\bar t^o)}{\sqrt{h(\bar t^o)}}\bar G(\bar k^o,\bar t^o)\mbox{ for }i=1,2,3.
\]
where $\bar M=(M_1,M_2,M_3)$, $\bar G=(G_1,G_2,G_3)$ and $\epsilon=\pm 1$. In the next lemma we will see that the value of  $\epsilon=1$ determines the sign that
appears in the offsetting construction. More precisely, in the proof of Proposition \ref{prop:ch4:ExtendableSolutions} we have seen that if $\bar y^o=P(\bar t^o)$, and
\[P_0(\bar t^o)h(\bar t^o)\bigl(P_2(\bar t^o)n_3(\bar t^o)-P_3(\bar t^o)n_2(\bar t^o)\bigr)\neq 0.\]
then it holds that
\[k_2^on_3(\bar t^o)-k_3^on_2(\bar t^o)\neq0 \mbox{ and }k_2^oP_3(\bar t^o)-k_3^oP_2(\bar t^o)\neq 0.\]
Furthermore, the point $\bar x^o$, constructed as follows
\begin{equation}\label{eq:ch4:LiftingFormulae}
\bar x^o=\dfrac{P_3(\bar t^o)n_2(\bar t^o)-P_2(\bar t^o)n_3(\bar t^o)}{-P_0(\bar t^o)G_1(\bar k^o,\bar t^o)}\bar k^o,
\end{equation}
 is the point in $\cO_{d^o}(\Sigma)\cap\cL_{\bar k^o}$ associated with $\bar y^o$. Thus, one has:
\[
\bar x^o=\bar y^o+\epsilon'\dfrac{d^o\nabla f(\bar y^o)}{\sqrt{h_{\operatorname{imp}(\bar y^o)}}}.
\]
where $\epsilon'=\pm 1$.
\end{Remark}
\vspace{3mm}
\begin{Lemma}\label{lem:ch4:SignOfLambdaAndOffsetting}
With the notation of Remark \ref{rem:ch4:SignOfLambdaAndOffsetting}, it holds that $\epsilon=\epsilon'$.
\end{Lemma}
\begin{proof}
From the Equations
\[
M_2(\bar k^o,\bar t^o)=\epsilon\dfrac{d^oP_0(\bar t^o)}{\sqrt{h(\bar t^o)}}G_2(\bar k^o,\bar t^o)
\mbox{ and }
M_3(\bar k^o,\bar t^o)=\epsilon\dfrac{d^oP_0(\bar t^o)}{\sqrt{h(\bar t^o)}}G_3(\bar k^o,\bar t^o),
\]
multiplying the first equation by $n_2(\bar t^o)$, the second by $n_3(\bar t^o)$ and adding the results, one has:
\[
-G_1(\bar k^o,\bar t^o)P_1(\bar t^o)-k^o_1(P_3n_2-P_2n_3)(\bar t^o)=\epsilon\, n_1(\bar t^o)\dfrac{d^oP_0(\bar t^o)}{\sqrt{h(\bar t^o)}}G_1(\bar k^o,\bar t^o)
\]
Using Equation \ref{eq:ch4:LiftingFormulae} in Remark \ref{rem:ch4:SignOfLambdaAndOffsetting}, this is:
\[
-G_1(\bar k^o,\bar t^o)P_1(\bar t^o)+x^o_1G_1(\bar k^o,\bar t^o)P_0(\bar t^o)=\epsilon\, n_1(\bar t^o)\dfrac{d^oP_0(\bar t^o)}{\sqrt{h(\bar t^o)}}G_1(\bar k^o,\bar t^o).
\]
Dividing by $G_1(\bar k^o,\bar t^o)P_0(\bar t^o)$:
\[
-\dfrac{P_1(\bar t^o)}{P_0(\bar t^o)}+x^o_1=\epsilon\,\dfrac{d^on_1(\bar t^o)}{\sqrt{h(\bar t^o)}},
\]
and finally
\[
x^o_1=\dfrac{P_1(\bar t^o)}{P_0(\bar t^o)}+\epsilon\, \dfrac{d^on_1(\bar t^o)}{\sqrt{h(\bar t^o)}}.
\]
Similar results are obtained for $x^o_2$ and $x^o_3$. Thus we have proved that $\epsilon'=\epsilon$.
\end{proof}

\subsection{Fake points}\label{subsec:ch4:FakePoints}

Using Proposition \ref{prop:ch4:ExtendableSolutions} (page \pageref{prop:ch4:ExtendableSolutions}) we can now define the set of fake points associated with this problem.
\vspace{3mm}
\begin{Definition}\label{def:ch4:FakePoints}
A point $\bar t^o\in\C^2$ is a {\sf fake point} if
\[
\left\{\begin{array}{c}
P_0(\bar t^o)h(\bar t^o)\bigl(P_2(\bar t^o)n_3(\bar t^o)-P_3(\bar t^o)n_2(\bar t^o)\bigr)=0\\
\mbox{ and }\\
P_0(\bar t^o)h(\bar t^o)\bigl(P_1(\bar t^o)n_3(\bar t^o)-P_3(\bar t^o)n_1(\bar t^o)\bigr)=0\\
\mbox{ and }\\
P_0(\bar t^o)h(\bar t^o)\bigl(P_1(\bar t^o)n_2(\bar t^o)-P_2(\bar t^o)n_1(\bar t^o)\bigr)=0
\end{array}\right.
\]
Equivalently,
\begin{equation}\label{eq:ch4:EquationForFakePoints}
P_0(\bar t^o)h(\bar t^o)=0\mbox{ or }(P_1(\bar t^o),P_2(\bar t^o),P_3(\bar t^o))\wedge\bar n(\bar t^o)=\bar 0
\end{equation}
The set of fake points will be {\sf denoted} by $\cF$.
\end{Definition}
\vspace{3mm}
\begin{Definition}\label{def:ch4:InvariantSolutionsSystemS3P}
Let $\Omega_0$ be as in Theorem \ref{thm:ch4:TheoreticalFoundation} (page \pageref{thm:ch4:TheoreticalFoundation}) and let $\Omega$ be any open subset of $\Omega_0$. {\sf The set of invariant solutions of ${{\mathfrak S}^P_3}(d,\bar k)$ w.r.t $\Omega$.} is defined as the set:
\[
\cI^P_3(\Omega)=\bigcap_{(d^o,\bar k^o)\in\Omega}{\Psi_3^P}(d^o,\bar k^o)
\]
\end{Definition}
\vspace{3mm}
\begin{Remark}
Note that if $\bar t^o\in\cF$, we do not assume that $\bar t^o\in{\Psi_3^P}(d^o,\bar k^o)$ for some $(d^o,\bar k^o)\in\C\times\C^3$.
\end{Remark}
\vspace{3mm}
 \nin We have introduced the fake points starting from the notion non-extendable solutions of $\mfS^P_3(d^o,\bar k^o)$. Another point of view is to define fake points as the {\em invariant solutions} of $\mfS^P_3(d,\bar k)$. First we will define what we mean by invariant in this context, and then we will show that, in a certain open subset of values $(d,\bar k)$, both notions actually coincide.

\nin To prove the equivalence between the notions of fake points and invariant points we need to further restrict the set of values of $(d,\bar k)$. The following lemma gives the required restrictions.
\vspace{3mm}
\begin{Lemma}\label{lem:ch4:ExcludingAdditional_dk_AfterTheoreticalFoundation}
Let $\Omega_0$ be the open set in Theorem \ref{thm:ch4:TheoreticalFoundation}. There exists an open non-empty $\Omega_1\subset\Omega_0$ such that if $(d^o,\bar k^o)\in\Omega_1$, then
\begin{enumerate}
 \item[(1)] $\bar k^o$ is not isotropic.
 \item[(2)] $d^o$ is not a critical distance of $\Sigma$ (see Corollary \ref{cor:ch1:PositiveDimensionComponentsOfBotVareInSpheres} in page \pageref{cor:ch1:PositiveDimensionComponentsOfBotVareInSpheres}).
 \item[(3)] The system
 \begin{equation}\label{sys:ch4:SystemOfMiPlusP0}
 \{P_0(\bar t)=M_1(\bar k^o,1,\bar t)=M_2(\bar k^o,1,\bar t)=M_3(\bar k^o,1,\bar t)=0\}
 \end{equation}
has no solutions unless $P_0(\bar t^o)=P_1(\bar t^o)=P_2(\bar t^o)=P_3(\bar t^o)=0$.
\end{enumerate}
\end{Lemma}
\begin{proof}
 \begin{enumerate}
 \item[]
 \item[(1)] Set $\Omega_1^1=\Omega_0\cap(\C\times\mfQ)$, where $\mfQ=\{\bar k^o/(k_1^o)^2+(k_2^o)^2+(k_3^o)^2=0\}$ is the cone of isotropy in $\bar k$.
 \item[(2)] Let $\Upsilon(\Sigma^{\bot})$ is the set of critical distances of $\Sigma$ (defined in page \pageref{cor:ch1:PositiveDimensionComponentsOfBotVareInSpheres}), and set  $\Omega_1^2=\Omega_1^1\cap(\Upsilon(\Sigma^{\bot})\times\C^3)$.
 \item[(3)] First we will show that the set of values $\bar k^o\neq\bar 0$ for which the System \ref{sys:ch4:SystemOfMiPlusP0} has a solution is contained in an at most two-dimensional closed subset $\mfR\subset\C^3$.  If $P_0$ is constant the result is trivial. Assuming that $P_0$ is not constant, let $\cC_0$ be the affine curve defined by $P_0$, and let $\cC_1$, $\cC_2$, $\cC_3$ be the  varieties defined by $P_1, P_2, P_3$ respectively.
Let $\cJ_{P_1}\subset\C[\bar k,\bar t,v]$ be the ideal defined as follows:
\[\cJ_{P_1}=<P_0,k_2P_3-k_3P_2,k_1P_3-k_3P_1,k_1P_2-k_2P_1,v P_1-1>,\]
and let $\bV(\cJ_{P_1})\subset\C^3\times\C^2\times\C$ be the solution set of this ideal. Consider the projections defined by:
\[
\begin{cases}
\pi_1(\bar k,\bar t,v)=\bar k\\
\pi_2(\bar k,\bar t,v)=\bar t
\end{cases}
\]
Let $A_0$ be an irreducible component of $\bV(\cJ_{P_1})$, and let $(\bar k^o,\bar t^o,v^o)\in A_0$. Then the points in $\pi_2^{-1}(\pi_2(\bar k^o,\bar t^o,v^o))$ are the solutions of the following system:
\[
\begin{cases}
\bar t=\bar t^o,\\
M_1(\bar k,1,\bar t^o)=M_2(\bar k,1,\bar t^o)=M_3(\bar k,1,\bar t^o)=0,\\
 v P_1(\bar t^o)-1=0
\end{cases}
\]
The dimension of the set of solutions is $1$. On the other hand, $\pi_2(A_0)\subset\cC_0$ implies that $\dim(\pi_2(A_0))\leq 1$. Thus, using Lemma \ref{lem:ch1:FiberDimension} (page \pageref{lem:ch1:FiberDimension}), one has that $\dim(\bV(\cJ_{P_1}))\leq 2$. Thus, $\dim(\pi_1\left(\bV(\cJ_{P_1})\right)^*)\leq 2$. Now, defining $\cJ_{P_2}$ and $\cJ_{P_3}$ in a similar way (that is, replacing the equation $v P_1(\bar t^o)-1=0$ by $v P_2(\bar t^o)-1=0$ and $v P_3(\bar t^o)-1=0$ respectively), we set:
\[\mfR=\bigcup_{i=1,2,3}\pi_1\left(\bV(\cJ_{P_i})\right)^*\]
Now let
$\Omega_1^3=\Omega_1^2\cap(\C\times\mfR)$.
\end{enumerate}
\nin The above construction shows that $\Omega_1=\Omega_1^3$ satisfies the required properties.
\end{proof}
	
\nin Now we can prove the announced equivalence between the notions of fake points and invariant points.
\vspace{3mm}
\begin{Proposition}\label{prop:ch4:FakePointsAndInvariantSolutionsCoincide} Let $\Omega_1$ be as in Lemma \ref{lem:ch4:ExcludingAdditional_dk_AfterTheoreticalFoundation} (page \pageref{lem:ch4:ExcludingAdditional_dk_AfterTheoreticalFoundation}). If $\Omega$ is a non-empty open subset of $\Omega_1$, then it holds that:
\[\cI^P_3(\Omega)=\displaystyle\cF\cap\left(\bigcup_{(d^o,\bar k^o)\in\Omega}{\Psi_3^P}(d^o,\bar k^o)\right).\]
\end{Proposition}
\nin {\em Proof.}\quad
Let $\bar t^o\in\cI^P_3(\Omega)$. Then $s_i(d^o,\bar k^o,\bar t^o)=0$ for $i=1,2,3$ and any $(d^o,\bar k^o)\in\Omega$.  Thus $\bar t^o\in\cup_{(d^o,\bar k^o)\in\Omega}{\Psi_3^P}(d^o,\bar k^o)$. Furthermore, considering $s_i$ as polynomials in $\C[\bar t][d,\bar k]$, it follows that $\bar t^o$ must be a solution of:
\[
\begin{cases}
h(\bar t)P_1(\bar t)=h(\bar t)P_2(\bar t)=h(\bar t)P_3(\bar t)=0\\
P_0(\bar t)n_1(\bar t)=P_0(\bar t)n_2(\bar t)=P_0(\bar t)n_3(\bar t)=0
\end{cases}
\]
It follows that $h(\bar t^o)P_0(\bar t^o)=0$, and so $\bar t^o\in\cF$. In fact, if we suppose $h(\bar t^o)P_0(\bar t^o)\neq 0$, then from $P_0(\bar t^o)\neq 0$ one gets $\bar n(\bar t^o)=0$, and so $h(\bar t^o)=0$, a contradiction.\\
Conversely, let  $\bar t^o\in\cF\cap\left(\bigcup_{(d^o,\bar k^o)\in\Omega}{\Psi_3^P}(d^o,\bar k^o)\right)$. Then:
\begin{enumerate}
 \item If $P_0(\bar t^o)=h(\bar t^o)=0$, then $s_i(d,\bar k,\bar t^o)=0$ identically in $(d,\bar k)$ for $i=1,2,3$, and so $\bar t^o\in\cI^P_3(\Omega)$.
 \item If $P_0(\bar t^o)\neq 0$ and $h(\bar t^o)=0$, then since $\bar t^o\in{\Psi_3^P}(d^o,\bar k^o)$ for some $(d^o,\bar k^o)\in\Omega$, one has the following two possibilities:
 \begin{enumerate}
 \item $\bar n(\bar t^o)$ is isotropic and parallel to  $\bar k^o$. This is impossible because of the construction of $\Omega_1$ (see Lemma \ref{lem:ch4:ExcludingAdditional_dk_AfterTheoreticalFoundation}(1), page \pageref{lem:ch4:ExcludingAdditional_dk_AfterTheoreticalFoundation}).
 \item $\bar n(\bar t^o)=\bar 0$. In this case, again $s_i(d,\bar k,\bar t^o)=0$ identically in $(d,\bar k)$ for $i=1,2,3$, and so $\bar t^o\in\cI^P_3(\Omega)$.
\end{enumerate}
 \item Let us suppose that $P_0(\bar t^o)=0$ and $h(\bar t^o)\neq 0$. Then, since $\bar t^o\in{\Psi_3^P}(d^o,\bar k^o)$ for some $(d^o,\bar k^o)\in\Omega$, one has that $\bar t^o$ is a solution of:
 \[P_0(\bar t)=0, \quad M_1(\bar t,\bar k^o)=M_2(\bar t,\bar k^o)=M_3(\bar t,\bar k^o)=0,\]
Thus, two cases are possible:
\begin{enumerate}
 \item $(P_1,P_2,P_3)(\bar t^o)=\bar 0$. In this case, $s_i(d,\bar k,\bar t^o)=0$ identically in $(d,\bar k)$ for $i=1,2,3$, and so $\bar t^o\in\cI^P_3(\Omega)$.
 \item $(P_1,P_2,P_3)(\bar t^o)$ is non-zero. This contradicts the construction of $\Omega_1$ in Lemma \ref{lem:ch4:ExcludingAdditional_dk_AfterTheoreticalFoundation}(3).
\end{enumerate}
\item Finally, let us suppose that $P_0(\bar t^o)h(\bar t^o)\neq 0$. Then
it follows that the point $P(\bar t^o)$ is well defined, and it belongs to $\Sigma_{\bot}^*$ (recall that $\Sigma_{\bot}$ was introduced in Lemma \ref{lem:Ch1:LineMeetsVarietyNormallyInProperClosedSet}, page \pageref{lem:Ch1:LineMeetsVarietyNormallyInProperClosedSet}). Thus $d^o$ would be one of the critical distances, and this contradicts the construction of $\Omega_0$ in Lemma \ref{lem:ch4:ExcludingAdditional_dk_AfterTheoreticalFoundation}(2). \qed
\end{enumerate}

\section{Total Degree Formula for Parametric Surfaces}\label{sec:ch4:TotalDegreeFormulaForRationalSurfaces}

According to Proposition \ref{prop:ch4:ExtendableSolutions} (page \pageref{prop:ch4:ExtendableSolutions}), if $(d^o,\bar k^o)\in\Omega_0$ it holds that
 \[m\delta=\#\left(\cA\cap{\Psi_3^P}(d^o,\bar k^o)\right).\]
Recall that $m$ is the tracing index of $P$, and $\delta$ is the total degree w.r.t $\bar x$ of the generic offset equation. Moreover, $\cA$ was introduced in Equation \ref{def:ch4:Set_A_OfExtensibleSolutions} (page \pageref{def:ch4:Set_A_OfExtensibleSolutions}), and ${\Psi_3^P}(d^o,\bar k^o)$ was also introduced in page \pageref{def:ch4:SolutionSet3}, as the solution set of System \ref{sys:ch4:AuxiliaryCurvesSystem} (page \pageref{sys:ch4:AuxiliaryCurvesSystem}). In this section, we will derive a formula for the total degree $\delta$, using the tools in Section \ref{sec:ch1:IntersectionCurvesResultants} (page \pageref{sec:ch1:IntersectionCurvesResultants}) to analyze the intersection   $\cA\cap{\Psi_3^P}(d^o,\bar k^o).$

\nin In order to do this:
\begin{itemize}
    \item in Subsection \ref{subsec:ProjectivizationParametrizationSurface} we will consider the projective closure  of the auxiliary curves introduced in the preceding section. This in turn, requires as a first step the projectivization of the parametrization $P$.  At the end of the subsection we introduce the Projective Auxiliary System \ref{sys:ch4:AuxiliaryCurvesSystem-ProjectiveAndPrimitive} (page \pageref{sys:ch4:AuxiliaryCurvesSystem-ProjectiveAndPrimitive}), which will play a key r\^ole in the degree formula.
    \item Subsection \ref{subsec:ch4:InvariantSolutionsOfS5}. (page \pageref{subsec:ch4:InvariantSolutionsOfS5}) is devoted to the study of the invariant solutions of the Projective Auxiliary System, connecting them with the corresponding affine notions in Section \ref{sec:ch4:AuxiliaryCurvesForRationalSurfaces}. A crucial step in our strategy concerns the multiplicity of intersection of the auxiliary curves at their non-invariant points of intersection.
    \item In Subsection \ref{subsec:ch4:MultiplicityOfIntersectionAtNon-FakePoints} (page \pageref{subsec:ch4:MultiplicityOfIntersectionAtNon-FakePoints}) we will prove (in Proposition \ref{prop:ch4:MultiplicityAtNonFakePoints}, page \pageref{prop:ch4:MultiplicityAtNonFakePoints}) that the value of that multiplicity of intersection has the required property for the use of generalized resultants (according to Lemma \ref{lem:ch1:GeneralizedResultants}, page \pageref{lem:ch1:GeneralizedResultants}) .
    \item After this is done, everything is ready for the proof of the degree formula, which is the topic  of Subsection \ref{subsec:ch4:DegreeFormula} (page \pageref{subsec:ch4:DegreeFormula}). The formula appears in Theorem \ref{thm:ch4:DegreeFormula} (page \pageref{thm:ch4:DegreeFormula}).
\end{itemize}

\subsection{Projectivization of the parametrization and auxiliary curves}
\label{subsec:ProjectivizationParametrizationSurface}

Let ${P}$ be the parametrization of $\Sigma$, introduced in Equation (\ref{eq:ch4:SurfaceAffineParametrization}). If we homogenize the components of $P$ w.r.t. a new variable $t_0$, multiplying both the numerators and denominators if necessary by a suitable power of $t_0$ we arrive at an expression of the form:
\begin{equation}\label{eq:ch4:SurfaceProjectiveParametrization}
P_h(\bar t_h)=
\left(
\dfrac{X(\bar t_h)}{W(\bar t_h)},
\dfrac{Y(\bar t_h)}{W(\bar t_h)},
\dfrac{Z(\bar t_h)}{W(\bar t_h)}
\right)
\end{equation}
where $\bar t_h=(t_0:t_1:t_2)$, and $X, Y, Z, W\in\C[\bar t_h]$ are homogeneous polynomials of the same degree $d_{P}$, for which $\gcd(X,Y,Z,W)=1$ holds. This $P_h$ will be called the {\sf projectivization} of ${P}$.
\vspace{3mm}
\begin{Remark}
Note that those projective values of $\bar t_h$ of the form $(0:a:b)$ correspond to points at infinity in the parameter plane.
\end{Remark}

\nin In Section \ref{subsec:ch4:SurfaceParametrizationsAndtheirAssociatedNormalVector} (page \pageref{subsec:ch4:SurfaceParametrizationsAndtheirAssociatedNormalVector}) we defined $\bar n=(n_1,n_2,n_3)$, the associated  normal vector to $P$. A similar construction, applied to $P_h$, leads to a normal vector $N=(N_1,N_2,N_3)$, where $N_i$ are homogeneous polynomials in $\bar t_h$ of the same degree, such that $\gcd(N_1,N_2,N_3)=1$. This vector $N$ will be called the {\sf associated homogeneous normal vector} to $P_h$.
\nin The homogeneous polynomial $H$ defined by
\[H(\bar t_h)=(N_1(\bar t))^2+(N_2(\bar t))^2+(N_3(\bar t))^2\]
is the {\sf parametric projective normal-hodograph\index{hodograph, projective parametric}\index{normal-hodograph, projective parametric}} of the parametrization $P_h$.
\vspace{3mm}
\begin{Remark}\label{rem:ch4:OneComponentof_N_isNotDivisibleByt0}
The polynomials $N_i$ are, up to multiplication by a power of $t_0$, the homogenization of the components of $\bar n$ w.r.t. $t_0$. However, since $\gcd(N_1,N_2,N_3)=1$, at least one of the components $N_i$ ($i=1,2,3$) is not divisible by $t_0$. Besides, note that if two components $N_i, N_j$, with $i\neq j$, are divisible by $t_0$, then $H$ is not.
\end{Remark}
\vspace{3mm}
\begin{Lemma}\label{lem:SurfaceParametrizationProperties2}
\begin{enumerate}
 \item[]
 \item If $W$ does not depend on $t_0$, then at least one of the polynomials $X, Y, Z$ must depend on $t_0$.
 \item If $W$ does not depend on $t_0$, and there is exactly one of the polynomials
 $X,Y,Z$ depending on $t_0$, then the surface is a cylinder with its axis parallel to the direction of the component with numerator depending on $t_0$.
\end{enumerate}
\end{Lemma}

\noindent{\em Proof}\vspace{-3mm}
\begin{enumerate}
 \item Otherwise, the rank of the jacobian matrix of $P$ would be less than two. To see this, let us suppose that $X,Y,Z,W$ depend only on $t_1,t_2$.
 Let $\partial_i P_h$ be the vector obtained as the partial derivative of $P_h$ w.r.t. $t_i$, that is;
 \[
 \partial_i P_h=\left(\frac{X_i W-XW_i}{W^2},\frac{Y_iW -Y W_i}{W^2},\frac{Z_iW-Z W_i}{W^2}\right)
 \]
 where $X_i,Y_i,Z_i,W_i$ denotes the partial derivative of $X,Y,Z,W$ w.r.t. $t_i$. Using Euler's formula, and taking into account that the polynomials $X,Y,Z,W$ have the same degree $n$, one has that $t_1 \partial_1 P_h=-t_2 \partial_2 P_h$. Substituting $t_0=1$, we see that the rank of the jacobian of $P$ would be less than 2.
\item Assume w.l.o.g. that $X,Y$ do not depend on $t_0$, but $Z$ does. The rational map
\[
\phi(\bar t)=\left(\dfrac{X(\bar t)}{W(\bar t)},\dfrac{Y(\bar t)}{W(\bar t)}\right)
\]
has rank one, because $X, Y, W$ are homogeneous polynomials in $\bar t$ of the same degree. Thus, $\phi$ parametrizes a curve $\cC$ in the $(y_1,y_2)$-plane. Let $\mbox{Cyl}(\cC)$ be the cylinder over $\cC$ with axis parallel to the $y_3$-axis. The points of the form $(\phi(\bar t^o),y_3^o)$, with  $W(\bar t^o)\neq 0$, are dense in $\mbox{Cyl}(\cC)$. Given one of these points, let $t^o_0$ be any solution of the equation (in $t_0$):
\[Z(\bar t^o,t_0)=y^o_3W(\bar t^o)\]
Then we have
\[P_h(\bar t^o,t^o_0)=(\phi(\bar t^o),y_3^o)\]
and so $P_h(\P^2)$ is dense in $\mbox{Cyl}(\cC)$. \hspace{8.5cm}$\Box$
\end{enumerate}

\nin Now we are ready to introduce the projective auxiliary polynomials. We consider the following system:
\begin{equation}\label{sys:ch4:AuxiliaryCurvesSystem-ProjectiveNoPrimitive}
\hspace{-5mm}\begin{minipage}{14cm}
\[\hspace{3mm}
{{\mathfrak S}^{P_h}_4}(d,\bar k)\equiv\begin{cases}
S_0(\bar k, \bar t_h):=k_1(YN_3-ZN_2)-k_2(XN_3-ZN_1) +k_3(XN_2-YN_1)
\\
S_1(d,\bar k, \bar t_h):=H(\bar t_h)(k_2Z-k_3Y)^2-d^2W(\bar t_h)^2(k_2N_3-k_3N_2)^2
\\
S_2(d,\bar k, \bar t_h):=H(\bar t_h)(k_1Z-k_3X)^2-d^2W(\bar t_h)^2(k_1N_3-k_3N_1)^2
\\
S_3(d,\bar k, \bar t_h):=H(\bar t_h)(k_1Y-k_2X)^2-d^2W(\bar t_h)^2(k_1N_2-k_2N_1)^2
\end{cases}
\]
\end{minipage}
\end{equation}

\nin As usual, for $(d^o,\bar k^o)\in\C^4$, we {\sf denote} by ${\Psi^{P_h}_4}(d^o,\bar k^o)$ the set of projective solutions of ${{\mathfrak S}^{P_h}_4}(d^o,\bar k^o)$. Our next goal is the analysis of the relation between ${\Psi^{P_h}_4}(d^o,\bar k^o)$ and $\Psi^P_3(d^o,\bar k^o)$  (the set of solutions of ${\mathfrak S_3^P}(d^o,\bar k^o)$, see Subsection \ref{subsec:ch4:EliminationAndAuxiliaryPolynomials}, page \pageref{subsec:ch4:EliminationAndAuxiliaryPolynomials}). In particular, an in order to obtain the degree formula, we will characterize those points in ${\Psi^{P_h}_4}(d^o,\bar k^o)$ that correspond to the points in $\cA\cap{\Psi_3^P}(d^o,\bar k^o)$.  In Proposition \ref{prop:ch4:FakePointsAndInvariantSolutionsCoincide} (page \pageref{prop:ch4:FakePointsAndInvariantSolutionsCoincide}) we have seen that the invariant solutions of ${\Psi_3^P}(d,\bar k)$ correspond to fake points. Thus, as a first step, we will characterize certain invariant solutions of ${{\mathfrak S}^{P_h}_4}(d,\bar k)$.
\vspace{3mm}
\begin{Lemma}\label{lem:ch4:ContentProjectiveAuxiliaryCurves}
Let $S=c_1S_1+c_2S_2+c_3S_3$. Then:
\[\Con_{(d,\bar k)}(S(\bar c, d,\bar k,\bar t_h))=\gcd(H,W^2).\]
\end{Lemma}
\begin{proof}
Since $S=c_1S_1+c_2S_2+c_3S_3$, one has:
\[\Con_{(d,\bar k)}(S)=
\gcd\left(\Con_{(d,\bar k)}(S_1),\Con_{(d,\bar k)}(S_2),\Con_{(d,\bar k)}(S_3)\right).\]
Now, considering $S_i$ for $i=1,2,3$ as polynomials in $\C[\bar t_h][d,\bar k]$ one has:
\[\Con_{(d,\bar k)}(S_1)=\gcd(HZ^2,HZY,HY^2,W^2N_2^2,W^2N_2N_3,W^2N_3^2).\]
That is,
\[\Con_{(d,\bar k)}(S_1)=\gcd(H\gcd(Y,Z)^2,W^2\gcd(N_2,N_3)).\]
Similarly,
\[\Con_{(d,\bar k)}(S_2)=\gcd(H\gcd(X,Z)^2,W^2\gcd(N_1,N_3)).\]
and
\[\Con_{(d,\bar k)}(S_3)=\gcd(H\gcd(X,Y)^2,W^2\gcd(N_1,N_2))\]
Taking into account that $\gcd(N_1,N_2,N_3)=1$ and $\gcd(X,Y,Z,W)=1$, one has
\[
\Con_{(d,\bar k)}(S)=
\gcd(H\gcd(X,Y,Z)^2,W^2)=\gcd(H,W^2).
\]
\end{proof}

\nin In order to use the above results, and to state the degree formula, we need to introduce some additional notation. We {\sf denote} by:
\begin{equation}\label{eq:ch4:ContentProjectiveAuxiliaryCurves}
Q_0(\bar t_h)=\Con_{\bar k}(S_0(\bar k,\bar t_h))\quad\mbox{ and }\quad
Q(\bar t_h)=\Con_{(d,\bar k)}(S(\bar c, d,\bar k,\bar t_h)).
\end{equation}
Observe that, by Lemma \ref{lem:ch4:ContentProjectiveAuxiliaryCurves}, $Q$ does not depend on $\bar c$, a fact that is reflected in our notation. Furthermore, note that:
\[Q_0(\bar t_h)=\gcd(YN_3-ZN_2,XN_3-ZN_1,XN_2-YN_1)\]
and
\[Q(\bar t_h)=\gcd(H,W^2).\]
We also denote by:
\begin{equation}\label{eq:ch4:PolynomialsTildeH_andTildeW}
\tilde H(\bar t_h)=\dfrac{H(\bar t_h)}{Q(\bar t_h)},\quad\quad
\tilde W(\bar t_h)=\dfrac{W^2(\bar t_h)}{Q(\bar t_h)},
\end{equation}
and
\begin{equation}\label{eq:ch4:PolynomialsU}
\begin{cases}
U_1(\bar t_h)=\dfrac{(YN_3-ZN_2)(\bar t_h)}{Q_0(\bar t_h)},\\[3mm]
U_2(\bar t_h)=\dfrac{(ZN_1-XN_3)(\bar t_h)}{Q_0(\bar t_h)},\\[3mm]
U_3(\bar t_h)=\dfrac{(XN_2-YN_1)(\bar t_h)}{Q_0(\bar t_h)}.
\end{cases}
\end{equation}
Thus, one has:
\[S_0(\bar k,\bar t_h)=Q_0(\bar t_h)(k_1U_1(\bar t_h)+k_2U_2(\bar t_h)+k_3U_3(\bar t_h)).\]
We denote as well:
\begin{equation}\label{eq:ch4:Polynomials_MhGh}
\begin{cases}
M_{h,1}(\bar k,\bar t_h)=k_2Z(\bar t_h)-k_3Y(\bar t_h),&G_{h,1}(\bar k,\bar t_h)=k_2N_3(\bar t_h)-k_3N_2(\bar t_h),\\[3mm]
M_{h,2}(\bar k,\bar t_h)=k_3X(\bar t_h)-k_1Z(\bar t_h),&G_{h,2}(\bar k,\bar t_h)=k_3N_1(\bar t_h)-k_1N_3(\bar t_h),\\[3mm]
M_{h,3}(\bar k,\bar t_h)=k_1Y(\bar t_h)-k_2X(\bar t_h),&G_{h,3}(\bar k,\bar t_h)=k_1N_2(\bar t_h)-k_2N_1(\bar t_h).
\end{cases}
\end{equation}
and so, for $i=1,2,3$,
\[S_i(d,\bar k,\bar t_h)=
Q(\bar t_h)\left({\tilde H}(\bar t_h)M_{h,i}^2(\bar k,\bar t_h)-d^2{\tilde W}(\bar t_h)G_{h,i}^2(\bar k,\bar t_h)\right)\]
We denote:
\begin{equation}\label{eq:ch4:ReducedAuxiliaryPolynomials}
T_0(\bar k,\bar t_h)=\dfrac{S_0(\bar k,\bar t_h)}{Q_0(\bar t_h)}
\end{equation}
and, for $i=1,2,3$,
\begin{equation}\label{eq:ch4:ReducedAuxiliaryPolynomials_2}
T_i(d,\bar k,\bar t_h)=\dfrac{S_i(d,\bar k,\bar t_h)}{Q(\bar t_h)},
\end{equation}
Finally, we denote:
\begin{equation}\label{eq:ch4:ReducedAuxiliaryPolynomials_3}
T(\bar c,d,\bar k,\bar t_h)=\dfrac{S(\bar c,d,\bar k,\bar t_h)}{Q(\bar t_h)}.
\end{equation}
Note that:
\[T(\bar c,d,\bar k,\bar t_h)=c_1T_1(d,\bar k,\bar t_h)+c_2T_2(d,\bar k,\bar t_h)+c_3T_3(d,\bar k,\bar t_h).\]

\nin With this notation we can introduce the system of equations that will play the central role in the degree formula:
\begin{equation}\label{sys:ch4:AuxiliaryCurvesSystem-ProjectiveAndPrimitive}
\hspace{-5mm}\begin{minipage}{14cm}
\[\hspace{3mm}
{{\mathfrak S}^{P_h}_5}(d,\bar k)\equiv\begin{cases}
{T_0}(\bar k, \bar t_h)=k_1U_1(\bar t_h)+k_2U_2(\bar t_h)+k_3U_3(\bar t_h)=0
\\[3mm]
{T_i}(d,\bar k, \bar t_h)={\tilde H}(\bar t_h)M_{h,i}^2(\bar k,\bar t_h)-d^2{\tilde W}(\bar t_h)G_{h,i}^2(\bar k,\bar t_h),
\\
\mbox{for }i=1,2,3.
\end{cases}
\]
\end{minipage}
\end{equation}
We will refer to this as the {\sf Projective Auxiliary System}.

\subsection{Invariant solutions of the projective auxiliary system}
\label{subsec:ch4:InvariantSolutionsOfS5}

\nin In passing from ${{\mathfrak S}^P_3}(d,\bar k)$ to ${{\mathfrak S}^{P_h}_4}(d,\bar k)$, and then to ${{\mathfrak S}^{P_h}_5}(d,\bar k)$, we have introduced additional solutions at infinity, in the space of parameters (that is, with $t_0=0$). The following results will show that, in a certain open subset of values of $(d,\bar k)$, these solutions at infinity are invariant w.r.t. $(d,\bar k)$. We start with some technical lemmas.
\vspace{3mm}
\begin{Lemma}\label{lem:ch4:U_iAndS_iNotIdenticallyZero}
There is always $i^o\in\{1,2,3\}$ such that $U_{i^o}(0,t_1,t_2)$ and $T_{i^o}(d,\bar k,0,t_1,t_2)$ are both not identically zero.
\end{Lemma}
\begin{proof}
First, let us prove that there are always $i,j\in\{1,2,3\}$ such that $t_0$ does not divide $T_i$ and $T_j$. Suppose, on the contrary that, for example $T_1(d,\bar k,0,t_1,t_2)\equiv 0$ and $T_2(d,\bar k,0,t_1,t_2)\equiv 0$. Considering $T_1$ and $T_2$ as polynomials in $\C[\bar t][d,\bar k]$, if $t_0$ divides $T_1$ and $T_2$ one concludes that $t_0$ must divide
\[{\tilde H}X, {\tilde H}Y, {\tilde H}Z, {\tilde W}N_1, {\tilde W}N_2\mbox{ and }{\tilde W}N_3.\]
If one assumes that $t_0$ divides ${\tilde W}$, then it does not divide ${\tilde H}$, because $\gcd({\tilde H},{\tilde W})=1$. Thus it divides $X$, $Y$ and $Z$. But this is again a contradiction, since $\gcd(X,Y,Z,W)=1$, and ${\tilde W}$ divides $W$. Thus, $t_0$ does not divide ${\tilde W}$. Then it must divide $N_1, N_2, N_3$. This is also a contradiction, since $\gcd(N_1,N_2,N_3)=1$. Therefore we can assume w.l.o.g. that e.g. $t_0$ does not divide $T_1$ and $T_2$. To finish the proof in this case we need to show that, if $t_0$ divides $T_3$, then it does not divide at least one of $U_1$ and $U_2$. The hypothesis that $t_0$ divides $T_3$ implies that it divides
\[{\tilde H}X, {\tilde H}Y, {\tilde W}N_1\mbox{ and }{\tilde W}N_2.\]
If $t_0$ divides ${\tilde W}$, again, it must divide $X$ and $Y$. Thus it does not divide $Z$.
Now, observe that $XU_1+YU_2+ZU_3=0$. Therefore, one concludes that $t_0$ divides $U_3$. Thus, $t_0$ does not divide at least one of $U_1$ and $U_2$, since $\gcd(U_1,U_2,U_3)=1$. If $t_0$ does not divide ${\tilde W}$, then it divides $N_1$ and $N_2$. Observing that $N_1U_1+N_2U_2+N_3U_3=0$, we again conclude that $t_0$ does not divide at least one of $U_1$ and $U_2$, since $\gcd(U_1,U_2,U_3)=1$.
\end{proof}
\vspace{3mm}
\begin{Lemma}\label{lem:ch4:FactoringT0}
Let $i^o\in\{1,2,3\}$ be such that $T_{i^o}(d,\bar k,0,t_1,t_2)$ and $U_{i^o}(0,t_1,t_2)$ are both not identically zero (see Lemma \ref{lem:ch4:U_iAndS_iNotIdenticallyZero}). Then
\[\gcd(T_0(\bar k,0,t_1,t_2),T_{i^o}(d,\bar k,0,t_1,t_2))\]
does not depend on $\bar k$ (it certainly does not depend on $d$).
\end{Lemma}
\begin{proof}
The claim follows observing that $T_0(\bar k,0,t_1,t_2)$ depends linearly on $k_{i^o}$, and $T_{i^o}(d,\bar k,0,t_1,t_2)$ does not.
\end{proof}
\nin In order to describe what we mean when we say that a solution is invariant w.r.t. $(d,\bar k)$, we make the following definition (recall Definition \ref{def:ch4:InvariantSolutionsSystemS3P}, page \pageref{def:ch4:InvariantSolutionsSystemS3P}):
\vspace{3mm}
\begin{Definition}\label{def:ch4:InvariantSolutionsSystemS5P}
Let $\Omega_1$ be as in Lemma \ref{lem:ch4:ExcludingAdditional_dk_AfterTheoreticalFoundation} (page \pageref{lem:ch4:ExcludingAdditional_dk_AfterTheoreticalFoundation}), and let $\Omega$ be a non-empty open subset of $\Omega_1$. {\sf The set of invariant solutions of ${{\mathfrak S}^{P_h}_5}(d,\bar k)$ w.r.t $\Omega$.} is defined as the set:
\[
\cI^{P_h}_5(\Omega)=\bigcap_{(d^o,\bar k^o)\in\Omega}{\Psi^{P_h}_5}(d^o,\bar k^o)
\]
\end{Definition}
\vspace{3mm}
\begin{Remark}\label{rem:ch4:DescriptionInvariantSolutionsS5}
\begin{enumerate}
 \item[]
 \item Considering $T_i$ (for $i=0,\ldots,3$) as polynomials in $\C[\bar t_h][d,\bar k]$, it is easy to see that $\cI^{P_h}_5(\Omega)$ is the set of solutions of:
 \begin{equation}\label{eq:ch4:EquationsInvariantSolutionsS5}
\begin{cases}
U_1(\bar t_h)=U_2(\bar t_h)=U_3(\bar t_h)=0,\\
(\tilde H\cdot X)(\bar t_h)=(\tilde H\cdot Y)(\bar t_h)=(\tilde H\cdot Z)(\bar t_h)=0,\\
(\tilde W\cdot N_1)(\bar t_h)=(\tilde W\cdot N_2)(\bar t_h)=(\tilde W\cdot N_3)(\bar t_h)=0.
\end{cases}
\end{equation}
 \item In particular, since $\gcd(U_1,U_2,U_3)=1$, the set  $\cI^{P_h}_5(\Omega)$ is always a finite set.
\end{enumerate}
\end{Remark}
\nin The following proposition shows that, restricting the values of $(d,\bar k)$ to a certain open set, we can ensure that all the solutions at infinity  of ${{\mathfrak S}^{P_h}_5}(d,\bar k)$ are invariant w.r.t. the particular choice of $(d,\bar k)$ in that open set.
\vspace{3mm}
\begin{Proposition}\label{prop:ch4:SystemS5HasOnlyInvariantSolutionsAtInfinity}
There exists an open non-empty subset $\Omega_2\subset\Omega_1$ (with $\Omega_1$ as in Lemma \ref{lem:ch4:ExcludingAdditional_dk_AfterTheoreticalFoundation}, page \pageref{lem:ch4:ExcludingAdditional_dk_AfterTheoreticalFoundation}), such that if $(d^o,\bar k^o)\in\Omega_2$, and $\bar t^o_h=(0:t^o_1:t^o_2)\in\Psi^{P_h}_5(d^o,\bar k^o)$,
then $\bar t^o_h\in\cI^{P_h}_5(\Omega_2)$.
\end{Proposition}
\begin{proof}
We know that $T_0(\bar k,0,\bar t)\not\equiv 0$. Suppose, in the first place, that $T_0(\bar k,0,\bar t)$ depends only in $\bar k$ and one of the variables $t_1, t_2$. Say, e.g., $T_0(\bar k,0,\bar t)=T_0^*(\bar k)t_1^{p}$ for some $p\in\N$.  This implies that, for any given $(d^o,\bar k^o)$ such that $T_0^*(\bar k^o)\neq 0$,  $(0:0:1)$ is the only possible point of $\Psi^{P_h}_5(d^o,\bar k^o)$ with $t_0=0$. Obviously,  if
\[T_i(d,\bar k,0,0,1)\equiv 0\mbox{ for }i=1,2,3,\]
then one may take $\Omega_2=\Omega_1\cap\{(d^o,\bar k^o)\in\C^4 / T_0^*(\bar k^o)\neq 0\}$, and the result is proved. On the other hand, if not all $T_i(d,\bar k,0,0,1)\equiv 0$, say w.l.o.g that
\[T_1(d,\bar k,0,0,1)\not\equiv 0\]
then we may take $\Omega_2=\Omega_1\cap\{(d^o,\bar k^o)/T_0^*(\bar k^o)T_1(d^o,\bar k^o,0,0,1)\neq 0\}$, and the result is proved.

\nin Thus, w.l.o.g. we can assume that $T_0(\bar k,0,\bar t)$ depends on both $t_1$ and $t_2$. Let $i^o\in\{1,2,3\}$ be such that $U_i(0,\bar t)\not\equiv 0$ and $T_i(d,\bar k,0,\bar t)\not\equiv 0$ (see Lemma \ref{lem:ch4:FactoringT0}, \pageref{lem:ch4:FactoringT0}). By
Lemma \ref{lem:ch4:U_iAndS_iNotIdenticallyZero} (page \pageref{lem:ch4:U_iAndS_iNotIdenticallyZero}) we know that this is the case at least for one value of $i^o$. Let us consider (see Lemma \ref{lem:ch4:FactoringT0}):
\[T^*_{i^o}(\bar t)=\gcd(T_0(\bar k,0,\bar t),T_{i^o}(d,\bar k,0,\bar t))\in\C[\bar t].\]
Note that $T^*_{i^o}$ is homogeneous in $\bar t$, and so, if it is not constant, it factors as:
\[T^*_{i^o}(\bar t)=\gamma\prod_{j=1}^p(\beta_jt_1-\alpha_jt_2)\]
for some $\gamma\in\C^{\times}$, and $(\alpha_j,\beta_j)\in\C^2, j=1,\ldots,p$.  For each point $(0:\alpha_j:\beta_j)$ we can repeat the construction that we did for $(0:0:1)$. Thus, one obtains a non-empty open set $\Omega_2^1\subset\Omega_1$ such that, if $(d^o,\bar k^o)\in\Omega_2^1$, and $T^*_{i^o}(\alpha_j,\beta_j)=0$, then either $(0:\alpha_j:\beta_j)\not\in\Psi^{P_h}_5(d^o,\bar k^o)$, or $(0:\alpha_j:\beta_j)\in\cI^{P_h}_5(\Omega_2^1)$.

\nin Let
\[
T'_0(\bar k,\bar t)=\dfrac{T_0(\bar k,0,\bar t)}{T^*_{i^o}(\bar t)},
\quad\mbox{ and }\quad
T'_{i^o}(d,\bar k,\bar t))=\dfrac{T_{i^o}(d,\bar k,0,\bar t)}{T^*_{i^o}(\bar t)}.
\]
Note that both $T'_0(\bar k,\bar t)$ and $T'_{i^o}(d,\bar k,\bar t)$ are homogeneous in $\bar t$, and by construction they have a trivial gcd. If we define:
\[
  \Gamma(d,\bar k,t_2)=
   \begin{cases}
   \Res_{t_1}(T'_0(\bar k,\bar t),T'_{i^o}(d,\bar k,\bar t))&\mbox{ if }\deg_{t_1}(T'_0(\bar k,\bar t))>0,\\[3mm]
   T'_0(\bar k,\bar t)&\mbox{ in other case.}
   \end{cases}
\]
Then $\Gamma$ is not identically zero, and since $T'_0$ and $T'_{i^o}$ are both homogeneous in $\bar t$, we have a factorization:
\[\Gamma(d,\bar k,t_2)=t_2^q\Gamma^*(d,\bar k)\]
for some $q\in\N$. Note also that, by construction, since $\gcd(T'_0,T'_{i^o})=1$, $t_2$ cannot divide both $T'_0$ and $T'_{i^o}$. In particular, since these polynomials are homogeneous in $\bar t$, one concludes that $\bar t^o=(1,0)$ is not a solution of
\[T'_0(\bar k,\bar t)=T'_{i^o}(d,\bar k,\bar t)=0.\]
We define
\[\Omega_2=\Omega^1_2\cap\{(d^o,\bar k^o)/\Gamma^*(d^o,\bar k^o)\neq 0\}\]
If $(d^o,\bar k^o)\in\Omega_2$, and $\bar t^o_h=(0:t^o_1:t^o_2)\in\Psi^{P_h}_5(d^o,\bar k^o)$, then either $T'_0(\bar k^o,\bar t^o)\neq 0$ or $T'_{i^o}(d^o,\bar k^o,\bar t^o)\neq 0$. In any case, one has $T^*_{i^o}(\bar t^o)=0$ (that is, $(0:t^o_1:t^o_2)=(0:\alpha_j:\beta_j)$ for some $j=1,\ldots,p$). The construction of $\Omega^1_2$ implies that  $(0:t^o_1:t^o_2)\in\cI^{P_h}_5(\Omega_2)$.
\end{proof}
\nin Let $\cA_h$ {\sf denote}\label{def:ch4:Set_Ah_OfProjectiveExtensibleSolutions} the set of values $\bar t^o_h\in\P^2$ such that (compare to the Definition of the set $\cA$ in Equation \ref{def:ch4:Set_A_OfExtensibleSolutions}, page \pageref{def:ch4:Set_A_OfExtensibleSolutions}):
\begin{equation}
\left\{\begin{array}{c}
t^o_0\,W(\bar t^o_h)H(\bar t^o_h)\bigl(Y(\bar t^o_h)N_3(\bar t^o_h)-Z(\bar t^o_h)N_2(\bar t^o_h)\bigr)\neq 0\\
\mbox{ or }\\
t^o_0\,W(\bar t^o_h)H(\bar t^o_h)\bigl(Z(\bar t^o_h)N_1(\bar t^o_h)-X(\bar t^o_h)N_3(\bar t^o_h)\bigr)\neq 0\\
\mbox{ or }\\
t^o_0\,W(\bar t^o_h)H(\bar t^o_h)\bigl(X(\bar t^o_h)N_1(\bar t^o_h)-Y(\bar t^o_h)N_1(\bar t^o_h)\bigr)\neq 0
\end{array}\right.
\end{equation}
Equivalently, $\cA_h$ consists in those points $\bar t^o_h\in\P^2$ such that:
\begin{equation}\label{eq:ch4:Set_Ah_OfProjectiveExtensibleSolutions}
t^o_0 W(\bar t^o_h)H(\bar t^o_h)\neq 0\mbox{ and }(X(\bar t^o_h),Y(\bar t^o_h),Z(\bar t^o_h))\wedge\bar N(\bar t^o_h)\neq\bar 0.
\end{equation}
We will see that the non-invariant solutions of $\Psi^{P_h}_5(d,\bar k)$ are points in $\cA_h$. Note that we are explicitly asking these points to be affine (recall Proposition \ref{prop:ch4:SystemS5HasOnlyInvariantSolutionsAtInfinity}, page \pageref{prop:ch4:SystemS5HasOnlyInvariantSolutionsAtInfinity})

\nin With this notation we are ready to state the main theorem about System ${{\mathfrak S}^{P_h}_5}(d,\bar k)$.
\vspace{3mm}
\begin{Theorem}\label{thm:ch4:RelationBetweenSystemS5AndSystemS4}
Let $\Omega_2$ be as in Proposition \ref{prop:ch4:SystemS5HasOnlyInvariantSolutionsAtInfinity}. If $(d^o,\bar k^o)\in\Omega_2$,
\[\bar t^o=(t^o_1,t^o_2)\in\cA\cap\Psi_3^P(d^o,\bar k^o)\Leftrightarrow \bar t^o_h=(1:t^o_1:t^o_2)\in\cA_h\cap\Psi_5^{P_h}(d^o,\bar k^o)\]
\end{Theorem}
\nin Recall that every point in $\cA_h$ is affine, see Equation \ref{def:ch4:Set_A_OfExtensibleSolutions} (page \pageref{def:ch4:Set_A_OfExtensibleSolutions}), for the Definition of $\cA$, and page \pageref{lem:ch4:AuxiliaryPolynomialsBelongToEliminationIdeal} for the definition of $\Psi_3^P(d^o,\bar k^o)$.
\begin{proof}
Let us prove that $\Rightarrow$ holds. If $\bar t^o\in\cA\cap\Psi_3^P(d^o,\bar k^o)$, then, e.g.,
\[P_0(\bar t^o)h(\bar t^o)(P_2n_3-P_3n_2)(\bar t^o)\neq 0.\]
Therefore,
\[W(\bar t^o_h)H(\bar t^o_h)(YN_3-ZN_2)(\bar t^o_h)\neq 0,\]
and so $\bar t^o_h\in\cA_h$. Besides, this last inequality implies that $Q_0(\bar t^o_h)Q(\bar t^o_h)\neq 0$. Since $\bar t^o\in\cA\cap\Psi_3^P(d^o,\bar k^o)$, by Proposition \ref{prop:ch4:ExtendableSolutions}(b) (page \pageref{prop:ch4:ExtendableSolutions}), one has that $(1:t^o_1:t^o_2)\in\Psi_4^{P_h}(d^o,\bar k^o)$. From Equations \ref{eq:ch4:ReducedAuxiliaryPolynomials} and \ref{eq:ch4:ReducedAuxiliaryPolynomials_2} (page \pageref{eq:ch4:ReducedAuxiliaryPolynomials}), and from $Q_0(\bar t^o_h)Q(\bar t^o_h)\neq 0$, one concludes that $\bar t^o_h=(1:t^o_1:t^o_2)\in\cA_h\cap\Psi_5^{P_h}(d^o,\bar k^o)$. Thus, $\Rightarrow$ is proved.\\
\nin The proof of $\Leftarrow$ is similar, simply reversing the implications.
\end{proof}
\vspace{3mm}
\begin{Remark}\label{rem:ch4:OffsetDegreeAsCardinal_Projective}
Let $\Omega_2$ be as in Proposition \ref{prop:ch4:SystemS5HasOnlyInvariantSolutionsAtInfinity}, and let $(d^o,\bar k^o)\in\Omega_2$. Then, Theorem \ref{thm:ch4:RelationBetweenSystemS5AndSystemS4} implies that:
 \[m\delta=\#\left(\cA\cap{\Psi_3^P}(d^o,\bar k^o)\right)=\#\left(\cA_h\cap{\Psi_5^{P_h}}(d^o,\bar k^o)\right).\]
\end{Remark}
\nin Theorem \ref{thm:ch4:RelationBetweenSystemS5AndSystemS4} establishes the link between $\cA\cap\Psi_3^P(d^o,\bar k^o)$ and $\cA_h\cap{\Psi_5^{P_h}}(d^o,\bar k^o)$ for a fixed $(d^o,\bar k^o)\in\Omega_2$. As we said before, the non-invariant solutions of $\Psi^{P_h}_5(d,\bar k)$ should be the points in $\cA_h$. As a first step, we have this result.
\vspace{3mm}
\begin{Proposition}\label{prop:ch4:NonInvariantSolutionsCanBeAvoided}
Let $\Omega_2$ be as in Proposition \ref{prop:ch4:SystemS5HasOnlyInvariantSolutionsAtInfinity} (page \pageref{prop:ch4:SystemS5HasOnlyInvariantSolutionsAtInfinity}).  If $\bar t^o_h\in\cA_h$, then $\bar t^o_h\not\in\cI^{P_h}_5(\Omega_2)$.
\end{Proposition}
\begin{proof}
Let us suppose that for every $(d^o,\bar k^o)\in\Omega_2$, one has $\bar t^o_h\in\left(\cA_h\cap{\Psi_5^{P_h}}(d^o,\bar k^o)\right)$. Then (recall Equation \ref{def:ch4:Set_A_OfExtensibleSolutions}, page \pageref{def:ch4:Set_A_OfExtensibleSolutions}), $\bar t^o_h$ is of the form $(1:t^o_1:t^o_2)$ and, by Theorem \ref{thm:ch4:RelationBetweenSystemS5AndSystemS4},
\[\bar t^o=(t^o_1,t^o_2)\in\bigcap_{(d^o,\bar k^o)\in\Omega_2}{\Psi_3^P}(d^o,\bar k^o)=\cI^P_3(\Omega_2).\]
Since $\Omega_2\subset\Omega_1$, Proposition \ref{prop:ch4:FakePointsAndInvariantSolutionsCoincide} (page \pageref{prop:ch4:FakePointsAndInvariantSolutionsCoincide}) applies, and one concludes that
$(t^o_1,t^o_2)\in\cF$. But then,
\[P_0(\bar t^o)h(\bar t^o)=0\mbox{ or }(P_1(\bar t^o),P_2(\bar t^o),P_3(\bar t^o))\wedge\bar n(\bar t^o)=\bar 0\]
This implies that:
\[W(\bar t^o_h)H(\bar t^o_h)=0\mbox{ or }(X(\bar t^o_h),Y(\bar t^o_h),Z(\bar t^o_h))\wedge\bar N(\bar t^o_h)=\bar 0,\]
contradicting $\bar t^o_h\in\cA_h$.
\end{proof}

\nin As we said before, we will prove the converse of this proposition. That is, if $\bar t^o_h\in\Psi_5^{P_h}(d^o,\bar k^o)$ for some $(d^o,\bar k^o)\in\Omega_2$, but $\bar t^o_h\not\in\cI^{P_h}_5(\Omega_2)$, then we would like to conclude that $\bar t^o_h\in\cA_h$. However, the open set $\Omega_2$ can be too large for this to hold. More precisely, the problem is caused by the solutions of $Q(\bar t_h)Q_0(\bar t_h)=0$ (when $Q$ and $Q_0$ are not equal $1$). These points do not belong to $\cA_h$. However, since $\cI_5^{P_h}(\Omega_2)$ is finite (see Remark \ref{rem:ch4:DescriptionInvariantSolutionsS5}, page \pageref{rem:ch4:DescriptionInvariantSolutionsS5}), most of the solutions $Q(\bar t_h)Q_0(\bar t_h)=0$  are not invariant.
Therefore, we need to impose some more restrictions in the values of $(d,\bar k)$. Note, however, that we have already dealt with the points at infinity; thus, we need only consider the affine solutions of $Q(\bar t_h)Q_0(\bar t_h)=0$. We do the necessary technical work in the following lemma. First, we introduce some notation for the affine versions of some polynomials. We {\sf denote:}
\[
\begin{cases}
q(\bar t)=Q(1,t_1,t_2), q_0(\bar t)=Q_0(1,t_1,t_2), \tilde w(\bar t)=\tilde W(1,t_1,t_2), \tilde h(\bar t)=\tilde H(1,t_1,t_2),\\
u_i(\bar t)=U_i(1,t_1,t_2)\mbox{ for }i=1,2,3.
\end{cases}
\]
and we consider a new auxiliary set of variables $\bar\rho=(\rho_1,\rho_2,\rho_3)$, in order to do the necessary Rabinowitsch's tricks.\\
\nin Let $\cG\subset\C\times\C^3\times\C^3\times\C^2$ be the set of solutions (in the variables $(d,\bar k,\bar\rho,\bar t)$) of this system of equations:
\begin{equation}\label{eq:ch4:SystemToExcludeZerosOfQQ}
\begin{cases}
q(\bar t)q_0(\bar t)=0\\
\tilde h(\bar t)M_i^2(\bar k,\bar t)-d^2\tilde w(\bar t)G_i^2(\bar k,\bar t)=0&\mbox{ for }i=1,2,3.\\
k_1u_1(\bar t)+k_2u_2(\bar t)+k_3u_3(\bar t)=0\\
\rho_1 \tilde w(\bar t)\tilde h(\bar t)-1=0\\
\prod_{i=1}^3(\rho_2 u_i(\bar t)-1)=0\\
\prod_{i=1}^3(\rho_3 n_i(\bar t)-1)=0
\end{cases}
\end{equation}
and consider the projection $\pi_1(d,\bar k,\lambda,\bar\rho,\bar t)=(d,\bar k)$.
\vspace{3mm}
\begin{Lemma} \label{lem:ch4:DimensionSetGlessThan4}
$\cG$ is empty or $\dim(\pi_1(\cG))\leq 3$.
\end{Lemma}
\begin{proof}
Let $\cG\neq\emptyset$. Then $q(\bar t)q_0(\bar t)$ is not constant. We will use Lemma \ref{lem:ch1:FiberDimension} (page \pageref{lem:ch1:FiberDimension}), to  prove that $\dim(\cG)\leq 3$. From this the result follows immediately.  Consider the projection $\pi_2(d,\bar k,\bar\rho,\bar t)=\bar t$. Clearly, $\pi_2(\cG)$ is contained in the affine curve defined by $q(\bar t)q_0(\bar t)=0$. Thus, $\dim(\pi_2(\cG))\leq 1$. Let $\bar t^o\in\pi_2(\cG)$. First, let us suppose that for all $(d^o,\bar k^o,\bar\rho^o,\bar t^o)\in\pi_2^{-1}(\bar t^o)$, one has $\bar k^o=\bar 0$. Then, if
$(d^o,\bar 0,\bar\rho^o,\bar t^o)\in\pi_2^{-1}(\bar t^o)$, $\rho_1^o, \rho_2^o$, and $\rho_3^o$ must be one of the finitely many solutions of the polynomial equations:
\[\rho_1 \tilde w(\bar t^o)\tilde h(\bar t^o)-1=0,\quad \prod_{i=1}^3(\rho_2 u_i(\bar t^o)-1)=0, \mbox{ and }\prod_{i=1}^3(\rho_3 n_i(\bar t^o)-1)=0.\]
The condition $\bar t^o\in\pi_2(\cG)$ implies that these equations can be solved. Note that, in this case, $(d^o,\bar 0,\bar\rho^o,\bar t^o)\in\pi_2^{-1}(\bar t^o)$ does not impose any condition on $d^o$. It follows that, in this case,
one has $\mu=\dim(\pi_2^{-1}(\bar t^o))=1$.

\nin Now, let us suppose that $(d^o,\bar k^o,\bar\rho^o,\bar t^o)\in\pi_2^{-1}(\bar t^o)$, with $\bar k^o\neq\bar 0$.
Then, by a similar argument to the proof of Proposition \ref{prop:ch4:ExtendableSolutions}(a) (page \pageref{prop:ch4:ExtendableSolutions}), and taking $\tilde w(\bar t^o)\tilde h(\bar t^o)\neq 0$ into account, there exists $\lambda^o\in\C^\times$ such that
 \[M_i(\bar k^o,\bar t^o)=\lambda^o G_i(\bar k^o,\bar t^o)\mbox{ for }i=1,2,3.\]
Thus, in this case $d^o$ must be a solution of:
\[\tilde h(\bar t^o)(\lambda^o)^2-(d^o)^2\tilde w(\bar t^o)=0.\]
Besides, there exists also $j^o\in\{1,2,3\}$ with $u_{j^o}(\bar  t^o)\neq 0$. Then, $\bar k^o$ must belong to the two-dimensional space defined by
\[k_1u_1(\bar t^o)+k_2u_2(\bar t^o)+k_3u_3(\bar t^o)=0.\]
Finally, $\rho_1^o, \rho_2^o$, and $\rho_3^o$ must be one of the finitely many solutions of the polynomial equations:
\[\rho_1 \tilde w(\bar t^o)\tilde h(\bar t^o)-1=0,\quad \prod_{i=1}^3(\rho_2 u_i(\bar t^o)-1)=0, \mbox{ and }\prod_{i=1}^3(\rho_3 n_i(\bar t^o)-1)=0.\]
The condition $\bar t^o\in\pi_2(\cG)$ implies that these equations can be solved. These remarks show that for every $\bar t^o\in\pi_2(\cG)$, one has $\mu=\dim(\pi_2^{-1}(\bar t^o))\leq 2$. Thus, using Lemma \ref{lem:ch1:FiberDimension}:
\[\dim(\cG)=\dim(\pi_2(\cG))+\mu\leq 1+2=3,\]
and the lemma is proved.
\end{proof}

\nin If $\bar t^o_h$ is such that
\[T_0(\bar k,\bar t^o_h)\equiv 0\mbox{ and }T_i(d,\bar k,\bar t^o_h)\equiv 0,\mbox{ for }i=1,2,3,\]
then $\bar t^o_h\in\cI^{P_h}_5(\Omega)$ for any choice of $\Omega$. However, if this is not the case, then sometimes we need to remove from $\Omega$ precisely those values $(d^o,\bar k^o)$ such that $\bar t^o_h\in\Psi^{P_h}_5(d^o,\bar k^o)$. In the proof of the following lemma we will need to do this several times. Thus we introduce the necessary notation.

\vspace{3mm}
\begin{Definition}\label{def:ch4:InvariantSetOfaPoint}
Let $\Omega\subset\C\times\C^3$ be non-empty and open. For $\bar t^o_h\in\P^2$ we define:
\[
\Omega^{inv}(\bar t^o_h)=
\begin{cases}
\Omega,\quad \mbox{ if }T_0(\bar k,\bar t^o_h)\equiv 0\mbox{ and }T_i(d,\bar k,\bar t^o_h)\equiv 0,\mbox{ for }i=1,2,3.\\[3mm]
\Omega\setminus\bigl\{(d^o,\bar k^o)/\bar t^o_h\in\Psi^{P_h}_5(d^o,\bar k^o)\bigr\}, \mbox{in other case.}
\end{cases}
\]
\end{Definition}
\vspace{3mm}
\begin{Remark}\label{rem:ch4:InvariantSetOfaPoint}
\begin{enumerate}
 \item[]
 \item[(1)] Note that if $\Omega\subset\Omega^{inv}(\bar t^o_h)$, then $\bar t^o_h\in\cI^{P_h}_5(\Omega)$.
 \item[(2)] Observe that $\Omega^{inv}(\bar t^o_h)\neq\emptyset$.
\end{enumerate}
\end{Remark}
\vspace{3mm}
\begin{Lemma}\label{lem:ch4:QQ_zero_impliesInvariant}
Let $\Omega_2$ be as in Proposition \ref{prop:ch4:SystemS5HasOnlyInvariantSolutionsAtInfinity}. There exists an open non-empty set $\Omega_3\subset\Omega_2$ such that the following hold:
\begin{enumerate}
 \item[(a)] If $\bar t^o_h\in\Psi_5^{P_h}(d^o,\bar k^o)$ for some $(d^o,\bar k^o)\in\Omega_3$, and $Q_0(\bar t^o_h)Q(\bar t^o_h)=0$, then $\bar t^o_h\in\cI^{P_h}_5(\Omega_3)$.
 \item[(b)] If $\bar t^o_h$ satisfies
 \[T_1(d,k,\bar t^o_h)=T_2(d,k,\bar t^o_h)=T_3(d,k,\bar t^o_h)=0\mbox{ identically in }(d,\bar k),\]
then $\bar t^o_h\in\cI^{P_h}_5(\Omega_3)$.
\end{enumerate}
\end{Lemma}
\begin{proof}
Let
\[A_0=\{\bar t^o_h\,|\,X(\bar t^o_h)=Y(\bar t^o_h)=Z(\bar t^o_h)=W(\bar t^o_h)=0\}.\]
Since $\gcd(X,Y,Z,W)=1$, one sees that $A_0$ is (empty or) a finite set.  Thus, if we define
(see Definition \ref{def:ch4:InvariantSetOfaPoint}):
\[\Omega^0_3=\Omega_2\cap\left(\bigcap_{\bar t^o_h\in A_0}\Omega^{inv}(\bar t^o_h)\right).\]
By Remark \ref{rem:ch4:InvariantSetOfaPoint}, $\Omega^0_3$ is an open non-empty set.
\nin Let
\[A_1=\{\bar t^o_h\,|\,N(\bar t^o_h)=\bar 0\},\]
where $N=(N_1,N_2,N_3)$. Recalling that $\gcd(N_1,N_2,N_3)=1$, $A_1$ is (empty or) a finite set.  We define:
\[\Omega^1_3=\Omega^0_3\cap\left(\bigcap_{\bar t^o_h\in A_1}\Omega^{inv}(\bar t^o_h)\right).\]
By Remark \ref{rem:ch4:InvariantSetOfaPoint}, $\Omega^1_3$ is an open non-empty set.\\
\nin Similarly,
since $\gcd(\tilde H,\tilde W)=1$,
the set
\[A_2=\{\bar t^o_h\,|\,\tilde H(\bar t^o_h)=\tilde W(\bar t^o_h)=0\}\]
is (empty or) finite. We define
\[\Omega^2_3=\Omega^1_3\cap\left(\bigcap_{\bar t^o_h\in A_2}\Omega^{inv}(\bar t^o_h)\right),\]
and $\Omega^2_3$ is an open non-empty set. Moreover, since $\gcd(U_1,U_2,U_3)=1$,
the set
\[A_3=\{\bar t^o_h\,|\,U_1(\bar t^o_h)=U_2(\bar t^o_h)=U_3(\bar t^o_h)=0\}\]
is (empty or) finite. We define
\[\Omega^3_3=\Omega^2_3\cap\left(\bigcap_{\bar t^o_h\in A_3}\Omega^{inv}(\bar t^o_h)\right),\]
and $\Omega^3_3$ is an open non-empty set. We define
\[\Omega^4_3=\Omega^3_3\setminus(\pi_1(\cG)^*),\]
where $\pi_1(\cG)$ is as in Lemma \ref{lem:ch4:DimensionSetGlessThan4} (page \pageref{lem:ch4:DimensionSetGlessThan4}), and, as usual, the asterisk denotes Zariski closure.

\nin Finally, since $T(\bar c,d,\bar k,\bar t)$ is primitive w.r.t. $(d,\bar k)$ (recall Equations \ref{eq:ch4:ContentProjectiveAuxiliaryCurves}, page \pageref{eq:ch4:ContentProjectiveAuxiliaryCurves}, and \ref{eq:ch4:ReducedAuxiliaryPolynomials_3}, page \pageref{eq:ch4:ReducedAuxiliaryPolynomials_3}), it follows that the set
\[A_4=\left\{\bar t^o_h\,|\,T_1(d,k,\bar t^o_h)=T_2(d,k,\bar t^o_h)=T_3(d,k,\bar t^o_h)=0\mbox{ identically in }(d,\bar k)\right\}\]
is (empty or) finite. We define
\[\Omega_3=\Omega^4_3\cap\left(\bigcap_{\bar t^o_h\in A_4}\Omega^{inv}(\bar t^o_h)\right).\]
Let us now suppose that $(d^o,\bar k^o)\in\Omega_3$ and $\bar t^o_h\in\Psi_5^{P_h}(d^o,\bar k^o)$, with $Q_0(\bar t^o_h)Q(\bar t^o_h)=0$. We will show that $\bar t^o_h\in\cI^{P_h}_5(\Omega_3)$. This will prove that statement (a) holds. If $\bar t^o_h$ is of the form $(0:t^o_1:t^o_2)$, by Proposition \ref{prop:ch4:SystemS5HasOnlyInvariantSolutionsAtInfinity}, $\bar t^o_h\in\cI^{P_h}_5(\Omega_3)$ holds trivially. Thus, in the rest of the proof we can assume w.l.o.g. that $\bar t^o_h$ is of the form $(1:t^o_1:t^o_2)$.

\nin If $\bar t^o_h\in\cup_{i=0,\ldots,3}A_i$, then we have $\Omega_3\subset\Omega^4_3\subset\Omega^{inv}(\bar t^o_h)$, and by Remark \ref{rem:ch4:InvariantSetOfaPoint}, $\bar t^o_h\in\cI^{P_h}_5(\Omega_3)$. So, let $\bar t^o_h\not\in\cup_{i=0,\ldots,3}A_i$. Then the following hold:
\begin{enumerate}
 \item[(0)] Since $\bar t^o_h\not\in A_0$, $P_i(\bar t^o)\neq 0$ for some $i=0,\ldots,3$ (recall that $\bar t^o_h=(1:\bar t^o)$).
 \item[(1)] Since $\bar t^o_h\not\in A_1$, $N(\bar t^o_h)\neq\bar 0$.
 \item[(2)] Since $\bar t^o_h\not\in A_2$, $\tilde H(\bar t^o_h)\neq 0$ or $\tilde W(\bar t^o_h)\neq 0$.
 \item[(3)] Since $\bar t^o_h\not\in A_3$, $U_i(\bar t^o_h)\neq 0$ for some $i=1,2,3$.
\end{enumerate}
 Let us show that (0) and (2) imply the following:
\begin{enumerate}
 \item[(4)] $\tilde H(\bar t^o_h)\tilde W(\bar t^o_h)\neq 0$.
\end{enumerate}
Indeed, if we suppose that $\tilde H(\bar t^o_h)=0$ but $\tilde W(\bar t^o_h)\neq 0$, then from $\bar t^o_h\in\Psi_5^{P_h}(d^o,\bar k^o)$ one concludes that $d^oG_i(\bar k^o,\bar t^o_h)=0$ for $i=1,2,3$. Since $d^o\neq 0$ and $\bar k^o\neq\bar 0$ in $\Omega_3$, one has that $N(\bar t^o_h)$ is isotropic and parallel to $\bar k^o$, contradicting Lemma \ref{lem:ch4:ExcludingAdditional_dk_AfterTheoreticalFoundation}(1) (page \pageref{lem:ch4:ExcludingAdditional_dk_AfterTheoreticalFoundation}). On the other hand, if we suppose $\tilde H(\bar t^o_h)\neq 0$ but $\tilde W(\bar t^o_h)=0$, then from $\bar t^o_h\in\Psi_5^{P_h}(d^o,\bar k^o)$ one concludes that
$d^oG_i(\bar k^o,\bar t^o_h)=0$ for $i=1,2,3$. Since $d^o\neq 0$, we conclude that $G_i(\bar k^o,\bar t^o_h)=0$ for $i=1,2,3$.  Thus, $\bar t^o$ is a solution of:
\[P_0(\bar t)=M_1(\bar k^o,1,\bar t)=M_2(\bar k^o,1,\bar t)=M_3(\bar k^o,1,\bar t)=0.\]
However, by (0), there exists $j^o\in\{0,1,2,3\}$ such that $P_j(\bar t^o)\neq 0$. Therefore, we get a contradiction with Lemma \ref{lem:ch4:ExcludingAdditional_dk_AfterTheoreticalFoundation}(3) (page \pageref{lem:ch4:ExcludingAdditional_dk_AfterTheoreticalFoundation}).

\nin From (1), (3), (4), and since $\bar t^o_h\in\Psi_5^{P_h}(d^o,\bar k^o)$ and  $Q_0(\bar t^o_h)Q(\bar t^o_h)=0$, it follows that $(d^o,\bar k^o,\bar t^o_h)$ can be extended to $(d^o,\bar k^o,\bar\rho^o,\bar t^o)\in\cG$. Thus, one has $(d^o,\bar k^o)\in\pi_1(\cG)$, contradicting the construction of $\Omega^3_3$.
This finishes the proof of statement (a).

\nin The proof of statement (b) is a consequence of the construction of $\Omega_3$ (in particular, see the construction of $A_4$); indeed, if $\bar t^o_h$ satisfies
 \[T_1(d,k,\bar t^o_h)=T_2(d,k,\bar t^o_h)=T_3(d,k,\bar t^o_h)=0\mbox{ identically in }(d,\bar k),\]
then
$\bar t^o_h\in A_4$. It follows that $\Omega_3\subset\Omega^{inv}(\bar t^o_h)$ and so $\bar t^o_h\in\cI^{P_h}_5(\Omega)$ (see Remark \ref{rem:ch4:InvariantSetOfaPoint}(1), page \pageref{rem:ch4:InvariantSetOfaPoint}).
\end{proof}
\nin Now, restricting the values of $(d,\bar k)$ to a new open set, we are ready to prove the announced converse of Proposition \ref{prop:ch4:NonInvariantSolutionsCanBeAvoided} (page \pageref{prop:ch4:NonInvariantSolutionsCanBeAvoided}).
\vspace{3mm}
\begin{Proposition}\label{prop:ch4:NonInvariantSolutionsCanBeAvoided_Part2}
Let $\Omega_{3}$ be as in Lemma \ref{lem:ch4:QQ_zero_impliesInvariant} (page \pageref{lem:ch4:QQ_zero_impliesInvariant}).
If $\bar t^o_h\in\Psi_5^{P_h}(d^o,\bar k^o)$ for some $(d^o,\bar k^o)\in\Omega_{3}$, but $\bar t^o_h\not\in\cI^{P_h}_5(\Omega_{3})$, then $\bar t^o_h\in\cA_h$.
\end{Proposition}
\begin{proof}
If $\bar t^o_h\not\in\cI^{P_h}_5(\Omega_3)$, then $t^o_0\neq 0$ (by Proposition \ref{prop:ch4:SystemS5HasOnlyInvariantSolutionsAtInfinity},  page \pageref{prop:ch4:SystemS5HasOnlyInvariantSolutionsAtInfinity}).
Let us write $\bar t^o_h=(1:t^o_1:t^o_2)$. Then, since $\bar t^o_h\in\Psi_5^{P_h}(d^o,\bar k^o)$, one has $\bar t^o\in\Psi_3^{P}(d^o,\bar k^o)$. Note also that, since $\bar t^o_h\not\in\cI^{P_h}_5(\Omega_{3})$, by Lemma \ref{lem:ch4:QQ_zero_impliesInvariant}, we must have $Q_0(\bar t^o_h)Q(\bar t^o_h)\neq 0$. If we suppose $\bar t^o_h\not\in\cA_h$, then $\bar t^o\not\in\cA$. Thus $\bar t^o\in\cF$, and by Proposition \ref{prop:ch4:FakePointsAndInvariantSolutionsCoincide} (page \pageref{prop:ch4:FakePointsAndInvariantSolutionsCoincide}), $\bar t^o\in\cI^P_3(\Omega_{3})$.  Taking Equation \ref{eq:ch4:ReducedAuxiliaryPolynomials_2} (page \pageref{eq:ch4:ReducedAuxiliaryPolynomials_3}) into account, and using $Q(\bar t^o_h)\neq 0$, we conclude that
\[T_1(d,k,\bar t^o_h)=T_2(d,k,\bar t^o_h)=T_3(d,k,\bar t^o_h)=0\mbox{ identically in }(d,\bar k).\]
Then. by Lemma \ref{lem:ch4:QQ_zero_impliesInvariant}(b) (page \pageref{lem:ch4:QQ_zero_impliesInvariant}), one has that $\bar t^o_h\in\cI^{P_h}_5(\Omega_{3})$. This is a contradiction, and so we obtain that $t^o_h\in\cA_h$.
\end{proof}

\subsection{Multiplicity of intersection at non-fake points}
\label{subsec:ch4:MultiplicityOfIntersectionAtNon-FakePoints}

The auxiliary polynomials $S_i$ (for $i=0,\ldots,3$) were introduced in Section \ref{sec:ch4:AuxiliaryCurvesForRationalSurfaces} (page \pageref{sec:ch4:AuxiliaryCurvesForRationalSurfaces}), in order to reduce the offset degree problem to a problem of intersection between planar curves. More precisely, the preceding results in this paper indicate that the offset degree problem can be reduced to an intersection problem between the planar curves defined by the auxiliary polynomials $T_i$. A crucial step in this reduction concerns the multiplicity of intersection of these curves at their non-invariant points of intersection. In this subsection we will prove that the value of that multiplicity of intersection is one (in Proposition \ref{prop:ch4:MultiplicityAtNonFakePoints}, page \pageref{prop:ch4:MultiplicityAtNonFakePoints}). We first introduce some notation for the curves involved in this problem.
\vspace{3mm}
\begin{Definition}\label{def:ch4:CurvesDefinedByAuxiliaryPolynomials}
Let $\Omega_0$ be as in Theorem \ref{thm:ch4:TheoreticalFoundation} (page \pageref{thm:ch4:TheoreticalFoundation}). If $(d^o,\bar k^o)\in\Omega_0$, and $T_i$ (for $i=0,\ldots,3$) are the polynomials introduced in Equations \ref{eq:ch4:ReducedAuxiliaryPolynomials} and \ref{eq:ch4:ReducedAuxiliaryPolynomials_2} (page \pageref{eq:ch4:ReducedAuxiliaryPolynomials}), we denote by $\cT^a_0(\bar k^o)$ (resp. $\cT^a_i(d^o,\bar k^o)$ for $i=1,2,3$) the affine algebraic set defined by the polynomial  $T_0(\bar k^o,1,\bar t)$ (resp. $T_i(d^o,\bar k^o,1,\bar t)$ for $i=1,2,3$).  Similarly, we denote by $\cT^h_0(\bar k^o)$ (resp. $\cT^h_i(d^o,\bar k^o)$ for $i=1,2,3$) the projective algebraic set defined by the polynomial $T_0(\bar k^o,\bar t_h)$ (resp. $T_i(d^o,\bar k^o,\bar t_h)$ for $i=1,2,3$).
\end{Definition}
\vspace{3mm}
\begin{Remark}
Note that the homogenization of the polynomials $T_0(\bar k^o,1,\bar t)$ and $T_i(d^o,\bar k^o,1,\bar t)$ w.r.t. $t_0$ does not necessarily coincide with $T_0(\bar k^o,\bar t_h)$ and $T_i(d^o,\bar k^o,\bar t_h)$. They may differ in a power of $t_0$. In particular, it is not necessarily true that $\overline{\cT^a_0(\bar k^o)}=\cT^h_0(\bar k^o)$ and $\overline{\cT^a_i(d^o,\bar k^o)}=\cT^h_i(d^o,\bar k^o)$ (the overline denotes projective closure, as usual). However, it holds that $\cT^h_i(d^o,\bar k^o)\cap\C^n=\cT^a_i(d^o,\bar k^o)$ and $\cT^h_0(d^o,\bar k^o)\cap\C^n=\cT^a_0(\bar k^o)$.
\end{Remark}
\vspace{3mm}
\begin{Proposition}\label{prop:ch4:MultiplicityAtNonFakePoints}
Let $\Omega_3$ be as in Lemma \ref{lem:ch4:QQ_zero_impliesInvariant} (page \pageref{lem:ch4:QQ_zero_impliesInvariant}). There exists a non-empty open  $\Omega_4\subset\Omega_3$, such that if $(d^o,\bar k^o)\in\Omega_4$, and $\bar t^o_h\in\cA_h\cap{\Psi_5^{P_h}}(d^o,\bar k^o)$, then:
\[\min_{i=1,2,3}\left(\mult_{\bar t^o}(\cT_o(\bar k^o),\cT_i(d^o,\bar k^o))\right)=1.\]
\end{Proposition}
\begin{proof}
Since $\bar t^o_h\in\cA_h\cap{\Psi_5^{P_h}}(d^o,\bar k^o)$, we can write  $\bar t^o_h=(1:t^o_1:t^o_2)$. Let $\bar t^o=(t^o_1,t^o_2)$. By Theorem \ref{thm:ch4:RelationBetweenSystemS5AndSystemS4} (page \pageref{thm:ch4:RelationBetweenSystemS5AndSystemS4}) we know that $\bar t^o=(t^o_1,t^o_2)\in\cA\cap\Psi_3^P(d^o,\bar k^o)$. W.l.o.g. we will suppose that
\[P_0(\bar t^o)h(\bar t^o)\bigl(P_2(\bar t^o)n_3(\bar t^o)-P_3(\bar t^o)n_2(\bar t^o)\bigr)\neq 0\]
(see the definition of the set $\cA$ in Equation \ref{def:ch4:Set_A_OfExtensibleSolutions}, page \pageref{def:ch4:Set_A_OfExtensibleSolutions}). In this case, it holds
(see Remark \ref{rem:ch4:SignOfLambdaAndOffsetting}, page \pageref{rem:ch4:SignOfLambdaAndOffsetting}) that
\[k_2^on_3(\bar t^o)-k_3^on_2(\bar t^o)\neq0 \mbox{ and }k_2^oP_3(\bar t^o)-k_3^oP_2(\bar t^o)\neq 0.\]
Furthermore, by Remark \ref{rem:ch4:RelationBetweenImplicitAndParametricNormalVector} (page \pageref{eq:ch4:RelationBetweenImplicitAndParametricNormalVectors}), one has
\begin{equation}\label{eq:ch4:ImplicitAndParamNormalVectors_Multip}
f_j(P(\bar t))=\dfrac{\beta(\bar t)}{P^{\mu}_0(\bar t)}n_j(\bar t)
\mbox{ for }j=1,2,3,
\end{equation}
and therefore
\begin{equation}\label{eq:ch4:ImplicitAndParamNormalVectors_Multip_h}
\sqrt{h_{\operatorname{imp}}(P(\bar t))}=\dfrac{\beta(\bar t)}{P^{\mu}_0(\bar t)}\sqrt{h(\bar t)}.
\end{equation}
with $\beta(\bar t^o)\neq 0$ (see Lemma \ref{lem:ch4:ExcludeBadParameterValues}, page \pageref{lem:ch4:ExcludeBadParameterValues}).

\nin For this case we will construct a non-empty open set  $\Omega_{4,1}\subset\Omega_3$ such that if $(d^o,\bar k^o)\in\Omega_{4,1}$, and $\bar t^o_h\in\cA_h\cap{\Psi_5^{P_h}}(d^o,\bar k^o)$, then:
\[\mult_{\bar t^o}\bigl(\cT_o(\bar k^o),\cT_1(d^o,\bar k^o)\bigr)=1.\]
If the second, respectively third, defining equation of $\cA$ is used, then analogous open subsets $\Omega_{4,2}$, respectively $\Omega_{4,3}$ can be constructed, and the corresponding  result for $\cT_2(d^o,\bar k^o)$, respectively $\cT_3(d^o,\bar k^o)$, is obtained. Finally, it suffices to take
\[\Omega_4=\Omega_{4,1}\cap\Omega_{4,2}\cap\Omega_{4,3}.\]

\nin The construction of $\Omega_{4,1}$ will proceed in several steps:
\begin{enumerate}
 \item[(1)] \nin By Proposition \ref{prop:ch4:ExtendableSolutions} (page \pageref{prop:ch4:ExtendableSolutions}), $(d^o,\bar k^o,\bar t^o)\in\pi_{(2,1)}({\Psi_2^P}(d^o,\bar k^o))$. Thus, by Theorem \ref{thm:ch4:TheoreticalFoundation} (page \pageref{thm:ch4:TheoreticalFoundation}), the point $\bar y^o=P(\bar t^o)$ is an affine, non normal-isotropic point of $\Sigma$, and it is associated with $\bar x^o\in\cL_{\bar k^o}\cap\cO_{d^o}(\Sigma)$, where $\bar x^o$ is a non normal-isotropic point of $\cO_{d^o}(\Sigma)$. Besides, since $(d^o,\bar k^o)\in\Omega_3\subset\Omega_0$, (see Remark \ref{rem:ch4:InOmega0GoodSpecialization}, page \pageref{rem:ch4:InOmega0GoodSpecialization}), $g(d^o,\bar x)$ is the defining polynomial of $\cO_{d^o}(\Sigma)$.

It follows that there is an open neighborhood $U^0$ of $(d^o,\bar y^o)$ (in the usual unitary topology of $\C\times\C^3$) such that the equation
\[
\left(\dfrac{\partial f}{\partial y_1}(\bar y)\right)^2+\left(\dfrac{\partial f}{\partial y_2}(\bar y)\right)^2+\left(\dfrac{\partial f}{\partial y_3}(\bar y)\right)^2=0
\]
has no solutions in $U^0$. Similarly, there is an open neighborhood $V^0$ of $(d^o,\bar x^o)$ (in the usual unitary topology of $\C\times\C^3$) such that the equation
\[
\left(\dfrac{\partial g}{\partial x_1}(d,\bar x)\right)^2+\left(\dfrac{\partial g}{\partial x_2}(d,\bar x)\right)^2+\left(\dfrac{\partial g}{\partial x_3}(d,\bar x)\right)^2=0
\]
has no solutions in $V^0$. Let us consider the map:
\[\varphi:U^0\to\C^3\]
defined by
\begin{equation}\label{eq:ch4:DefinitionVarphiMap_Multiplicity}
\bar\varphi(d,\bar y)=(\varphi_1(d,\bar y),\varphi_2(d,\bar y),\varphi_3(d,\bar y))=\bar y\pm d\dfrac{\nabla f(\bar y)}{\sqrt{h_{\operatorname{imp}}(\bar y)}}
\end{equation}
We assume w.l.o.g. that the $+$ sign in this expression is chosen so that $\bar\varphi(d^o,\bar y^o)=\bar x^o$; our discussion does not depend on this choice of sign in this expression, as will be shown below. According to Remark \ref{rem:ch4:SignOfLambdaAndOffsetting} and Lemma \ref{lem:ch4:SignOfLambdaAndOffsetting} (page \pageref{lem:ch4:SignOfLambdaAndOffsetting}), this implies that:
 \begin{equation}\label{eq:ch4:ChoiceOfEpsilonInMultiplicityProof}
\bar M(\bar k^o,\bar t^o)=\epsilon\dfrac{d^oP_0(\bar t^o)}{\sqrt{h(\bar t^o)}}\bar G(\bar k^o,\bar t^o)\mbox{ for }i=1,2,3.
\end{equation}
We will use Equation \ref{eq:ch4:ChoiceOfEpsilonInMultiplicityProof} later in the proof. Since $\bar y^o$ is not normal-isotropic in $\Sigma$, it follows that $\bar\varphi$ is analytic in $U^0$. Furthermore, we consider the  map
\[\bar\eta:V^0\to\C^3\]
defined by:
\[
\bar\eta(d,\bar x)=\bar x+d\dfrac{\nabla_{\bar x}g(d,\bar x)}{\|\nabla_{\bar x}g(d,\bar x)\|}.
\]
Here $\nabla_{\bar x}$ refers to the gradient computed w.r.t. $\bar x$; that is:
\[
\nabla_{\bar x}g(d,\bar x)=\left(\dfrac{\partial g}{\partial x_1}(d,\bar x),\dfrac{\partial g}{\partial x_2}(d,\bar x),\dfrac{\partial g}{\partial x_3}(d,\bar x)\right).
\]
In the definition of $\bar\eta$, w.l.o.g. the sign $+$ is chosen so that $\bar\eta(d^o,\bar x^o)=\bar y^o$. Then, since $\bar x^o$ is non normal-isotropic in  $\cO_{d^o}(\Sigma)$, it follows that $\bar\eta$ is analytic in $V^o$. Thus, there are open neighborhoods $U^1$ of $(d^o,\bar y^o)$ and $V^1$ of  $(d^o,\bar x^o)$ (in the unitary topology of $\C\times\C^3$), such that $\bar\varphi$ is an analytic isomorphism between $U^1$ and $V^1$, with inverse given by $\bar\eta$. We can assume w.l.o.g. that $\|\nabla f(\bar y)\|\neq 0$ holds in $U^1$, and $\|\nabla_{d,\bar x} g(\bar x)\|\neq 0$ holds in $V^1$. Note also that if $(d^o,\bar y^1)\in U^1$, with $\bar y^1\in\Sigma$, then $\bar\varphi(d^o,\bar y^1)\in\cO_{d^o}(\Sigma)$. It follows that the map $\bar\varphi_{d^o}$, obtained by restricting $\bar\varphi$ to $d=d^o$, induces an isomorphism:
\[d\bar\varphi_{d^o}:T_{\bar y^o}(\Sigma)\to T_{\bar x^o}(\cO_{d^o}(\Sigma))\]
where $T_{\bar y^o}(\Sigma)$ is the tangent plane to $\Sigma$ at $\bar y^o$, and $T_{\bar x^o}(\cO_{d^o}(\Sigma))$ is the tangent plane to $\cO_{d^o}((\Sigma)$ at $\bar x^o$.

\nin Since $(d^o,\bar k^o)\in\Omega_3\subset\Omega_0$, we have $\bar t^o\in\Upsilon_1$, with $\Upsilon_1$ as in Lemma \ref{lem:ch4:PropertiesSurfaceParametrization}, page \pageref{lem:ch4:PropertiesSurfaceParametrization} (see the construction of $\Omega^4_0$ in the proof of Theorem \ref{thm:ch4:TheoreticalFoundation}, \pageref{thm:ch4:TheoreticalFoundation}). Thus, the jacobian $\dfrac{\partial P}{\partial\bar t}(\bar t^o)$ has rank two. It follows that $P$ induces an isomorphism:
\[dP:T_{\bar t^o}(\C^2)\to T_{\bar y^o}(\Sigma)\]
where $T_{\bar t^o}(\C^2)$ is the tangent plane to $\C^2$ at $\bar t^o$. Therefore, the map defined by
\begin{equation}\label{eq:ch4:DefinitionNu_Multiplicity}
\bar\nu_{d^o}(\bar t)=\bar\varphi_{d^o}(P(\bar t))=\bar\varphi(d^o,P(\bar t))
\end{equation}
induces an isomorphism $d\bar\nu_{d^o}$ between $T_{\bar t^o}(\C^2)$ and $T_{\bar x^o}(\cO_{d^o}(\Sigma))$.

\item[(2)] Consider the following polynomials in $\C[d,\bar k,\rho]$:
\[
\begin{cases}
K(d,\bar k,\rho)=k_1\dfrac{\partial g}{\partial x_1}(d,\rho\bar k)+k_2\dfrac{\partial g}{\partial x_2}(d,\rho\bar k)+k_3\dfrac{\partial g}{\partial x_3}(d,\rho\bar k)\\
\tilde g(d,\bar k,\rho)=g(d,\rho\bar k)
\end{cases}
\]
and let
\[\Theta(d,\bar k)=\Res_{\rho}(K(d,\bar k,\rho),\tilde g(d,\bar k,\rho)).\]
Let us show that this resultant does not vanish identically.  If it does, then there are $A,B_1,B_2\in\C[d,\bar k,\rho]$, with $\deg_{\rho}(A(d,\bar k,\rho))>0$, such that
\[
\begin{cases}
K(d,\bar k,\rho)=A(d,\bar k,\rho)B_1(d,\bar k,\rho),\\
\tilde g(d,\bar k,\rho)=A(d,\bar k,\rho)B_2(d,\bar k,\rho).
\end{cases}
\]
Then $g(d,\rho\bar k)=A(d,\bar k,\rho)B_2(d,\bar k,\rho)$, and $\deg_{\bar k}(A(d,\bar k,\rho))>0$ (because $\tilde g$ cannot have a non constant factor in $\C[d,\rho]$). Thus, setting
$\rho=1$ and $\bar k=\bar x$, one has $g(d,\bar x)=A(d,\bar x,1)B_2(d,\bar x,1)$. It follows (see Remark \ref{rem:ch1:GenericOffsetEqSqfreeAndHasAtMostTwoFactors}(1), page \pageref{rem:ch1:GenericOffsetEqSqfreeAndHasAtMostTwoFactors}) that if $\tilde A(d,\bar x)$ is any irreducible factor of $A(d,\bar x,1)$, then $\tilde A(d,\bar x)$  defines an irreducible component $\cM$ of the generic offset, such that
\[x_1\dfrac{\partial g}{\partial x_1}(d,\bar x)+x_2\dfrac{\partial g}{\partial x_2}(d,\bar x)+x_3\dfrac{\partial g}{\partial x_3}(d,\bar x)=0\]
holds identically on $\cM$. Besides, for an open set of points $\bar x^o\in\cM$, one has $\nabla_{\bar x}g(d,\bar x^o)=\nabla_{\bar x}\tilde A(d,\bar x^o)$. Thus the above equation implies that
\[x_1\dfrac{\partial\tilde A}{\partial x_1}(d,\bar x)+x_2\dfrac{\partial\tilde A}{\partial x_2}(d,\bar x)+x_3\dfrac{\partial\tilde A}{\partial x_3}(d,\bar x)=0\]
holds identically in $\cM$. Therefore, since $\tilde A$ is irreducible, we get
\[x_1\dfrac{\partial\tilde A}{\partial x_1}(d,\bar x)+x_2\dfrac{\partial\tilde A}{\partial x_2}(d,\bar x)+x_3\dfrac{\partial\tilde A}{\partial x_3}(d,\bar x)=\kappa^o\tilde A(d,\bar x)\]
for some constant $\kappa^o$. This implies that the polynomial $\tilde A(d,\bar x)$ is homogeneous w.r.t. $\bar x$, and it follows that, for any value $d^o\not\in\Delta$ (with $\Delta$ as in Corollary \ref{cor:ch1:BadDistancesFiniteSet}, page \pageref{cor:ch1:BadDistancesFiniteSet}), $\cO_{d^o}(\Sigma)$, has a homogeneous component. This implies that $\bar 0\in\cO_{d^o}(\Sigma)$ for $d^o\not\in\Delta$, which is a contradiction with our hypothesis (see Remark \ref{rem:ch4:NotInfinitelyManyOffsetsThroughOrigin}(1), page \pageref{rem:ch4:NotInfinitelyManyOffsetsThroughOrigin}). Thus, $\Theta(d,\bar k)$ is not constant. Let us define $\Omega_{4,1}^1=\Omega_3\setminus\{(d^o,\bar k^o)/\Theta(d^o,\bar k^o)=0\}$.

\item[(3)] Let us consider the following polynomials in $\C[d, \bar k,\bar x]$
\begin{equation}\label{eq:ch4:SigmaPolynomials}
\begin{cases}
\sigma_0(d,\bar k,\bar x)=\det(\bar k,\bar x,\nabla_{\bar x}g(d,\bar x))\\
\sigma_1(d,\bar k,\bar x)=k_2x_3-k_3x_2\\
\end{cases}
\end{equation}
Let $\Omega_{4,1}^2\subset\Omega_{4,1}^1$ be such that, for $(d^o,\bar k^o)\in\Omega_{4,1}^2$, these polynomials are non identically zero (note that $\sigma_0$ and $\sigma_1$ are both homogeneous w.r.t. $\bar k$). Therefore, for $(d^o,\bar k^o)\in\Omega_{4,1}^2$, and for $i=0,1$, $\sigma_i(d^o,\bar k^o,\bar x)$ defines a surface $\Sigma_i(d^o,\bar k^o)$. From
\[\nabla_{\bar x}\sigma_1(\bar k,\bar x)=(0,-k_3,k_2)\]
one has
\[\nabla_{\bar x}g\wedge\nabla_{\bar x}\sigma_1=\left(k_2\dfrac{\partial g}{\partial x_2}+k_3\dfrac{\partial g}{\partial x_3},-k_2\dfrac{\partial g}{\partial x_1},-k_3\dfrac{\partial g}{\partial x_1}\right).\]
Let $(d^o,\bar k^o)\in\Omega_{4,1}^2$ and $\bar x^o\in\cO_{d^o}(\Sigma)\cap\cL_{\bar k^o}$. We will show that
\begin{equation}\label{eq:ch4:GradientPlaneAndOffsetCrossProduct}
\nabla_{\bar x}g(d^o,\bar x^o)\wedge\nabla_{\bar x}\sigma_1(\bar k^o,\bar x^o)\neq\bar 0.
\end{equation}
First, note that there is $\rho^o\in\C$ such that $\bar x^o=\rho^o\bar k^o$. If $\dfrac{\partial g}{\partial x_1}(d^o,\bar x^o)\neq 0$, then since $k_i\neq 0$ for $i=1,2,3$, the result follows. Thus, let $\dfrac{\partial g}{\partial x_1}(d^o,\bar x^o)=0$. If we suppose that $\nabla_{\bar x}g(d^o,\bar x^o)\wedge\nabla_{\bar x}\sigma_1(\bar k^o,\bar x^o)=\bar 0$, then
\[
k^o_2\dfrac{\partial g}{\partial x_2}(d^o,\bar x^o)+k^o_3\dfrac{\partial g}{\partial x_3}(d^o,\bar x^o)=0.
\]
Thus, one obtains
\[
\begin{cases}
g(d^o,\rho^o\bar k^o)=0\\
k^o_1\dfrac{\partial g}{\partial x_1}(d^o,\rho^o\bar k^o)+k^o_2\dfrac{\partial g}{\partial x_2}(d^o,\rho^o\bar k^o)+k^o_3\dfrac{\partial g}{\partial x_3}(d^o,\rho^o\bar k^o)=0
\end{cases}
\]
and it follows that $\Theta(d^o,\bar k^o)=0$ (with $\Theta$ as in step (2) of the proof), contradicting the construction of $\Omega_{4,1}^1$. Thus, Equation \ref{eq:ch4:GradientPlaneAndOffsetCrossProduct} is proved.

\nin We will prove the analogous result for $\sigma_0$. That is, for $(d^o,\bar k^o)\in\Omega_{4,1}^2$ and $\bar x^o\in\cO_{d^o}(\Sigma)\cap\cL_{\bar k^o}$, we will show that:
\begin{equation}\label{eq:ch4:GradientSigma0AndOffsetCrossProduct}
\nabla_{\bar x}g(d^o,\bar x^o)\wedge\nabla_{\bar x}\sigma_0(d^o,\bar k^o,\bar x^o)\neq\bar 0.
\end{equation}
From
\[\sigma_0(d,\bar k,\bar x)=\det(\bar k,\bar x,\nabla_{\bar x}g(d,\bar x))\]
and applying the derivation properties of determinants, one has, e.g.
\[
\dfrac{\partial\sigma_0}{\partial x_1}(d,\bar k,\bar x)=
\left|\begin{array}{ccc}
k_1&k_2&k_3\\
1&0&0\\
\partial_1g&\partial_2g&\partial_3g
\end{array}\right|
+
\left|\begin{array}{ccc}
k_1&k_2&k_3\\
x_1&x_2&x_3\\
\partial_{1,1}g&\partial_{2,1}g&\partial_{3,1}g
\end{array}\right|
\]
where $\partial_ig=\dfrac{\partial g}{\partial x_i}$ and $\partial_{i,j}g=\dfrac{\partial^2 g}{\partial x_i\partial x_j}$ for $i,j\in\{1,2,3\}$. Let as before, $\bar x^o=\rho^o\bar k^o$ for some $\rho^o\in\C$. Then:
\[
\dfrac{\partial\sigma_0}{\partial x_1}(d^o,\bar k^o,\bar x^o)=
\left|\begin{array}{ccc}
k^o_1&k^o_2&k^o_3\\
1&0&0\\
\partial_1g&\partial_2g&\partial_3g
\end{array}\right|
+
\left|\begin{array}{ccc}
k^o_1&k^o_2&k^o_3\\
\rho^o k^o_1&\rho^o k^o_2&\rho^o k^o_3\\
\partial_{1,1}g&\partial_{2,1}g&\partial_{3,1}g
\end{array}\right|
\]
with all the partial derivatives evaluated at $(d^o,\bar x^o)$. Since the second determinant in the above equation vanishes, one concludes that
\[
\dfrac{\partial\sigma_0}{\partial x_1}(d^o,\bar k^o,\bar x^o)=
k^o_3\dfrac{\partial g}{\partial x_2}(d^o,\bar x^o)-k^o_2\dfrac{\partial g}{\partial x_3}(d^o,\bar x^o).
\]
Similar results are obtained for the other two partial derivatives, leading to:
\[
\nabla_{\bar x}\sigma_0(d^o,\bar k^o,\bar x^o)=\nabla_{\bar x}g(d^o,\bar x^o)\wedge\bar k^o.
\]
Therefore,
\[
\nabla_{\bar x}g(d^o,\bar x^o)\wedge\nabla_{\bar x}\sigma_0(d^o,\bar k^o,\bar x^o)=
\nabla_{\bar x}g(d^o,\bar x^o)\wedge(\nabla_{\bar x}g(d^o,\bar x^o)\wedge\bar k^o).
\]
Note that $\nabla_{\bar x}g(d^o,\bar x^o)\wedge\bar k^o\neq \bar 0$ because, by construction, $\cL_{\bar k^o}$ is not normal to $\cO_{d^o}(\Sigma)$ at $\bar x^o$. If we suppose that
\[
\nabla_{\bar x}g(d^o,\bar x^o)\wedge\nabla_{\bar x}\sigma_0(d^o,\bar k^o,\bar x^o)=\bar 0,
\]
then the vectors $\nabla_{\bar x}g(d^o,\bar x^o)$ and $\nabla_{\bar x}g(d^o,\bar x^o)\wedge\bar k^o$ are parallel and perpendicular to each other. However, if two vectors are parallel and perpendicular, and one of them is not zero, then the other one must be isotropic.
One concludes that $\|\nabla_{\bar x}g(d^o,\bar x^o)\|=0$. This is a contradiction (see step (1) of the proof); therefore, Equation \ref{eq:ch4:GradientSigma0AndOffsetCrossProduct} is proved.

\nin From Equations \ref{eq:ch4:GradientPlaneAndOffsetCrossProduct} and \ref{eq:ch4:GradientSigma0AndOffsetCrossProduct}, and using Theorem 9 in \cite{Cox1997} (page 480), we conclude that $\bar x^o$ is a regular point in $\Sigma_i(d^o,\bar k^o)\cap\cO_{d^o}(\Sigma)$ (for $i=0,1$). Besides, $\bar x^o$ belongs to a unique one-dimensional component of $\Sigma_i(d^o,\bar k^o)\cap\cO_{d^o}(\Sigma)$. For $i=0,1$, let $\cC_{i}(d^o,\bar k^o)$ be the one-dimensional component of $\Sigma_i(d^o,\bar k^o)\cap\cO_{d^o}(\Sigma)$ containing $\bar x^o$.

\item[(4)] The non-zero vector
\[\bar v_i(d^o,\bar k^o,\bar x^o)=\nabla_{\bar x}g(d^o,\bar x^o)\wedge\nabla_{\bar x}\sigma_i(\bar k^o,\bar x^o),\quad (i=0,1)\]
obtained in step (3) of the proof, is a tangent vector to $\cC_{i}(d^o,\bar k^o)$ at $\bar x^o$. We will show that
\begin{equation}\label{eq:ch4:CrossProductOfTangentVectors}
\bar v_0(d^o,\bar k^o,\bar x^o)\wedge\bar v_1(d^o,\bar k^o,\bar x^o)\neq\bar 0.
\end{equation}
It holds that
\[
\begin{array}{l}
\bar v_0(d^o,\bar k^o,\bar x^o)\wedge\bar v_1(d^o,\bar k^o,\bar x^o)=\\
-(k^o_3\dfrac{\partial g}{\partial x_2}-k^o_2\dfrac{\partial g}{\partial x_3})\cdot
\left(k^o_1\dfrac{\partial g}{\partial x_1}+k^o_2\dfrac{\partial g}{\partial x_2}+k^o_3\dfrac{\partial g}{\partial x_3}\right)\cdot
\nabla_{\bar x}g(d^o,\bar x^o),
\end{array}
\]
with all the partial derivatives evaluated at $(d^o,\bar x^o)$.  Since
\[\|\nabla f(\bar y^o)\|\cdot\|\nabla_{\bar x}g(d^o,\bar x^o)\|\neq 0,\]
by the fundamental property of the offset (Proposition \ref{prop:ch1:FundamentalPropertyOffsets}, page \pageref{prop:ch1:FundamentalPropertyOffsets}), there is some $\kappa^o\in\C^\times$ such that
\[\nabla_{\bar x}g(d^o,\bar x^o)=\kappa^o\nabla f(\bar y^o).\]
Then, using Equation \ref{eq:ch4:ImplicitAndParamNormalVectors_Multip} (page \pageref{eq:ch4:ImplicitAndParamNormalVectors_Multip}) one has (see Remark \ref{rem:ch4:SignOfLambdaAndOffsetting}, page \pageref{rem:ch4:SignOfLambdaAndOffsetting}) that
\[k^o_3\dfrac{\partial g}{\partial x_2}-k^o_2\dfrac{\partial g}{\partial x_3}=
\kappa^o(k^o_3f_2(\bar y^o)-k^o_2f_3(\bar y^o))=\kappa^o\dfrac{\beta(\bar t^o)}{P^{\mu}_0(\bar t^o)}(k_3^on_2(\bar t^o)-k_2^on_3(\bar t^o))\neq 0.\]
Besides, in step (2) of the proof we have already seen that
\[\left(k^o_1\dfrac{\partial g}{\partial x_1}+k^o_2\dfrac{\partial g}{\partial x_2}+k^o_3\dfrac{\partial g}{\partial x_3}\right)\neq 0.\]
Thus, the proof of Equation \ref{eq:ch4:CrossProductOfTangentVectors} is finished.

\item[(5)]  From Equation \ref{eq:ch4:ChoiceOfEpsilonInMultiplicityProof} (page \pageref{eq:ch4:ChoiceOfEpsilonInMultiplicityProof}) one has
\[
M_1\bar k^o,\bar t^o)=\epsilon\dfrac{d^oP_0(\bar t^o)}{\sqrt{h(\bar t^o)}}G_1(\bar k^o,\bar t^o).
\]
Therefore:
\begin{equation}
\sqrt{h(\bar t^o)} (k^o_2P_3(\bar t^o)-k^o_3P_2(\bar t^o))=d^oP_0(\bar t^o)(k_2(\bar t^o)n_3(\bar t^o)-k_3(\bar t^o)n_2(\bar t^o)).
\end{equation}
Multiplying by $\dfrac{\beta(\bar t^o)}{P^{\mu+1}_0(\bar t^o)}$, it holds that
\[
\dfrac{\beta(\bar t^o)\sqrt{h(\bar t^o)}}{P^{\mu}_0(\bar t^o)}\left(k^o_2\dfrac{P_3(\bar t^o)}{P_0(\bar t^o)}-k^o_3\dfrac{P_2(\bar t^o)}{P_0(\bar t^o)}\right)=
d^o\left(k_2(\bar t^o)\dfrac{\beta(\bar t^o)n_3(\bar t^o)}{P^{\mu}_0(\bar t^o)}-k_3(\bar t^o)\dfrac{\beta(\bar t^o)n_2(\bar t^o)}{P^{\mu}_0(\bar t^o)}\right).
\]
Using Equation \ref{eq:ch4:ImplicitAndParamNormalVectors_Multip} (page \pageref{eq:ch4:ImplicitAndParamNormalVectors_Multip}), one obtains (recall that $\bar y^o=P(\bar t^o)$):
\begin{equation}\label{eq:ch4:Multiplicity_BranchPLusSign}
\sqrt{h_{\operatorname{imp}}(\bar y^o)}(k^o_2y^o_3-k^o_3y^o_2)-d^o(k^o_2f^o_3(\bar y^o)-k^o_3f_2(\bar y^o))=0.
\end{equation}
Note also that, since $\sqrt{h(\bar t^o)} (k^o_2P_3(\bar t^o)-k^o_3P_2(\bar t^o))\neq 0$ we also have
\begin{equation}\label{eq:ch4:Multiplicity_BranchMinusSign}
\sqrt{h_{\operatorname{imp}}(\bar y^o)}(k^o_2y^o_3-k^o_3y^o_2)+d^o(k^o_2f^o_3(\bar y^o)-k^o_3f_2(\bar y^o))\neq 0.
\end{equation}
Observe that, if the sign $\epsilon=-1$ is used in the offsetting construction (see step (1) of the proof), the results in Equations \ref{eq:ch4:Multiplicity_BranchPLusSign} and \ref{eq:ch4:Multiplicity_BranchMinusSign} are reversed.

\nin Recall (see Equation \ref{sys:ch4:AuxiliaryCurvesSystem}, page \pageref{sys:ch4:AuxiliaryCurvesSystem}) that the auxiliary polynomial $s_1$ is given by:
\[s_1(d,\bar k, \bar t)=h(\bar t)(k_2P_3-k_3P_2)^2-d^2P_0(\bar t)^2(k_2n_3-k_3n_2)^2.\]
Thus, one has:
\[\begin{array}{l}
\dfrac{\beta^2(\bar t)}{P^{2\mu+2}_0(\bar t)}s_1(d,\bar k, \bar t)=\\[5mm]
\dfrac{\beta^2(\bar t)}{P^{2\mu}_0(\bar t)}
\left(k_2\dfrac{P_3(\bar t)}{P_0(\bar t)}-k_3\dfrac{P_2(\bar t)}{P_0(\bar t)}\right)^2-
d^2\left(k_2\dfrac{\beta(\bar t)n_3(\bar t)}{P^{\mu}_0(\bar t)}-k_3\dfrac{\beta(\bar t)n_2(\bar t)}{P^{\mu}_0(\bar t)}\right)^2.
\end{array}
\]
And substituting $\bar y=P(\bar t)$ in $\dfrac{\beta^2(\bar t)}{P^{2\mu+2}_0(\bar t)}s_1(d,\bar k, \bar t)$, one obtains:
\begin{equation}\label{eq:ch4:AuxiliaryCurveImplicitVersion}
\dfrac{\beta^2(\bar t)}{P^{2\mu+2}_0(\bar t)}s_1(d,\bar k, \bar t)=
h_{\operatorname{imp}}(\bar y)(k_2y_3-k_3y_2)^2-d^2(k_2f_3(\bar y)-k_3f_2(\bar y))^2.
\end{equation}
Let us consider the polynomial $\sigma'_1\in\C[d,\bar k,\bar y]$ defined by
\[\sigma'_1(d,\bar k,\bar y)=h_{\operatorname{imp}}(\bar y)(k_2y_3-k_3y_2)^2-d^2(k_2f_3(\bar y)-k_3f_2(\bar y))^2,\]
and let $\Sigma'_1(d^o,\bar k^o)\subset\C^3$ be the algebraic closed set defined by the equation $\sigma'_1(d^o,\bar k^o,\bar y)=0$.
Let $\bar\tau=(\tau^1,\tau^2)$, and let $\cQ_1(\bar\tau)$ be a place of $\cT^a_1(d^o,\bar k^o)$ centered at $\bar t^o$. We assume that $\cQ_1(\bar 0)=\bar t^o$.
Since $T_1(d^o,\bar k^o,1,\cQ_1(\bar\tau))=0$ identically in $\bar\tau$, from Equation \ref{eq:ch4:ReducedAuxiliaryPolynomials_2} (page \pageref{eq:ch4:ReducedAuxiliaryPolynomials}) it follows that
\[s_1(d^o,\bar k^o,\cQ_1(\bar\tau))=S_1(d^o,\bar k^o,1,\cQ_1(\bar\tau))=0\]
identically in $\bar\tau$. Thus, from Equation \ref{eq:ch4:AuxiliaryCurveImplicitVersion} (recall that $\bar y=P(\bar t)$ in the lhs of Equation \ref{eq:ch4:AuxiliaryCurveImplicitVersion}) one has:
\[\sigma'_1(d^o,\bar k^o,P(\cQ_1(\bar\tau)))=0\]
identically in $\bar\tau$. Note that:
\[
\sigma'_1(d,\bar k^o,\bar y)=\sigma'_{1,+}(d,\bar k^o,\bar y)\sigma'_{1,-}(d,\bar k^o,\bar y)
\]
with
\[
\begin{cases}
\sigma'_{1,+}(d,\bar k^o,\bar y)=\sqrt{h_{\operatorname{imp}}(\bar y)}(k^o_2y_3-k^o_3y_2)+d(k^o_2f_3(\bar y)-k^o_3f_2(\bar y)),\\[3mm]
\sigma'_{1,-}(d,\bar k^o,\bar y)=\sqrt{h_{\operatorname{imp}}(\bar y)}(k^o_2y_3-k^o_3y_2)-d(k^o_2f_3(\bar y)-k^o_3f_2(\bar y)).
\end{cases}
\]
The functions $\sigma'_{1,+}(d,\bar k^o,\bar y)$ and $\sigma'_{1,-}(d,\bar k^o,\bar y)$ are analytic in the neighborhood $U^1$ of $(d^o,\bar x^o)$ introduced in step (1) of the proof. Therefore:
\[\sigma'_{1,+}(d^o,\bar k^o,\cQ_1(\bar\tau))\sigma'_{1,-}(d^o,\bar k^o,\cQ_1(\bar\tau))=0,\]
identically in $\bar\tau$. However, evaluating at $\bar\tau=\bar 0$, and taking Equations \ref{eq:ch4:Multiplicity_BranchPLusSign} and \ref{eq:ch4:Multiplicity_BranchMinusSign} into account, one sees that
\[\sigma'_{1,+}(d,\bar k^o,\bar y^o)\neq 0,\mbox{ while }\sigma'_{1,-}(d,\bar k^o,\bar y^o)=0.\]
By the analytic character of the functions, one concludes that
\[\sigma'_{1,-}(d^o,\bar k^o,\cQ_1(\bar\tau))=0,\]
identically in $\bar\tau$. Dividing by $\sqrt{h_{\operatorname{imp}}(\bar y)}$, this relation implies that:
{\small
\[
k^o_2\left(\dfrac{P_3(\cQ_1(\bar\tau))}{P_0(\cQ_1(\bar\tau))}+d^o\dfrac{f_3(P(\cQ_1(\bar\tau)))}{\sqrt{h_{\operatorname{imp}}(P_3(\cQ_1(\bar\tau)))}}\right)
-
k^o_3\left(\dfrac{P_2(\cQ_1(\bar\tau))}{P_0(\cQ_1(\bar\tau))}+d^o\dfrac{f_2(P(\cQ_1(\bar\tau)))}{\sqrt{h_{\operatorname{imp}}(P_2(\cQ_1(\bar\tau)))}}\right)=0.
\]}
\nin That is,
\[k^o_3\varphi_2(d^o,P(\cQ_1(\bar\tau)))-k^o_2\varphi_3(d^o,P(\cQ_1(\bar\tau)))=0\]
identically in $\bar\tau$, where $\bar\varphi=(\varphi_2,\varphi_2,\varphi_3)$  was defined in step (1) of the proof (see Equation \ref{eq:ch4:DefinitionVarphiMap_Multiplicity}, page \pageref{eq:ch4:DefinitionVarphiMap_Multiplicity}). With the notation introduced in step (3) of the proof (see Equation \ref{eq:ch4:SigmaPolynomials}, page \pageref{eq:ch4:SigmaPolynomials}), this is
\[
\sigma_1(d^o,\bar\varphi(d^o,P(\cQ_1(\bar\tau))))=0,
\]
identically in $\bar\tau$. This implies that if $\cB_1$ is the branch of $\cT^a_1(d^o,\bar k^o)$ at $\bar t^o$ determined by $\cQ_1(\bar\tau)$, then
\[\bar\nu_{d^o}(\cB_1)\subset\cC_1(d^o,\bar k^o),\]
where $\bar\nu_{d^o}$ was defined in Equation \ref{eq:ch4:DefinitionNu_Multiplicity} (page \pageref{eq:ch4:DefinitionNu_Multiplicity}), and  $\cC_1(d^o,\bar k^o)$ was introduced at the end of step (3) of the proof.

\item[(6)] Let $\cQ_0(\bar\tau)$ be a place of $\cT^a_0(\bar k^o)$ centered at $\bar t^o$. We assume that $\cQ_0(\bar 0)=\bar t^o$. Since $T_0(\bar k^o,1,\cQ_0(\bar\tau))=0$ identically in $\bar\tau$, from Equation \ref{eq:ch4:ReducedAuxiliaryPolynomials} (page \pageref{eq:ch4:ReducedAuxiliaryPolynomials}) it follows that
\[s_0(\bar k^o,\cQ_0(\bar\tau))=S_0(\bar k^o,1,\cQ_0(\bar\tau))=0,\]
identically in $\bar\tau$. That is,
\begin{equation}\label{eq:ch4:S0InAPLace_Multiplicity}
s_0(\bar k^o,\cQ_0(\bar\tau))=\det\left(
\begin{array}{ccc}
k^o_1&k^o_2&k^o_3\\
P_1(\cQ_0(\bar\tau))&P_2(\cQ_0(\bar\tau))&P_3(\cQ_0(\bar\tau))\\
n_1(\cQ_0(\bar\tau))&n_2(\cQ_0(\bar\tau))&n_3(\cQ_0(\bar\tau))
\end{array}
\right),
\end{equation}
identically in $\bar\tau$. Multiplying this by $\dfrac{\beta(\cQ_0(\bar\tau))}{P^{\mu+1}_0(\cQ_0(\bar\tau))}$, one has:
\[
\det\left(
\begin{array}{ccc}
k^o_1&k^o_2&k^o_3\\[3mm]
\dfrac{P_1(\cQ_0(\bar\tau))}{P_0(\bar\tau)}&\dfrac{P_2(\cQ_0(\bar\tau))}{{P_0(\cQ_0(\bar\tau))}}&\dfrac{P_3(\cQ_0(\bar\tau))}{P_0(\cQ_0(\bar\tau))}\\[5mm]
\dfrac{\beta(\cQ_0(\bar\tau))n_1(\cQ_0(\bar\tau))}{P_0^{\mu}(\cQ_0(\bar\tau))}&\dfrac{\beta(\bar\tau)n_2(\cQ_0(\bar\tau))}{P_0^{\mu}(\cQ_0(\bar\tau))}&
\dfrac{\beta(\bar\tau)n_3(\cQ_0(\bar\tau))}{P_0^{\mu}(\cQ_0(\bar\tau))}
\end{array}
\right)=0,
\]
identically in $\bar\tau$. Using Equation \ref{eq:ch4:ImplicitAndParamNormalVectors_Multip} (page \pageref{eq:ch4:ImplicitAndParamNormalVectors_Multip}), this implies that:
\[\det(\bar k^o,P(\cQ_0(\bar\tau)),\nabla f(P(\cQ_0(\bar\tau))))=0.\]
Since
\[
\bar\varphi(d^o,P(\cQ_0(\bar\tau)))=P(\cQ_0(\bar\tau))\pm d^o\dfrac{\nabla f(P(\cQ_0(\bar\tau)))}{\sqrt{h_{\operatorname{imp}}(\bar P(\cQ_0(\bar\tau)))}}
\]
and the second term in the sum is parallel to $\nabla f(P(\cQ_0(\bar\tau)))$, we have
\[\det(\bar k^o,\bar\varphi(d^o,P(\cQ_0(\bar\tau))),\nabla f(P(\cQ_0(\bar\tau))))=0.\]
Besides, by the fundamental property of the offset (Proposition \ref{prop:ch1:FundamentalPropertyOffsets}, page \pageref{prop:ch1:FundamentalPropertyOffsets}), and the construction in step (1) of the proof, the vectors
\[\nabla f(y)\mbox{ and }\nabla_{\bar x}g(d^o,\bar\varphi(d^o,\bar y))\]
are parallel for every value of $(d^o,\bar y)$ in $V^1$. It follows that
\[\det(\bar k^o,\bar\varphi(d^o,P(\cQ_0(\bar\tau))),\nabla_{\bar x}g(d^o,\bar\varphi(d^o,P(\cQ_0(\bar\tau)))))=0,\]
identically in $\bar\tau$. Recalling the definition of $\sigma_0$ in Equation \ref{eq:ch4:SigmaPolynomials} (page \pageref{eq:ch4:SigmaPolynomials}), this implies
that
\[\sigma_0(d^o,\bar k^o,\bar\varphi(d^o,P(\cQ_0(\bar\tau))))=0,\]
identically in $\bar\tau$. It follows that, if $\cB_0$ is the branch of $\cT^a_0(\bar k^o)$ at $\bar t^o$ determined by $\cQ_0(\bar\tau)$, then
\[\bar\nu_{d^o}(\cB_0)\subset\cC_0(d^o,\bar k^o),\]
where $\bar\nu_{d^o}$ was defined in Equation \ref{eq:ch4:DefinitionNu_Multiplicity} (page \pageref{eq:ch4:DefinitionNu_Multiplicity}), and  $\cC_0(d^o,\bar k^o)$ was introduced at the end of step (3) of the proof.
\end{enumerate}
Now we can finish the proof of the proposition. In steps (5) and (6) of the proof we have shown that any branch at $\bar t^o$ of the curves $\cT^a_0(\bar k^o)$ or $\cT^a_1(d^o,\bar k^o)$ is mapped by $\bar\nu_{d^o}$ respectively into the curve $\cC_1(d^o,\bar k^o)$ or $\cC_0(d^o,\bar k^o)$ (these curves are constructed in step (3)).  Since $d\bar\nu_{d^o}$ is an isomorphism of vector spaces (see step (1)), it follows that:
\begin{itemize}
 \item By the results in step (3), there is only one branch at $\bar t^o$ of each of the curves $\cT^a_0(\bar k^o)$ and $\cT^a_1(d^o,\bar k^o)$. Besides, since the rank of the Jacobian matrix (and therefore, the condition in \cite{Cox1997}, Theorem 9, page 480) is preserved under analytic isomorphisms, the unique branch of each the curves $\cT^a_0(\bar k^o)$ and $\cT^a_1(d^o,\bar k^o)$  passing through $\bar t^o$ is regular at that point.
 \item By the results in step (4), if $\ell_1$ and $\ell_0$ are the two tangent lines of these two branches, then $\ell_1$ and $\ell_0$ are different.
\end{itemize}
Then
\[\mult_{\bar t^o}(\cT_o,\cT_1)=1\]
follows from Theorem 5.10 in \cite{Walker1950} (page 114).
\end{proof}

\subsection{The degree formula}
\label{subsec:ch4:DegreeFormula}

\nin Before the statement of the degree formula we need to introduce some more notation and a technical lemma. Let
\[
R(\bar c,d,\bar k,\bar t)=\operatorname{Res}_{t_0}\left(T_0(\bar k,\bar t_h),T(\bar c,d,\bar k,\bar t_h)\right)
\]
(for the definition of $T_0$ and $T$ see Equations \ref{eq:ch4:ReducedAuxiliaryPolynomials} and \ref{eq:ch4:ReducedAuxiliaryPolynomials_3}, in page \pageref{eq:ch4:ReducedAuxiliaryPolynomials}). Then $R$ factors as follows:
\[
R(\bar c,d,\bar k,\bar t)=
N_1(d,\bar k,\bar t)M_3(\bar c,d,\bar k,\bar t)
\]
where $N_1(d,\bar k,\bar t)=\operatorname{Con}_{\bar c}(R)$ and
$M_3(\bar c,d,\bar k,\bar t)=\operatorname{PP}_{\bar c}(R)$.

\nin Besides, $N_1$ factors as follows:
\[N_1(d,\bar k,\bar t)=M_1(\bar t)M_2(d,\bar k,\bar t)\]
where $M_1(\bar t)=\operatorname{Con}_{(d,\bar k)}(N_1)$ and $M_2(d,\bar k,\bar t)=\operatorname{PP}_{(d,\bar k)}(N_1)$.  Thus, one has
\[
R(\bar c,d,\bar k,\bar t)=
M_1(\bar t)M_2(d,\bar k,\bar t)M_3(\bar c,d,\bar k,\bar t)
\]
and
\[M_2(d,\bar k,\bar t)=\operatorname{PP}_{(d,\bar k)}(\operatorname{Con}_{\bar c}(R)).\]
Note that $M_1, M_2$ and $M_3$ are homogeneous polynomials in $\bar t=(t_1,t_2)$.

\nin The following lemma deals with the specialization of the resultant $R(\bar c,d,\bar k,\bar t)$. More precisely, for $(d^o,\bar k^o)\in\C\times\C^3$ we denote:
\[
T^{\bar k^o}_0(\bar t_ h)=T_0(d^o,\bar k^o,\bar t_ h),
\quad T^{(d^o,\bar k^o)}(\bar c,\bar t_ h)=T(\bar c,d^o,\bar k^o,\bar t_ h).
\]
and
\[
R^{(d^o,\bar k^o)}(\bar c,\bar t)=
\operatorname{Res}_{t_0}\left(T^{\bar k^o}_0(\bar t_h),T^{(d^o,\bar k^o)}(\bar c,\bar t_h)\right)
\]
\vspace{3mm}
\begin{Lemma}\label{lem:ch4:ResultantSpecializationForDegreeFormula}
Let $\Omega_4$ be as in Proposition \ref{prop:ch4:MultiplicityAtNonFakePoints} (page \pageref{prop:ch4:MultiplicityAtNonFakePoints}). There exists a non-empty open $\Omega_5$,  such that for $(d^o,\bar k^o)\in\Omega_5$ the following hold:
\begin{enumerate}
 \item[(a)] The resultant $R(\bar c,d,\bar k,\bar t)$ specializes properly:
 \[
 R^{(d^o,\bar k^o)}(\bar c,\bar t)=R(\bar c,d^o,\bar k^o,\bar t)=
 M_1(\bar t)M_2(d^o,\bar k^o,\bar t)M_3(\bar c,d^o,\bar k^o,\bar t).
 \]
  \item[(b)] The content w.r.t $\bar c$ also specializes properly:
 \[\operatorname{Con}_{\bar c}(R^{(d^o,\bar k^o)})(\bar t)=M_1(\bar t)M_2(d^o,\bar k^o,\bar t).\]
 \item[(c)] the coprimality of $M_1$ and $M_2$ is invariant under specialization:
\[\gcd\left(M_1(\bar t),M_2(d^o,\bar k^o,\bar t)\right)=1\]
\end{enumerate}
\end{Lemma}
\begin{proof}
\nin For (a), consider $T_0$ and $T$ as polynomials in $\C[\bar c,d,\bar k,\bar t][t_0]$. Let $\operatorname{lc}(T_0)(\bar k,\bar t)$, (resp.  $\operatorname{lc}(T)(\bar c,d,\bar k,\bar t)$) be a leading coefficient w.r.t. $t_0$ of $T_0$ (resp. $T$). Take $A_1(\bar k)$ (resp. $B_1(d,\bar k)$) to be the coefficient of a term of $\operatorname{lc}(T_0)(\bar k,\bar t)$ (resp. $\operatorname{lc}(T)(d,\bar k,\bar t)$) of degree equal to $\deg_{\bar t}(\operatorname{lc}(T_0)(\bar k,\bar t))$ (resp. $\deg_{\{\bar c,\bar t\}}(\operatorname{lc}(T)(\bar k,\bar t))$). Now, if $A_1(\bar k^o)B_1(d^o,\bar k^o)\neq 0$, then (a) holds. Thus, set
\[\Omega^1_5=\Omega_4\cap\{(d^o,\bar k^o)/A_1(\bar k^o)B_1(d^o,\bar k^o)\neq 0\}.\]
For (b), we know that $M_3(\bar c,d,\bar k,\bar t)$ is primitive w.r.t. $\bar c$. If $M_3(\bar c,d,\bar k,\bar t)$ (considered as a polynomial in $\C[d,\bar k,\bar t][\bar c]$) has only one term, then its coefficient w.r.t. $\bar c$ must be constant, and so $M_3$ remains primitive under specialization of $(d,\bar k)$. Suppose, on the other hand, that $M_3(\bar c,d,\bar k,\bar t)$ has more than one term, and let:
\[M_{3,1}(d,\bar k,\bar t),\ldots,M_{3,\rho}(d,\bar k,\bar t)\]
be an (arbitrary) ordering of its non-zero coefficients w.r.t. $\bar c$.
Let
$\Gamma_1(d,\bar k,\bar t)=M_{3,1}(d,\bar k,\bar t)$, and for $j=2,\ldots,\rho$ let
\[
\Gamma_j(d,\bar k,\bar t)=\gcd\left(M_{3,j}(d,\bar k,\bar t),\Gamma_{j-1}(d,\bar k,\bar t)\right)
\]
Note that, for $j=1,\ldots,\rho$, the $M_{3,j}$ are homogeneous in $\bar t $ of the same degree. Thus, the  $\Gamma_j$ are either homogeneous in $\bar t$, or they only depend on $(d,\bar k)$.

\nin Since $M_3$ is primitive w.r.t. $\bar c$, let $j^o$ be the first index value in $1,\ldots,\rho$ for which $\Gamma_{j^o}(d,\bar k,\bar t)=1$. If $j^o=1$, then $M_{3,1}(d,\bar k,\bar t)$ is a constant, and in this case it is obvious that $M_3$ remains primitive under specialization of $(d,\bar k)$. If $j^o>1$,  we consider:
\[
\Res_{t_1}(M_{3,j^o}(d,\bar k,\bar t),\Gamma_{j^o-1}(d,\bar k,\bar t))
\]
This resultant is not identically zero, because we have assumed that $\Gamma_{j^o-1}(d,\bar k,\bar t)=1$. Since the involved polynomials are homogeneous in $\bar t$, this resultant is of the form $t_2^p\Phi(d,\bar k)$ for some $p\in\N$ and some $\Phi\in\C[d,\bar k]$. Now, because of the construction, if $\Phi(d^o,\bar k^o)\neq 0$, the specialization $M_3(\bar c,d^o,\bar k^o,\bar t)$ is primitive w.r.t. $\bar c$. Thus, set:
\[\Omega^2_5=\Omega^1_5\cap\{(d^o,\bar k^o)/\Phi(d^o,\bar k^o)\neq 0\}.\]
For (c) we use a similar construction. If either $M_1$ or $M_2$ do not depend on $\bar t_h$, the result is trivial. Otherwise, $M_1$ and $M_2$ are both homogeneous polynomials in $\bar t$, so the resultant
\[
\Res_{t_1}(M_1(\bar t),M_2(d,\bar k,\bar t))
\]
is of the form $t_2^{\tilde p}{\tilde \Phi}(d,\bar k)$ for some $\tilde p\in\N$ and some $\tilde \Phi_1\in\C[d,\bar k]$. Thus, if $\tilde\Phi_1(d^o,\bar k^o)\neq 0$, then $M_1(\bar t)$ and $M_2(d^o,\bar k^o,\bar t)$ do not have common factors of positive degree in $t_1$. A similar construction can be carried out w.r.t. $t_2$, obtaining a certain $\tilde\Phi_2$. Thus, set:
\[\Omega^3_5=\Omega^2_5\cap\{(d^o,\bar k^o)/\tilde\Phi_1(d^o,\bar k^o)\tilde\Phi_2(d^o,\bar k^o)\neq 0\}.\]
The construction shows that the lemma holds for $\Omega_5=\Omega^3_5$.
\end{proof}

\nin Finally, we are ready to state and prove the degree formula.
\vspace{3mm}
\begin{Theorem}[\sf Total Degree Formula for the Offset of a Parametric Surface]
\label{thm:ch4:DegreeFormula}
Let $T_0$ and $T$ be as in Equations \ref{eq:ch4:ReducedAuxiliaryPolynomials} and \ref{eq:ch4:ReducedAuxiliaryPolynomials_3} (page \pageref{eq:ch4:ReducedAuxiliaryPolynomials}). Then:
\[
\fbox{$
m\cdot\delta=
\deg_{\{\bar t\}}\left(
\operatorname{PP}_{(d,\bar k)}\left(
\operatorname{Con}_{\bar c}\left(
\operatorname{Res}_{t_0}\left(T_0(\bar k,\bar t_h),T(\bar c,d,\bar k,\bar t_h)\right)\right)\right)\right)=
\deg_{\bar t}(M_2(d,\bar k,\bar t))
$}
\]
where $m$ is the tracing index of $P$ (see Remark \ref{rem:ch4:TracingIndex_m}, page \pageref{rem:ch4:TracingIndex_m}), and if $g(d,\bar x)$ is the generic offset polynomial of $\Sigma$, then $\delta=\deg_{\bar x}(g(d,\bar x))$.
\end{Theorem}
\begin{proof}\label{proof:ch4:DegreeFormula}
Recall that (see Remark \ref{rem:ch4:OffsetDegreeAsCardinal_Projective}, page \pageref{rem:ch4:OffsetDegreeAsCardinal_Projective}), since $\Omega_3\subset\Omega_2$, if $(d^o,\bar k^o)\in\Omega_3$ it holds that
\[m\delta=\#\left(\cA_h\cap{\Psi_5^{P_h}}(d^o,\bar k^o)\right).\]
Thus, to prove the theorem it suffices to show that for any of these $(d^o,\bar k^o)$, it holds that
\[
\deg_{\bar t}\left(M_2(d,\bar k,\bar t)\right)=\#\left(\cA_h\cap{\Psi_5^{P_h}}(d^o,\bar k^o)\right).
\]
In order to do this, we will specialize at $(d^o,\bar k^o)$. More specifically, we will show that there is an open non-empty subset  $\Omega_6\subset\Omega_5$, such that, if $(d^o,\bar k^o)\in\Omega_6$, then $\deg_{\bar t}\left(M_2(d^o,\bar k^o,\bar t)\right)$ equals $\#\left(\cA_h\cap{\Psi_5^{P_h}}(d^o,\bar k^o)\right)$.

\nin Let $\Omega_5$ be as in Lemma \ref{lem:ch4:ResultantSpecializationForDegreeFormula} (page \pageref{lem:ch4:ResultantSpecializationForDegreeFormula}), and let $(d^o,\bar k^o)\in\Omega_5$.
Note that  $M_1(\bar t)$ and $M_2(d^o,\bar k^o,\bar t)$ both factor as product of linear factors.
there exists $\gamma\in\C$ such that $(\gamma:\alpha:\beta)\in{\Psi_5^P}(d^o,\bar k^o)$. Let us see that if $M_1(\bar t^o)=0$, with $\bar t^o=(t^o_1,t^o_2)$, and there is $t^o_0$ such that $(t^o_0:t^o_1:t^o_2)\in{\Psi_5^P}(d^o,\bar k^o)$, then $\bar t^o_h\not\in\cA_h$. In fact, if $t^o_0=0$, the result follows from  Proposition \ref{prop:ch4:SystemS5HasOnlyInvariantSolutionsAtInfinity} (page \pageref{prop:ch4:SystemS5HasOnlyInvariantSolutionsAtInfinity}) and Proposition \ref{prop:ch4:NonInvariantSolutionsCanBeAvoided} (page \pageref{prop:ch4:NonInvariantSolutionsCanBeAvoided}). Thus, w.l.o.g we suppose that $\bar t^o_h=(1:t^o_1:t^o_2)$, with $\bar t^o_h\in\cA_h$. Then using Proposition \ref{prop:ch4:NonInvariantSolutionsCanBeAvoided} (page \pageref{prop:ch4:NonInvariantSolutionsCanBeAvoided}), we get $\bar t^o_h\not\in\cI^{P_h}_5(\Omega_2)$. This is a contradiction, since $M_1(\bar t)$ does not depend on $(d,\bar k)$.

\nin We will now show that there is an open set $\Omega_6\subset\Omega_5$ such that if $(d^o,\bar k^o)\in\Omega_6$ and $M_2(d^o,\bar k^o,\bar t^o)=0$, then
there is $t^o_0$ such that $\bar t^o_h=(t^o_0,\bar t^o)\in{\Psi_5^P}(d^o,\bar k^o)$. This follows from Lemma \ref{lem:ch1:GeneralizedResultants}, page \pageref{lem:ch1:GeneralizedResultants}. Let us define:
\[\begin{cases}
\Omega_6^1=\Omega_5\cap\{(d^o,\bar k^o)\,|\,\deg_{\bar t_0}(T_i(d^o,\bar k^o,\bar t_h))=\deg_{\bar t_0}(T_i(d,\bar k,\bar t_h))\mbox{ for }i=1,2,3\}
\\
\Omega_6^2=\Omega_6^1\cap\{(d^o,\bar k^o)\,|\,\deg_{\bar t_h}(T_i(d^o,\bar k^o,\bar t_h))=\deg_{\bar t_h}(T_i(d,\bar k,\bar t_h))\mbox{ for }i=1,2,3\}
\\
\Omega_6^3=\Omega_6^2\cap\{(d^o,\bar k^o)\,|\,\gcd(T_1(d^o,\bar k^o,\bar t_h),T_2(d^o,\bar k^o,\bar t_h),T_3(d^o,\bar k^o,\bar t_h))=1\}
\end{cases}\]
The sets $\Omega_6^1$ and $\Omega_6^2$ are open and non-empty because they are defined by the non-vanishing of the corresponding leading coefficients. The fact that $\Omega_6^3$ is open and non-empty follows from a similar argument to the proof of Lemma \ref{lem:ch4:ResultantSpecializationForDegreeFormula}(c) (page \pageref{lem:ch4:ResultantSpecializationForDegreeFormula}).
Finally, take  $\Omega_6=\Omega_6^3$. Then, (i), (ii) and (iii) in Lemma \ref{lem:ch1:GeneralizedResultants} hold because of the construction of $\Omega_6^i$ for $i=1,2,3$, respectively. And also
$$\operatorname{lc}_{t_0}(T_0)(\bar t^o)\cdot\operatorname{lc}_{t_0}(T)(\bar c,\bar t^o)\neq 0$$
holds because of the construction of $\Omega_5^1$ in Lemma \ref{lem:ch4:ResultantSpecializationForDegreeFormula} (page \pageref{lem:ch4:ResultantSpecializationForDegreeFormula}), and  because $\Omega_6\subset\Omega_5$.

\nin Let $(d^o,\bar k^o)\in\Omega_6$. If $\bar t^o_h\in\cA\cap\Psi^{P_h}_5(d^o,\bar k^o)$, then $M_1(\bar t^o)M_2(d^o,\bar k^o,\bar t^o)=0$. Since we have seen that $M_1(\bar t^o)\neq 0$, one concludes that  $M_2(d^o,\bar k^o,\bar t^o)=0$. Conversely, let $\bar t^o$ be such that $M_2(d^o,\bar k^o,\bar t^o)=0$. Then, by the construction of $\Omega_6$, there is $t^o_0$ such that $(t^o_0:\bar t^o)\in{\Psi_5^P}(d^o,\bar k^o)$. Let us see that $\bar t^o_h\in\cA_h$. If $\bar t^o_h\in\cI^{P_h}_5(\Omega_2)$, then  because of the invariance, $M_1(\bar t^o)=0$, and this contradicts Lemma \ref{lem:ch4:ResultantSpecializationForDegreeFormula}(c) (page \pageref{lem:ch4:ResultantSpecializationForDegreeFormula}). Thus, $\bar t^o_h\not\in\cI^{P_h}_5(\Omega_2)$, and by Proposition \ref{prop:ch4:NonInvariantSolutionsCanBeAvoided_Part2} (page \pageref{prop:ch4:NonInvariantSolutionsCanBeAvoided_Part2}), one has $\bar t^o_h\in\cA_h$.

\nin Thus, we have shown that for each of the factors of $M_2(d^o,\bar k^o,\bar t)$ there is a point $\bar t^o_h\in\cA_h\cap{\Psi_5^{P_h}}(d^o,\bar k^o)$ such that $M_2(d^o,\bar k^o,\bar t^o)=0$, and conversely. Let $L^{(\alpha,\beta)}(\bar  t)=\beta t_1-\alpha t_2$ be one of these factors of $M_2(d^o,\bar k^o,\bar t)$, and let $\cL^{(\alpha,\beta)}$ the line defined by the equation $L^{(\alpha,\beta)}(\bar  t)=0$. By Lemma \ref{lem:ch4:ResultantSpecializationForDegreeFormula}(c) (page \pageref{lem:ch4:ResultantSpecializationForDegreeFormula}), one has
\begin{equation}\label{eq:ch4:CardinalEqualityInDegreeFormulaProof}
\#\left(\cL^{(\alpha,\beta)}\cap\cA_h\cap{\Psi_5^{P_h}}(d^o,\bar k^o)\right)=\#\bigl(\cL^{(\alpha,\beta)}\cap{\Psi_5^{P_h}}(d^o,\bar k^o)\bigr).
\end{equation}
If we define
\[p(\alpha,\beta)=\#\bigl(\cL^{(\alpha,\beta)}\cap{\Psi_5^{P_h}}(d^o,\bar k^o)\bigr),\]
then we will show that $L^{(\alpha,\beta)}(\bar  t)$ appears in $M_2(d^o,\bar k^o,\bar t)$ with exponent equal to $p(\alpha,\beta)$. From this it will follow that:
\[
\#\left(\cA_h\cap{\Psi_5^{P_h}}(d^o,\bar k^o)\right)=\sum_{(\alpha,\beta)}\#\bigl(\cL^{(\alpha,\beta)}\cap{\Psi_5^{P_h}}(d^o,\bar k^o)\bigr)=
\sum_{(\alpha,\beta)}p(\alpha,\beta)=\deg_{\bar t}\left(M_2(d,\bar k,\bar t)\right),
\]
and this will conclude the proof of the theorem.

\nin To prove our claim, note that $\cA_h\cap{\Psi_5^{P_h}}(d^o,\bar k^o)$ is a finite set, and by Proposition \ref{prop:ch4:MultiplicityAtNonFakePoints} (page \pageref{prop:ch4:MultiplicityAtNonFakePoints}), if $\bar t^o_h\in\cA_h\cap{\Psi_5^{P_h}}(d^o,\bar k^o)$, then:
\begin{equation}\label{eq:ch4:MultiplicityOfIntersection}
\min_{i=1,2,3}\left(\mult_{\bar t^o}(\cT_o(\bar k^o),\cT_i(d^o,\bar k^o))\right)=1.
\end{equation}
Recall that
\[
T^{\bar k^o}_0(\bar t_ h)=T_0(d^o,\bar k^o,\bar t_ h)\]
and
\[
\quad T^{(d^o,\bar k^o)}(\bar c,\bar t_ h)=T(\bar c,d^o,\bar k^o,\bar t_ h)
=c_1T_1(d^o,\bar k^o,\bar t_h)+c_2T_2(d^o,\bar k^o,\bar t_h)+c_3T_3(d^o,\bar k^o,\bar t_h).\]
For $\bar c^o\in{\C}^3$, let $\cT^{(\bar c^o,d^o,\bar k^o)}$ be the algebraic closed subset of $\P^2$ defined by the equation $T^{(d^o,\bar k^o)}(\bar c^o,\bar t_ h)=0$. Note that there is an open set of values $\bar c^o$ for which $\cT^{(\bar c^o,d^o,\bar k^o)}$ is indeed a curve. Let us see that there is an open subset $A(\bar t^o_h)\subset{\C}^3$, such that if $c^o\in A(\bar t^o_h)$, then
\[\mult_{\bar t^o_h}({\cT}^{\bar k^o}_0,\cT^{(\bar c^o,d^o,\bar k^o)})=1.\]
To prove this, let $\cP(\bar\tau)$ be a place of ${\cT}^{\bar k^o}_0$ at $\bar t^o_h$. Then, by Equation \ref{eq:ch4:MultiplicityOfIntersection} the order of the power series
$T^{(d^o,\bar k^o)}(\bar c,\cP(\bar\tau))$ is one. From this, one sees that it suffices to take $A(\bar t^o_h)$ to be the open set of values $\bar c^o$ for which the order of $T^{(d^o,\bar k^o)}(\bar c,\cP(\bar\tau))$ does not increase.

\nin Let now
\[\bar c^o\in\bigcap_{\bar t^o_h\in\cA_h\cap{\Psi_5^{P_h}}(d^o,\bar k^o)}A(\bar t^o).\]
Applying Lemma \ref{lem:ch1:MultiplicityUsingResultants} (page \pageref{lem:ch1:MultiplicityUsingResultants}) to the curves ${\cT}^{\bar k^o}_0$ and $\cT^{(\bar c^o,d^o,\bar k^o)}$, and the line $\cL^{(\alpha,\beta)}$, one concludes that the factor $\beta t_1-\alpha_j t_2$ appears in $\operatorname{Res}_{t_0}\left(T^{\bar k^o}_0(\bar t_h),T^{(d^o,\bar k^o)}(\bar c^o,\bar t_h)\right)$ with exponent equal to:
\[\sum_{\bar t^o_h\in\cL^{(\alpha,\beta)}\cap\cA_h\cap{\Psi_5^{P_h}}(d^o,\bar k^o)}
\mult_{\bar t^o_h}({\cT}^{\bar k^o}_0,\cT^{(\bar c^o,d^o,\bar k^o)})=\#\bigl(\cL^{(\alpha,\beta)}\cap\cA_h\cap{\Psi_5^{P_h}}(d^o,\bar k^o)\bigr).\]
Taking Equation \ref{eq:ch4:CardinalEqualityInDegreeFormulaProof} into account, this finishes the proof of our claim, and of the theorem.
\end{proof}

\nin We will finish this section with some examples, illustrating the use of the degree formula in Theorem \ref{thm:ch4:DegreeFormula} (page \pageref{thm:ch4:DegreeFormula}). The implicit equations in these examples have been obtained with the Computer Algebra System \cocoa (see \cite{CocoaSystem}).
\vspace{3mm}
\begin{Example}\label{exam:ch4:HypParabNonSym}
Let $\Sigma$ be the surface (a hyperbolic paraboloid) with implicit equation
\[y_3-y_1^2+\dfrac{y_2^2}{4}=0.\]
A rational --in fact polynomial-- parametrization of $\Sigma$ is given by:
\[P(t_1,t_2)=(t_1,2t_2,t_1^2-t_2^2).\]
From the form of its two first components, it is clear that this is a proper parametrization. This surface and its offset at $d^o=1$ are illustrated in Figure \ref{fig:HypParabNonSym}.
\begin{figure}
\centering
\includegraphics[height=8cm]{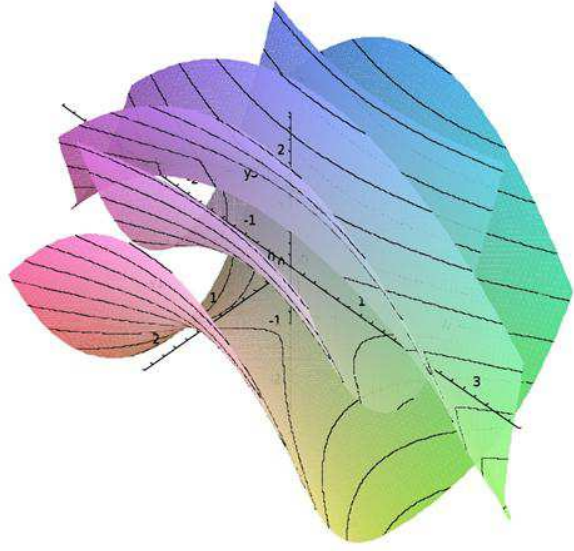}
\caption{Hyperbolic paraboloid and one of its offsets}
\label{fig:HypParabNonSym}
\end{figure}
The homogeneous associated normal vector is
\[N(\bar t_h)=(-2t_1,t_2,t_0).\]
\nin Then the auxiliary curves are:\\[3mm]
\nin $T_0(\bar t_h)=2\,k_1\,t_2\,t_0^{2}-k_1\,t_2\,t_1^{2}+k_1\,t_2^{3}-t_1\,t_0^{2}k_2+5\,t_1\,t_0\,k_3\,t_2-2\,k_2\,t_1^{3}+2\,t_1\,k_2\,t_2^{2},$\\[3mm]
\nin$T_1(\bar t_h)=4\,t_1^{6}k_2^{2}-7\,t_1^{4}k_2^{2}t_2^{2}-16\,t_1^{4}k_2\,k_3\,t_2\,t_0+2\,t_1^{2}k_2^{2}t_2^{4}+12\,t_1^{2}k_2\,t_2^{3}k_3\,t_0+16\,t_1^{2}k_3^{2}t_2^{2}t_0^{2}+t_2^{6}k_2^{2}+4\,t_2^{5}k_2\,k_3\,t_0+4\,t_2^{4}k_3^{2}t_0^{2}+t_0^{2}k_2^{2}t_1^{4}-2\,t_0^{2}k_2^{2}t_1^{2}t_2^{2}-4\,t_0^{3}k_2\,t_1^{2}k_3\,t_2+t_0^{2}k_2^{2}t_2^{4}+4\,t_0^{3}k_2\,t_2^{3}k_3+4\,t_0^{4}k_3^{2}t_2^{2}-{d}^{2}t_0^{6}k_2^{2}+2\,{d}^{2}t_0^{5}k_2\,k_3\,t_2-{d}^{2}t_0^{4}k_3^{2}t_2^{2},$\\[3mm]
\nin$T_2(\bar t_h)=4\,t_1^{4}t_0^{2}k_3^{2}-8\,t_1^{5}t_0\,k_3\,k_1+6\,t_1^{3}t_0\,k_3\,k_1\,t_2^{2}+4\,t_1^{6}k_1^{2}-7\,t_1^{4}k_1^{2}t_2^{2}+2\,t_1^{2}k_1^{2}t_2^{4}+t_1^{2}k_3^{2}t_2^{2}t_0^{2}+2\,t_2^{4}t_1\,t_0\,k_3\,k_1+t_2^{6}k_1^{2}+t_0^{4}t_1^{2}k_3^{2}-2\,t_0^{3}t_1^{3}k_3\,k_1+2\,t_0^{3}t_1\,k_3\,k_1\,t_2^{2}+t_0^{2}k_1^{2}t_1^{4}-2\,t_0^{2}k_1^{2}t_1^{2}t_2^{2}+t_0^{2}k_1^{2}t_2^{4}-4\,{d}^{2}t_0^{4}t_1^{2}k_3^{2}-4\,{d}^{2}t_0^{5}k_3\,t_1\,k_1-{d}^{2}t_0^{6}k_1^{2},$\\[3mm]
\nin$T_3(\bar t_h)=t_0^{2} (4\,k_2^{2}t_1^{4}-16\,t_1^{3}k_2\,k_1\,t_2+16\,k_1^{2}t_1^{2}t_2^{2}+k_2^{2}t_1^{2}t_2^{2}-4\,t_2^{3}k_2\,t_1\,k_1+4\,k_1^{2}t_2^{4}+t_0^{2}k_2^{2}t_1^{2}-4\,k_2\,t_1\,t_0^{2}k_1\,t_2+4\,k_1^{2}t_2^{2}t_0^{2}-4\,t_0^{2}{d}^{2}k_2^{2}t_1^{2}-4\,t_0^{2}{d}^{2}t_1\,k_2\,k_1\,t_2-t_0^{2}{d}^{2}k_1^{2}t_2^{2}).
$\\[3mm]
Denoting, as in the degree formula,
\[
R(\bar c,d,\bar k,\bar t)=\operatorname{Res}_{t_0}\left(T_0(\bar k,\bar t_h),T(\bar c,d,\bar k,\bar t_h)\right),
\]
one has:

\nin{$
R(\bar c,d,\bar k,\bar t)=
 ( t_1-t_2) ^{2} ( t_1+t_2) ^{2}( 4\,t_2^{4}c_3^{2}k_1^{4}+t_2^{4}c_2^{2}k_1^{4}+t_2^{4}c_1^{2}k_2^{4}+t_1^{2}t_2^{2}c_3^{2}k_1^{2}k_2^{2}+4\,t_1^{2}t_2^{2}c_1\,c_3\,k_1^{2}k_2^{2}+4\,t_1^{2}t_2^{2}c_2\,c_3\,k_1^{2}k_2^{2}+4\,t_1^{2}t_2^{2}c_2\,c_3\,k_1^{4}+4\,t_1^{2}t_2^{2}c_1\,c_3\,k_2^{4}+4\,t_1^{2}t_2^{2}c_1\,c_2\,k_2^{2}k_3^{2}+4\,t_1^{2}t_2^{2}c_1\,c_2\,k_1^{2}k_3^{2}+17\,t_1^{2}t_2^{2}c_1\,c_3\,k_2^{2}k_3^{2}+17\,t_1^{2}t_2^{2}c_2\,c_3\,k_1^{2}k_3^{2}-4\,t_1^{2}t_2^{2}c_1^{2}k_2^{2}k_3^{2}+17\,t_1^{2}c_2\,c_1\,k_3^{4}t_2^{2}-4\,t_1^{2}t_2^{2}c_2^{2}k_1^{2}k_3^{2}-6\,t_1\,t_2^{3}c_2\,c_3\,k_1^{3}k_2-6\,t_1\,t_2^{3}c_1\,c_3\,k_1\,k_2^{3}+12\,t_1\,t_2^{3}c_1\,c_3\,k_1\,k_2\,k_3^{2}-12\,t_1\,t_2^{3}c_1\,c_2\,k_1\,k_2\,k_3^{2}+12\,t_1\,t_2^{3}c_3^{2}k_1^{3}k_2+2\,t_2^{4}c_1\,c_2\,k_1^{2}k_2^{2}-4\,t_2^{4}c_1\,c_3\,k_1^{2}k_2^{2}-4\,t_2^{4}c_2\,c_3\,k_1^{4}-4\,t_1^{4}c_1\,c_3\,k_2^{4}-2\,t_1^{2}t_2^{2}c_2^{2}k_1^{4}-2\,t_1^{2}t_2^{2}c_1^{2}k_2^{4}+4\,t_1^{4}c_2^{2}k_3^{4}-4\,t_1^{4}c_1\,c_2\,k_2^{2}k_3^{2}+4\,t_1^{4}c_2^{2}k_1^{2}k_3^{2}+2\,t_1^{4}k_1^{2}k_2^{2}c_2\,c_1-4\,t_1^{4}c_2\,c_3\,k_1^{2}k_2^{2}+8\,t_1^{4}c_2\,c_3\,k_2^{2}k_3^{2}-12\,t_1^{3}t_2\,c_3^{2}k_1\,k_2^{3}+6\,t_1^{3}t_2\,c_1\,c_3\,k_1\,k_2^{3}+12\,t_1^{3}t_2\,c_1\,c_2\,k_1\,k_2\,k_3^{2}-12\,t_1^{3}t_2\,c_2\,c_3\,k_1\,k_2\,k_3^{2}+6\,t_1^{3}t_2\,c_2\,c_3\,k_1^{3}k_2-4\,t_1^{2}t_2^{2}c_1\,c_2\,k_1^{2}k_2^{2}-4\,t_2^{4}c_1\,c_2\,k_1^{2}k_3^{2}+8\,t_2^{4}c_1\,c_3\,k_1^{2}k_3^{2}+4\,c_1^{2}t_2^{4}k_3^{4}+4\,t_2^{4}c_1^{2}k_2^{2}k_3^{2}+t_1^{4}k_2^{4}c_1^{2}+4\,t_1^{4}c_3^{2}k_2^{4}+t_1^{4}k_1^{4}c_2^{2})\,\cdot\,( -48\,{d}^{2}t_1^{6}t_2^{4}k_2^{6}+64\,{d}^{2}t_1^{2}t_2^{8}k_1^{6}-72\,{d}^{2}t_1^{4}t_2^{6}k_1^{6}-128\,{d}^{4}t_1^{8}t_2^{2}k_2^{6}+64\,{d}^{4}t_1^{6}t_2^{4}k_2^{6}+{d}^{4}t_1^{4}t_2^{6}k_1^{6}-2\,{d}^{4}t_1^{2}t_2^{8}k_1^{6}+25\,t_1^{6}k_2^{4}t_2^{4}k_3^{2}+4\,k_2^{6}t_1^{10}-78\,t_1^{5}t_2^{5}k_1\,k_2^{5}-154\,t_1^{7}t_2^{3}k_1\,k_2^{5}+{d}^{4}k_1^{6}t_2^{10}+64\,{d}^{4}k_2^{6}t_1^{10}+8\,{d}^{2}k_1^{6}t_2^{10}+32\,{d}^{2}k_2^{6}t_1^{10}+136\,{d}^{2}t_1^{5}t_2^{5}k_1\,k_2^{5}+88\,{d}^{2}t_1^{7}t_2^{3}k_1\,k_2^{5}+320\,{d}^{2}t_1^{8}t_2^{2}k_1^{2}k_2^{4}+440\,{d}^{2}t_1^{4}t_2^{6}k_1^{4}k_2^{2}-100\,{d}^{2}t_1^{6}k_2^{2}k_3^{2}t_2^{4}k_1^{2}-170\,{d}^{2}t_1^{3}t_2^{7}k_1^{3}k_2^{3}-110\,{d}^{2}t_1^{5}t_2^{5}k_1^{3}k_2^{3}+280\,{d}^{2}t_1^{7}t_2^{3}k_1^{3}k_2^{3}+1200\,{d}^{2}t_1^{7}k_2^{3}k_3^{2}t_2^{3}k_1-400\,{d}^{2}t_1^{8}k_2^{4}k_3^{2}t_2^{2}+300\,{d}^{2}t_1^{5}k_2^{3}k_3^{2}t_2^{5}k_1-1200\,{d}^{2}t_1^{5}k_2\,k_3^{2}t_2^{5}k_1^{3}-25\,{d}^{2}t_1^{4}k_2^{2}k_3^{2}t_2^{6}k_1^{2}-400\,{d}^{2}t_1^{4}k_3^{2}t_2^{6}k_1^{4}-100\,{d}^{2}t_1^{2}k_3^{2}t_2^{8}k_1^{4}-370\,{d}^{2}t_1^{6}t_2^{4}k_1^{4}k_2^{2}-344\,{d}^{2}t_1^{5}t_2^{5}k_1^{5}k_2+16\,{d}^{2}t_1\,t_2^{9}k_1^{5}k_2+328\,{d}^{2}t_1^{3}t_2^{7}k_1^{5}k_2-300\,{d}^{2}t_1^{3}k_2\,k_3^{2}t_2^{7}k_1^{3}-70\,{d}^{2}t_1^{2}t_2^{8}k_1^{4}k_2^{2}+465\,t_1^{8}t_2^{2}k_1^{2}k_2^{4}-120\,{d}^{4}t_1^{4}t_2^{6}k_1^{4}k_2^{2}+160\,{d}^{4}t_1^{3}t_2^{7}k_1^{3}k_2^{3}-24\,{d}^{4}t_1^{3}t_2^{7}k_1^{5}k_2+60\,{d}^{4}t_1^{2}t_2^{8}k_1^{4}k_2^{2}+12\,{d}^{4}t_1\,t_2^{9}k_1^{5}k_2+2480\,t_1^{4}t_2^{6}k_1^{4}k_2^{2}+2400\,t_1^{6}k_2^{2}k_3^{2}t_2^{4}k_1^{2}-440\,t_1^{3}t_2^{7}k_1^{3}k_2^{3}-1920\,t_1^{5}t_2^{5}k_1^{3}k_2^{3}-1640\,t_1^{7}t_2^{3}k_1^{3}k_2^{3}-800\,t_1^{7}k_2^{3}k_3^{2}t_2^{3}k_1+100\,t_1^{8}k_2^{4}k_3^{2}t_2^{2}+16\,{d}^{2}t_1^{8}t_2^{2}k_2^{6}-200\,t_1^{5}k_2^{3}k_3^{2}t_2^{5}k_1-3200\,t_1^{5}k_2\,k_3^{2}t_2^{5}k_1^{3}+600\,t_1^{4}k_2^{2}k_3^{2}t_2^{6}k_1^{2}+1600\,t_1^{4}k_3^{2}t_2^{6}k_1^{4}+400\,t_1^{2}k_3^{2}t_2^{8}k_1^{4}+3160\,t_1^{6}t_2^{4}k_1^{4}k_2^{2}-3168\,t_1^{5}t_2^{5}k_1^{5}k_2-128\,t_1\,t_2^{9}k_1^{5}k_2-1504\,t_1^{3}t_2^{7}k_1^{5}k_2-800\,t_1^{3}k_2\,k_3^{2}t_2^{7}k_1^{3}+12\,t_1^{8}t_2^{2}k_2^{6}+9\,t_1^{6}t_2^{4}k_2^{6}+360\,t_1^{2}t_2^{8}k_1^{4}k_2^{2}-224\,{d}^{2}t_1^{9}t_2\,k_1\,k_2^{5}+20\,{d}^{2}t_1^{4}t_2^{6}k_1^{2}k_2^{4}-340\,{d}^{2}t_1^{6}t_2^{4}k_1^{2}k_2^{4}+192\,{d}^{4}t_1^{9}t_2\,k_1\,k_2^{5}+240\,{d}^{4}t_1^{8}t_2^{2}k_1^{2}k_2^{4}-384\,{d}^{4}t_1^{7}t_2^{3}k_1\,k_2^{5}+160\,{d}^{4}t_1^{7}t_2^{3}k_1^{3}k_2^{3}-480\,{d}^{4}t_1^{6}t_2^{4}k_1^{2}k_2^{4}+60\,{d}^{4}t_1^{6}t_2^{4}k_1^{4}k_2^{2}+192\,{d}^{4}t_1^{5}t_2^{5}k_1\,k_2^{5}-320\,{d}^{4}t_1^{5}t_2^{5}k_1^{3}k_2^{3}+12\,{d}^{4}t_1^{5}t_2^{5}k_1^{5}k_2+240\,{d}^{4}t_1^{4}t_2^{6}k_1^{2}k_2^{4}+16\,k_1^{6}t_2^{10}+288\,t_1^{2}t_2^{8}k_1^{6}+1296\,t_1^{4}t_2^{6}k_1^{6}-68\,t_1^{9}t_2\,k_1\,k_2^{5}-100\,{d}^{2}t_1^{6}k_2^{4}t_2^{4}k_3^{2}+265\,t_1^{4}t_2^{6}k_1^{2}k_2^{4}+770\,t_1^{6}t_2^{4}k_1^{2}k_2^{4}).
$}\\

\nin From this expression it is easy to check that $\operatorname{PP}_{(d,\bar k)}\left(
\operatorname{Con}_{\bar c}\left(
\operatorname{Res}_{t_0}\left(T_0(\bar k,\bar t_h),T(\bar c,d,\bar k,\bar t_h)\right)\right)\right)$ is the last factor in the above expression, and so:
\[
\deg_{\bar t}\left(\operatorname{PP}_{(d,\bar k)}\left(
\operatorname{Con}_{\bar c}\left(
\operatorname{Res}_{t_0}\left(T_0(\bar k,\bar t_h),T(\bar c,d,\bar k,\bar t_h)\right)\right)\right)\right)=10.
\]
Using Theorem \ref{thm:ch4:DegreeFormula} one has that the total offset degree is $\delta=10$. In fact, in this case, using elimination techniques, it is possible to check this result, computing the generic offset polynomial:\\

\nin{\small $
g(d,\bar x)=
-256 x_1^{10}-640 x_1^8 x_2^2-256 x_1^8 x_3^2+1408 x_1^8 d^2-400 x_1^6 x_2^4-384 x_1^6 x_2^2 x_3^2+3232 x_1^6 x_2^2 d^2+1152 x_1^6 x_3^2 d^2-3088 x_1^6 d^4+80 x_1^4 x_2^6-16 x_1^4 x_2^4 x_3^2+2448 x_1^4 x_2^4 d^2+1696 x_1^4 x_2^2 x_3^2 d^2-5904 x_1^4 x_2^2 d^4-1936 x_1^4 x_3^2 d^4+3376 x_1^4 d^6+80 x_1^2 x_2^8+96 x_1^2 x_2^6 x_3^2+832 x_1^2 x_2^6 d^2+736 x_1^2 x_2^4 x_3^2 d^2-3744 x_1^2 x_2^4 d^4-2272 x_1^2 x_2^2 x_3^2 d^4+4672 x_1^2 x_2^2 d^6+1440 x_1^2 x_3^2 d^6-1840 x_1^2 d^8-16 x_2^{10}-16 x_2^8 x_3^2+208 x_2^8 d^2+192 x_2^6 x_3^2 d^2-928 x_2^6 d^4-736 x_2^4 x_3^2 d^4+1696 x_2^4 d^6+960 x_2^2 x_3^2 d^6-1360 x_2^2 d^8-400 x_3^2 d^8+400 d^{10}+3200 x_1^8 x_3-320 x_1^6 x_2^2 x_3+3072 x_1^6 x_3^3-12608 x_1^6 x_3 d^2-1560 x_1^4 x_2^4 x_3-2944 x_1^4 x_2^2 x_3^3+496 x_1^4 x_2^2 x_3 d^2-9088 x_1^4 x_3^3 d^2+18088 x_1^4 x_3 d^4+1640 x_1^2 x_2^6 x_3+1696 x_1^2 x_2^4 x_3^3+7016 x_1^2 x_2^4 x_3 d^2+5184 x_1^2 x_2^2 x_3^3 d^2-2184 x_1^2 x_2^2 x_3 d^4+8480 x_1^2 x_3^3 d^4-11080 x_1^2 x_3 d^6-320 x_2^8 x_3-288 x_2^6 x_3^3+2912 x_2^6 x_3 d^2+2272 x_2^4 x_3^3 d^2-7072 x_2^4 x_3 d^4-3680 x_2^2 x_3^3 d^4+2080 x_2^2 x_3 d^6-2400 x_3^3 d^6+2400 x_3 d^8+2544 x_1^8-9144 x_1^6 x_2^2-10752 x_1^6 x_3^2-10520 x_1^6 d^2+4479 x_1^4 x_2^4-6976 x_1^4 x_2^2 x_3^2+25770 x_1^4 x_2^2 d^2-11776 x_1^4 x_3^4+21568 x_1^4 x_3^2 d^2+16583 x_1^4 d^4-684 x_1^2 x_2^6+10304 x_1^2 x_2^4 x_3^2+2700 x_1^2 x_2^4 d^2+9088 x_1^2 x_2^2 x_3^4+25056 x_1^2 x_2^2 x_3^2 d^2-23444 x_1^2 x_2^2 d^4+14720 x_1^2 x_3^4 d^2-7840 x_1^2 x_3^2 d^4-11980 x_1^2 d^6+24 x_2^8-2472 x_2^6 x_3^2+160 x_2^6 d^2-1936 x_2^4 x_3^4+14488 x_2^4 x_3^2 d^2-2752 x_2^4 d^4+8480 x_2^2 x_3^4 d^2-15160 x_2^2 x_3^2 d^4+6080 x_2^2 d^6-400 x_3^4 d^4-3000 x_3^2 d^6+3400 d^8-19008 x_1^6 x_3+25896 x_1^4 x_2^2 x_3+3328 x_1^4 x_3^3+44072 x_1^4 x_3 d^2-6534 x_1^2 x_2^4 x_3+23616 x_1^2 x_2^2 x_3^3+2484 x_1^2 x_2^2 x_3 d^2+15360 x_1^2 x_3^5+15680 x_1^2 x_3^3 d^2-31790 x_1^2 x_3 d^4+312 x_2^6 x_3-9112 x_2^4 x_3^3+3112 x_2^4 x_3 d^2-5760 x_2^2 x_3^5+31120 x_2^2 x_3^3 d^2-22360 x_2^2 x_3 d^4+9600 x_3^5 d^2-17400 x_3^3 d^4+7800 x_3 d^6-6240 x_1^6+3360 x_1^4 x_2^2+29984 x_1^4 x_3^2+15816 x_1^4 d^2-510 x_1^2 x_2^4-18472 x_1^2 x_2^2 x_3^2+11022 x_1^2 x_2^2 d^2+14080 x_1^2 x_3^4+8440 x_1^2 x_3^2 d^2-16520 x_1^2 d^4+15 x_2^6+1319 x_2^4 x_3^2-669 x_2^4 d^2-15680 x_2^2 x_3^4+15010 x_2^2 x_3^2 d^2+1045 x_2^2 d^4-6400 x_3^6+27200 x_3^4 d^2-28825 x_3^2 d^4+8025 d^6+14880 x_1^4 x_3-4800 x_1^2 x_2^2 x_3-13120 x_1^2 x_3^3+14320 x_1^2 x_3 d^2+270 x_2^4 x_3+1720 x_2^2 x_3^3-4270 x_2^2 x_3 d^2-9600 x_3^5+17400 x_3^3 d^2-7800 x_3 d^4-400 x_1^4+200 x_1^2 x_2^2-11040 x_1^2 x_3^2+7840 x_1^2 d^2-25 x_2^4+1440 x_2^2 x_3^2-140 x_2^2 d^2-400 x_3^4-3000 x_3^2 d^2+3400 d^4+800 x_1^2 x_3-200 x_2^2 x_3+2400 x_3^3-2400 x_3 d^2-400 x_3^2+400 d^2
$}\\[3mm]
\nin This is, as predicted by our formula, a polynomial of degree $10$ in $\bar x$.
\end{Example}
\begin{center}
\rule{2cm}{0.5pt}
\quad\\
\end{center}
\vspace{3mm}
\begin{Example}\label{exam:ch4:NonProperCircularParaboloid}
To illustrate the behavior of the degree formula in the case of non-proper parametrizations, let us consider the surface $\Sigma$ defined by the parametrization:
\[P(t_1,t_2)=(t_1^3,t_2,t_1^6+t_2^2).\]
This is a parametrization of the circular paraboloid with implicit equation $y_3=y_1^2+y_2^2$; the parametrization $P$ has been obtained by replacing $t_1$ with $t_1^3$ in the usual proper parametrization $\tilde P$ of $\Sigma$, which is given by:
\[\tilde P(t_1,t_2)=(t_1,t_2,t_1^2+t_2^2).\]
Thus, the tracing index of the parametrization $P$ in this example is $\mu=3$.
Computing with $P$ we obtain the following associated normal vector:
\[N(\bar t_h)=(-2t_1^3,-2t_0^2t_2,t_0^3).\]
\nin Then the auxiliary curves are:\\[3mm]
\nin $T_0(\bar t_h)=
-k_2\,t_1^{3}+k_1\,t_0^{2}{t_2}
$\\[3mm]
\nin $T_1(\bar t_h)=
4\,t_1^{18}{k_2}^{2}+12\,t_1^{12}{k_2}^{2}t_2^{2}t_0^{4}-8\,t_1^{12}k_2\,k_3\,{\it t2}\,t_0^{5}+12\,t_1^{6}{k_2}^{2}t_2^{4}t_0^{8}-16\,t_1^{6}k_2\,t_2^{3}t_0^{9}k_3+4\,t_1^{6}{k_3}^{2}t_2^{2}t_0^{10}+4\,t_2^{6}t_0^{12}{k_2}^{2}-8\,t_2^{5}t_0^{13}k_2\,k_3+4\,t_2^{4}t_0^{14}{k_3}^{2}+t_0^{6}{k_2}^{2}t_1^{12}+2\,t_0^{10}{k_2}^{2}t_1^{6}t_2^{2}-2\,t_0^{11}k_2\,t_1^{6}k_3\,{t_2}+t_0^{14}{k_2}^{2}t_2^{4}-2\,t_0^{15}k_2\,t_2^{3}k_3+t_0^{16}{k_3}^{2}t_2^{2}-{d}^{2}t_0^{18}{k_2}^{2}-4\,{d}^{2}t_0^{17}k_2\,k_3\,{\it t2}-4\,{d}^{2}t_0^{16}{k_3}^{2}t_2^{2}
$\\[3mm]
\nin $T_2(\bar t_h)=
4\,t_1^{12}{k_3}^{2}t_0^{6}-8\,t_1^{15}k_3\,t_0^{3}k_1-16\,t_1^{9}k_3\,t_0^{7}k_1\,t_2^{2}+4\,t_1^{18}{k_1}^{2}+12\,t_1^{12}{k_1}^{2}t_2^{2}t_0^{4}+12\,t_1^{6}{k_1}^{2}t_2^{4}t_0^{8}+4\,t_1^{6}{k_3}^{2}t_2^{2}t_0^{10}-8\,t_2^{4}t_0^{11}k_3\,t_1^{3}k_1+4\,t_2^{6}t_0^{12}{k_1}^{2}+t_0^{12}{k_3}^{2}t_1^{6}-2\,t_0^{9}k_3\,t_1^{9}k_1-2\,t_0^{13}k_3\,t_1^{3}k_1\,t_2^{2}+t_0^{6}{k_1}^{2}t_1^{12}+2\,t_0^{10}{k_1}^{2}t_1^{6}t_2^{2}+t_0^{14}{k_1}^{2}t_2^{4}-4\,{d}^{2}t_0^{12}t_1^{6}{k_3}^{2}-4\,{d}^{2}t_0^{15}k_3\,{{\it t1
}}^{3}k_1-{d}^{2}t_0^{18}{k_1}^{2}
$\\[3mm]
\nin $T_3(\bar t_h)=
t_0^{6} ( k_2\,t_1^{3}-k_1\,t_0^{2}{t_2} ) ^{2} ( 4\,t_1^{6}+t_0^{6}+4\,t_2^{2}t_0^{4}-4\,{d}^{2}t_0^{6})
$\\[3mm]
Denoting, as in the degree formula,
\[
R(\bar c,d,\bar k,\bar t)=\operatorname{Res}_{t_0}\left(T_0(\bar k,\bar t_h),T(\bar c,d,\bar k,\bar t_h)\right),
\]
one has:\\
$
R(\bar c,d,\bar k,\bar t)=
( k_1^{2}a_2+k_2^{2}a_1) ^{2}t_1^{36}( -8\,k_2^{12}t_1^{12}k_1^{2}t_2^{6}k_3^{4}{d}^{2}+16\,k_1^{18}t_2^{18}+16\,k_2^{12}k_1^{6}t_2^{18}+96\,k_2^{10}k_1^{8}t_2^{18}+240\,k_2^{8}k_1^{10}t_2^{18}+320\,k_2^{6}k_1^{12}t_2^{18}+240\,k_1^{14}t_2^{18}k_2^{4}+96\,k_1^{16}t_2^{18}k_2^{2}+t_1^{18}{d}^{4}k_2^{18}+t_1^{6}k_1^{4}t_2^{12}k_2^{14}+8\,t_1^{3}k_1^{5}t_2^{15}k_2^{13}+4\,t_1^{6}k_1^{6}t_2^{12}k_2^{12}+40\,t_1^{3}k_1^{7}{{\it t2}
}^{15}k_2^{11}+6\,t_1^{6}k_1^{8}t_2^{12}k_2^{10}+80\,t_1^{3}k_1^{9}t_2^{15}k_2^{9}+80\,t_1^{3}k_1^{11}t_2^{15}k_2^{7}+4\,k_2^{8}t_1^{6}k_1^{10}t_2^{12}+40\,k_2^{5}t_1^{3}k_1^{13}t_2^{15}+k_2^{6}t_1^{6}k_1^{12}t_2^{12}+8\,k_2^{3}t_1^{3}k_1^{15}t_2^{15}-32\,k_2^{6}t_1^{6}k_1^{10}t_2^{12}{d}^{2}k_3^{2}+16\,k_2^{10}t_1^{6}k_1^{4}t_2^{12}k_3^{4}+32\,k_2^{8}t_1^{6}k_1^{6}t_2^{12}k_3^{4}-192\,k_2^{7}t_1^{3}k_1^{9}t_2^{15}k_3^{2}-128\,k_2^{5}t_1^{3}k_1^{11}t_2^{15}k_3^{2}+16\,k_2^{6}t_1^{6}k_1^{8}t_2^{12}k_3^{4}-32\,k_2^{3}t_1^{3}k_1^{13}t_2^{15}k_3^{2}-4\,k_2^{11}t_1^{9}k_1^{5}t_2^{9}k_3^{2}+8\,k_2^{11}t_1^{9}k_1^{3}t_2^{9}k_3^{4}-2\,k_2^{9}t_1^{9}k_1^{7}t_2^{9}k_3^{2}+8\,k_2^{9}t_1^{9}k_1^{5}t_2^{9}k_3^{4}-48\,k_2^{8}t_1^{6}k_1^{8}t_2^{12}k_3^{2}-16\,k_2^{6}t_1^{6}k_1^{10}t_2^{12}k_3^{2}-32\,k_2^{12}k_1^{4}t_2^{12}t_1^{6}{d}^{2}{{\it
k3}}^{2}-32\,k_2^{11}k_1^{5}t_2^{15}t_1^{3}k_3^{2}-128\,k_2^{9}k_1^{7}t_2^{15}t_1^{3}k_3^{2}+16\,k_2^{12}t_1^{12}k_1^{2}t_2^{6}{d}^{4}k_3^{4}-144\,k_2^{11}t_1^{9}k_1^{5}t_2^{9}{d}^{2}k_3^{2}-32\,k_2^{11}t_1^{9}k_1^{3}t_2^{9}{d}^{2}k_3^{4}-96\,k_2^{10}t_1^{6}k_1^{6}t_2^{12}{d}^{2}k_3^{2}-2\,t_1^{12}{d}^{2}k_2^{16}k_1^{2}t_2^{6}-2\,t_1^{15}{d}^{2}k_2^{15}k_1\,t_2^{3}k_3^{2}-8\,t_1^{9}{d}^{2}k_2^{15}k_1^{3}t_2^{9}-8\,t_1^{15}{d}^{4}k_2^{15}k_1\,t_2^{3}k_3^{2}-4\,t_1^{12}{d}^{2}k_2^{14}k_1^{4}t_2^{6}-24\,t_1^{12}{d}^{2}k_2^{14}{{
\it k1}}^{2}t_2^{6}k_3^{2}-24\,t_1^{9}{d}^{2}k_2^{13}k_1^{5}t_2^{9}-2\,t_1^{12}{d}^{2}k_2^{12}k_1^{6}t_2^{6}-24\,t_1^{12}{d}^{2}k_2^{12}k_1^{4}t_2^{6}k_3^{2}-24\,t_1^{9}{d}^{2}k_2^{11}k_1^{7}t_2^{9}-8\,t_1^{9}{d}^{2}k_2^{9}k_1^{9}t_2^{9}-2\,t_1^{9}k_1^{3}t_2^{9}k_2^{13}k_3^{2}-72\,t_1^{9}k_1^{3}t_2^{9}k_2^{13}{d}^{2}k_3^{2}-16\,t_1^{6}k_1^{4}t_2^{12}k_2^{12}k_3^{2}-48\,t_1^{6}k_1^{6}t_2^{12}k_2^{10}k_3^{2}+k_2^{12}t_1^{12}k_1^{2}t_2^{6}k_3^{4}-72\,k_2^{9}t_1^{9}k_1^{7}t_2^{9}{d}^{2}k_3^{2}-32\,k_2^{9}t_1^{9}k_1^{5}t_2^{9}{d}^{2}k_3^{4}-96\,k_2^{8}t_1^{6}k_1^{8}t_2^{12}{d}^{2}k_3^{2}).
$\\[3mm]
\nin From this expression it is easy to check that
\[
\deg_{\bar t}\left(\operatorname{PP}_{(d,\bar k)}\left(
\operatorname{Con}_{\bar c}\left(
\operatorname{Res}_{t_0}\left(T_0(\bar k,\bar t_h),T(\bar c,d,\bar k,\bar t_h)\right)\right)\right)\right)=18.
\]
This agrees with the expected result $\mu\cdot\delta=3\cdot 6=18.$
\end{Example}
\begin{center}
\rule{2cm}{0.5pt}
\quad\\
\end{center}
\vspace{3mm}
\begin{Example}\label{exam:ch4:WhitneyUmbrella}
Let $\Sigma$ be the surface (Whitney Umbrella) with implicit equation $y_1^2 - y_2^2y_3=0$. A proper rational parametrization of $\Sigma$ is given by:
\[P(t_1,t_2)=(t_1t_2,t_2,t_1).\]
This surface is illustrated in Figure \ref{fig:WhitneyUmbrella}.
\begin{figure}
\centering
\includegraphics[height=8cm]{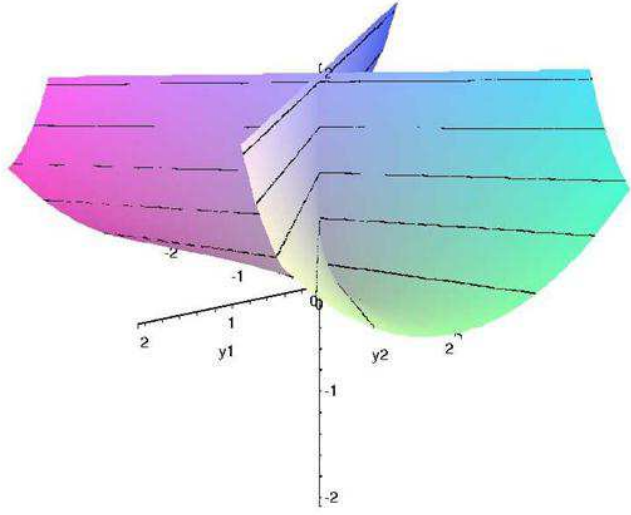}
\caption{The Whitney Umbrella}
\label{fig:WhitneyUmbrella}
\end{figure}
The homogeneous associated normal vector is
\[N(\bar t_h)=(2t_1t_2,-2t_1^2,-t_0t_2).\]
\nin Then the auxiliary curves are:\\[3mm]
\nin $T_0(\bar t_h)=-k_1\,t_0^{2}t_2^{2}+2\,k_1\,t_1^{4}+t_1\,t_0^{2}k_2\,t_2-2\,t_1^{3}t_0\,k_3+2\,t_1^{3}t_2\,k_2-2\,t_1\,t_2^{2}k_3\,t_0$\\[3mm]
\nin$T_1(\bar t_h)=4\,t_1^{6}t_2^{2}k_2^{2}-8\,t_1^{4}t_2^{3}k_2\,k_3\,t_0+4\,t_1^{2}t_2^{4}k_3^{2}t_0^{2}+4\,t_1^{8}k_2^{2}-8\,t_1^{6}k_2\,k_3\,t_0\,t_2+4\,t_1^{4}k_3^{2}t_0^{2}t_2^{2}+t_0^{2}t_2^{2}k_2^{2}t_1^{4}-2\,t_0^{3}t_2^{3}k_2\,t_1^{2}k_3+t_0^{4}t_2^{4}k_3^{2}-{d}^{2}t_2^{6}k_2^{2}t_0^{2}+4\,{d}^{2}t_2^{5}k_2\,t_1^{2}k_3\,t_0-4\,{d}^{2}t_2^{4}k_3^{2}t_1^{4}$\\[3mm]
\nin$T_2(\bar t_h)=4\,t_1^{4}k_3^{2}t_0^{2}t_2^{2}-8\,t_1^{5}t_2^{2}k_3\,t_0\,k_1+4\,t_1^{6}t_2^{2}k_1^{2}+4\,t_1^{6}k_3^{2}t_0^{2}-8\,t_1^{7}k_3\,t_0\,k_1+4\,t_1^{8}k_1^{2}+t_0^{4}t_2^{2}k_3^{2}t_1^{2}-2\,t_0^{3}t_2^{2}k_3\,t_1^{3}k_1+t_0^{2}t_2^{2}k_1^{2}t_1^{4}-4\,{d}^{2}t_2^{6}t_1^{2}k_3^{2}-4\,{d}^{2}t_2^{6}t_1\,k_3\,k_1\,t_0-{d}^{2}t_2^{6}k_1^{2}t_0^{2}$\\[3mm]
\nin$T_3(\bar t_h)=4\,t_0^{2}t_2^{2}k_2^{2}t_1^{4}-8\,t_1^{3}t_2^{3}k_2\,t_0^{2}k_1+4\,t_1^{2}t_2^{4}k_1^{2}t_0^{2}+4\,t_1^{6}k_2^{2}t_0^{2}-8\,t_1^{5}k_2\,t_0^{2}k_1\,t_2+4\,t_0^{2}t_2^{2}k_1^{2}t_1^{4}+t_0^{4}t_2^{2}k_2^{2}t_1^{2}-2\,t_0^{4}t_2^{3}k_2\,t_1\,k_1+t_0^{4}t_2^{4}k_1^{2}-4\,{d}^{2}t_2^{6}k_2^{2}t_1^{2}-8\,{d}^{2}t_2^{5}t_1^{3}k_2\,k_1-4\,{d}^{2}t_2^{4}k_1^{2}t_1^{4}
$\\[3mm]
Denoting, as in the degree formula,
\[
R(\bar c,d,\bar k,\bar t)=\operatorname{Res}_{t_0}\left(T_0(\bar k,\bar t_h),T(\bar c,d,\bar k,\bar t_h)\right),
\]
one has:\\

\nin{\small $
R(\bar c,d,\bar k,\bar t)=4\,t_2^{2}\,t_1^{4}\,( -28\,t_1^{6}t_2^{8}k_3^{2}{d}^{2}k_1^{2}+4\,k_2^{4}t_1^{8}t_2^{6}{d}^{2}+k_1^{4}{d}^{4}t_2^{12}t_1^{2}-6\,{d}^{2}k_1^{4}t_2^{8}t_1^{6}+k_2^{2}k_1^{2}{d}^{4}t_2^{14}-4\,k_1^{4}{d}^{2}t_2^{10}t_1^{4}+12\,k_1^{4}t_2^{6}t_1^{8}-20\,t_1^{4}t_2^{10}{d}^{2}k_1^{2}k_3^{2}-4\,k_2^{2}k_1^{2}{d}^{4}t_2^{12}t_1^{2}-4\,k_2^{2}t_1^{8}t_2^{6}k_3^{2}{d}^{2}-24\,k_2\,t_1^{11}t_2^{3}k_1\,k_3^{2}-2\,k_2\,{d}^{4}k_1^{3}t_2^{11}t_1^{3}+4\,k_2\,k_1^{3}{d}^{2}t_2^{9}t_1^{5}+2\,k_2\,k_1^{3}{d}^{4}t_2^{13}t_1+16\,k_2\,t_1^{5}t_2^{9}{d}^{2}k_3^{2}k_1-8\,k_2\,t_1^{13}t_2\,k_3^{2}k_1-8\,k_2\,t_1^{7}t_2^{7}k_3^{2}k_1+32\,k_2\,t_1^{7}t_2^{7}k_3^{2}{d}^{2}k_1+16\,k_2\,{d}^{2}k_1^{3}t_2^{7}t_1^{7}+16\,k_2\,k_3^{2}k_1\,{d}^{2}t_2^{5}t_1^{9}-24\,k_2\,k_3^{2}k_1\,t_2^{5}t_1^{9}-12\,t_1^{8}t_2^{6}k_3^{2}{d}^{2}k_1^{2}-4\,t_1^{10}t_2^{4}k_3^{4}{d}^{2}+4\,t_1^{12}t_2^{2}k_1^{2}k_3^{2}+4\,t_1^{6}t_2^{8}k_3^{2}k_1^{2}-16\,t_1^{8}t_2^{6}k_3^{4}{d}^{2}-4\,t_1^{2}t_2^{12}{d}^{2}k_3^{4}+12\,t_1^{8}t_2^{6}k_3^{2}k_1^{2}+12\,t_1^{10}t_2^{4}k_3^{2}k_1^{2}-6\,k_2^{3}k_1\,t_2^{5}t_1^{9}+12\,t_1^{10}t_2^{4}k_3^{2}k_2^{2}+4\,t_1^{8}t_2^{6}k_3^{2}k_2^{2}-16\,t_1^{4}t_2^{10}k_3^{4}{d}^{2}+4\,t_1^{4}t_2^{10}k_3^{2}{d}^{2}k_2^{2}-24\,t_1^{6}t_2^{8}k_3^{4}{d}^{2}+t_1^{2}t_2^{12}{d}^{4}k_2^{4}+2\,t_1^{6}t_2^{8}{d}^{2}k_2^{4}+4\,t_1^{6}t_2^{8}k_3^{2}{d}^{2}k_2^{2}-20\,k_2^{3}k_1\,t_2\,t_1^{13}+4\,k_2^{4}t_1^{14}+9\,k_1^{4}t_2^{4}t_1^{10}-22\,k_2^{3}k_1\,t_2^{3}t_1^{11}+37\,k_2^{2}k_1^{2}t_2^{2}t_1^{12}+44\,k_2^{2}k_1^{2}t_2^{4}t_1^{10}+10\,k_2^{2}k_1^{2}{d}^{2}t_2^{10}t_1^{4}+13\,k_2^{2}k_1^{2}t_2^{6}t_1^{8}-12\,k_2\,k_1^{3}t_2^{7}t_1^{7}-30\,k_2\,k_1^{3}t_2^{3}t_1^{11}-38\,k_2\,k_1^{3}t_2^{5}t_1^{9}-4\,k_2\,{d}^{2}k_1^{3}t_2^{11}t_1^{3}+t_1^{10}t_2^{4}k_2^{4}-4\,t_1^{2}t_2^{12}{d}^{2}k_1^{2}k_3^{2}+4\,k_1^{4}t_2^{8}t_1^{6}+4\,k_2^{2}t_1^{14}k_3^{2}+4\,k_2^{4}t_1^{12}t_2^{2}+12\,t_1^{12}t_2^{2}k_3^{2}k_2^{2}-8\,k_2^{3}k_1\,{d}^{2}t_2^{9}t_1^{5}-12\,k_2^{3}k_1\,{d}^{2}t_2^{7}t_1^{7}+2\,k_2^{3}k_1\,{d}^{4}t_2^{11}t_1^{3}+4\,k_2^{3}k_1\,{d}^{2}t_2^{5}t_1^{9}-2\,k_2^{3}k_1\,{d}^{4}t_2^{13}t_1+8\,k_2^{2}k_1^{2}{d}^{2}t_2^{8}t_1^{6}-14\,k_2^{2}k_1^{2}{d}^{2}t_2^{6}t_1^{8}-4\,k_2^{2}t_1^{10}t_2^{4}{d}^{2}k_3^{2}+k_2^{2}k_1^{2}{d}^{4}t_2^{10}t_1^{4})\,\cdot\,(8\,c_3\,c_1\,k_1\,k_3^{2}k_2\,t_1\,t_2^{3}-4\,c_3\,c_1\,k_1\,t_1^{3}k_2^{3}t_2-8\,c_3\,c_1\,k_1\,t_1^{3}k_3^{2}k_2\,t_2+4\,c_3\,c_1\,k_3^{2}k_2^{2}t_2^{4}-4\,c_3\,c_1\,k_2^{4}t_1^{2}t_2^{2}+4\,c_3\,c_1\,t_1^{4}k_2^{2}k_3^{2}-4\,c_1^{2}k_3^{2}k_2^{2}t_1^{2}t_2^{2}+8\,c_3^{2}k_1^{3}t_1\,k_2\,t_2^{3}-8\,c_3^{2}k_1^{3}t_1^{3}t_2\,k_2+4\,c_3^{2}k_2^{2}k_1^{2}t_2^{4}+4\,c_3^{2}k_1^{2}t_1^{4}k_2^{2}-16\,c_3^{2}k_1^{2}t_2^{2}k_2^{2}t_1^{2}-8\,c_3^{2}k_1\,k_2^{3}t_1\,t_2^{3}+8\,c_3^{2}k_1\,t_1^{3}k_2^{3}t_2+4\,{{c_3}}^{2}k_2^{4}t_1^{2}t_2^{2}+c_2^{2}t_1^{2}t_2^{2}k_1^{4}+4\,c_2^{2}t_2^{2}t_1^{2}k_3^{4}+4\,c_2\,{c_1}\,k_3^{4}t_2^{4}+4\,c_2\,c_1\,k_3^{4}t_1^{4}+4\,c_3\,c_2\,t_1^{2}t_2^{2}k_1^{4}+4\,c_3\,c_2\,k_1^{3}t_1\,k_2\,t_2^{3}-4\,c_3\,c_2\,k_1^{3}t_1^{3}t_2\,k_2+4\,c_3\,c_2\,k_1^{2}t_2^{4}k_3^{2}-4\,c_3\,c_2\,k_1^{2}t_2^{2}k_2^{2}t_1^{2}+4\,c_3\,c_2\,k_1^{2}t_1^{4}k_3^{2}-8\,c_3\,c_2\,k_1\,k_3^{2}k_2\,t_1\,t_2^{3}+8\,c_3\,c_2\,k_1\,t_1^{3}k_3^{2}k_2\,t_2+8\,c_3\,c_2\,k_3^{2}k_2^{2}t_1^{2}t_2^{2}+4\,c_2^{2}k_1^{2}t_1^{2}k_3^{2}t_2^{2}+4\,c_2\,c_1\,k_1^{2}t_1^{2}k_3^{2}t_2^{2}+2\,c_2\,c_1\,k_1^{2}t_2^{2}k_2^{2}t_1^{2}+8\,c_2\,c_1\,k_1\,k_3^{2}k_2\,t_1\,t_2^{3}-8\,c_2\,c_1\,k_1\,t_1^{3}k_3^{2}k_2\,t_2-4\,c_2\,{c_1}\,k_3^{2}k_2^{2}t_1^{2}t_2^{2}+c_1^{2}k_2^{4}t_1^{2}t_2^{2}+4\,c_1^{2}t_2^{2}t_1^{2}k_3^{4}+4\,{{c_3}}^{2}t_1^{2}t_2^{2}k_1^{4}+8\,c_3\,c_1\,k_1^{2}t_1^{2}k_3^{2}t_2^{2}+4\,c_3\,c_1\,k_1^{2}t_2^{2}k_2^{2}t_1^{2}+4\,c_3\,c_1\,k_1\,k_2^{3}t_1\,t_2^{3}).
$}

\nin From this expression it is easy to check that
\[
\deg_{\bar t}\left(\operatorname{PP}_{(d,\bar k)}\left(
\operatorname{Con}_{\bar c}\left(
\operatorname{Res}_{t_0}\left(T_0(\bar k,\bar t_h),T(\bar c,d,\bar k,\bar t_h)\right)\right)\right)\right)=14,
\]
and then, using Theorem \ref{thm:ch4:DegreeFormula} one concludes that the total offset degree in $\bar x$ is $\delta=14.$

\nin In fact, in this case, using elimination techniques, it is possible to check this result,  computing the generic offset polynomial (see Appendix \ref{Ap3-ComplementsToSomeProofs}, page \pageref{Ap3-ComplementsToSomeProofs-Whitney}). This is indeed a polynomial of degree $14$ in $\bar x$.
\end{Example}
\begin{center}
\rule{2cm}{0.5pt}
\quad\\
\end{center}

\section*{Appendix: Computational Complements}\label{Ap3-ComplementsToSomeProofs}
\subsection*{Coefficients $c^{(i)}_j$ in Lemma \ref{lem:ch4:AuxiliaryPolynomialsBelongToEliminationIdeal}
(page \pageref{lem:ch4:AuxiliaryPolynomialsBelongToEliminationIdeal}).}
\addcontentsline{toc}{section}{Appendix: Computational Complements}

The polynomials $s_i$ can be expressed as follows:
\[
s_i=c^{(i)}_1\,b^{P}+c^{(i)}_2\,\normal^{P}_{(1,2)}+c^{(i)}_3\,\normal^{P}_{(1,3)}+c^{(i)}_4\,\normal^{P}_{(2,3)}
+c^{(i)}_5\,w^{P}+c^{(i)}_6\,\ell_1+c^{(i)}_7\,\ell_2+c^{(i)}_8\,\ell_3
\]
where $c^{(i)}_j\in\C[d,\bar k,l,r,\bar t,\bar x]$ for $i=0,\ldots,3$ , $j=1,\ldots,8$ are the following polynomials:

\nin $c^{(0)}_1=0$

\nin $c^{(0)}_2=k_3$

\nin $c^{(0)}_3=-k_2$

\nin $c^{(0)}_4=k_1$

\nin $c^{(0)}_5=P_0\,n_2\,k_3-P_0\,n_3\,k_2$

\nin $c^{(0)}_6=-P_0\,n_1\,k_3+P_0\,n_3\,k_1$

\nin $c^{(0)}_7=P_0\,n_1\,k_2-P_0\,n_2\,k_1$

\nin $c^{(0)}_8=0$

\nin $c^{(1)}_1=-2\,P_0\,n_1^3\,n_2\,k_1\,k_2\,r+P_0\,n_1^2\,n_2^2\,k_1^2\,r-2\,P_0\,n_1\,n_2^3\,k_1\,k_2\,r-2\,P_0\,n_1\,n_2\,n_3^2\,k_1\,k_2\,r+P_0\,n_2^4\,k_1^2\,r+P_0\,n_2^2\,n_3^2\,k_1^2\,r+n_1^2\,k_2^2$

\nin $c^{(1)}_2=2\,P_1\,P_0\,n_1^3\,k_1\,k_2\,r-P_1\,P_0\,n_1^2\,n_2\,k_1^2\,r+2\,P_1\,P_0\,n_1\,n_2^2\,k_1\,k_2\,r-P_1\,P_0\,n_2^3\,k_1^2\,r+P_1\,n_2\,k_2^2-P_2\,P_0\,n_1^3\,k_1^2\,r-P_2\,P_0\,n_1\,n_2^2\,k_1^2\,r+P_2\,n_1\,k_2^2-2\,P_2\,n_2\,k_1\,k_2+2\,P_3\,P_0\,n_1^2\,n_3\,k_1\,k_2\,r+P_0^2\,n_1^3\,k_1^2\,r\,x_2-2\,P_0^2\,n_1^3\,k_1\,k_2\,r\,x_1+P_0^2\,n_1^2\,n_2\,k_1^2\,r\,x_1-2\,P_0^2\,n_1^2\,n_3\,k_1\,k_2\,r\,x_3+P_0^2\,n_1\,n_2^2\,k_1^2\,r\,x_2-2\,P_0^2\,n_1\,n_2^2\,k_1\,k_2\,r\,x_1+P_0^2\,n_2^3\,k_1^2\,r\,x_1-P_0\,n_1\,k_2^2\,x_2+P_0\,n_2\,k_2^2\,x_1$

\nin $c^{(1)}_3=2\,P_1\,P_0\,n_1\,n_2\,n_3\,k_1\,k_2\,r-P_1\,P_0\,n_2^2\,n_3\,k_1^2\,r+P_1\,n_3\,k_2^2-2\,P_2\,P_0\,n_1^2\,n_3\,k_1\,k_2\,r-P_3\,P_0\,n_1\,n_2^2\,k_1^2\,r+P_3\,n_1\,k_2^2-2\,P_3\,n_2\,k_1\,k_2+2\,P_0^2\,n_1^2\,n_3\,k_1\,k_2\,r\,x_2+P_0^2\,n_1\,n_2^2\,k_1^2\,r\,x_3-2\,P_0^2\,n_1\,n_2\,n_3\,k_1\,k_2\,r\,x_1+P_0^2\,n_2^2\,n_3\,k_1^2\,r\,x_1-P_0\,n_1\,k_2\,k_3\,x_2+2\,P_0\,n_2\,k_2\,k_3\,x_1-P_0\,n_3\,k_2^2\,x_1$

\nin $c^{(1)}_4=-2\,P_1\,n_3\,k_1\,k_2+P_2\,P_0\,n_1^2\,n_3\,k_1^2\,r+P_2\,n_3\,k_1^2+P_3\,P_0\,n_1^2\,n_2\,k_1^2\,r+P_3\,n_2\,k_1^2-P_0^2\,n_1^2\,n_2\,k_1^2\,r\,x_3-P_0^2\,n_1^2\,n_3\,k_1^2\,r\,x_2-P_0\,n_2\,k_1\,k_3\,x_1+P_0\,n_3\,k_1\,k_2\,x_1$

\nin $c^{(1)}_5=2\,P_1\,P_0\,n_1^2\,k_2^2-2\,P_1\,P_0\,n_2\,n_3\,k_2\,k_3+2\,P_1\,P_0\,n_3^2\,k_2^2-2\,P_2\,P_0\,n_1^2\,k_1\,k_2+2\,P_2\,P_0\,n_1\,n_2\,k_2^2-2\,P_2\,P_0\,n_2^2\,k_1\,k_2+P_2\,P_0\,n_2\,n_3\,k_1\,k_3-P_2\,P_0\,n_3^2\,k_1\,k_2+2\,P_3\,P_0\,n_1\,n_2\,k_2\,k_3-P_3\,P_0\,n_2^2\,k_1\,k_3-P_3\,P_0\,n_2\,n_3\,k_1\,k_2+P_0^2\,n_1^2\,k_1\,k_2\,x_2-P_0^2\,n_1^2\,k_2^2\,x_1-2\,P_0^2\,n_1\,n_2\,k_2^2\,x_2-2\,P_0^2\,n_1\,n_2\,k_3^2\,x_2+P_0^2\,n_2^2\,k_1\,k_2\,x_2+P_0^2\,n_2^2\,k_1\,k_3\,x_3+P_0^2\,n_2^2\,k_2^2\,x_1+2\,P_0^2\,n_2\,n_3\,k_2\,k_3\,x_1-P_0^2\,n_3^2\,k_2^2\,x_1$

\nin $c^{(1)}_6=-2\,P_1\,P_0\,n_1^2\,k_1\,k_2+P_1\,P_0\,n_1\,n_3\,k_2\,k_3-2\,P_1\,P_0\,n_3^2\,k_1\,k_2+2\,P_2\,P_0\,n_1^2\,k_1^2-2\,P_2\,P_0\,n_1\,n_2\,k_1\,k_2+2\,P_2\,P_0\,n_2^2\,k_1^2+P_2\,P_0\,n_3^2\,k_1^2-P_3\,P_0\,n_1^2\,k_2\,k_3+P_3\,P_0\,n_2\,n_3\,k_1^2-P_0^2\,n_1^2\,k_1^2\,x_2+P_0^2\,n_1^2\,k_1\,k_2\,x_1+P_0^2\,n_1^2\,k_2\,k_3\,x_3+2\,P_0^2\,n_1\,n_2\,k_1\,k_2\,x_2-2\,P_0^2\,n_1\,n_2\,k_1\,k_3\,x_3+2\,P_0^2\,n_1\,n_2\,k_3^2\,x_1-P_0^2\,n_1\,n_3\,k_2\,k_3\,x_1-P_0^2\,n_2^2\,k_1^2\,x_2-P_0^2\,n_2^2\,k_1\,k_2\,x_1-P_0^2\,n_2\,n_3\,k_1\,k_3\,x_1+P_0^2\,n_3^2\,k_1\,k_2\,x_1$

\nin $c^{(1)}_7=-P_1\,P_0\,n_1\,n_3\,k_2^2+2\,P_1\,P_0\,n_2\,n_3\,k_1\,k_2-P_2\,P_0\,n_2\,n_3\,k_1^2+P_3\,P_0\,n_1^2\,k_2^2-2\,P_3\,P_0\,n_1\,n_2\,k_1\,k_2+P_3\,P_0\,n_2^2\,k_1^2-P_0^2\,n_1^2\,k_2^2\,x_3+2\,P_0^2\,n_1\,n_2\,k_1\,k_2\,x_3+2\,P_0^2\,n_1\,n_2\,k_1\,k_3\,x_2-2\,P_0^2\,n_1\,n_2\,k_2\,k_3\,x_1+P_0^2\,n_1\,n_3\,k_2^2\,x_1-P_0^2\,n_2^2\,k_1^2\,x_3-P_0^2\,n_2\,n_3\,k_1\,k_2\,x_1$

\nin $c^{(1)}_8=2\,P_1\,P_2\,n_1^2\,k_1\,k_2-2\,P_1\,P_0\,n_1^2\,k_1\,k_2\,x_2-P_2^2\,n_1^2\,k_1^2+2\,P_2^2\,n_1\,n_2\,k_1\,k_2-P_2^2\,n_2^2\,k_1^2+2\,P_2\,P_0\,n_1^2\,k_1^2\,x_2-2\,P_2\,P_0\,n_1^2\,k_1\,k_2\,x_1-4\,P_2\,P_0\,n_1\,n_2\,k_1\,k_2\,x_2+2\,P_2\,P_0\,n_2^2\,k_1^2\,x_2+2\,P_3^2\,n_1\,n_2\,k_1\,k_2-P_3^2\,n_2^2\,k_1^2-4\,P_3\,P_0\,n_1\,n_2\,k_1\,k_2\,x_3+2\,P_3\,P_0\,n_2^2\,k_1^2\,x_3-P_0^2\,n_1^2\,k_1^2\,x_2^2+2\,P_0^2\,n_1^2\,k_1\,k_2\,x_1\,x_2-2\,P_0^2\,n_1\,n_2\,k_1\,k_2\,d^2+2\,P_0^2\,n_1\,n_2\,k_1\,k_2\,x_2^2+2\,P_0^2\,n_1\,n_2\,k_1\,k_2\,x_3^2+P_0^2\,n_2^2\,k_1^2\,d^2-P_0^2\,n_2^2\,k_1^2\,x_2^2-P_0^2\,n_2^2\,k_1^2\,x_3^2$

\nin $c^{(2)}_1=P_0\,n_1^2\,n_2^2\,k_3^2\,r-2\,P_0\,n_1^2\,n_2\,n_3\,k_2\,k_3\,r+P_0\,n_1^2\,n_3^2\,k_2^2\,r+P_0\,n_2^4\,k_3^2\,r-2\,P_0\,n_2^3\,n_3\,k_2\,k_3\,r+P_0\,n_2^2\,n_3^2\,k_2^2\,r+P_0\,n_2^2\,n_3^2\,k_3^2\,r-2\,P_0\,n_2\,n_3^3\,k_2\,k_3\,r+P_0\,n_3^4\,k_2^2\,r$

\nin $c^{(2)}_2=-P_1\,P_0\,n_1^2\,n_2\,k_3^2\,r+2\,P_1\,P_0\,n_1^2\,n_3\,k_2\,k_3\,r-P_1\,P_0\,n_2^3\,k_3^2\,r+2\,P_1\,P_0\,n_2^2\,n_3\,k_2\,k_3\,r-P_2\,P_0\,n_1^3\,k_3^2\,r-P_2\,P_0\,n_1\,n_2^2\,k_3^2\,r+P_0^2\,n_1^3\,k_3^2\,r\,x_2+P_0^2\,n_1^2\,n_2\,k_3^2\,r\,x_1-2\,P_0^2\,n_1^2\,n_3\,k_2\,k_3\,r\,x_1+P_0^2\,n_1\,n_2^2\,k_3^2\,r\,x_2+P_0^2\,n_2^3\,k_3^2\,r\,x_1-2\,P_0^2\,n_2^2\,n_3\,k_2\,k_3\,r\,x_1$

\nin $c^{(2)}_3=-P_1\,P_0\,n_1^2\,n_3\,k_2^2\,r-P_1\,P_0\,n_2^2\,n_3\,k_2^2\,r-P_1\,P_0\,n_2^2\,n_3\,k_3^2\,r+2\,P_1\,P_0\,n_2\,n_3^2\,k_2\,k_3\,r-P_1\,P_0\,n_3^3\,k_2^2\,r+2\,P_2\,P_0\,n_1^3\,k_2\,k_3\,r+2\,P_2\,P_0\,n_1\,n_2^2\,k_2\,k_3\,r-P_3\,P_0\,n_1^3\,k_2^2\,r-P_3\,P_0\,n_1\,n_2^2\,k_2^2\,r-P_3\,P_0\,n_1\,n_2^2\,k_3^2\,r+2\,P_3\,P_0\,n_1\,n_2\,n_3\,k_2\,k_3\,r-P_3\,P_0\,n_1\,n_3^2\,k_2^2\,r+P_0^2\,n_1^3\,k_2^2\,r\,x_3-2\,P_0^2\,n_1^3\,k_2\,k_3\,r\,x_2+P_0^2\,n_1^2\,n_3\,k_2^2\,r\,x_1+P_0^2\,n_1\,n_2^2\,k_2^2\,r\,x_3-2\,P_0^2\,n_1\,n_2^2\,k_2\,k_3\,r\,x_2+P_0^2\,n_1\,n_2^2\,k_3^2\,r\,x_3-2\,P_0^2\,n_1\,n_2\,n_3\,k_2\,k_3\,r\,x_3+P_0^2\,n_1\,n_3^2\,k_2^2\,r\,x_3+P_0^2\,n_2^2\,n_3\,k_2^2\,r\,x_1+P_0^2\,n_2^2\,n_3\,k_3^2\,r\,x_1-2\,P_0^2\,n_2\,n_3^2\,k_2\,k_3\,r\,x_1+P_0^2\,n_3^3\,k_2^2\,r\,x_1$

\nin $c^{(2)}_4=2\,P_2\,P_0\,n_1^2\,n_2\,k_2\,k_3\,r-P_2\,P_0\,n_1^2\,n_3\,k_2^2\,r+P_2\,P_0\,n_1^2\,n_3\,k_3^2\,r+2\,P_2\,P_0\,n_2^3\,k_2\,k_3\,r-P_2\,P_0\,n_2^2\,n_3\,k_2^2\,r+2\,P_2\,P_0\,n_2\,n_3^2\,k_2\,k_3\,r-P_2\,P_0\,n_3^3\,k_2^2\,r+P_2\,n_3\,k_3^2-P_3\,P_0\,n_1^2\,n_2\,k_2^2\,r+P_3\,P_0\,n_1^2\,n_2\,k_3^2\,r-2\,P_3\,P_0\,n_1^2\,n_3\,k_2\,k_3\,r-P_3\,P_0\,n_2^3\,k_2^2\,r-P_3\,P_0\,n_2\,n_3^2\,k_2^2\,r+P_3\,n_2\,k_3^2-2\,P_3\,n_3\,k_2\,k_3+P_0^2\,n_1^2\,n_2\,k_2^2\,r\,x_3-2\,P_0^2\,n_1^2\,n_2\,k_2\,k_3\,r\,x_2-P_0^2\,n_1^2\,n_2\,k_3^2\,r\,x_3+P_0^2\,n_1^2\,n_3\,k_2^2\,r\,x_2+2\,P_0^2\,n_1^2\,n_3\,k_2\,k_3\,r\,x_3-P_0^2\,n_1^2\,n_3\,k_3^2\,r\,x_2+P_0^2\,n_2^3\,k_2^2\,r\,x_3-2\,P_0^2\,n_2^3\,k_2\,k_3\,r\,x_2+P_0^2\,n_2^2\,n_3\,k_2^2\,r\,x_2+P_0^2\,n_2\,n_3^2\,k_2^2\,r\,x_3-2\,P_0^2\,n_2\,n_3^2\,k_2\,k_3\,r\,x_2+P_0^2\,n_3^3\,k_2^2\,r\,x_2-P_0\,n_2\,k_3^2\,x_3+P_0\,n_3\,k_3^2\,x_2$

\nin
$c^{(2)}_5=0$

\nin $c^{(2)}_6=2\,P_2\,P_0\,n_1^2\,k_3^2+2\,P_2\,P_0\,n_2^2\,k_3^2-2\,P_3\,P_0\,n_1^2\,k_2\,k_3-2\,P_3\,P_0\,n_2^2\,k_2\,k_3+2\,P_3\,P_0\,n_2\,n_3\,k_3^2-2\,P_3\,P_0\,n_3^2\,k_2\,k_3+P_0^2\,n_1^2\,k_2\,k_3\,x_3-P_0^2\,n_1^2\,k_3^2\,x_2+P_0^2\,n_2^2\,k_2\,k_3\,x_3-P_0^2\,n_2^2\,k_3^2\,x_2-2\,P_0^2\,n_2\,n_3\,k_3^2\,x_3+P_0^2\,n_3^2\,k_2\,k_3\,x_3+P_0^2\,n_3^2\,k_3^2\,x_2$

\nin $c^{(2)}_7=-2\,P_2\,P_0\,n_1^2\,k_2\,k_3-2\,P_2\,P_0\,n_2^2\,k_2\,k_3+2\,P_3\,P_0\,n_1^2\,k_2^2+2\,P_3\,P_0\,n_2^2\,k_2^2-2\,P_3\,P_0\,n_2\,n_3\,k_2\,k_3+2\,P_3\,P_0\,n_3^2\,k_2^2-P_0^2\,n_1^2\,k_2^2\,x_3+P_0^2\,n_1^2\,k_2\,k_3\,x_2-P_0^2\,n_2^2\,k_2^2\,x_3+P_0^2\,n_2^2\,k_2\,k_3\,x_2+2\,P_0^2\,n_2\,n_3\,k_2\,k_3\,x_3-P_0^2\,n_3^2\,k_2^2\,x_3-P_0^2\,n_3^2\,k_2\,k_3\,x_2$

\nin $c^{(2)}_8=-P_2^2\,n_1^2\,k_3^2-P_2^2\,n_2^2\,k_3^2+2\,P_2\,P_3\,n_1^2\,k_2\,k_3+2\,P_2\,P_3\,n_2^2\,k_2\,k_3-2\,P_2\,P_0\,n_1^2\,k_2\,k_3\,x_3+2\,P_2\,P_0\,n_1^2\,k_3^2\,x_2-2\,P_2\,P_0\,n_2^2\,k_2\,k_3\,x_3+2\,P_2\,P_0\,n_2^2\,k_3^2\,x_2-P_3^2\,n_1^2\,k_2^2-P_3^2\,n_2^2\,k_2^2-P_3^2\,n_2^2\,k_3^2+2\,P_3^2\,n_2\,n_3\,k_2\,k_3-P_3^2\,n_3^2\,k_2^2+2\,P_3\,P_0\,n_1^2\,k_2^2\,x_3-2\,P_3\,P_0\,n_1^2\,k_2\,k_3\,x_2+2\,P_3\,P_0\,n_2^2\,k_2^2\,x_3-2\,P_3\,P_0\,n_2^2\,k_2\,k_3\,x_2+2\,P_3\,P_0\,n_2^2\,k_3^2\,x_3-4\,P_3\,P_0\,n_2\,n_3\,k_2\,k_3\,x_3+2\,P_3\,P_0\,n_3^2\,k_2^2\,x_3-P_0^2\,n_1^2\,k_2^2\,x_3^2+2\,P_0^2\,n_1^2\,k_2\,k_3\,x_2\,x_3-P_0^2\,n_1^2\,k_3^2\,x_2^2-P_0^2\,n_2^2\,k_2^2\,x_3^2+2\,P_0^2\,n_2^2\,k_2\,k_3\,x_2\,x_3+P_0^2\,n_2^2\,k_3^2\,d^2-P_0^2\,n_2^2\,k_3^2\,x_2^2-P_0^2\,n_2^2\,k_3^2\,x_3^2-2\,P_0^2\,n_2\,n_3\,k_2\,k_3\,d^2+2\,P_0^2\,n_2\,n_3\,k_2\,k_3\,x_3^2+P_0^2\,n_3^2\,k_2^2\,d^2-P_0^2\,n_3^2\,k_2^2\,x_3^2$

\nin $c^{(3)}_1=-2\,P_0\,n_1^3\,n_3\,k_1\,k_3\,r+P_0\,n_1^2\,n_3^2\,k_1^2\,r-2\,P_0\,n_1\,n_2^2\,n_3\,k_1\,k_3\,r-2\,P_0\,n_1\,n_3^3\,k_1\,k_3\,r+P_0\,n_2^2\,n_3^2\,k_1^2\,r+P_0\,n_3^4\,k_1^2\,r+n_1^2\,k_3^2$

\nin $c^{(3)}_2=2\,P_1\,P_0\,n_1\,n_2\,n_3\,k_1\,k_3\,r+P_1\,n_2\,k_3^2-2\,P_2\,P_0\,n_1^2\,n_3\,k_1\,k_3\,r+P_2\,n_1\,k_3^2-2\,P_3\,P_0\,n_1^2\,n_2\,k_1\,k_3\,r-2\,P_3\,n_2\,k_1\,k_3+2\,P_0^2\,n_1^2\,n_2\,k_1\,k_3\,r\,x_3+2\,P_0^2\,n_1^2\,n_3\,k_1\,k_3\,r\,x_2-2\,P_0^2\,n_1\,n_2\,n_3\,k_1\,k_3\,r\,x_1-P_0\,n_1\,k_3^2\,x_2+P_0\,n_2\,k_3^2\,x_1$

\nin $c^{(3)}_3=2\,P_1\,P_0\,n_1^3\,k_1\,k_3\,r-P_1\,P_0\,n_1^2\,n_3\,k_1^2\,r+2\,P_1\,P_0\,n_1\,n_3^2\,k_1\,k_3\,r-P_1\,P_0\,n_2^2\,n_3\,k_1^2\,r-P_1\,P_0\,n_3^3\,k_1^2\,r+P_1\,n_3\,k_3^2+4\,P_2\,P_0\,n_1^2\,n_2\,k_1\,k_3\,r-P_3\,P_0\,n_1^3\,k_1^2\,r-P_3\,P_0\,n_1\,n_2^2\,k_1^2\,r-P_3\,P_0\,n_1\,n_3^2\,k_1^2\,r+P_3\,n_1\,k_3^2-2\,P_3\,n_3\,k_1\,k_3+P_0^2\,n_1^3\,k_1^2\,r\,x_3-2\,P_0^2\,n_1^3\,k_1\,k_3\,r\,x_1-4\,P_0^2\,n_1^2\,n_2\,k_1\,k_3\,r\,x_2+P_0^2\,n_1^2\,n_3\,k_1^2\,r\,x_1+P_0^2\,n_1\,n_2^2\,k_1^2\,r\,x_3+P_0^2\,n_1\,n_3^2\,k_1^2\,r\,x_3-2\,P_0^2\,n_1\,n_3^2\,k_1\,k_3\,r\,x_1+P_0^2\,n_2^2\,n_3\,k_1^2\,r\,x_1+P_0^2\,n_3^3\,k_1^2\,r\,x_1-P_0\,n_1\,k_3^2\,x_3+P_0\,n_3\,k_3^2\,x_1$

\nin $c^{(3)}_4=-P_2\,P_0\,n_1^2\,n_3\,k_1^2\,r+2\,P_2\,P_0\,n_1\,n_2^2\,k_1\,k_3\,r+2\,P_2\,P_0\,n_1\,n_3^2\,k_1\,k_3\,r-P_2\,P_0\,n_2^2\,n_3\,k_1^2\,r-P_2\,P_0\,n_3^3\,k_1^2\,r-P_3\,P_0\,n_1^2\,n_2\,k_1^2\,r-P_3\,P_0\,n_2^3\,k_1^2\,r-P_3\,P_0\,n_2\,n_3^2\,k_1^2\,r+P_0^2\,n_1^2\,n_2\,k_1^2\,r\,x_3+P_0^2\,n_1^2\,n_3\,k_1^2\,r\,x_2-2\,P_0^2\,n_1\,n_2^2\,k_1\,k_3\,r\,x_2-2\,P_0^2\,n_1\,n_3^2\,k_1\,k_3\,r\,x_2+P_0^2\,n_2^3\,k_1^2\,r\,x_3+P_0^2\,n_2^2\,n_3\,k_1^2\,r\,x_2+P_0^2\,n_2\,n_3^2\,k_1^2\,r\,x_3+P_0^2\,n_3^3\,k_1^2\,r\,x_2$

\nin $c^{(3)}_5=2\,P_1\,P_0\,n_1^2\,k_3^2+2\,P_2\,P_0\,n_1\,n_2\,k_3^2-2\,P_3\,P_0\,n_1^2\,k_1\,k_3+2\,P_3\,P_0\,n_1\,n_3\,k_3^2-2\,P_3\,P_0\,n_2^2\,k_1\,k_3-2\,P_3\,P_0\,n_3^2\,k_1\,k_3+P_0^2\,n_1^2\,k_1\,k_3\,x_3-P_0^2\,n_1^2\,k_3^2\,x_1-2\,P_0^2\,n_1\,n_2\,k_3^2\,x_2-2\,P_0^2\,n_1\,n_3\,k_3^2\,x_3+P_0^2\,n_2^2\,k_1\,k_3\,x_3+P_0^2\,n_2^2\,k_3^2\,x_1+P_0^2\,n_3^2\,k_1\,k_3\,x_3+P_0^2\,n_3^2\,k_3^2\,x_1$

\nin
$c^{(3)}_6=0$

\nin $c^{(3)}_7=-2\,P_1\,P_0\,n_1^2\,k_1\,k_3-2\,P_2\,P_0\,n_1\,n_2\,k_1\,k_3+2\,P_3\,P_0\,n_1^2\,k_1^2-2\,P_3\,P_0\,n_1\,n_3\,k_1\,k_3+2\,P_3\,P_0\,n_2^2\,k_1^2+2\,P_3\,P_0\,n_3^2\,k_1^2-P_0^2\,n_1^2\,k_1^2\,x_3+P_0^2\,n_1^2\,k_1\,k_3\,x_1+2\,P_0^2\,n_1\,n_2\,k_1\,k_3\,x_2+2\,P_0^2\,n_1\,n_3\,k_1\,k_3\,x_3-P_0^2\,n_2^2\,k_1^2\,x_3-P_0^2\,n_2^2\,k_1\,k_3\,x_1-P_0^2\,n_3^2\,k_1^2\,x_3-P_0^2\,n_3^2\,k_1\,k_3\,x_1$

\nin $c^{(3)}_8=2\,P_1\,P_3\,n_1^2\,k_1\,k_3-2\,P_1\,P_0\,n_1^2\,k_1\,k_3\,x_3+2\,P_2\,P_3\,n_1\,n_2\,k_1\,k_3-2\,P_2\,P_0\,n_1\,n_2\,k_1\,k_3\,x_3-P_3^2\,n_1^2\,k_1^2+2\,P_3^2\,n_1\,n_3\,k_1\,k_3-P_3^2\,n_2^2\,k_1^2-P_3^2\,n_3^2\,k_1^2+2\,P_3\,P_0\,n_1^2\,k_1^2\,x_3-2\,P_3\,P_0\,n_1^2\,k_1\,k_3\,x_1-2\,P_3\,P_0\,n_1\,n_2\,k_1\,k_3\,x_2-4\,P_3\,P_0\,n_1\,n_3\,k_1\,k_3\,x_3+2\,P_3\,P_0\,n_2^2\,k_1^2\,x_3+2\,P_3\,P_0\,n_3^2\,k_1^2\,x_3-P_0^2\,n_1^2\,k_1^2\,x_3^2+2\,P_0^2\,n_1^2\,k_1\,k_3\,x_1\,x_3+2\,P_0^2\,n_1\,n_2\,k_1\,k_3\,x_2\,x_3-2\,P_0^2\,n_1\,n_3\,k_1\,k_3\,d^2+2\,P_0^2\,n_1\,n_3\,k_1\,k_3\,x_3^2-P_0^2\,n_2^2\,k_1^2\,x_3^2+P_0^2\,n_3^2\,k_1^2\,d^2-P_0^2\,n_3^2\,k_1^2\,x_3^2$

\subsection*{Generic Offset Polynomial for Example \ref{exam:ch4:WhitneyUmbrella} (page \pageref{exam:ch4:WhitneyUmbrella})}\label{Ap3-ComplementsToSomeProofs-Whitney}

\nin {$ g(d,x_1,x_2,x_3)=-16{\,}x_1^2{\,}x_2^{10}{\,}x_3^2+16{\,}x_1^2{\,}x_2^{10}{\,}d^2-48{\,}x_1^2{\,}x_2^8{\,}x_3^4+128{\,}x_1^2{\,}x_2^8{\,}x_3^2{\,}d^2-80{\,}x_1^2{\,}x_2^8{\,}d^4-48{\,}x_1^2{\,}x_2^6{\,}x_3^6+240{\,}x_1^2{\,}x_2^6{\,}x_3^4{\,}d^2-352{\,}x_1^2{\,}x_2^6{\,}x_3^2{\,}d^4+160{\,}x_1^2{\,}x_2^6{\,}d^6-16{\,}x_1^2{\,}x_2^4{\,}x_3^8+160{\,}x_1^2{\,}x_2^4{\,}x_3^6{\,}d^2-432{\,}x_1^2{\,}x_2^4{\,}x_3^4{\,}d^4+448{\,}x_1^2{\,}x_2^4{\,}x_3^2{\,}d^6-160{\,}x_1^2{\,}x_2^4{\,}d^8+32{\,}x_1^2{\,}x_2^2{\,}x_3^8{\,}d^2-176{\,}x_1^2{\,}x_2^2{\,}x_3^6{\,}d^4+336{\,}x_1^2{\,}x_2^2{\,}x_3^4{\,}d^6-272{\,}x_1^2{\,}x_2^2{\,}x_3^2{\,}d^8+80{\,}x_1^2{\,}x_2^2{\,}d^{10}-16{\,}x_1^2{\,}x_3^8{\,}d^4+64{\,}x_1^2{\,}x_3^6{\,}d^6-96{\,}x_1^2{\,}x_3^4{\,}d^8+64{\,}x_1^2{\,}x_3^2{\,}d^{10}-16{\,}x_1^2{\,}d^{12}-16{\,}x_2^{12}{\,}x_3^2+16{\,}x_2^{12}{\,}d^2-64{\,}x_2^{10}{\,}x_3^4+160{\,}x_2^{10}{\,}x_3^2{\,}d^2-96{\,}x_2^{10}{\,}d^4-96{\,}x_2^8{\,}x_3^6+416{\,}x_2^8{\,}x_3^4{\,}d^2-560{\,}x_2^8{\,}x_3^2{\,}d^4+240{\,}x_2^8{\,}d^6-64{\,}x_2^6{\,}x_3^8+448{\,}x_2^6{\,}x_3^6{\,}d^2-1024{\,}x_2^6{\,}x_3^4{\,}d^4+960{\,}x_2^6{\,}x_3^2{\,}d^6-320{\,}x_2^6{\,}d^8-16{\,}x_2^4{\,}x_3^{10}+208{\,}x_2^4{\,}x_3^8{\,}d^2-768{\,}x_2^4{\,}x_3^6{\,}d^4+1216{\,}x_2^4{\,}x_3^4{\,}d^6-880{\,}x_2^4{\,}x_3^2{\,}d^8+240{\,}x_2^4{\,}d^{10}+32{\,}x_2^2{\,}x_3^{10}{\,}d^2-224{\,}x_2^2{\,}x_3^8{\,}d^4+576{\,}x_2^2{\,}x_3^6{\,}d^6-704{\,}x_2^2{\,}x_3^4{\,}d^8+416{\,}x_2^2{\,}x_3^2{\,}d^{10}-96{\,}x_2^2{\,}d^{12}-16{\,}x_3^{10}{\,}d^4+80{\,}x_3^8{\,}d^6-160{\,}x_3^6{\,}d^8+160{\,}x_3^4{\,}d^{10}-80{\,}x_3^2{\,}d^{12}+16{\,}d^{14}+32{\,}x_1^4{\,}x_2^8{\,}x_3+744{\,}x_1^4{\,}x_2^6{\,}x_3^3-808{\,}x_1^4{\,}x_2^6{\,}x_3{\,}d^2-120{\,}x_1^4{\,}x_2^4{\,}x_3^5-1080{\,}x_1^4{\,}x_2^4{\,}x_3^3{\,}d^2+1232{\,}x_1^4{\,}x_2^4{\,}x_3{\,}d^4+32{\,}x_1^4{\,}x_2^2{\,}x_3^7+408{\,}x_1^4{\,}x_2^2{\,}x_3^5{\,}d^2-272{\,}x_1^4{\,}x_2^2{\,}x_3^3{\,}d^4-168{\,}x_1^4{\,}x_2^2{\,}x_3{\,}d^6+32{\,}x_1^4{\,}x_3^7{\,}d^2-352{\,}x_1^4{\,}x_3^5{\,}d^4+608{\,}x_1^4{\,}x_3^3{\,}d^6-288{\,}x_1^4{\,}x_3{\,}d^8+32{\,}x_1^2{\,}x_2^{10}{\,}x_3+968{\,}x_1^2{\,}x_2^8{\,}x_3^3-1032{\,}x_1^2{\,}x_2^8{\,}x_3{\,}d^2+720{\,}x_1^2{\,}x_2^6{\,}x_3^5-3176{\,}x_1^2{\,}x_2^6{\,}x_3^3{\,}d^2+2488{\,}x_1^2{\,}x_2^6{\,}x_3{\,}d^4-184{\,}x_1^2{\,}x_2^4{\,}x_3^7-392{\,}x_1^2{\,}x_2^4{\,}x_3^5{\,}d^2+2168{\,}x_1^2{\,}x_2^4{\,}x_3^3{\,}d^4-1592{\,}x_1^2{\,}x_2^4{\,}x_3{\,}d^6+32{\,}x_1^2{\,}x_2^2{\,}x_3^9+632{\,}x_1^2{\,}x_2^2{\,}x_3^7{\,}d^2-1672{\,}x_1^2{\,}x_2^2{\,}x_3^5{\,}d^4+1320{\,}x_1^2{\,}x_2^2{\,}x_3^3{\,}d^6-312{\,}x_1^2{\,}x_2^2{\,}x_3{\,}d^8+32{\,}x_1^2{\,}x_3^9{\,}d^2-512{\,}x_1^2{\,}x_3^7{\,}d^4+1344{\,}x_1^2{\,}x_3^5{\,}d^6-1280{\,}x_1^2{\,}x_3^3{\,}d^8+416{\,}x_1^2{\,}x_3{\,}d^{10}+160{\,}x_2^{10}{\,}x_3^3-160{\,}x_2^{10}{\,}x_3{\,}d^2+224{\,}x_2^8{\,}x_3^5-736{\,}x_2^8{\,}x_3^3{\,}d^2+512{\,}x_2^8{\,}x_3{\,}d^4-32{\,}x_2^6{\,}x_3^7-256{\,}x_2^6{\,}x_3^5{\,}d^2+736{\,}x_2^6{\,}x_3^3{\,}d^4-448{\,}x_2^6{\,}x_3{\,}d^6-96{\,}x_2^4{\,}x_3^9+544{\,}x_2^4{\,}x_3^7{\,}d^2-928{\,}x_2^4{\,}x_3^5{\,}d^4+608{\,}x_2^4{\,}x_3^3{\,}d^6-128{\,}x_2^4{\,}x_3{\,}d^8+224{\,}x_2^2{\,}x_3^9{\,}d^2-1024{\,}x_2^2{\,}x_3^7{\,}d^4+1728{\,}x_2^2{\,}x_3^5{\,}d^6-1280{\,}x_2^2{\,}x_3^3{\,}d^8+352{\,}x_2^2{\,}x_3{\,}d^{10}-128{\,}x_3^9{\,}d^4+512{\,}x_3^7{\,}d^6-768{\,}x_3^5{\,}d^8+512{\,}x_3^3{\,}d^{10}-128{\,}x_3{\,}d^{12}-16{\,}x_1^6{\,}x_2^6-2073{\,}x_1^6{\,}x_2^4{\,}x_3^2+873{\,}x_1^6{\,}x_2^4{\,}d^2+384{\,}x_1^6{\,}x_2^2{\,}x_3^4-324{\,}x_1^6{\,}x_2^2{\,}x_3^2{\,}d^2+2052{\,}x_1^6{\,}x_2^2{\,}d^4-16{\,}x_1^6{\,}x_3^6+504{\,}x_1^6{\,}x_3^4{\,}d^2-1728{\,}x_1^6{\,}x_3^2{\,}d^4+216{\,}x_1^6{\,}d^6-16{\,}x_1^4{\,}x_2^8-2797{\,}x_1^4{\,}x_2^6{\,}x_3^2+1277{\,}x_1^4{\,}x_2^6{\,}d^2-2529{\,}x_1^4{\,}x_2^4{\,}x_3^4+2602{\,}x_1^4{\,}x_2^4{\,}x_3^2{\,}d^2+2615{\,}x_1^4{\,}x_2^4{\,}d^4+560{\,}x_1^4{\,}x_2^2{\,}x_3^6+980{\,}x_1^4{\,}x_2^2{\,}x_3^4{\,}d^2+552{\,}x_1^4{\,}x_2^2{\,}x_3^2{\,}d^4-3372{\,}x_1^4{\,}x_2^2{\,}d^6-16{\,}x_1^4{\,}x_3^8+744{\,}x_1^4{\,}x_3^6{\,}d^2-3992{\,}x_1^4{\,}x_3^4{\,}d^4+3768{\,}x_1^4{\,}x_3^2{\,}d^6-504{\,}x_1^4{\,}d^8-620{\,}x_1^2{\,}x_2^8{\,}x_3^2+364{\,}x_1^2{\,}x_2^8{\,}d^2-1060{\,}x_1^2{\,}x_2^6{\,}x_3^4+252{\,}x_1^2{\,}x_2^6{\,}x_3^2{\,}d^2+1256{\,}x_1^2{\,}x_2^6{\,}d^4-824{\,}x_1^2{\,}x_2^4{\,}x_3^6+1436{\,}x_1^2{\,}x_2^4{\,}x_3^4{\,}d^2+2448{\,}x_1^2{\,}x_2^4{\,}x_3^2{\,}d^4-3252{\,}x_1^2{\,}x_2^4{\,}d^6+192{\,}x_1^2{\,}x_2^2{\,}x_3^8+1952{\,}x_1^2{\,}x_2^2{\,}x_3^6{\,}d^2-4032{\,}x_1^2{\,}x_2^2{\,}x_3^4{\,}d^4+608{\,}x_1^2{\,}x_2^2{\,}x_3^2{\,}d^6+1280{\,}x_1^2{\,}x_2^2{\,}d^8+224{\,}x_1^2{\,}x_3^8{\,}d^2-2432{\,}x_1^2{\,}x_3^6{\,}d^4+4544{\,}x_1^2{\,}x_3^4{\,}d^6-2688{\,}x_1^2{\,}x_3^2{\,}d^8+352{\,}x_1^2{\,}d^{10}+8{\,}x_2^{10}{\,}x_3^2-8{\,}x_2^{10}{\,}d^2-480{\,}x_2^8{\,}x_3^4+296{\,}x_2^8{\,}x_3^2{\,}d^2+184{\,}x_2^8{\,}d^4+296{\,}x_2^6{\,}x_3^6-8{\,}x_2^6{\,}x_3^4{\,}d^2+152{\,}x_2^6{\,}x_3^2{\,}d^4-440{\,}x_2^6{\,}d^6-240{\,}x_2^4{\,}x_3^8+104{\,}x_2^4{\,}x_3^6{\,}d^2+424{\,}x_2^4{\,}x_3^4{\,}d^4-584{\,}x_2^4{\,}x_3^2{\,}d^6+296{\,}x_2^4{\,}d^8+672{\,}x_2^2{\,}x_3^8{\,}d^2-1792{\,}x_2^2{\,}x_3^6{\,}d^4+1600{\,}x_2^2{\,}x_3^4{\,}d^6-512{\,}x_2^2{\,}x_3^2{\,}d^8+32{\,}x_2^2{\,}d^{10}-448{\,}x_3^8{\,}d^4+1408{\,}x_3^6{\,}d^6-1536{\,}x_3^4{\,}d^8+640{\,}x_3^2{\,}d^{10}-64{\,}d^{12}+2106{\,}x_1^8{\,}x_2^2{\,}x_3-216{\,}x_1^8{\,}x_3^3+1944{\,}x_1^8{\,}x_3{\,}d^2+2946{\,}x_1^6{\,}x_2^4{\,}x_3+3282{\,}x_1^6{\,}x_2^2{\,}x_3^3+54{\,}x_1^6{\,}x_2^2{\,}x_3{\,}d^2-312{\,}x_1^6{\,}x_3^5+4176{\,}x_1^6{\,}x_3^3{\,}d^2-5400{\,}x_1^6{\,}x_3{\,}d^4+760{\,}x_1^4{\,}x_2^6{\,}x_3+800{\,}x_1^4{\,}x_2^4{\,}x_3^3+56{\,}x_1^4{\,}x_2^4{\,}x_3{\,}d^2+1744{\,}x_1^4{\,}x_2^2{\,}x_3^5+2048{\,}x_1^4{\,}x_2^2{\,}x_3^3{\,}d^2-3696{\,}x_1^4{\,}x_2^2{\,}x_3{\,}d^4-96{\,}x_1^4{\,}x_3^7+2784{\,}x_1^4{\,}x_3^5{\,}d^2-8992{\,}x_1^4{\,}x_3^3{\,}d^4+5280{\,}x_1^4{\,}x_3{\,}d^6-16{\,}x_1^2{\,}x_2^8{\,}x_3+1362{\,}x_1^2{\,}x_2^6{\,}x_3^3-546{\,}x_1^2{\,}x_2^6{\,}x_3{\,}d^2-1560{\,}x_1^2{\,}x_2^4{\,}x_3^5+2880{\,}x_1^2{\,}x_2^4{\,}x_3^3{\,}d^2-2472{\,}x_1^2{\,}x_2^4{\,}x_3{\,}d^4+480{\,}x_1^2{\,}x_2^2{\,}x_3^7+2120{\,}x_1^2{\,}x_2^2{\,}x_3^5{\,}d^2-3408{\,}x_1^2{\,}x_2^2{\,}x_3^3{\,}d^4+1192{\,}x_1^2{\,}x_2^2{\,}x_3{\,}d^6+672{\,}x_1^2{\,}x_3^7{\,}d^2-5088{\,}x_1^2{\,}x_3^5{\,}d^4+6624{\,}x_1^2{\,}x_3^3{\,}d^6-2208{\,}x_1^2{\,}x_3{\,}d^8-72{\,}x_2^8{\,}x_3^3+72{\,}x_2^8{\,}x_3{\,}d^2+440{\,}x_2^6{\,}x_3^5+592{\,}x_2^6{\,}x_3^3{\,}d^2-1032{\,}x_2^6{\,}x_3{\,}d^4-320{\,}x_2^4{\,}x_3^7-864{\,}x_2^4{\,}x_3^5{\,}d^2+896{\,}x_2^4{\,}x_3^3{\,}d^4+288{\,}x_2^4{\,}x_3{\,}d^6+1120{\,}x_2^2{\,}x_3^7{\,}d^2-1440{\,}x_2^2{\,}x_3^5{\,}d^4+32{\,}x_2^2{\,}x_3^3{\,}d^6+288{\,}x_2^2{\,}x_3{\,}d^8-896{\,}x_3^7{\,}d^4+2176{\,}x_3^5{\,}d^6-1664{\,}x_3^3{\,}d^8+384{\,}x_3{\,}d^{10}-729{\,}x_1^{10}-1053{\,}x_1^8{\,}x_2^2-1377{\,}x_1^8{\,}x_3^2+2673{\,}x_1^8{\,}d^2-300{\,}x_1^6{\,}x_2^4+684{\,}x_1^6{\,}x_2^2{\,}x_3^2+3132{\,}x_1^6{\,}x_2^2{\,}d^2-888{\,}x_1^6{\,}x_3^4+5616{\,}x_1^6{\,}x_3^2{\,}d^2-3672{\,}x_1^6{\,}d^4+8{\,}x_1^4{\,}x_2^6-1500{\,}x_1^4{\,}x_2^4{\,}x_3^2+1072{\,}x_1^4{\,}x_2^4{\,}d^2+2232{\,}x_1^4{\,}x_2^2{\,}x_3^4+456{\,}x_1^4{\,}x_2^2{\,}x_3^2{\,}d^2-2576{\,}x_1^4{\,}x_2^2{\,}d^4-240{\,}x_1^4{\,}x_3^6+4384{\,}x_1^4{\,}x_3^4{\,}d^2-8272{\,}x_1^4{\,}x_3^2{\,}d^4+2336{\,}x_1^4{\,}d^6+276{\,}x_1^2{\,}x_2^6{\,}x_3^2-164{\,}x_1^2{\,}x_2^6{\,}d^2-1336{\,}x_1^2{\,}x_2^4{\,}x_3^4+1048{\,}x_1^2{\,}x_2^4{\,}x_3^2{\,}d^2-1024{\,}x_1^2{\,}x_2^4{\,}d^4+640{\,}x_1^2{\,}x_2^2{\,}x_3^6+480{\,}x_1^2{\,}x_2^2{\,}x_3^4{\,}d^2-112{\,}x_1^2{\,}x_2^2{\,}x_3^2{\,}d^4+272{\,}x_1^2{\,}x_2^2{\,}d^6+1120{\,}x_1^2{\,}x_3^6{\,}d^2-5760{\,}x_1^2{\,}x_3^4{\,}d^4+5088{\,}x_1^2{\,}x_3^2{\,}d^6-704{\,}x_1^2{\,}d^8-1{\,}x_2^8{\,}x_3^2+1{\,}x_2^8{\,}d^2+184{\,}x_2^6{\,}x_3^4-112{\,}x_2^6{\,}x_3^2{\,}d^2-72{\,}x_2^6{\,}d^4-240{\,}x_2^4{\,}x_3^6-896{\,}x_2^4{\,}x_3^4{\,}d^2+656{\,}x_2^4{\,}x_3^2{\,}d^4+480{\,}x_2^4{\,}d^6+1120{\,}x_2^2{\,}x_3^6{\,}d^2-480{\,}x_2^2{\,}x_3^4{\,}d^4-864{\,}x_2^2{\,}x_3^2{\,}d^6+224{\,}x_2^2{\,}d^8-1120{\,}x_3^6{\,}d^4+2080{\,}x_3^4{\,}d^6-1056{\,}x_3^2{\,}d^8+96{\,}d^{10}-648{\,}x_1^8{\,}x_3+834{\,}x_1^6{\,}x_2^2{\,}x_3-968{\,}x_1^6{\,}x_3^3+2664{\,}x_1^6{\,}x_3{\,}d^2-336{\,}x_1^4{\,}x_2^4{\,}x_3+1352{\,}x_1^4{\,}x_2^2{\,}x_3^3-1408{\,}x_1^4{\,}x_2^2{\,}x_3{\,}d^2-320{\,}x_1^4{\,}x_3^5+3456{\,}x_1^4{\,}x_3^3{\,}d^2-3776{\,}x_1^4{\,}x_3{\,}d^4+2{\,}x_1^2{\,}x_2^6{\,}x_3-464{\,}x_1^2{\,}x_2^4{\,}x_3^3+176{\,}x_1^2{\,}x_2^4{\,}x_3{\,}d^2+480{\,}x_1^2{\,}x_2^2{\,}x_3^5-568{\,}x_1^2{\,}x_2^2{\,}x_3^3{\,}d^2+1176{\,}x_1^2{\,}x_2^2{\,}x_3{\,}d^4+1120{\,}x_1^2{\,}x_3^5{\,}d^2-3776{\,}x_1^2{\,}x_3^3{\,}d^4+2144{\,}x_1^2{\,}x_3{\,}d^6+8{\,}x_2^6{\,}x_3^3-8{\,}x_2^6{\,}x_3{\,}d^2-96{\,}x_2^4{\,}x_3^5-256{\,}x_2^4{\,}x_3^3{\,}d^2+352{\,}x_2^4{\,}x_3{\,}d^4+672{\,}x_2^2{\,}x_3^5{\,}d^2-64{\,}x_2^2{\,}x_3^3{\,}d^4-608{\,}x_2^2{\,}x_3{\,}d^6-896{\,}x_3^5{\,}d^4+1280{\,}x_3^3{\,}d^6-384{\,}x_3{\,}d^8-216{\,}x_1^8+132{\,}x_1^6{\,}x_2^2-456{\,}x_1^6{\,}x_3^2+720{\,}x_1^6{\,}d^2-1{\,}x_1^4{\,}x_2^4+376{\,}x_1^4{\,}x_2^2{\,}x_3^2-388{\,}x_1^4{\,}x_2^2{\,}d^2-240{\,}x_1^4{\,}x_3^4+1416{\,}x_1^4{\,}x_3^2{\,}d^2-856{\,}x_1^4{\,}d^4-32{\,}x_1^2{\,}x_2^4{\,}x_3^2+20{\,}x_1^2{\,}x_2^4{\,}d^2+192{\,}x_1^2{\,}x_2^2{\,}x_3^4-288{\,}x_1^2{\,}x_2^2{\,}x_3^2{\,}d^2+416{\,}x_1^2{\,}x_2^2{\,}d^4+672{\,}x_1^2{\,}x_3^4{\,}d^2-1472{\,}x_1^2{\,}x_3^2{\,}d^4+416{\,}x_1^2{\,}d^6-16{\,}x_2^4{\,}x_3^4+8{\,}x_2^4{\,}x_3^2{\,}d^2+8{\,}x_2^4{\,}d^4+224{\,}x_2^2{\,}x_3^4{\,}d^2-64{\,}x_2^2{\,}x_3^2{\,}d^4-160{\,}x_2^2{\,}d^6-448{\,}x_3^4{\,}d^4+512{\,}x_3^2{\,}d^6-64{\,}d^8-96{\,}x_1^6{\,}x_3+40{\,}x_1^4{\,}x_2^2{\,}x_3-96{\,}x_1^4{\,}x_3^3+320{\,}x_1^4{\,}x_3{\,}d^2+32{\,}x_1^2{\,}x_2^2{\,}x_3^3-8{\,}x_1^2{\,}x_2^2{\,}x_3{\,}d^2+224{\,}x_1^2{\,}x_3^3{\,}d^2-352{\,}x_1^2{\,}x_3{\,}d^4+32{\,}x_2^2{\,}x_3^3{\,}d^2-32{\,}x_2^2{\,}x_3{\,}d^4-128{\,}x_3^3{\,}d^4+128{\,}x_3{\,}d^6-16{\,}x_1^6-16{\,}x_1^4{\,}x_3^2+48{\,}x_1^4{\,}d^2+32{\,}x_1^2{\,}x_3^2{\,}d^2-48{\,}x_1^2{\,}d^4-16{\,}x_3^2{\,}d^4+16{\,}d^6$.
}

\hrule
\paragraph{Contact Info:}\quad\\
\noindent{\small \rm
Departamento de Matem\'aticas, Facultad de Ciencias, Universidad de Alcal\'a\\
Ap. de Correos 20,E-28871 Alcal\'a de Henares (Madrid),SPAIN\\}
\href{mailto:fernando.sansegundo@uah.es}{fernando.sansegundo@uah.es } (corresponding author), \href{mailto:rafael.sendra@uah.es}{rafael.sendra@uah.es}

\paragraph{Funding:} This work has been partially supported by the research project {\sf MTM2008-04699-C03-01 ``Variedades paramétricas: algoritmos y aplicaciones''}, Ministerio de Ciencia e Innovación, Spain.

 \bibliographystyle{plain}
 \bibliography{Bibliography}

\begin{thebibliography}{10}

\bibitem{Alcazar2007}
J.G. Alc{\'a}zar.
\newblock {\em Effective Algorithms for the Study of the Topology of Algebraic
  Varieties, and Applications}.
\newblock PhD thesis, Universidad de Alcal{\'a}, 2007.

\bibitem{alcazar2008}
J.G. Alc{\'a}zar.
\newblock {\em Good global behavior of offsets to plane algebraic curves}.
\newblock {\em \em Journal of Symbolic Computation}, 43(9):659--680, 2008.

\bibitem{alcazar2008localsurfaces}
J.G. Alc{\'a}zar.
\newblock {\em Good local behavior of offsets to rational regular algebraic
  surfaces}.
\newblock {\em \em Journal of Symbolic Computation}, 43(12):845--857, 2008.

\bibitem{alcazarLocalCurves2009}
J.G. Alc{\'a}zar.
\newblock {\em Good Local Behavior of Offsets to Implicit Algebraic Curves}.
\newblock {\em \em Mathematics in Computer Science}, to appear.

\bibitem{alcazar2007local}
J.G. Alc{\'a}zar and J.R. Sendra.
\newblock {\em Local shape of offsets to algebraic curves}.
\newblock {\em \em Journal of Symbolic Computation}, 42(3):338--351, 2007.

\bibitem{anton2005offset}
F.~Anton, I.~Emiris, B.~Mourrain, and M.~Teillaud.
\newblock {\em The offset to an algebraic curve and an application to conics}.
\newblock {\em \em Lecture Notes in Computer Science}, 3480/2005:683--696,
  2005.

\bibitem{Arrondo1997}
E.~Arrondo, J.~Sendra, and J.R. Sendra.
\newblock {\em Parametric Generalized Offsets to Hypersurfaces}.
\newblock {\em {\em Journal of Symbolic Computation}}, 23(2-3):267--285, 1997.

\bibitem{Arrondo1999}
E.~Arrondo, J.~Sendra, and J.R. Sendra.
\newblock {\em Genus formula for generalized offset curves}.
\newblock {\em \em Journal of Pure and Applied Algebra}, 136(3):199--209, 1999.

\bibitem{Brieskorn1986}
E.~Brieskorn and H.~Kn{\"o}rrer.
\newblock {\em Plane algebraic curves}.
\newblock Birkh{\"a}user Verlag, Basel, 1986.

\bibitem{Castelnuovo1939}
G.~Castelnuovo.
\newblock {\em Sulle Superficie di Genere Zero}.
\newblock {\em Memorie Scelte}, pages 307--334, 1939.

\bibitem{CocoaSystem}
{CoCoA}Team.
\newblock {{\hbox{\rm C\kern-.13em o\kern-.07em C\kern-.13em o\kern-.15em A}}}:
  a system for doing {C}omputations in {C}ommutative {A}lgebra.
\newblock Available at \/ {\tt http://cocoa.dima.unige.it}.

\bibitem{Cox1997}
D.A. Cox, J.B. Little, and D.~O'Shea.
\newblock {\em Ideals, Varieties, and Algorithms: An Introduction to
  Computational Algebraic Geometry and Commutative Algebra}.
\newblock Springer, 2nd edition, 1997.

\bibitem{Farouki1990}
R.T. Farouki and C.A. Neff.
\newblock {\em Algebraic properties of plane offset curves}.
\newblock {\em \em Computer Aided Geometric Design}, 7:101--127, 1990.

\bibitem{Farouki1990a}
R.T. Farouki and C.A. Neff.
\newblock {\em Analytic properties of plane offset curves}.
\newblock {\em \em Computer Aided Geometric Design}, 7(1-4):83--99, 1990.

\bibitem{GPS}
G.-M. Greuel, G.~Pfister, and H.~Sch{\"o}nemann.
\newblock {\sc Singular} 3.1.0 --- {A} computer algebra system for polynomial
  computations.
\newblock 2009.
\newblock http://www.singular.uni-kl.de.

\bibitem{Harris1992}
J.~Harris.
\newblock {\em Algebraic Geometry: A First Course}.
\newblock Springer, 1992.

\bibitem{Hoffmann1989}
Ch.M. Hoffmann.
\newblock {\em Geometric and Solid Modelling: An Introduction}.
\newblock Morgan Kaufmann, San Mateo (CA), 1989.

\bibitem{hoschek1993fundamentals}
J.~Hoschek, D.~Lasser, and L.L. Schumaker.
\newblock {\em Fundamentals of Computer-Aided Geometric Design}.
\newblock AK Peters, Ltd., 1993.

\bibitem{Leibniz1692}
G.W. Leibniz.
\newblock {\em Generalia de natura linearum, anguloque contactus et osculi
  provocationibis aliisque cognatis et eorum usibus nonnullis}.
\newblock {\em \em Acta Eruditorum}, 1692.

\bibitem{loria1902}
G.~Loria and F.~Sch{\"u}tte.
\newblock {\em {Spezielle algebraische und transscendente ebene Kurven: Theorie
  und Geschichte}}.
\newblock BG Teubner, 1902.

\bibitem{Lu1995}
W.~L{\"u}.
\newblock {\em Offset-rational parametric plane curves}.
\newblock {\em \em Computer Aided Geometric Design}, 12(6):601--616, 1995.

\bibitem{Lu95TR}
W.~L{\"u}.
\newblock {\em Rational Parametrizations of Quadrics and their Offsets}.
\newblock {\em \em Technical Reports, Institut für Geometrie, Technische
  Universität Wien}, 24, 1995.

\bibitem{patrikalakis2002sic}
N.M. Patrikalakis and T.~Maekawa.
\newblock {\em Shape Interrogation for Computer Aided Design and
  Manufacturing}.
\newblock Springer, 2002.

\bibitem{perezdiaz2006ppr}
S.~P{\'e}rez-D{\'\i}az.
\newblock {\em On the problem of proper reparametrization for rational curves
  and surfaces}.
\newblock {\em {\em Computer Aided Geometric Design}}, 23(4):307--323, 2006.

\bibitem{Perez-Diaz2002}
S.~P{\'e}rez-D{\'\i}az, J.~Schicho, and J.R. Sendra.
\newblock {\em Properness and Inversion of Rational Parametrizations of
  Surfaces}.
\newblock {\em \em Applicable Algebra in Engineering, Communication and
  Computing}, 13(1):29--51, 2002.

\bibitem{SendraPerez2004JPAA-DegreeSurfaceParam}
S~P{\'e}rez-D{\'\i}az and J.R. Sendra.
\newblock {\em Computation of the degree of rational surface parametrizations}.
\newblock {\em Journal of Pure and Applied Algebra}, 193(1-3):99--121, 2004.

\bibitem{Pottmann1995}
H.~Pottmann.
\newblock {\em Rational curves and surfaces with rational offsets}.
\newblock {\em \em Computer Aided Geometric Design}, 12(2):175--192, 1995.

\bibitem{pottmann1996rational}
H.~Pottmann, W.~L{\"u}, and B.~Ravani.
\newblock {\em Rational ruled surfaces and their offsets}.
\newblock {\em \em Graphical Models and Image Processing}, 58(6):544--552,
  1996.

\bibitem{pottmann1998laguerre}
H.~Pottmann and M.~Peternell.
\newblock {\em Applications of Laguerre geometry in CAGD}.
\newblock {\em \em Computer Aided Geometric Design}, 15(2):165--186, 1998.

\bibitem{Roman2005}
S.~Roman.
\newblock {\em Advanced Linear Algebra}.
\newblock {\em \em Graduate Texts in Mathematics, vol. 135. Springer, New
  York}, 2005.

\bibitem{salmon-treatise}
G.~Salmon.
\newblock {\em {A Treatise on Higher Plane Curves, 1879}}.

\bibitem{PhdSanSegundo2010}
F.~San~Segundo.
\newblock {\em Effective Algorithms for the Study of the Degree of Algebraic
  Varieties in Offsetting Processes}.
\newblock PhD thesis, Universidad de Alcal\'a, 2010.

\bibitem{SS05}
F.~San~Segundo and J.R. Sendra.
\newblock {\em Degree formulae for offset curves}.
\newblock {\em \em J. Pure Appl. Algebra}, 195(3):301--335, 2005.

\bibitem{SS06}
F.~San~Segundo and J.R. Sendra.
\newblock {\em Partial Degree Formulae for Plane Offset Curves}.
\newblock {\em \em Arxiv preprint math.AG/0609137}, 2006.

\bibitem{SSS09}
F.~San~Segundo and J.R. Sendra.
\newblock {\em Partial Degree formulae for Plane Offset Curves}.
\newblock {\em \em Journal of Symbolic Computation}, 44(6):635--654, 2009.

\bibitem{schicho1998rps}
J.~Schicho.
\newblock {\em Rational Parametrization of Surfaces}.
\newblock {\em \em Journal of Symbolic Computation}, 26(1):1--29, 1998.

\bibitem{Sendra1999}
J.~Sendra.
\newblock {\em Algoritmos efectivos para la manipulaci{\'o}n de offsets de
  hipersuperficies}.
\newblock PhD thesis, Universidad Polit{\'e}cnica de Madrid, 1999.

\bibitem{Sendra2000}
J.R. Sendra and J.~Sendra.
\newblock {\em Algebraic analysis of offsets to hypersurfaces}.
\newblock {\em \em Mathematische Zeitschrift}, 234(4):697--719, 2000.

\bibitem{Sendra2000a}
J.R. Sendra and J.~Sendra.
\newblock {\em Rationality Analysis and Direct Parametrization of Generalized
  Offsets to Quadrics}.
\newblock {\em \em Applicable Algebra in Engineering, Communication and
  Computing}, 11(2):111--139, 2000.

\bibitem{Sendra2007}
J.R. Sendra, F.~Winkler, and S.~P{\'e}rez-D{\'\i}az.
\newblock {\em Rational Algebraic Curves---A Computer Algebra Approach}.
\newblock Springer-Verlag, Heidelberg, 2007.

\bibitem{Shafarevich1994}
I.R. Shafarevich.
\newblock {\em Basic Algebraic Geometry, Vol. 1}.
\newblock Springer-Verlag, 1994.

\bibitem{Snapper1971}
E.~Snapper and R.J. Troyer.
\newblock {\em Metric affine geometry}.
\newblock Academic Press, New York, 1971.

\bibitem{Walker1950}
R.J. Walker.
\newblock {\em Algebraic curves}.
\newblock Princeton, New Jersey, 1950.

\end{thebibliography}
\addcontentsline{toc}{section}{Bibliography}
\label{Bibliografia}

\end{document}